\documentclass[letterpaper,11pt,oneside,reqno]{amsbook}
\usepackage{bm}
\usepackage{mathrsfs}
\usepackage{amsfonts,amsmath, amssymb,amsthm,amscd,stmaryrd,centernot,mathtools,placeins}
\usepackage{pgfplots}
\pgfplotsset{compat=1.17}
\usepackage{tikz}
\usepackage{siunitx}
\usepackage[height=9.6in,width=5.95in]{geometry}
\usepackage[mpexclude,DIV13]{typearea}
\usepackage{verbatim}
\usepackage{hyperref}
\usepackage{graphicx}
\usepackage[latin1]{inputenc}
\usepackage{latexsym}
\usepackage{lscape}
\usepackage{epsfig}
\include{bibtex}

\usepackage{setspace}

\usepackage{etoolbox}

\usepackage{wrapfig}

\BeforeBeginEnvironment{theorem}{\vspace{-\parskip}}
\AfterEndEnvironment{theorem}{\vspace{-\parskip}}
\BeforeBeginEnvironment{lemma}{\vspace{-\parskip}}
\AfterEndEnvironment{lemma}{\vspace{-\parskip}}
\BeforeBeginEnvironment{definition}{\vspace{-\parskip}}
\AfterEndEnvironment{definition}{\vspace{-\parskip}}

\def\fff#1{&{{\pageref{#1}}}\cr}
\def\hfff#1{\label{#1}}

\makeatletter
\pretocmd{\@starttoc}{\begingroup\setlength{\parskip}{0pt}}{}{}
\apptocmd{\@starttoc}{\endgroup}{}{}
\makeatother

\usepackage{parskip}
\newcommand{\e}[0]{\epsilon}

\newcommand{\PP}{\ensuremath{\mathbb{P}}}

\newcommand{\N}{\ensuremath{\mathbb{N}}}
\newcommand{\R}{\ensuremath{\mathbb{R}}}

\newcommand{\Z}{\ensuremath{\mathbb{Z}}}

\newcommand{\E}[0]{\mathbb{E}}

\newcommand{\xmac}{\xi}
\newcommand{\xord}{x}
\newcommand{\ximac}{\varphi}
\newcommand{\battlefield}{v}

\newcommand{\Dslope}{Q}   
\newcommand{\Dsup}{Q_0}   

\newcommand{\imac}{h}
\newcommand{\hmac}{R}
\newcommand{\hmacprime}{H}

\newcommand{\udummy}{r}
\newcommand{\rdummy}{r}

\newtheorem{theorem}{Theorem}[chapter]
\numberwithin{figure}{chapter}

\newtheorem{lemma}[theorem]{Lemma}
\newtheorem{proposition}[theorem]{Proposition}
\newtheorem{corollary}[theorem]{Corollary}

\numberwithin{section}{chapter}

\theoremstyle{definition}

\theoremstyle{definition}

\theoremstyle{definition}
\newtheorem{definition}[theorem]{Definition}

\theoremstyle{definition}

\theoremstyle{definition}
\newtheorem{conjecture}[theorem]{Conjecture}

\def\lora{\longrightarrow}

\newcommand{\dd}{\, {\rm d}}

\newcommand{\mc}{\mathcal}

\newcommand{\rindep}{u}
\newcommand{\xaux}{z}

\newcommand{\rhomac}{\rho}

\newcommand{\incbij}{{\rm IB}}

\newcommand{\tlp}{\mathrm{TLP}}
\newcommand{\tlpkp}{$\mathrm{TLP}(\kappa,\rho)$}
\newcommand{\pfkp}{$\mathrm{PF}(\kappa,\rho)$}
\newcommand{\pfkpspace}{$\mathrm{PF}(\kappa,\rho)$ }
\newcommand{\bbrho}{$\mathrm{BB}(\rho)$}
\newcommand{\bbrhospace}{$\mathrm{BB}(\rho)$ }
\newcommand{\bbrhokappa}{$\mathrm{BB}_\kappa(\rho)$}
\newcommand{\bbrhokappaspace}{$\mathrm{BB}_\kappa(\rho)$ } 
\newcommand{\tlpkpspace}{$\mathrm{TLP}(\kappa,\rho)$ }
\newcommand{\scaledtlpkp}{$\mathrm{ScTLP}(\kappa,\rho)$}
\newcommand{\scaledtlpkpspace}{$\mathrm{ScTLP}(\kappa,\rho)$ }
\newcommand{\lambdamax}{\lambda_{\rm max}}
\newcommand{\lambdamaxkp}{\lambda_{\rm max}(\kappa,\rho)}

\newcommand{\vmac}{v}

\newcommand{\nwithzero}{\N}
\newcommand{\nwozero}{\N_+}

\newcommand{\deltam}{\Delta m}
\newcommand{\deltan}{\Delta n}

\newcommand{\trackS}{\mathsf{S}}

\newcommand{\macnot}{in view of the product notation} 
\newcommand{\macprod}{\prod_{i=0}^k \big( \seemac_i(x) - 1 \big) \, \, \,  \textrm{for $k \in \Z$}}
\newcommand{\macnprod}{x \prod_{i=0}^k \big( \deemac_i(x) - 1 \big) \, \, \,  \textrm{for $k \in \Z$}}
\newcommand{\macab}{$\adefault(x),\bdefault(x):\Z \to (0,\infty)$}
\newcommand{\macmdiff}{\mdefault_{i+1}(x) - \mdefault_{i-1}(x)} 
\newcommand{\macndiff}{\ndefault_{i-1}(x) - \ndefault_{i+1}(x)}
\newcommand{\dragmapright}{\drag^{i\rightarrow}:\tis \to \tis}
\newcommand{\dragmapleft}{\drag^{i\leftarrow}:\tis \to \tis}
\newcommand{\qhalf}{q_j/2}  
\newcommand{\macmid}{m_{-\infty} \pgameplay{S_-}{S_+}{i} \big( \sigma_{i+1} = \infty, E \big)} 
\newcommand{\hproduct}{h_{k+1}^{-1} \cdots h_{-1}^{-1}}
\newcommand{\minanum}{x \, \sum_{k \in \Z} \, \, \prod_{i=0}^k \big( \deemac_i(x) - 1 \big)}
\newcommand{\quadstand}{\big(\asta_i(x),\bst_i(x),\mst_i(x),\nst_i(x): i \in \Z \big)} 
\newcommand{\fmc}{finite mean costs}
\newcommand{\fmcspace}{finite mean costs }
\newcommand{\lessmin}{less than their minimum}
\newcommand{\dragbound}{\egameplay{\drag^{\leftarrow i}(S_-)}{S_+}{i}[P_-] > \egameplay{S_-}{S_+}{i}[P_-]}
\newcommand{\nas}{new altered scenario }
\newcommand{\hatseta}{\hat{S}_-(1,i) = b_i + \eta}
\newcommand{\mupperbound}{\pgameplay{S_-}{S_+}{i}(E_-)\cdot m_{-\infty} +  \pgameplay{S_-}{S_+}{i}(E_+) \cdot m_\infty}
\newcommand{\fpcrhs}{\vert k - \imac \vert \leq \ell \, , \, \vert k-\imac \vert + \ell  \in 2\nwithzero} 
\newcommand{\sminusrhs}{\E_{S_-,S_+}^\imac [P_-]} 
\newcommand{\devMina}{deviating for Mina}
\newcommand{\rvv}{real-valued variables }
\newcommand{\stoppingtime}{the stopping time $\inf \big\{ \ell \in \nwozero : X_\ell = i+1 \big\}$}
\newcommand{\erhs}{\egameplay{S_-}{S_+}{i} \big[ \egameplay{S_-}{S_+}{X(\sigma_{i+1})}[P_+] \big]}
\newcommand{\tworatios}{the two ratios on the right-hand side }

\newcommand{\typnot}{typical notation for a time-invariant strategy pair}
\newcommand{\lawandexpectation}{the law and expectation of play from $i \in \Z$ under $(S_-,S_+) \in \mc{S}^2$}
\newcommand{\escapes}{events of escape $\vert X \vert \to \infty$, left escape $X \to -\infty$ and right escape $X \to \infty$} 
\newcommand{\idlemac}{the idle zone, for $(S_-,S_+) \in \tis^2$, where neither player offers a stake}
\newcommand{\macmmm}{the Mina margin map, sending $\phi_0$ to $\tfrac{n_{-\infty} - n_\infty}{m_\infty - m_{-\infty}}$}
\newcommand{\counterpart}{the counterpart map  $\phi_0 \mapsto \tfrac{n_{-j-1} - n_{k+1}}{m_{k+1} - m_{-j-1}}$  for trail $\llbracket -j-1,k+1 \rrbracket$}

\newcommand{\macindex}{will index four real-valued }

\newcommand{\abmndef}{\adefault(x),\bdefault(x),\mdefault(x),\ndefault(x)}

\newcommand{\macform}{in the form $\big\{ *^{\rm def}_i(x): i \in \Z \big\}$ }

\newcommand{\centralratio}{{\rm CenRatio}}
\newcommand{\tis}{\mc{S}_0}

\newcommand\pgameplay[3]{\PP^{#3}_{#1,#2}}
\newcommand\egameplay[3]{\E^{#3}_{#1,#2}}

\newcommand{\seemac}{c}

\newcommand{\deemac}{d}
\newcommand{\macphi}{\phi}

\newcommand{\macD}{L}

\newcommand{\adefault}{a^{\rm def}}
\newcommand{\bdefault}{b^{\rm def}}
\newcommand{\mdefault}{m^{\rm def}}
\newcommand{\ndefault}{n^{\rm def}}
\newcommand{\asta}{a^{\rm st}}
\newcommand{\bst}{b^{\rm st}}
\newcommand{\mst}{m^{\rm st}}
\newcommand{\nst}{n^{\rm st}}

\newcommand{\minammkp}{\mathcal{M}_{\kappa,\rho}}

\newcommand{\abmnmac}{\textrm{ABMN}}
\newcommand{\abmnkpmac}{\textrm{ABMN}($\kappa,\rho$)}
\newcommand{\mnmac}{\textrm{MN}}

\newcommand{\abmnmacspace}{\textrm{ABMN} }
\newcommand{\abmnkpmacspace}{\textrm{ABMN}($\kappa,\rho$) }
\newcommand{\mnmacspace}{\textrm{MN} }

\newcommand{\rocket}{\textrm{Rocket}}

\newcommand{\drag}{\textrm{Drag}}
\newcommand{\macleft}{\textrm{Left}}
\newcommand{\macright}{\textrm{Right}}

\newcommand{\pimac}{\PP^i}

\newcommand{\sdev}{S_-^{\mathrm{dev}}}
\newcommand{\sdevalt}{S_-^{\mathrm{dev}}[\mathrm{alt}]}

\newcommand{\pidev}{\PP^\imac_{S_-^{\mathrm{dev}},S_+}}
\newcommand{\eidev}{\E^\imac_{S_-^{\mathrm{dev}},S_+}}
\newcommand{\pidevalt}{\PP^\imac_{S_-^{\mathrm{dev}}[\mathrm{alt}],S_+}}
\newcommand{\eidevalt}{\E^\imac_{S_-^{\mathrm{dev}}[\mathrm{alt}],S_+}}

\newcommand{\trap}{\mathsf{Trap}}
\newcommand{\proxytrap}{\mathsf{ProxyTrap}}

\makeatletter
\@namedef{subjclassname@2020}{\textup{2020} Mathematics Subject Classification}
\makeatother

\begin{document}

\title[From tug-of-war to Brownian Boost]{From tug-of-war to Brownian Boost: \\ explicit ODE solutions for player-funded stochastic-differential games}

\author[A. Hammond]{Alan Hammond}
\address{A. Hammond\\
  Departments of Mathematics and Statistics\\
 U.C. Berkeley \\
  899 Evans Hall \\
  Berkeley, CA, 94720-3840 \\
  U.S.A.}
  \email{alanmh@berkeley.edu}
\keywords{Dynamic contests, Hamilton-Jacobi-Bellman equations, infinite-turn  non-zero-sum two-player player-funded stake-governed random-turn games,     resource allocation, stochastic differential games, strategic move evaluation, tug-of-war,  Tullock contests.}

\begin{abstract} 
Brownian Boost is a one-parameter family of stochastic differential games played on the real line 
 in which two players
continuously expend resources in an effort to influence the drift 
 of a randomly diffusing point particle~$X$, with victory determined only at time infinity. Each player chooses a spending
rate over time, incurring cumulative cost, and receives a terminal reward if the particle
ultimately escapes in her preferred direction. We characterise and explicitly compute the
time-homogeneous Markov-perfect Nash equilibria of Brownian Boost, showing that the
derivatives of the players' equilibrium payoff functions solve a pair of coupled first-order
non-linear ODE.

Brownian Boost arises as a fine-mesh, high-noise limit of a two-parameter family of discrete,
player-funded tug-of-war games, one instance of which was studied in~\cite{LostPennies}. 
We analyse the discrete games and Brownian Boost, finding them  to exemplify key features studied in the economics literature of tug-of-war initiated by~\cite{HarrisVickers87}: a battlefield region where players spend heavily; 
  spending rates that decay rapidly but asymmetrically in distance to the battlefield; and an effect of discouragement that renders equilibria fragile under asymmetric perturbation of incentive. 
  
  Tug-of-war has a parallel mathematical literature derived from~\cite{PSSW09}, which solved the scaled fair-coin game in a Euclidean domain via the infinity Laplacian PDE.  By offering an analytic solution to Brownian Boost---a game modelling strategic interaction and
resource allocation---we aim to draw the attention of probabilists and analysts to
mathematically novel tug-of-war games whose delicate analytic structure
admits economically meaningful interpretation. 
\end{abstract}

\medskip

\maketitle
\frontmatter
\chapter*{Preface}

Fix $\rho \in (0,\infty)$.  Mina and Maxine play 
$\rho$-Brownian Boost. 
A point-particle counter $X:[0,\infty) \to \R$\hfff{X} starts at the origin.
Left to its own devices, $X$ is  Brownian; but it is  equipped with a motor that may  impute a drift, left or right, of magnitude at most one. At any given time $t \geq 0$, Maxine and Mina stake money at respective non-negative rates that for convenience we denote by $a(t)$ and $b(t)$, though in principle players' decisions may depend on the entire counter history until the present time. The counter evolves according to the stochastic differential equation 
 \begin{equation}\label{e.countersde}
  \dd X_t =   \frac{a(t)^\rho - b(t)^\rho}{a(t)^\rho + b(t)^\rho} \dd t  + \dd W_t \, ,
 \end{equation}
 where $W:[0,\infty) \to \R$ is standard Brownian motion (and $X(0) = 0$).
 That is to say,  the instantaneously imputed drift equals  $\frac{a(t)^\rho - b(t)^\rho}{a(t)^\rho + b(t)^\rho} \in [-1,1]$, with Mina or Maxine securing drift closer to~$-1$ or~$+1$ by dispersing funds more rapidly, and thus dominating more easily if $\rho$ is higher.

The victor of Brownian Boost is decided only at time infinity. It is Mina when $X(t) \to -\infty$ and Maxine when $X(t) \to \infty$.
The victor alone receives a terminal reward: one for Maxine and $\lambda \in (0.\infty)$ for Mina, so that the given pair $(\rho,\lambda) \in (0,\infty)^2$ specifies the game data.

Players pay for their  stakes, 
each accruing a running cost given by the integral over all positive time of her chosen stake rate.
Each player's reserves are great enough as to permit unfettered access to funds, but the accrued cost for a given player is deducted from her (possibly zero) terminal reward.
In this way, stake decisions are governed not by spending limits but by the cost of expenditure relative to terminal reward.

{\bf } This stochastic differential game is a cousin of tug-of-war games that have been the object of study in a long strand of economics research in dynamic contest theory originating with Harris and Vickers~\cite{HarrisVickers87}. 
Two reactions may seem natural to the reader.
\begin{itemize}
\item {\em Is the game not unplayable?} Brownian Boost is an unusual game from an economic standpoint, because no reward is offered according to how the counter travels over any finite period of time.
The choice is natural, forcing  any structure in equilibrium  to arise endogenously from strategic interaction, rather than from an indelicate preference for some region enforced by design. Yet the emergence of such structure is hardly guaranteed.
The real line is a featureless terrain on which to seek to obtain a local geographic advantage by irrevocably dispensing a valued resource. Why would a player spend at all? A player who has outbid her opponent for a long time will likely have induced near-uniform motion at speed close to one, carrying the counter far from the origin in her preferred direction. But what benefit compensates her spent funds when the present location is indistinguishable from any other, according to the translation invariance of the real line? Surely the funds would have better gone unspent.
\item {\em Discouragement for the less incentivized.} If $\lambda > 1$, then Mina has a greater incentive to win. She may plausibly meet any given stake rate offered by Maxine by proportionally outstaking her, even by a slight margin. The counter would then drift in her favour at a uniform rate, assuring her victory. Maxine would recognize this, accept her certain defeat, and save her running cost; Mina could then win with a very small cost. It would seem then that any asymmetry in incentive is incompatible with non-trivial structure.  
\end{itemize} 

As we will see, the first objection is false. The conclusion reached by the  second is correct for Brownian Boost but false for counterpart discrete games, 
though in a weakened guise it retains some accuracy.

A Tullock contest~\cite{Tullock} with exponent~$\rho \in (0,\infty)$ is a single-round competition in which a player who stakes $a \in (0,\infty)$ against the opponent's $b \in (0,\infty)$ wins with probability $\frac{a^\rho}{a^\rho + b^\rho}$. A natural discrete stake-governed tug-of-war game may be played on a finite integer interval: each round is a Tullock contest that results in a movement of the counter one step to the left or the right, until one or other end of the gameboard is reached. This game is a close cousin of Harris and Vickers', in which two firms invest resources in research over time in an ongoing competition to secure a patent.

In this article, we study a family of such tug-of-war games, 
on a gameboard given by the set of integers 
rather than on a finite interval therein. 
These Trail of Lost Pennies games
  serve as natural discrete counterparts to Brownian Boost.
  The family has two parameters: first, the Tullock exponent $\rho$; 
  second, 
  a noise parameter $\kappa \in (0,1)$ that  governs the strength of mixing of stake-based moves with independent $\pm 1$ moves decided by a fair coin flip.
    
  These games  offer a means to rigorously interpret and study Brownian Boost, and they provide a cohesive framework in which to examine natural hypotheses such as the two mooted above, for both the continuum model and the discrete counterparts. 
  
  We exhibit Brownian Boost, for $\rho \in (0,1]$, as a fine-mesh, high-noise scaling limit of the discrete games. In this representation, we obtain a classification of the Markov perfect Nash equilibria of Brownian Boost, describing them explicitly in terms of a coupled ODE pair. 
 In the discrete games, we find a rich collection of behaviour at Nash equilibrium, exhibiting  robustness under perturbations of incentive from equality that, while  positive, is highly fragile.  
This behaviour illustrates and quantifies economically significant effects such as discouragement, and the presence of battlefield regions wherein players commit resources intensively and away from which these commitments decay rapidly but asymmetrically.


\setcounter{tocdepth}{1}
\tableofcontents

\mainmatter
\chapter{Introduction}

This article has two principal aims: to rigorously analyse Brownian Boost, and to develop and solve a broad class of stake-governed tug-of-war games whose equilibria illuminate themes in dynamic contest theory such as the discouragement effect and robustness (or fragility) under asymmetric perturbations of terminal incentives.

The two aims are pursued in tandem. We introduce a two-parameter family of discrete tug-of-war games on~$\Z$, the Trail of Lost Pennies, one instance of which was studied in~\cite{LostPennies}. A natural solution concept is Markov-perfect Nash equilibrium, under which players seek to maximise mean net receipts, acting solely on the basis of present counter location. 
  We will classify explicitly the resulting equilibria in these discrete games, over a broad swathe of parameter space.   
  
Brownian Boost is then exhibited as a fine-mesh, high-noise scaling limit of the Trail of Lost Pennies. This representation provides a rigorous framework for the stochastic differential game and circumvents the difficulties associated with instantaneous feedback, while the accompanying scaling of the discrete game solutions yields an explicit ODE-based characterization of Brownian Boost equilibria. 
Together, the discrete models and their continuum limit form a fairly broad menagerie of explicitly solved models with which to analyse themes in dynamic contest theory.

\section{Brownian Boost: setup and signposts}\label{s.sdg}

We develop our presentation of the stochastic differential game at the heart of this study, and signpost one of the principal inferences we will reach,
concerning gameplay under the Brownian Boost equilibria classified in the article. 
 
\subsection{Game setup}\label{s.gamesetup}
We introduce some notation to support the description of Brownian Boost offered in the preface. The game $\rho$-Brownian Boost \bbrhospace\hfff{bbrho} is indexed by the parameter $\rho \in (0,\infty)$, which governs the stake-dependent evolution of the counter, together with a further parameter~$\lambda \in (0,\infty)$, which determines the terminal payoffs. Recall that $a(t)$ and $b(t)$ are the rates at which Maxine and Mina spend at time~$t \geq 0$. 
These stakes are raised to the $\rho$\textsuperscript{th} power
to specify the {\em boosts} offered by the players at time~$t$.  
 The instantaneous drift imputed to the counter appearing in the SDE~(\ref{e.countersde}) equals $2p(t)-1$, where $p(t) = \tfrac{a(t)^\rho}{a(t)^\rho+b(t)^\rho}$ is the proportion of the present total boost that is due to Maxine; in this way, the drift interpolates linearly in the proportion $p(t) \in [0,1]$ between the leftmost and rightmost values~$-1$ and~$+1$. 
  
The game ends at time infinity with  {\em left escape} if $E_- := \big\{ X(t) \to - \infty \big\}$ occurs and with {\em right escape} if $E_+  := \big\{ X(t) \to \infty \big\}$ does; the limits are in high~$t$.
Each player receives a terminal payment according to which event occurs. In principle, four parameters are needed, and we may write Maxine's terminal receipt $T_+$ in the form $m_\infty {\bf 1}_{E_+} + m_{-\infty} {\bf 1}_{E_-}$ and 
Mina's $T_-$ as 
$n_\infty {\bf 1}_{E_+} 
+ n_{-\infty} {\bf 1}_{E_-}$.
On the vector $(m_{-\infty},m_\infty,n_{-\infty},n_\infty) \in \R^4$ we impose $m_{-\infty} < m_\infty$ and $n_\infty < n_{-\infty}$
in order that $E_-$ and $E_+$ may rightly be viewed as events of victory for Mina and Maxine. In practice, three trivial symmetries serve to reduce this boundary-data parameter space to a one-dimensional subspace. Suitable pre-game payments to the players permit us to set the losing receipts $m_{-\infty}$ and $n_\infty$ equal to zero. And, since the drift in~(\ref{e.countersde}) is homogeneous of degree zero in the stake pair~$(a,b)$, currency revaluation then permits Maxine's victory receipt to equal one. The remaining free parameter is Mina's winning receipt, which we set equal to the given value $\lambda \in (0,\infty)$\hfff{lambda}. To wit, we take    
  $T_+ = {\bf 1}_{E_+}$ and $T_- = \lambda \cdot {\bf 1}_{E_-}$, as in the preface. 
 
  Maxine's running cost is $R_+ = \int_0^\infty a(t) {\rm d} t$; Mina's,  $R_- = \int_0^\infty b(t) {\rm d} t$. 
   Maxine and Mina's net total payoffs from playing the game are thus equal to $P_+ := T_+ - R_+$ and $P_- := T_- - R_-$.

 Each player seeks to maximize her expected net receipt, $\E \, P_-$ or $\E \, P_+$.
 Each chooses a strategy that determines her stake-rate $a(t)$ or $b(t)$ at every time~$t \geq 0$ in a manner adapted to the game's history until the present moment.
 A Nash equilibrium is a strategy-pair from which neither player would unilaterally choose to deviate given the declared aim of maximizing mean receipt.
 

This description of admissible strategies is imprecise, and we will shortly indicate more about the framework in which we interpret and analyse Brownian Boost. For now, it is useful 
 to explain a basic aspect which is needed to explain our analysis of the game.
In playing \bbrho, a player may in principle draw on a broad range of strategies determined by game history in choosing her stakes. We will restrict attention to a narrower class that includes all viable options according to an intuitively appealing principle akin to the Markov property: 
in the history of gameplay until a given moment, the one piece of data that should be determinative for deciding the stake rate is the present counter position~$X_t$.
Stake pairs that meet this condition are {\em time-homogeneous} and  {\em Markov perfect}, and in shorthand we will call them~{\em time-invariant}. By focussing on such pairs,  we reinterpret the stakes specified in the preface
as profiles $a,b:\R \to [0,\infty)$,
with $a(x)$ and $b(x)$ denoting the rate at which Maxine and Mina stake at any moment $t \geq 0$ for which $X_t = x$. In this way, the Brownian Boost equilibria that we seek to classify are the time-invariant Nash equilibria (or TINE). Each TINE is a map $(a,b): \R \to [0,\infty)$ whose domain is the gameboard~$\R$. A given time-invariant strategy pair, including any TINE, induces gameplay satisfying the SDE~(\ref{e.countersde}) recorded in the preface where the drift term takes the form 
 $\frac{a(X_t)^\rho - b(X_t)^\rho}{a(X_t)^\rho + b(X_t)^\rho}$.

Moscarini and Smith~\cite{MoscariniSmith2007} study a model whose specification is close to the game we consider; formally, it is obtained from Brownian Boost by setting $\rho = 1/2$ and omitting the denominator in the drift term in~(\ref{e.countersde}). For games played on finite intervals, the symmetric Markov-perfect equilibria are characterized via Hamilton--Jacobi--Bellman arguments. We will discuss~\cite{MoscariniSmith2007} and its relation to our work in Section~\ref{s.bbrho} after giving a heuristic derivation of the HJB equations for our system.

\subsection{Some signposts for  Brownian Boost results}\label{s.solvingbb}  
  
  In the preface, we stated two premises concerning Brownian Boost: that it is unplayable because of the translation invariance of its real-line gameboard; and that any incentive asymmetry disrupts non-trivial solution structure. Considering this pair of notions may offer a useful route towards understanding Brownian Boost and its mathematical structure. 
  Later in the introduction, we will introduce the discrete Trail of Lost Pennies games and explicitly classify their TINE; and then scale them to find a counterpart classification for Brownian Boost. Doing so will naturally involve some elaboration
  in notation and statement.
  Here, we offer some signposts to the form of the Brownian Boost equilibria we will exhibit. While what we indicate is merely a fragment of the solution structure, this piece of the picture is enough to permit us to evaluate the validity of the two premises. 
  
  Our analytic presentation of these equilibria is governed by two coupled ODE, the \bbrho-ODE pair~(\ref{bbode}), whose solutions describe (up to sign) the spatial derivatives of the players' value functions; 
  the corresponding  equilibrium stake profiles (or TINE) and the induced counter dynamics are then expressed as simple functionals of this solution pair. 
  Here we describe precisely only one aspect, the dynamics, which admits a simple description.

  The set of TINE in \bbrhospace will be classified when $\rho \in (0,1]$. For each such $\rho$, this collection is a one-dimensional space, invariant under real shifts, and indexed by a `battlefield' value in~$\R$. From the equilibrium $(a,b): \R \to (0,\infty)$ with battlefield value zero, all other TINE are obtained in the form $\big( a ( \bullet  - \battlefield  ), b( \bullet  - \battlefield ) \big)$ as the battlefield value $\battlefield$ varies over~$\R$. We defer recording the form of $(a,b)$, and instead indicate the gameplay that results when the players follow this zero-indexed TINE.
  This gameplay $X:[0,\infty) \to \R$ solves
  \begin{equation}\label{e.xsde}
{\rm d}X_t =  R_\rho(X_t) \, {\rm d}t + {\rm d} W_t \, , \, \, \, \textrm{with $W_t$ standard Brownian motion} \, ,
\end{equation}
an SDE in whose drift term $R_\rho: \R \to (-1,1)$ is given by $R_\rho(u) =  \frac{1-J(u)}{1 + J(u)}$ where $J:\R \to (0,\infty)$ is a decreasing bijection that is the unique solution of the ODE 
$$
\frac{{\rm d} J(u)}{{\rm d} u} \, = \,  -8 \rho^2  \frac{J(u)^2}{(1 + J(u))^2}
$$  
with $J(0) = 1$. The function $R_\rho$ has qualitative similarities to the $(2/\pi)$-multiple of the arctangent: 
it odd, smooth, and satisfies  $R_\rho(u) \to \pm 1$ as $u \to \pm \infty$.

Brownian Boost is a playable game when $\rho \in (0,1]$, in the sense that solutions exist for a standard notion of equilibrium. The first premise advanced in the preface is false because the translation-invariance of the real line is manifest in the space of TINE; this symmetry breaks for any given element of this set. Indeed, the deferred explicit form of the zero-indexed equilibrium $(a,b):\R \to (0,\infty)^2$ comprises two functions of unit order in a compact neighbourhood of the origin, with exponential decay at higher distances. When this TINE is played, both players vie for control of the counter near the origin, rapidly curtailing stakes at further distances. The origin (and its locale) is a battlefield where the outcome is typically decided. 

What of the second premise, that no asymmetry in incentive is tolerated in Brownian Boost? The solution $J$ of the above initial value problem is straightforwardly seen to satisfy $J(-u)=1/J(u)$; we will obtain this equality in the proof of Proposition~\ref{p.sfacts}(1,2). Consequently, the zero-indexed TINE gameplay SDE has a drift term that is an odd function of spatial location (as we stated above).
This gameplay symmetry is in fact incompatible with any asymmetry in terminal receipt. Which is to say: the value $\lambda$ of Mina's terminal receipt in the event of her victory is equal to one. 
This conclusion holds for the zero-battlefield TINE for $\rho \in (0,1]$; and it holds for any other TINE by translation invariance of the TINE set. 
So the second premise is validated for Brownian Boost: no asymmetry in incentive is permissible. 

With these signposts indicated, we are ready to turn to specify the form and solution of the Trail of Lost Pennies games that underpin our Brownian Boost analysis.
The game studied in~\cite{LostPennies} is the special case~$\tlp(1,1)$ of a two-parameter family \tlpkpspace that we consider. 
Our study of these games was prompted by 
a referee of~\cite{LostPennies} who asked whether mixing the game in question with random-walk noise changes the dynamics.
There are two independent directions in the broadening we consider: a noise parameter $\kappa\in(0,1)$ that mixes stake-governed motion with fair coin flips, per the referee's suggestion; and an  exponent $\rho\in(0,\infty)$ that specifies the Tullock contest on which each turn is modelled. 
In this way, the one-dimensional collection of Brownian Boost games, parameterised by $\rho$, is rigorously analysed as 
a regularization of the discrete games. This is achieved by taking  a fine-mesh, high-noise limit $\kappa\searrow 0$.

These discrete games are akin to random-turn tug-of-war games that have been considered since the 1980s in the economics literature of dynamic contests; indeed, 
 beyond their role in interpreting Brownian Boost, the Trail of Lost Pennies games form a rich source for examining prevalent themes in that literature. 
It may thus be instructive to first review this vein of research; and any such review for mathematicians naturally raises the prospect of comparison with a parallel but oddly disjoint tug-of-war research vein in probability and PDE that dates from the 2000s.

In the next section, we indeed offer a review of these two literatures. Three further sections then complete the introduction. The first addresses the setup and solution of the discrete games; the second, how they scale to Brownian Boost, and the resulting classification of equilibria for the continuum game.

While the discouragement effect may be said to be absolute in Brownian Boost---all TINE have $\lambda =1$---the counterpart question for the discrete games has a nuanced answer
which develops that found for $\tlp(1,1)$ in~\cite{LostPennies}:  asymmetric perturbation is highly fragile but not infinitely so, with the degree of this fragility depending on the parameter pair~$(\kappa,\rho)$ in non-trivial and in some cases very surprising ways.
 In the final section of the introduction, 
we examine how the explicit forms for equilibria in the discrete and continuous games cast light on themes in the economics literature such as discouragement, explaining how numerical evidence
supports such  conclusions about fragility under incentive asymmetry.



\section{Tug of war, in economics and mathematics}\label{s.tugofwar}

In 1987, Harris and Vickers~\cite{HarrisVickers87}  introduced a model of a pair of competing firms who spend on research in a race to secure a patent. The principal features they sought to capture were the uncertainty in how effort leads to progress, and the strategic interaction of the competitors as the race unfolds.  In a model they called {\em tug-of-war}, the race is comprised of a sequence of rounds, at each of which a firm expends research effort at a chosen rate, with higher rates improving its odds for the round. Victory for a firm in a given round brings its aim one step closer, and puts its rival's aim one step further away. The race stops when one firm secures the patent and is rewarded with a prize; the opposing firm receives a lesser reward, and both firms must deduct the costs of their respective cumulative research efforts to compute their net receipts. (We will call games with such rules {\em player-funded}.)

In 2009, Peres, Schramm, Sheffield and Wilson~\cite{PSSW09} studied a class of random-turn games, which they also named tug-of-war. Played on a discrete graph~$G=(V,E)$ with boundary~$B$, or in a domain~$D$ in Euclidean space, the game begins with a counter located at a vertex in~$V$ or at an interior point of~$D$. At each turn, a fair coin is flipped and the turn victor moves the counter to a location of his choosing: an adjacent vertex in the discrete setting; and, in the continuous one, a point in $D$ at distance at most $\e$ away, where $\e > 0$ is a parameter fixed for the game. On the boundary $B$ or~$\partial D$ is specified a real-valued payment function $f$.
The game ends when the counter arrives in the boundary with a payment from one player to the other given by the evaluation of $f$ at the terminal counter location. 
 In the discrete setting, the game value~$h(v)$ expressed as a function of starting location~$v$ is the extension of~$f$ that satisfies 
 $h(v)  =   \big( \max_{u \sim v} h(u)   +   \min_{u \sim v} h(u) \big)/2$, the minimum and maximum over neighbours reflecting the choices made when playing from~$v$. The equation is an $\infty$-version of the mean value property in which only the two extremes contribute to the average. In the Euclidean setting, the infinity-harmonic extension of~$f$ to~$D$ 
 is the viscosity solution~$h:D \to \R$ of the infinity Laplace equation 
 $\sum_{i,j} \partial_{x_i} h \, \partial_{x_i x_j} h \,  \partial_{x_j} h = 0$ subject to $h \big\vert_B = f$, whose second derivative in the gradient direction vanishes.
In~\cite{PSSW09}, it is proved the value of tug-of-war 
played on~$D$
converges in the low-$\e$ limit to this extension.

These two seminal contributions each initiated a wave of interest in their respective domains. 

\vspace{-1mm}

\subsection{The economics vein}
The relationship between research allocation and contest outcome is dominant in the economics literature, with works from~\cite{HarrisVickers87} onwards examining the premise that firms contest intensely at a certain pivot location (where the principal battle may be said to take place), with effort that is rapidly decaying away from this location in an asymmetric sense, so that the player close to securing the patent continues to invest an effort that while small exceeds the opposing firm's. The {\em discouragement effect} is another prevalent theme: if one firm will be more rewarded in obtaining the patent, it may plan greater research effort, so that the other, knowing this, may make little, leaving the more incentivized firm in the happy position of winning at little cost.  

One rule to model a single round in player-funded tug-of-war is 
a Tullock contest~\cite{Tullock}.  As we discussed in the preface, this is a single-stage game in which player $A$ stakes $x \in [0,\infty)$ and player $B$, $y \in [0,\infty)$, the contest won by $A$ with probability $\tfrac{x^\rho}{x^\rho + y^\rho}$, where   $\rho \in (0,\infty)$ is the Tullock exponent. 
This family of win probabilities (or contest success functions) is a natural
class of expressions that are homogeneous of degree zero. 
Jia~\cite{Jia2008} has shown that the contest success function must in fact take the Tullock form when finitely many players compete in a contest in which the victor is the player with the highest realized performance, defined as the product of her stake and an independent random factor with a law shared across players, provided that the
vector of win probabilities conditional on  the stake vector satisfies Luce's choice axiom, which encodes the independence of irrelevant alternatives.
From an axiomatic standpoint, Skaperdas~\cite{Skaperdas1996} shows that basic structural requirements---scale invariance, monotonicity, and a  {\em contest}-level independence-of-irrelevant-alternatives principle---again single out the Tullock form.

When $\rho \to \infty$ in the Tullock contest,  all-pay auctions are obtained, in which the higher staking player wins. Player-funded tug-of-war, including the role of battlefields and discouragement, has been studied~\cite{KonradKovenock05,AgastyaMcAfee,Konrad2012} on finite integer intervals with the 
all-pay auction rule used to decide turn victor and in variants~\cite{Hafner16,Hafner17} where a firm is composed of several individuals who are responsible for different payments.
The player-funded game has been studied with the majoritarian objective in which the patent is awarded to the firm who first achieves a certain number of turn victories,
as a model of the premise that early expenditure is decisive, in~\cite{KlumppPolborn}; with intermediate prizes~\cite{KonradKovenock09}; and with discounting~\cite{Gelder} viewed as a dissipator of momentum for the leading player. Two phases of play---no site revisits, then tug of war---occur in a more general graphical framework studied in~\cite{ET19}.    

A separate thread in the economics literature concerns a form of stake-governed tug-of-war where, rather than pay from their own savings, players finance their stakes from a budget allocated to them as part of the game design.  See~\cite{KlumpKonradSolomon} for an analysis with the majoritarian objective, and~\cite{Klumpp} for finite integer intervals. In~\cite{HP2022}, a leisurely or lazy version of the game is studied on a class of trees, with connections drawn to constant-bias tug-of-war.

\subsection{Tug-of-war in PDE and probability}
   As \cite{Manfredi2012} surveys,  the game theory connection identified in~\cite{PSSW09} has attracted a lot of attention from PDE specialists.   
  New boundary rules for  $\e$ tug-of-war led to more regular game value functions in~\cite{ArmstrongSmart2012}. Heavy-tailed moves connect to the infinity fractional Laplacian in~\cite{CCF2012}.
  A noisy version of the game has been considered, in which the counter makes a random displacement of magnitude $c \e$ at the end of each turn.  
The $p$-Laplacian~\cite{Lindqvist} interpolates, as $p$ ranges over $(1,\infty)$, between the total-variation
($p=1$) operator---whose level-set evolution is motion by mean curvature;
see~\cite{KohnSerfaty} for  a game-theoretic approach---and the
infinity Laplacian obtained in the high-$p$ limit.  
In~\cite{PeresSheffield}, the value of the noisy game is shown to converge to a $p$-harmonic extension of boundary data, for $p$ suitably chosen as a function of $c$: the survey~\cite{Lewicka}
 takes this perspective as central.  A variant of this game has been used to study $p$-Laplacian obstacle problems~\cite{LewickaMarta2017}.
The abundant PDE connections of tug-of-war are reviewed in the book~\cite{BlancRossi}. 

\subsection{Weaving together the two research strands}  As of 2025,~\cite{HarrisVickers87} and~\cite{PSSW09} have both garnered over five hundred citations, with no article citing both until~\cite{HP2022,LostPennies}.
Despite the thematic similarities and coincidence of names in the economists' and mathematicians' tug-of-war, the two veins of research appear to have developed quite independently for decades.
The economists' work treats much more developed random decision rules for turn victory than the mathematicians' trivial fair-coin (or constant-bias~\cite{PPS10}) versions, but the mathematicians' studies have a much richer geometric flavour.

Weaving together the two strands is therefore a natural aim, but one should note important structural distinctions, between or within the strands.
One  concerns the asymptotic role of individual turns.
Player-funded tug-of-war has an intrinsically discrete character, with players, even on long integer-interval gameboards, committing significant resources only in a bounded window around a pivot or battlefield location.
By contrast, in $\varepsilon$-tug-of-war on Euclidean domains, individual turns have asymptotically negligible weight, and analytic connections emerge, via PDE.

A second distinction, for the stake-governed games, concerns the source of funding.
In player-funded games, strategic expenditure is endogenous, while in budget-allocated games the spending constraint is imposed by the game design.
Brownian Boost occupies an intermediate position.
It is player-funded, yet shares the asymptotic feature that individual actions carry vanishing weight; its analytic structure is governed by coupled ODE rather than PDE.

When budgets are allocated and play takes place on a fine mesh in a Euclidean domain, the infinity-harmonic structure of value functions identified in discrete settings in~\cite{HP2022} may plausibly give rise to PDE connections, with the optimal stake proportion (shared by the players) determined formulaically by the associated infinity-harmonic value function.

In this sense, the present work and~\cite{HP2022} outline two coherent routes by which economically motivated contest models may merit the attention of analysts and probabilists.



\section{The Trail of Lost Pennies}

In three subsections, 
we specify the discrete games; derive a coupled
system of four equations indexed by~$\Z$ that govern equilibrium values and stakes; and present explicit solution formulas for this system.

\subsection{Specifying the discrete games}\label{s.tlp}

Let $(\kappa,\rho) \in (0,1] \times (0,\infty)$. In brief, the game \tlpkpspace is player-funded tug-of-war on~$\Z$ with turns decided with probability~$\kappa$ by a Tullock contest of exponent~$\rho$, and otherwise by a fair coin flip.

More thoroughly:~\tlpkpspace\hfff{tlpkp} is also specified by  a quadruple $(m_{-\infty},m_\infty,n_{-\infty},n_\infty) \in \R^4$
that satisfies $m_{-\infty} < m_\infty$ and $n_\infty < n_{-\infty}$, and an integer starting location $\ell \in \Z$.
 The counter~$X$ makes~$\pm 1$ moves at each turn, starting at $X(0) = \ell$.
At the start of the $(k+1)$\textsuperscript{st} turn, for $k \in \N$ (including zero), the counter locations, given by $X$ on the integer interval~$\llbracket 0,k \rrbracket$\hfff{intint}, form the history, including the present counter location~$X(k)$. The turn begins with a request for a non-negative stake from each player: say $S_-(k)$ for Mina and $S_+(k)$ for Maxine. The stakes are collected and held in reserve. The umpire now tosses a coin whose sides are marked  {\em stake} and {\em flip} that lands stake with probability~$\kappa$\hfff{kappa}. When the coin lands, the umpire announces suitably `the turn is stake' or `the turn is flip'.    

If the turn is stake,  a coin is tossed that lands heads with probability~$\frac{S_+(k)^\rho}{S_-(k)^\rho + S_+(k)^\rho}$ determined by the $\rho$\textsuperscript{th}\hfff{rho} stake powers. Should neither player offer a positive stake, a fair coin is used. If the coin lands heads, Maxine wins the turn; tails, and Mina does.
If the turn is flip, the coin used is fair. 

The turn victor moves the counter one unit to the left or the right, so that the value of~$X(k+1)$ is recorded. Our specification will make it clear that it is always in Mina's interest to move left and in Maxine's to move right, and we encode these choices in the rules.

The game is being played on~$\Z$ and is necessarily of infinite duration. Its victor is Maxine if the counter evolution $X: \N \to \Z$\hfff{counter}
 satisfies the right-escape event $E_+ := \big\{ X(n) \to \infty \big\}$; and it is Mina if left-escape $E_- :=  \big\{ X(n) \to - \infty \big\}$ occurs. When escape~$E := E_- \cup E_+$\hfff{escape} fails to occur, the game is called {\em unfinished}. 

When Maxine wins a game of \tlpkp, she receives a terminal payment of~$m_\infty$, while Mina receives~$n_\infty$. When Mina wins, she receives~$n_{-\infty}$ and Maxine,~$m_{-\infty}$. Note that the pair of bounds on the boundary data quadruple serve to enforce the preference of Mina to play left and Maxine right.
When the game is unfinished, the terminal payment to Maxine is~$m_*$ and to Mina it is~$n_*$, where $m_*$ and $n_*$ are fixed real numbers that satisfy $m_* < m_{-\infty}$ and $n_* < n_\infty$: outcomes worse than losing the game, for both players.

Players are unrestricted in their choice of stake at each turn, but each must pay all of her stakes from her own funds. As such, Maxine and Mina accrue running costs 
\begin{equation}\label{e.runningcosts}
C_+ = \sum_{k=0}^\infty S_+(k) \, \, \, \,  \textrm{and}  \, \, \, \, C_- = \sum_{k=0}^\infty S_-(k) \, , 
\end{equation}
where $S_+(k)$ and $S_-(k)$ are their stakes at the  $(k+1)$\textsuperscript{st} turn. These costs are deducted from terminal payments to compute a player's overall net receipt. That is, writing $T_+$ and $T_-$ for the terminal payments, the net receipts for Maxine\hfff{receiptmaxine}  and Mina\hfff{receiptmina} are equal to 
\begin{equation}\label{e.receipt}
P_+ = T_+ - C_+ \, \, \, \, \textrm{and} \, \, \, \,  P_- = T_- - C_- \, .
\end{equation} 
The decisions players face in a game of \tlpkpspace are how much to stake at each turn. In formulating a suitable space of strategies from which the players may choose,
we seek to restrict the space so as to unburden notation while ensuring that players may choose from all plausibly appealing options.

For $k \in \N$, write $\Lambda_k$\hfff{pathspace} for the space of $k$-length paths $\psi: \llbracket 0,k \rrbracket \to \Z$ such that $\vert \psi(\ell + 1) - \psi(\ell) \vert = 1$ for $\ell \in \llbracket 0, k-1 \rrbracket$;
set  $\Lambda = \bigcup_{k =0}^\infty \Lambda_k$. 
Let~$\mc{S}$\hfff{strategy} denote the space of maps $S: \Lambda \to (0,\infty)$. The element $S$ is a deterministic strategy that dictates a stake of $S\big( X \big\vert_{\llbracket 0, k \rrbracket} \big)$ at the $(k+1)$\textsuperscript{st}
turn. In this way, a player decides how much to stake in light of the counter's history $X(0),\cdots$ up to its present location~$X(k)$.

The information permitted is a little limited, but in fact most of the strategies needed for our study make do with even less. 
For time-homogeneous Markov-perfect strategies, the only pertinent data in the record~$X: \llbracket 0,k \rrbracket \to \Z$ available at the outset of the  $(k+1)$\textsuperscript{st} turn is the present counter location~$X(k)$.  As in Section~\ref{s.analytic}, we call any such strategy~$S$, namely one whose value on every path is determined by the path's terminal value, {\em time-invariant}; and write  $\mc{S}_0$\hfff{timeinvariant} for the space of these strategies. 
 When Mina and Maxine play the respective elements of a time-invariant strategy pair $(S_-,S_+) \in \mc{S}_0^2$\hfff{typicalnotation}, we will abusively denote the  pair $(b,a)$, for $a,b:\Z \to [0,\infty)$
 given by
 \begin{equation}\label{e.ba}
  \textrm{$a_i \, = \, S_+(\psi)$ \, and \,  $b_i \, = \, S_-(\psi)$ \, for any $i \in \N$ and  $\psi \in \Lambda_i$} \, . 
\end{equation}

For $(S_-,S_+) \in \mc{S}^2$,
the law of gameplay in \tlpkpspace given $X(0)=\ell$ governed by the strategy pair~$(S_-,S_+)$ will be denoted~$\pgameplay{S_-}{S_+}{\ell}$, with~$\egameplay{S_-}{S_+}{\ell}[\cdot]$\hfff{lawexpect}
  
the corresponding expectation. Note also that the usage  $(S_-,S_+) \in \mc{S}^2$ entails a conflict where the stake offered under $S_-$ at the  $(k+1)$\textsuperscript{st} turn, which is formally 
$S_-\big( X \big\vert_{\llbracket 0, k \rrbracket} \big)$, is referred to simply as $S_-(k)$ in~(\ref{e.runningcosts}). We will continue with the simpler usage in most instances since there is little prospect of confusion.

The pair  $(S_-,S_+) \in \mc{S}^2$ is a Nash equilibrium if 
$$
 \egameplay{S_-}{S_+}{\ell} [P_+] \geq  \egameplay{S_-}{S}{\ell} [P_+]  \, \, \, \,
\textrm{and} \, \, \, \,  \egameplay{S_-}{S_+}{\ell} [P_-]  \geq  \egameplay{S}{S_+}{\ell} [P_-] 
$$
 for all $S \in \mc{S}$ and $\ell \in \Z$.

Let $\mc{N}_{\kappa,\rho} = \mc{N}_{\kappa,\rho}(m_{-\infty},m_\infty,n_{-\infty},n_\infty) \subset \mc{S}^2$ denote the space\hfff{nash} 
of Nash equilibria. 
Under a time-invariant Nash equilibrium\hfff{tine}, which is an element~$(S_-,S_+)$ of $\tis^2$ that satisfies the displayed condition,
neither player would gain in expectation by a unilateral deviation in strategy, including by deviation to strategies in~$\mc{S}$ that are not time-invariant.

\medskip

\medskip

\subsection{Time-invariant Nash equilibria and \abmnkpmacspace solutions}

 \begin{definition}\label{d.quadruple}
\vspace{-2mm}
 For $(S_-,S_+) \in \mc{S}_0^2$, set\hfff{mini}
  $m_i = \egameplay{S_-}{S_+}{i} [ P_+]$ and  $n_i = \egameplay{S_-}{S_+}{i}  [P_-]$ for $i \in \Z$. 
 The values $a_i$ and $b_i$\hfff{aibi} are determined by~(\ref{e.ba}).
  Thus to each time-invariant strategy pair $(S_-,S_+)$ we associate a quadruple $(a,b,m,n):\Z \to [0,\infty)^2 \times \R^2$, and conversely any such quadruple determines~$(S_-,S_+)$.
 \end{definition}

We will record differences of elements in the $m$- and $n$-sequences by setting $m_{i,j} = m_j - m_i$ and $n_{j,i} = n_i - n_j$\hfff{mndiffernces}
whenever  $i,j \in \Z \cup \{ -\infty ,\infty\}$ satisfy~$i < j$.   The $m$-sequence is always increasing and the $n$-sequence decreasing; thus $m_{i,j}$ and $n_{j,i}$ 
are non-negative whenever $i < j$, and in our usage of this notation the pair-index order will always increase for $m$ and decrease for~$n$.

 \begin{definition}\label{d.abmn} 
 Let $(\kappa,\rho) \in (0,1] \times (0,\infty)$. 
The ABMN($\kappa,\rhomac$) system on~$\Z$
 is the set of equations in the four variables $(a_i, b_i,  m_i,n_i) \in (0,\infty)^2 \times \R^2$, indexed by~$i \in \Z$, 
\begin{eqnarray*}
  2 \big(a_i^\rhomac + b_i^\rhomac \big) (m_i + a_i) & = &  \Big( a_i^\rhomac(1-\kappa) + b_i^\rhomac(1+\kappa) \Big) m_{i-1} +   \Big( a_i^\rhomac (1+\kappa) + b_i^\rhomac (1-\kappa) \Big) m_{i+1}  \\
 2 \big(a_i^\rhomac  + b_i^\rhomac  \big) (n_i + b_i) &  = &   \Big( a_i^\rhomac (1-\kappa) + b_i^\rhomac (1+\kappa) \Big) n_{i-1} +   \Big( a_i^\rhomac (1+\kappa) + b_i^\rhomac (1-\kappa) \Big) n_{i+1}    \\
 \big( a_i^\rho + b_i^\rho \big)^2  &  = & 
  \rho \kappa \, a_i^{\rho -1} b_i^\rho  m_{i-1,i+1}  \\ 
  \big( a_i^\rho + b_i^\rho \big)^2  &  = & 
  \rho \kappa \, a_i^\rho b_i^{\rho-1}  n_{i+1,i-1}  \, .
  \end{eqnarray*}
 where $i$ ranges over $\Z$. We will call the respective equations \abmnmac($i$) for $i \in \{1,2,3,4\}$.
 \abmnmac(3,4) would require a convention to interpret for $\rho \in (0,1)$ were one of $a_i$ or $b_i$ to vanish, but note that, by definition, we take every $a$- and $b$-value to be positive.
 
 The space of solutions~$(a,b,m,n):\Z \to (0,\infty)^2 \times \R^2$\hfff{positiveabmn} will be denoted \abmnkpmac.
  An element is said to have boundary data $(m_{-\infty},m_\infty,n_{-\infty},n_\infty)$\hfff{boundarydata}
 when 
 \begin{equation}\label{e.boundarydata}
 \lim_{k \to \infty} m_{-k} = m_{-\infty} \, \, \, , \, \, \,
 \lim_{k \to \infty} m_k = m_\infty \, \, \, , \, \, \, 
 \lim_{k \to \infty} n_{-k} = n_{-\infty} \,\,\,\,
 \textrm{and}
 \,\,\,\, \lim_{k \to \infty} n_k = n_\infty \, . 
 \end{equation}
 \end{definition}
On this data, we will impose that 
\begin{equation}\label{e.quadruple}
m_{-\infty}<m_\infty \, \, \, \,  \textrm{and} \, \, \, \,   n_\infty < n_{-\infty} \, .
\end{equation}
The next result states the basic relationship between the trail game and \abmnmac: a time-invariant strategy pair is a Nash equilibrium if and only if it is the $(b,a)$-projection of an element of \abmnkpmac. Note that the assertion is made only under the condition that $\rho \leq 1$. 
 \begin{theorem}\label{t.nashabmn}
 Let  $(\kappa,\rho) \in (0,1]^2$, and let $(m_{-\infty},m_\infty,n_{-\infty},n_\infty) \in \R^4$ satisfy~(\ref{e.quadruple}). 
 \begin{enumerate}
 \item
 Suppose that $(S_-,S_+) \in \mc{S}_0^2$ is an element of ~$\mc{N}_{\kappa,\rho}(m_{-\infty},m_\infty,n_{-\infty},n_\infty)$. 
 The quadruple $\big\{ (a_i,b_i,m_i,n_i): i \in \Z \big\}$ associated to~$(S_-,S_+)$ by Definition~\ref{d.quadruple}  is an element of \abmnkpmac, with boundary data $(m_{-\infty},m_\infty,n_{-\infty},n_\infty)$.
   \item Conversely, if $\big\{ (a_i,b_i,m_i,n_i) \in (0,\infty)^2 \times \R^2 : i \in \Z \big\}$ with boundary  data $(m_{-\infty},m_\infty,n_{-\infty},n_\infty)$ belongs to \abmnkpmac, then the associated pair
  $(S_-,S_+) \in \mc{S}_0^2$ 
    lies in $\mc{N}_{\kappa,\rho}(m_{-\infty},m_\infty,n_{-\infty},n_\infty)$.
  \end{enumerate}
 \end{theorem}

The next two results state basic aspects of how boundary data determines whether the \abmnmacspace system is solvable. When operating with \abmnkpmac, without regard to the game~\tlpkp, we need typically demand only that the pair~$(\kappa,\rho)$ satisfy a {\em weaker} condition than membership of the box~$(0,1]^2$. This condition takes the form $(\kappa,\rho) \in W$, where we set\hfff{weakregion}
\begin{equation}\label{e.weakregion}
W \, = \, \big\{ (\kappa,\rho) \in (0,1] \times (0,\infty):  \rho^2 \kappa \leq 1 \big\} \, .
\end{equation}
 The hypothesis $(\kappa,\rho) \in W$ will be recalled from time to time in our study of \abmnkpmac, but in practice it is  {\em almost always in force}.

 \begin{theorem}\label{t.abmnpositive}
 Let $(\kappa,\rho) \in W$ and $(a,b,m,n) \in$ \abmnkpmac.
 \begin{enumerate}
 \item  For $i \in \Z$, $m_{i+1} > m_i$ and $n_i > n_{i+1}$.
 \item The boundary conditions satisfy $\infty > m_\infty > m_{-\infty} > -\infty$ and $\infty > n_{-\infty} > n_\infty > -\infty$.
 \end{enumerate}
 \end{theorem}

 The {\em Mina margin}\hfff{minamargin} of a solution $(a,b,m,n) \in$ \abmnkpmacspace 
 is set equal to $\frac{n_{\infty,-\infty}}{m_{-\infty,\infty}}$.
 This real-valued quantity has a fundamental role to play in determining whether the \abmnkpmacspace  system can be solved, as we now see. 
 
\begin{definition}\label{d.lambdamax}
For $(\kappa,\rho) \in W$, set\hfff{lambdamax}
$$
 \lambdamax(\kappa,\rho) \, = \,  \sup \, \left\{  \frac{n_{\infty,-\infty}}{m_{-\infty,\infty}}: (a,b,m,n) \in \abmnmac(\kappa,\rho) \right\} \, .
$$
\end{definition}

 \begin{theorem}\label{t.minamarginvalues}
The function  $(\kappa,\rho) \to \lambdamax(\kappa,\rho)$ maps $W$ to $[1,\infty)$. Let $(\kappa,\rho) \in W$, and consider 
 $(m_{-\infty},m_\infty,n_{-\infty},n_\infty) \in \R^4$ with $m_{-\infty} < m_\infty$ and  $n_\infty < n_{-\infty}$.
 An element of \abmnkpmacspace exists with this boundary data quadruple if and only if   $\frac{n_{\infty,-\infty}}{m_{-\infty,\infty}} \in \big[\lambdamax(\kappa,\rho)^{-1},\lambdamax(\kappa,\rho)\big]$.
\end{theorem}

  \subsection{Explicit \abmnmacspace solutions}\label{s.solvingabmn}
 \mbox{}\medskip
   \subsubsection{Ingredients for solving \abmnmac}\label{s.basicfunctions}   
 Some basic functions are needed in preparation for an explicit solution of the \abmnkpmacspace equations. 
 \begin{definition}\label{d.fourfunctions}
We define four real-valued functions $\gamma,\delta,\phi_0,\phi_1$\hfff{fourfunctions}
of the triple $(\kappa,\rhomac,\beta) \in 
W \times (0,\infty)$, where the trail game parameters $(\kappa,\rhomac)$ are now accompanied by $\beta \in (0,\infty)$. These are
\begin{eqnarray}
 \gamma(\kappa,\rhomac,\beta)  & = & \frac{(1-\kappa)\beta^{2\rhomac}  + 2(1-\rhomac\kappa)\beta^\rhomac + 1 +\kappa}{2 (1+\beta^\rhomac)^2}  \, , \label{e.cbeta.rho} \\
\delta(\kappa,\rhomac,\beta) & = & \frac{(1-\kappa)\beta^{2\rhomac} + 2(1+\rhomac\kappa) \beta^\rhomac + 1 + \kappa}{2(1+\beta^\rhomac)^2} \, ,  \nonumber \\ 
 \phi_0(\kappa,\rhomac,\beta) & = &  \frac{\beta \Big( (1-\kappa)\beta^{2\rhomac} + 2(1+\kappa\rhomac)\beta^\rhomac  + 1 +\kappa\Big)}{ (1-\kappa)\beta^{2\rhomac}  + 2(1-\kappa\rhomac)  \beta^\rhomac + 1+\kappa}  \label{e.phizero.rho} 
\end{eqnarray}
and
\begin{equation}\label{e.phione.rho}
\phi_1(\kappa,\rhomac,\beta) \, = \, \frac{\beta \Big( (1+\kappa)\beta^{2\rhomac} + 2(1-\kappa\rhomac) \beta^\rhomac + 1 -\kappa \Big)}{(1+\kappa)\beta^{2\rhomac} + 2(1+\rhomac\kappa)\beta^\rhomac + 1 -  \kappa} \, . 
\end{equation}
\end{definition}
The four functions are positive, because our minimal hypothesis, that $(\kappa,\rho)$  belongs to the set~$W$ specified in~(\ref{e.weakregion}), implies that $\kappa$ and $\kappa \rho$ are at most one, and  every displayed coefficient is then non-negative. 

The map $s$ defined by $s(\phi_0) = \phi_1$, and its forward and backward iterates,  are also fundamental in solving the \abmnmacspace system.
\begin{definition}\label{d.scd}
Let $(\kappa,\rho) \in W$.
\begin{enumerate}
\item As we will show in Lemma~\ref{l.incphi}, $\phi_0(\kappa,\rhomac,\bullet)$ and $\phi_1(\kappa,\rhomac,\bullet)$ are  increasing bijections on $(0,\infty)$. Consequently, the map that sends $\phi_0 \in (0,\infty)$ to $\phi_1$ is well defined. We label this function~$s:(0,\infty) \to (0,\infty)$\hfff{smap}. Thus, for any given $x \in (0,\infty)$, $s(x) = \phi_1(\kappa,\rhomac,\beta)$ for the unique value of $\beta \in (0,\infty)$ for which $\phi_0(\kappa,\rhomac,\beta) = x$. 
\item We further define functions $c,d:(0,\infty) \to (0,\infty)$ by taking $c = 1/\gamma$ and $d = 1/\delta$, with the argument of $c$ and $d$ being $x = \phi_0$ in the same sense as above. Which is to say, we set
$c(x) = 1/\gamma(\kappa,\rhomac,\beta)$ and $d(x) = 1/\delta(\kappa,\rhomac,\beta)$ , the right-hand sides specified by Definition~\ref{d.fourfunctions}, with the value of 
 $\beta \in (0,\infty)$ being the unique choice such that $\phi_0(\kappa,\rhomac,\beta) = x$.   
\end{enumerate}
\end{definition}
Here is notation for the two-sided $s$-orbit.
\begin{definition}\label{d.stabc}
Let $s_{-1}:(0,\infty) \to (0,\infty)$ denote the inverse of~$s$.
Define functions $s_i:(0,\infty) \to (0,\infty)$ indexed by $i \in \Z$. First set $s_0(x) = x$ for $x \in (0,\infty)$. 
Then iteratively specify forward and backward orbits, $s_i(x) = s \big( s_{i-1}(x) \big)$ and  $s_{-i}(x) = s_{-1} \big( s_{-(i-1)}(x) \big)$ 
for $i \in \nwozero$ and $x \in (0,\infty)$.

Set $\seemac_j,\deemac_j:(0,\infty) \to (0,\infty)$, $j \in \Z$, via $\seemac_j(x) =  \seemac (s_j(x))$ and  $\deemac_j(x) =  \deemac (s_j(x))$.
\end{definition}
As we will see in Proposition~\ref{p.sminusone}, the inverse map
$s_{-1}(x)$ is equal to $1/s(1/x)$. 
  
  \medskip

\subsubsection{The solution formulas}\label{s.solutionformulas}

Here we present an explicit form for all members of \abmnkpmac.
The boundary condition
\begin{equation}\label{e.boundarycondition}
 \textrm{$(m_{-\infty},m_\infty,n_{-\infty},n_\infty) \in \R^4$ satisfies $m_{-\infty}<m_\infty$ and  $n_\infty < n_{-\infty}$} \, . 
 \end{equation}
 We may, and will, harmlessly suppose that $m_{-\infty} = 0$ and $n_\infty = 0$, conditions that correspond to zero terminal payment for a player who loses a game of \tlpkp. Indeed, the transformation $(m_i,n_i,a_i,b_i) \to (m_i + \psi,n_i + \zeta,a_i,b_i)$, $i \in \Z$, for arbitrary $(\psi,\zeta) \in \R^2$, maps the \abmnkpmacspace solution space to itself. With  $m_{-\infty} = n_\infty = 0$ set, a further trivial symmetry is manifest via  dilation of real quadruples by an arbitrary positive real: this transformation is a revaluation of currency that also maps the solution space to itself. 

Given the four parameters in~(\ref{e.boundarycondition}) and the three noted symmetries, we may expect the reduced solution space to be parametrized by one free parameter. What is a natural choice for this? We propose two, one local, the other global.  
 For $(a,b,m,n) \in$  \abmnkpmac, the local choice is the  {\em central ratio}
\hfff{centralratio} 
$\centralratio$, which we set to be $\frac{n_{0,-1}}{m_{-1,0}}$.
The global choice is  the solution's {\em Mina margin} which, recall,
is defined to be $\tfrac{n_{\infty,-\infty}}{m_{-\infty,\infty}}$, or $\tfrac{n_{-\infty}}{m_\infty}$ given our assumptions.  Both of these parameters offer measures of relative incentive: the central ratio
expresses the average gain involved for Mina in winning a turn at which the counter starts at the origin relative to the average loss for Maxine due to the same outcome; while the Mina margin measures the terminal reward for Mina relative to this reward for Maxine.
The local choice is useful for describing explicit formulas for solutions.
The global choice is less useful as a parameter, because it does not bijectively index solutions up to symmetry, but this global statistic is important for understanding the game-theoretic consequences of the form of the solutions.

In summary, then, taking $m_{-\infty} = n_\infty = 0$, and expressing the choice of currency valuation by means of the parameter $m_{-1,0} = m_0 - m_{-1} \in (0,\infty)$, we will express our explicit solutions by working with the local choice of the remaining free parameter: we will set $n_{0,-1}/m_{-1,0}$ equal to a given value $x \in (0,\infty)$.

\begin{definition}\label{d.zdef}
For a sequence $h$, write as usual $\prod_{i=0}^k h_i = h_0 \cdots h_k$ for $k \in \N$. A device extends this notation to negative $k \in \Z$: we set 
$$\prod_{i=0}^k h_i \, = \, \begin{cases}
  \, 1  &  \text{for $k=-1$} \\
 \,  \hproduct  &  \text{for $k \leq -2$} \, .
\end{cases}
$$
Let $x \in (0,\infty)$. This parameter \macindex sequences 
$$
\abmndef :\Z \to (0,\infty)
$$ 
which we denote \macform for $* \in \{a,b,m,n\}$. The resulting~$(a,b,m,n)$ is a normalized or `default' quadruple that (as we will state shortly) solves the \abmnmacspace system. 

We first specify $\mdefault(x):\Z \to \R$. This increasing sequence is given by
\begin{equation}\label{e.mdef}
\mdefault_{-\infty}(x) = 0 \, , \, \, \, \, \textrm{and} \, \, \, \, \mdefault_{k+1}(x)- \mdefault_k(x) \, = \,  \kappa \macprod \, ,
\end{equation}
in the notation of Definition~\ref{d.stabc}.
Note that $\mdefault_0(x) - \mdefault_{-1}(x) = \kappa$ \macnot.

The decreasing sequence   $\ndefault(x):\Z \to \R$ satisfies 
$$
\ndefault_\infty(x) = 0 \, , \, \, \, \,   \textrm{and} \, \, \, \, \ndefault_k(x)- \ndefault_{k+1}(x)  \, = \, \macnprod \, . 
$$
Note that $\ndefault_{-1}(x) - \ndefault_0(x) = \kappa x$.

To specify \macab, we set 
$$
M_i(x) = \macmdiff
  \, \, \, \, \textrm{and} \, \, \, \, N_i(x) = \macndiff
$$ 
for $i \in \Z$. We further write 
$$
\adefault_i(x) = \frac{\kappa \rho \, M_i(x)^{1+\rho} N_i(x)^\rho}{\big(M_i(x)^\rho+N_i(x)^\rho\big)^2} \, \, \, \, \textrm{and} \, \, \, \, 
\bdefault_i(x) = \frac{\kappa \rho  \,  M_i(x)^\rho N_i(x)^{1+\rho}}{\big(M_i(x)^\rho+N_i(x)^\rho\big)^2} \, .
$$
\end{definition}

\begin{theorem}\label{t.defexplicit}
Let $(\kappa,\rho) \in W$ and $x \in (0,\infty)$. A quadruple sequence  $(a,b,m,n):\Z \to \R^4$
is an element of  \abmnkpmacspace  satisfying $m_{-\infty} = n_\infty = 0$
and  $\centralratio = x$ if and only if $(a,b,m,n)$ is the dilation by some factor $\mu \in (0,\infty)$ of the  sequence  $\big(  \big( \adefault_i(x),\bdefault_i(x),\mdefault_i(x),\ndefault_i(x) \big): i \in \Z \big)$  specified in Definition~\ref{d.zdef}.
The value $m_{-1,0} = m_0 - m_{-1}$ of the solution is equal to~$\mu \kappa$.
\end{theorem}

We distinguish two choices\hfff{defaultstandard} of currency revaluation for solutions with central ratio~$x$. The default solution has $\mu = 1$.
The other choice is $\mu = \mdefault_\infty(x)^{-1}$ where $\mdefault_\infty(x) = \kappa \sum_{k \in \Z} \prod_{i=0}^k \big( c_i(x) - 1 \big)$ is Maxine's default prize.
This solution~$(a,b,m,n)$  
satisfies the boundary condition $(m_{-\infty},n_\infty) = (0,1)$ and $m_\infty = 1$, which is a natural specification for game data. We label the solution accordingly.
\begin{definition}\label{d.standard}
Let $x \in (0,\infty)$. The unique element of \abmnkpmacspace with $\centralratio = \frac{n_{0,-1}}{m_{-1,0}}$ equal to~$x$ and 
  $(m_{-\infty},m_\infty,n_\infty) = (0,1,0)$ is called {\em standard}. We denote it  
  $$
  \big(\asta_i(\kappa,\rho,x),\bst_i(\kappa,\rho,x),\mst_i(\kappa,\rho,x),\nst_i(\kappa,\rho,x): i \in \Z \big) \, , 
  $$
  omitting the $\kappa$ and $\rho$ arguments when the context is clear.  
\end{definition}
The default and standard normalizations may appear to diverge as $\kappa \searrow 0$, but in fact the sum $\sum_{k \in \Z} \prod_{i=0}^k \big( c_i(x) - 1 \big)$ is $\Theta(\kappa^{-1})$, making the conversion factor bounded.

{\em Remark.} The representation of solutions in Theorem~\ref{t.defexplicit} is governed by orbits of the $(\kappa,\rho)$-parameterised map $s:(0,\infty) \to (0,\infty)$. For generic $(\kappa,\rho) \in (0,1]^2$, there is no explicit form for $s: \phi_0 \mapsto \phi_1$. 
In the case $(\kappa,\rho) = (1,1)$ analysed in~\cite{LostPennies}, 
$\phi_0 = \beta(2\beta+1)$ and $\phi_1 = \beta^2/(\beta +2)$. Since $\phi_0$ and $\phi_1$ appear linearly in coefficients of quadratic equations in the~$\beta$-variable, $s$ has an explicit form, given in~\cite[Definition~$2.18$]{LostPennies}. When $\kappa \in (0,1)$ and $\rho =1$, $\phi_0$ and $\phi_1$ appear in coefficients of cubic equations in~$\beta$, and $s$ may be expressed as a rational function of the unique positive root of a cubic polynomial.  
The generic inexplicitness of~$s$ may disconcert at first, but its implications for this study have been limited to the use of approximate root solving in the numerical investigation of \abmnkpmacspace elements.

\section{Brownian Boost and the high-noise limit}\label{s.bbhighnoise}
 
 We now return to Brownian Boost. 
 Our analysis of the game operates by regularizing it via a fine-mesh high-noise scaling of the  Trail of Lost Pennies.
 In the first of three subsections, we advocate this discretization as a natural means of rigorous analysis of Brownian Boost. 
 In the second, we develop the analytic framework for \bbrhospace equilibria and their properties that was glimpsed in Section~\ref{s.solvingbb}. In the third, we further the presentation of our principal conclusions about \bbrhospace by stating Theorem~\ref{t.lowkappasde}, which describes equilibrium stakes and gameplay in the high-noise regime via this analytic framework.

\subsection{The scaled high-noise trail game as a regularized Brownian Boost}

The space of strategies in \bbrhospace may in principle be chosen to permit decisions on stake rates that are determined by the history of counter evolution and stake profiles up to the present time. 
That said, anomalous outcomes arising from joint adoption of such strategies as `I'll stake twice what she just staked' must be excluded.
In the Elliott-Kalton formalism~\cite{ElliottKalton} of stochastic differential games (as it applies in the non-zero-sum case), upper and lower value functions for each player are specified in terms of pairs of non-anticipating strategies in which one or other player is given first access to information at the instant it arises. When the two values coincide, for both players, they encode the expected payoffs achievable under these non-anticipative strategies.

We do not seek to implement this approach, and instead study a concrete feedback-safe regularization of \bbrho.
An effective time-delay on feedback is implemented by insisting that players commit to stakes for short periods. For $\kappa > 0$, consider a variant game \bbrhokappaspace specified by iterative construction of the counter evolution $X:[\tau_i,\tau_{i+1}] \to \R$ for an increasing sequence of stopping times~$\tau_i$ such that $\tau_0 = 0$ and  $X(\tau_i) \in \kappa \, \Z$ for $i \in \N$. At time~$\tau_i$, Maxine and Mina declare stake rates $a(i)$ and $b(i)$ and spend at these rates during $[\tau_i,\tau_{i+1}]$, with $X$ on this interval given by setting the drift $d_i$ equal to $\frac{a(i)^\rho - b(i)^\rho}{a(i)^\rho + b(i)^\rho}$ solving 
 $\dd X_t = \dd B_t + d_i \dd t $
 from the already constructed starting point $X(\tau_i)$. Set $\tau_{i+1} = \inf  \, \{ \, t > \tau_i : \vert X(t) - X(\tau_i) \vert = \kappa \, \}$.
The $i$\textsuperscript{th} turn of \bbrhokappaspace is called {\em positive} or {\em negative} according to the sign in $X(\tau_{i+1}) = X(\tau_i) \pm \kappa$.
Given the value of $d_i$, the probability $p_i$ that the $i$\textsuperscript{th} turn is positive equals $u(0)$, where 
 \(u\) solves the boundary value problem
$\tfrac12 u''(x) + d_i u'(x) = 0$ with $u(-\kappa)=0$  and $u(\kappa)=1$.  We have then that 
$$
u(x) \;=\; \frac{1 - e^{-2d_i(x+\kappa)}}{1 - e^{-4 d_i \kappa}} \, ,
$$
 so that 
$$
p_i =  
\frac{1}{1+e^{-2d_i\kappa}}  
\, = \, 1/2 +    \frac{a(i)^\rho - b(i)^\rho}{2(a(i)^\rho + b(i)^\rho)} \kappa +  
 \kappa^3 O(1) \, .
$$
We may compare the games \bbrhokappaspace and \tlpkp. When the stake-pair $(a,b)$ is offered at a turn in the latter game, Maxine's win probability equals
 $$
 (1-\kappa)\cdot \tfrac12 + \kappa \cdot \frac{a^\rho}{a^\rho+b^\rho} \,  = \,  1/2 +    \frac{a^\rho - b^\rho}{2(a^\rho + b^\rho)} \kappa \, ,
 $$ 
 the left-hand summands contributed by the turn being flip or stake.
 Maxine's turn-win probabilities coincide to order $O(\kappa^3)$ in the two games. If we code $\pm$-valued sequences indexed by~$\N$  according to whether turns in \bbrhokappaspace are positive or negative, and do likewise in an evident way for \tlpkp, then we see that any given stake-pair sequence, when played in one or other game, gives rise to very similar laws on $\pm$ sequences:
since the per-turn Bernoulli success probabilities differ by $O(\kappa^3)$, the first disagreement between the coupled sequences has mean $O(\kappa^{-3})$ and occurs at a much later time than the $\kappa^{-2}$-scale on which the counter in \bbrhokappaspace has moved a unit order.

Counter displacement at a turn in  \tlpkpspace  has magnitude one, but in  \bbrhokappa, it has magnitude~$\kappa$. And 
while a player in \tlpkpspace simply surrenders her stake at each turn, the counterpart cost in \bbrhokappaspace also involves the duration for which she spends at the committed rate.
For example, the mean running cost for Mina at a turn where she commits to~$b$
equals $\big( \kappa^2 + O(\kappa^4) \big) \cdot b$ where the prefactor is the mean turn duration, which is exactly $\kappa^2$ in the driftless case, with the $O(\kappa^4)$ error 
enough (by a short omitted computation)
to accommodate the drift of magnitude at most one.

As such, \bbrhokappaspace may be more closely compared to a scaled version \scaledtlpkpspace of the Trail of Lost Pennies. The scaled game operates by  the rules of \tlpkpspace with two changes: it is played on $\kappa \, \Z$, not $\Z$; and  the running costs~(\ref{e.runningcosts}) that enter the net receipt formulas~(\ref{e.receipt}) now include $\kappa^2$-prefactors,  $C_\pm = \kappa^2 \sum_{k=0}^\infty S_\pm(k)$.

The effect of these changes  is to put the turn-by-turn counter locations in \bbrhokappaspace and \scaledtlpkpspace on an equal footing, while ensuring consistent units for measuring trail game and Brownian Boost stakes. 
 As a result, for any given strategy pair, the turn-win sequence in the vanishing-$\kappa$ limit is practically indistinguishable between \bbrhokappaspace and \tlpkp, and when compared to \scaledtlpkp, this agreement is accompanied by asymptotically equal mean net receipts for the players and by asymptotically close counter evolutions. In this sense, the status of \bbrhokappaspace as a natural instant-feedback-safe  surrogate for Brownian Boost passes to the scaled trail game  \scaledtlpkp, due to the match in  both payoff structure and gameplay dynamics.

It is natural to pose the problem of determining the `domain of attraction' of discretized approximant games for \bbrho. 
Adapting the methods of~\cite{ElliottKalton1974,EvansSouganidis}, \cite{FlemingSouganidis} addresses this type of question for stochastic differential games of zero-sum; while
 the strongly non-anticipatory framework in~\cite{BCR} is
adapted to the non-zero-sum case.
   Implementing a framework such as~\cite{BCR}'s  rigorously for \bbrhospace would require careful handling of the infinite horizon, non-zero-sum payoffs, and the absence of discounting in Brownian Boost.
Instead we choose  the concrete discretization~\tlpkpspace in the limit of low~$\kappa$ (or high flip noise) as the rigorous point of contact with \bbrho.

\subsection{Analytic formulation and solutions for Brownian Boost equilibria}\label{s.analytic}
Now we present the ODE system, including its solutions and some important properties, that governs our characterization of \bbrhospace equilibria.

Here is the ODE pair that will be shown to govern time-invariant Nash equilibria in \bbrho.
\begin{definition}\label{d.depair}
Let $\rho \in (0,\infty)$.  A pair of differentiable functions $f,g:\R \to (0,\infty)$\hfff{fg} is called a $\rho$-Brownian Boost ODE pair\hfff{bbode} if it satisfies at every point on the real line
\begin{eqnarray}
 2 \rho f^{1+\rho} g^\rho & = & \, \, \, \, \, \, \big( f^{2\rho}- g^{2\rho} \big)f \, + \, \tfrac{1}{2} f' \big(f^\rho + g^\rho \big)^2 \label{e.fg.rho} \, , \\
 2 \rho f^\rho g^{1+\rho} & = & - \, \big( f^{2\rho}- g^{2\rho} \big)g \, - \, \tfrac{1}{2} g' \big(f^\rho + g^\rho \big)^2 \, . \nonumber
\end{eqnarray}
A pair of functions $f,g:\R \to (0,\infty)$ is called {\em default} if $f(0) = 1$ and $g(0) > 0$.
\end{definition}
In our game-theoretic interpretation of the pair $(f,g)$, the domain variable is a location on the real-line gameboard, and $f$ and $-g$ are spatial derivatives of expected net receipt at equilibrium for Maxine and Mina, viewed as functions of initial counter location. The minus sign attached to~$g$ renders this function positive, reflecting that Mina plays to the left and improves her position by doing so.

The ODE pair arises as the coupled system of Hamilton--Jacobi--Bellman~[HJB] equations associated to the non-zero-sum game \bbrho.
In Chapter~\ref{c.brownianboost}, we will explain this connection with a simple but non-rigorous argument.
Using Markovian balance equations and stability under momentary perturbation of stake by a given player, the argument yields necessary conditions for a stake-profile pair $(a,b):\R \to [0,\infty)^2$ to constitute a Nash equilibrium.
Associated to $(a,b)$ are value functions $m,n:\R \to [0,\infty)$, representing the players' expected net receipts as a function of the initial counter location.
Under a differentiability hypothesis, these satisfy $m'>0$ and $n'<0$, and the resulting HJB conditions identify $(f,g)=(m',-n')$ as 
solutions of~(\ref{e.fg.rho}). 
(In the theory of stochastic differential games, formal derivations of HJB equations would suppose sufficient differentiability; but, in contrast to \bbrho, value functions often do not enjoy that regularity, and are instead exhibited rigorously as viscosity solutions~\cite{CrandallLions}: see~\cite{FlemingSouganidis}
and~\cite{BensoussanFrehse} respectively for zero- and non-zero-sum treatments.)

In the preceding section, we identified explicitly the time-invariant Nash equilibria of the discrete Trail of Lost Pennies games.
In the next subsection, we will state how, when $\rho \in (0,1]$,  in the high-noise fine-mesh limit these discrete equilibria converge to Brownian Boost equilibria, and that the limiting equilibrium objects satisfy the ODE pair in Definition~\ref{d.depair}.
Thus, for $\rho \in (0,1]$, the ODE system provides the analytic framework in which time-invariant Nash equilibria in \bbrhospace are classified.

For now, the prospect of such a characterization may provoke the question,
how to solve the above pair of equations? We record the answer next, noting that currency revaluation permits us to consider only default solutions.
Our analytic deductions hold whenever $\rho \in (0,\infty)$, even if the game-theoretic meaning of the \bbrhospace ODE pair is unsettled for~$\rho > 1$.

\begin{definition}\label{d.fg}
For $\rho,\xmac \in (0,\infty)$,
let $S_\rho(\xmac,\bullet):\R \to (0,\infty)$\hfff{srho} denote the unique solution to the differential equation
\[
 \frac{{\rm d}}{{\rm d}u} S_\rho(\xmac,u) \,=\, - \frac{8\rho \, S_\rho(\xmac,u)^{1+\rho}}{\big(1+S_\rho(\xmac,u)^\rho\big)^2} \, ,
 \qquad S_\rho(\xmac,0)=\xmac \, .
\]
Associate to this solution the pair of functions $f_\rho(\xmac,\bullet),g_\rho(\xmac,\bullet): \R \to (0,\infty)$\hfff{fgrho} by means of
\[
f_\rho(\xmac,\rindep) \,=\, \exp\!\left\{ 2 \int_0^{\rindep} \Bigg(
 1 - \frac{2}{\big(1+S_\rho(\xmac,\udummy)^\rho\big)^2}\,\Big(1+(1-\rho)S_\rho(\xmac,\udummy)^\rho \Big)  \Bigg)\, {\rm d}\udummy \right\}
\]
and
$$
g_\rho(\xmac,\rindep) \,=\, \xmac \cdot  \exp\!\left\{ -2 \int_0^{\rindep} \Bigg(
 1 - \frac{2}{\big(1+S_\rho(\xmac,\udummy)^\rho\big)^2}\, \Big((1-\rho)S_\rho(\xmac,\udummy)^\rho + S_\rho(\xmac,\udummy)^{2\rho}\Big) \Bigg)\, {\rm d}\udummy \right\}
$$
for $\rindep \in \R$.
When $\rindep <0$, the integrals are specified in the usual way: $\int_a^b h = - \int_b^a h$ for $a > b$.
\end{definition}
The parameter $\xmac \in (0,\infty)$ acts as an index for the space of flows  $S_\rho(\xmac,\bullet):\R \to (0,\infty)$.
We reserve the symbol~$\xmac$ for this usage, calling it the {\em flow index}\hfff{flowindex}. This index parameterises default solutions of the \bbrhospace ODE pair as follows. 
\begin{theorem}\label{t.fg}
Let $\rho \in (0,\infty)$. The space of default solutions to the system~(\ref{e.fg.rho}) is equal to
$$
\Big\{ \big(  f_\rho(\xmac,\bullet), g_\rho(\xmac,\bullet)  \big): \R \to (0,\infty)^2 \Big\} \, ,
$$
where the index runs over $\xmac \in (0,\infty)$. For each $\xmac$, we have $g_\rho(\xmac,\bullet)=  f_\rho(\xmac,\bullet) S_\rho(\xmac,\bullet)$.
\end{theorem}
In the discrete context, we interpreted the central ratio 
$n_{0,-1}/m_{-1,0}$
 as a measure of local relative incentive for a turn at which the counter is at the origin.
Recall that this parameter acts as the domain variable~$x$ that  parameterises the explicit forms of  \abmnkpmacspace elements recorded in Subsection~\ref{s.solutionformulas}.   A counterpart role for Brownian Boost is played by 
the flow index $\xmac = S_\rho(\xi,0)$, since this quantity  is equal to the ratio  $g_\rho(\xmac,0)/f_\rho(\xmac,0)$  of  the spatial derivative (in absolute value) of expected net receipts evaluated at~$0$.

Given a solution $(f,g)$ of the ODE pair, 
 the associated time-invariant stake profile~$(a,b)$ is recovered explicitly.
\begin{definition}\label{d.arhobrho}
For $\xmac\in(0,\infty)$, define functions $a_\rho(\xmac,\bullet),b_\rho(\xmac,\bullet):\R\to(0,\infty)$\hfff{abrho} by
\[
a_\rho(\xmac,u)
=
2\rho\,
\frac{f_\rho(\xmac,u)^{1+\rho}\,g_\rho(\xmac,u)^\rho}
{\big(f_\rho(\xmac,u)^\rho+g_\rho(\xmac,u)^\rho\big)^2},
\qquad
b_\rho(\xmac,u)
=
2\rho\,
\frac{f_\rho(\xmac,u)^\rho\,g_\rho(\xmac,u)^{1+\rho}}
{\big(f_\rho(\xmac,u)^\rho+g_\rho(\xmac,u)^\rho\big)^2}.
\]
\end{definition}

The default solutions of Definition~\ref{d.depair}
are normalized by a derivative condition imposed at the origin, in correspondence with the scaling symmetry given by currency revaluation. 
For any solution arising from Definition~\ref{d.fg}, however, the associated global quantities
\begin{equation}\label{e.mndef}
m_\rho(\xmac,\infty)
:=
\int_{-\infty}^{\infty} f_\rho(\xmac,r)\,{\rm d}r \, ,
\qquad
n_\rho(\xmac,-\infty)
:=
\int_{-\infty}^{\infty} g_\rho(\xmac,r)\,{\rm d}r
\end{equation}
are finite and positive.
Interpreted as integrals of spatial derivatives of expected net receipts,
these quantities represent the players' terminal rewards in the event of victory.
As such, these quantities are attributes of \bbrhospace game data, and, given our convention for specifying the continuum game, it is natural  to use a normalization, counterpart to the discrete case in Definition~\ref{d.standard}, for which Maxine's winning receipt equals one.
\begin{definition}\label{d.bbstandard}
A solution of the \bbrho-ODE pair $(f,g):\R \to (0,\infty)$ is called {\em standard}\hfff{bbstandard} if $\int_\R f(r) \, {\rm d} r =1$.
\end{definition}
Up to trivial symmetries, the family
\[
\big\{(f_\rho(\xmac,\bullet),g_\rho(\xmac,\bullet)):\xmac\in(0,\infty)\big\}
\]
exhausts all solutions of the $\rho$-Brownian Boost ODE pair. 
The switch from default to standard normalization may be effected by dilating the displayed pair by a suitable positive $\xmac$-dependent factor.


Our analysis of  the \bbrho-ODE system 
will show that in fact the integral expressions specifying $m_\rho(\xmac,\infty)$ and $n_\rho(\xmac,-\infty)$ in~(\ref{e.mndef}) coincide for any given $\xmac \in (0,\infty)$.
This equality is maintained by dilation and thus holds for any normalization of solutions. Consequently,  $\int_\R g(r) \, {\rm d} r =1$ for any standard solution. 

\subsection{Scaled gameplay in the high-noise trail game}

Here we substantiate that for $\rho \in (0,1]$ the 
time-invariant equilibria of \bbrhospace are given by the prescription in Definitions~\ref{d.fg} and~\ref{d.arhobrho}, with a result showing that this description captures (in the limit of low~$\kappa$) all time-invariant stake-profiles and gameplay in the scaled trail game~\scaledtlpkp.

For $(\kappa,\rho,x) \in (0,1]^2 \times (0,\infty)$, recall that the default solution 
$$  \big(\adefault_i(\kappa,\rho,x),\bdefault_i(\kappa,\rho,x),\mdefault_i(\kappa,\rho,x),\ndefault_i(\kappa,\rho,x): i \in \Z \big) 
$$
 is the unique element of \abmnkpmacspace with $\phi_0 = x$ and $m_{-1,0} = \kappa$ (as well as $m_{-\infty} = n_\infty = 0$).
 In view of Theorems~\ref{t.nashabmn} and ~\ref{t.defexplicit}, time-invariant Nash equilibria in \tlpkpspace  are characterized up to the trivial symmetries by these solutions, and  we use them to express our result Theorem~\ref{t.lowkappasde}.

The result has two parts. In its first, we see that stake profiles in \tlpkp, when multiplied by~$\kappa^{-2}$,  mimic profiles arising from Brownian Boost ODE pairs.
In the second, gameplay in \tlpkpspace is scaled as $\kappa \searrow 0$, sped up by a factor of $\kappa^{-2}$. The resulting SDE weak limit is counter evolution in $\rho$-Brownian Boost played at the time-invariant Nash equilibrium (which is unique up to a real shift indexed by the battlefield value, which we take to be zero). The drift coefficient has a simple expression in terms of the ODE solution~$S_\rho$ in Definition~\ref{d.fg}.

The scaling factors cohere with the transform of \tlpkpspace to \scaledtlpkp, which squeezes space and time by respective factors of~$\kappa$ and~$\kappa^2$.

\begin{theorem}\label{t.lowkappasde}
Let $(\kappa,\rho) \in (0,1]^2$ and $\xmac \in  (0,\infty)$.
\begin{enumerate}
\item As $\kappa \searrow 0$,
\begin{align}
\kappa^{-2} \adefault_{\lfloor \kappa^{-1} \rindep\rfloor}(\kappa,\rho,\xmac)
& \, = \, a_\rho(\xmac,\rindep) \;\big(1+\ O(\kappa) \big) \, , \\
\kappa^{-2} \bdefault_{\lfloor \kappa^{-1} \rindep\rfloor}(\kappa,\rho,\xmac)
& \, = \,  b_\rho(\xmac,\rindep)
\;\big(1+\ O(\kappa) \big) \, ,
\end{align}
with the implicit constant in the $O$-terms being uniform in $(\rho,\xmac,\rindep) \in (0,1] \times K \times K$ for compact $K \subset (0,\infty)$. 
\item For $y \in \R$,  let $X_{\kappa,\rho}(y,\bullet):\N \to \Z$ denote the evolution of the counter with $X(0) = \lfloor y \rfloor$
 under the time-invariant Nash equilibrium of battlefield index zero in the game~$\tlp(\kappa,\rho)$.

For $z \in \R$, consider  the scaled process
$$
[0,\infty) \to \R: t \to \kappa \, X_{\kappa,\rho}(\kappa^{-1}z,\kappa^{-2}t) \, ,
$$
whose domain of definition is enlarged from $\kappa^2  \N$ to $[0,\infty)$ by interpolation.

Equip the space $\mc{C}$ of continuous functions $f:[0,\infty) \to \R$ with the topology of uniform convergence on compact intervals.
As \(\kappa \searrow 0\), this process
converges weakly in~$\mc{C}$ to the unique solution $Z:[0,\infty) \to \R$ of the stochastic differential equation
\[
{\rm d}Z_t =  R_\rho(Z_t) \, {\rm d}t + {\rm d} W_t ,
\]
with $Z_0 = z$,
where \(W_t\) is standard Brownian motion. The function $R_\rho:\R \to (-1,1)$ is given by
 $R_\rho(u) = \frac{1 - S_\rho(1,u)^\rho}{1 + S_\rho(1,u)^\rho}$ or equivalently  $\frac{1 - S_1(1, \rho^2 u)}{1 + S_1(1, \rho^2 u)}$.  
 It has asymptotics 
 $$
  R_\rho(u) = 1 - \frac{1}{4 \rho^2 u} + O(u^{-2}) \, \, \, \textrm{as $u \to \infty$} \, \, \, \, \textrm{and} \, \, \,\, R_\rho(u) = -1 + \frac{1}{4 \rho^2 \vert u \vert} + O(u^{-2})  \,  , \,  \textrm{as $u \to -\infty$} \, .
 $$
 \end{enumerate}
\end{theorem}
{\em Remarks: (1).} In introducing the theorem, we referred to its description of \bbrhospace TINE as being up to trivial symmetries. To be precise, our convention for terminal receipts in \bbrhospace is $\lambda \cdot {\bf 1}_{E_-}$ for Mina and ${\bf 1}_{E_+}$ for Maxine. Since 
$m_\rho(\xmac,\infty)$
and 
$n_\rho(\xmac,-\infty)$
in~(\ref{e.mndef})
are equal for any given $\xmac \in (0,\infty)$, we have $\lambda=1$ for
all \bbrhospace TINE. The correct normalization of such TINE is the standard one, with the $\xmac$-indexed TINE taking the form  
$\R \to (0,\infty)^2: u \to Q(\xmac)^{-1} \big( a_\rho(\xmac,u),  b_\rho(\xmac,u) \big)$, with normalization factor 
$Q(\xmac)$ equal to the common value of
$m_\rho(\xmac,\infty)$
and 
$n_\rho(\xmac,-\infty)$.

{\em (2).} The function $J(u)$ in the earlier signpost~(\ref{e.xsde}) equals $S_1(1,\rho^2 u)$. The form for $J'$ recorded there is given by taking $\rho=1$ in Definition~\ref{d.fg}
with a linear change of variable: see Lemma~\ref{l.de}.

\section{Robustness of inferences: the discouragement effect and asymmetric decay}\label{s.robust}

Here we examine the implications of the games \bbrhospace and \tlpkpspace for some of the principal themes in dynamic contest theory seen in the economics literature:
how rapidly and asymmetrically stakes decay away from a battlefield at which they concentrate; and the degree to which a less incentivized player may be discouraged from staking, permitting her opponent to win the contest at little cost.

\subsection{Brownian Boost equilibrium behaviour}

We present properties of the Brownian Boost equilibrium objects from Section~\ref{s.analytic} pertinent for comparing to dynamic contest theory.

Default solutions of the \bbrhospace ODE pair are indexed by the flow index $\xmac \in (0,\infty)$, and take the form 
$\big( f_\rho(\xmac,\bullet),g_\rho(\xmac,\bullet) \big):\R \to (0,\infty)^2$
as described in Theorem~\ref{t.fg}.
In the signposts Section~\ref{s.solvingbb}, we indicated that the space of time-invariant Nash equilibria
is invariant under spatial shifts, and is naturally indexed by a real
{\em battlefield} value~$\battlefield$.
This translation invariance acts by shifting the spatial argument of the equilibrium profiles.
The time-invariant Nash equilibria associated to default solutions
are given by the stake-profile pairs $\big( a_\rho(\xmac,\bullet),b_\rho(\xmac,\bullet) \big)$,
as specified in Definition~\ref{d.arhobrho}.
The flow-index and battlefield-value descriptions are reconciled in Section~\ref{s.odeproofs} 
where a one-to-one correspondence $(0,\infty)\longleftrightarrow\R$ between the two indices is presented.


When $\xmac=1$, the associated battlefield value $\battlefield$ is the origin.
It is enough to study this choice.
The next result does so, capturing the rapid decay and the asymmetry that are characteristic of all solutions' behaviour.
Before stating it, we note that the equality $m_\rho(1,\infty) = n_\rho(1,-\infty)$ noted at the end of Section~\ref{s.analytic} has a direct game-theoretic interpretation: 
since these quantities represent terminal rewards, time-invariant equilibria exist  
in Brownian Boost only when player incentives are exactly balanced, that is, when $\lambda=1$ in the Brownian Boost setup from~Section~\ref{s.sdg}.
\begin{proposition}\label{p.fgab}
For $\rho\in(0,\infty)$, the functions
$f_\rho(1,u)$, $g_\rho(1,u)$, $a_\rho(1,u)$ and $b_\rho(1,u)$
satisfy, as $u\to\infty$,
\[
f_\rho(1,u)=u^{\zeta_f} e^{-2u}\,\Theta(1),\qquad
g_\rho(1,u)=u^{\zeta_g} e^{-2u}\,\Theta(1),
\]
and
\[
a_\rho(1,u)=u^{\zeta_a} e^{-2u}\,\Theta(1),\qquad
b_\rho(1,u)=u^{\zeta_b} e^{-2u}\,\Theta(1),
\]
where
\[
 \zeta_f = \frac{1+\rho}{2\rho^2}, \qquad
 \zeta_g = \frac{1-\rho}{2\rho^2}, \qquad
 \zeta_a = \frac{1+\rho}{2\rho^2}-1, \qquad
 \zeta_b = \frac{1-\rho}{2\rho^2}-1 .
\]

As $u\to-\infty$, the corresponding asymptotics are
\[
f_\rho(1,u)=|u|^{\zeta_g} e^{-2|u|}\,\Theta(1),\qquad
g_\rho(1,u)=|u|^{\zeta_f} e^{-2|u|}\,\Theta(1),
\]
and
\[
a_\rho(1,u)=|u|^{\zeta_b} e^{-2|u|}\,\Theta(1),\qquad
b_\rho(1,u)=|u|^{\zeta_a} e^{-2|u|}\,\Theta(1).
\]

The $\Theta(1)$ factors are uniformly bounded away from zero and infinity
for $\rho$ in any compact subset of~$(0,\infty)$.
\end{proposition}

\begin{figure}[htbp]
\centering
\includegraphics[width=0.75\textwidth]{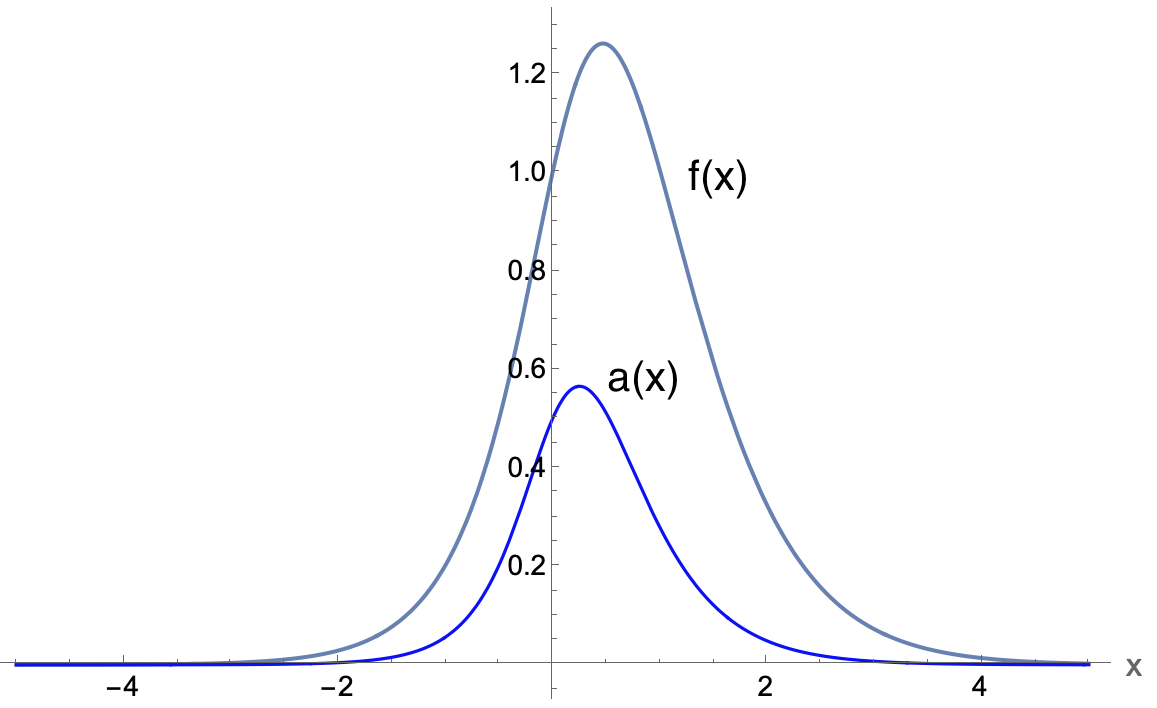}
\caption{The curves $a= a_1(1,\bullet)$ and  $f = f_1(1,\bullet)$, as specified by Definitions~\ref{d.fg} and~\ref{d.arhobrho}; $b$ and $g$ are given by reflecting in the vertical axis.
Maxine's stake profile~$a$ takes maximum value $0.57$ at $u = 0.25$ to two decimal places.}\label{f.odepair}
\end{figure}

When $\xmac=1$, players stake at unit order in a bounded neighbourhood of the origin. (As Figure~\ref{f.odepair} depicts in the case~$\rho =1$, a player spends most as she begins to lead.)
All four functions, including the stake profiles $a_\rho(1,\cdot)$ and $b_\rho(1,\cdot)$, decay rapidly, as $e^{-2\vert u \vert}$, at high distances from the origin.
There is a more modest but clear asymmetry in the rate of decay, manifest in the values of the power-law exponents.
When $u \gg 0$, the presumptively leading player, Maxine, is staking at normalized rate $u^{\zeta_a}$, above the analogous level of $u^{\zeta_b}$ for Mina; and this circumstance is swapped in the opposite regime $u \ll 0$.
As $\zeta_a > \zeta_b$, so $\zeta_f > \zeta_g$: $1 \gg m' \gg -n' > 0$ when $u \gg 0$; integrating on $[u,\infty)$, the leading player's shortfall in expectation relative to her winning terminal receipt is seen to exceed the opponent's excess over her losing terminal receipt.
This imbalance reflects the effort of expenditure that the leading player must exert---small in an absolute sense, but large relative to the opponent's---in order to convert a likely victory.

\subsection{Fixed-parameter \abmnkpmacspace asymptotics, and asymmetric decay}

Harris and Vickers~\cite{HarrisVickers87} enquire `whether the leader in a race makes greater
efforts than a follower' and `whether efforts are greatest when the competitors are neck-and-neck'.
The 2012  review~\cite{Konrad2012} of dynamic contests surveys how the discouragement effect (the subject of the next subsection)
`may cause violent conflict in early rounds, but may also lead to long periods of peaceful interaction'.

These themes are apparent in~\bbrhospace from Proposition~\ref{p.fgab}, wherein the flow index choice $\xmac=1$
corresponds to the battlefield value $\battlefield$ being the origin. Stake-profiles $\Theta(1)$ in a compact neighbourhood containing~$0$, while satisfying
\begin{equation}\label{e.abstake}
 a_{-u} \ll b_{-u} \ll a_{-u}/b_{-u} \ll 1 \, \, \, \, \textrm{for} \, \, \, \, u \gg 0  :
\end{equation}
in negative territory, where Mina leads, stakes have fallen exponentially, the more so for Maxine, though the decay in
stake ratio has a more modest polynomial rate. `Battlefield Cyl Fog' (cut your losses, foot on gas) is a mnemonic for the premise~(\ref{e.abstake}), the phrases descriptive of the trailing and leading player's respective approach far from the battlefield.

In~\cite{LostPennies}, {\rm TLP}$(1,1)$ was studied: the Trail of Lost Pennies without flip moves whose turn outcomes are decided by the simple $\frac{a}{a+b}$ `lottery' rule. The Battlefield Cyl Fog was verified (in a manner we will recall shortly).
In Theorem~\ref{t.abmn}, we present fixed-parameter asymptotics for \abmnkpmacspace elements throughout the region $(\kappa,\rho) \in (0,1]^2$ in which such elements describe Nash equilibria in \tlpkpspace according to Theorem~\ref{t.nashabmn}. This result permits us to
interrogate the validity of this premise~(\ref{e.abstake}) more broadly, via a two-dimensional family of models.

The $\phi$-sequence now defined\hfff{phi} will permit us to identify the battlefield in the definition that follows.
\begin{definition}\label{d.deltai}
Let $(a,b,m,n) \in$ \abmnkpmac. For $i \in \Z$, set
$\macphi_i = \frac{n_{i,i-1}}{m_{i-1,i}}$.
\end{definition}
This sequence is an important part of our apparatus for computing \abmnkpmacspace elements. Note that $\phi_0 = n_{0,-1}/m_{-1,0}$ equals the central ratio.  Being a measure of single-turn relative incentive for play at the origin, it is a discrete counterpart to the flow index~$\xmac$ in Brownian Boost.

Note that the quantities $\phi_0$ and $\phi_1$ now have two meanings: as sequence elements for an \abmnkpmacspace element, and as functions $\beta \to \phi_i(\kappa,\rho,\beta)$ in Definition~\ref{d.fourfunctions}. The coincidence is intentional, with the choice $\beta = n_{1,-1}/m_{-1,1}$ reconciling the objects, as we will see in solving the \abmnmacspace system in Chapter~\ref{c.abmnformfixed}. The meaning of $\phi_0$ and $\phi_1$ will, we hope, be clear from context: specified by the formulas when $\beta \in (0,\infty)$ is a free parameter, and as sequence elements when an \abmnkpmacspace element is given.

For any parameter pair $(\kappa,\rho)$ belonging to the region $W$ in~(\ref{e.weakregion}),
the orbit $\{ s_i(x): i \in \Z \}$ specified in Definition~\ref{d.stabc} will be shown to be decreasing, for any $x \in (0,\infty)$.
As we will substantiate in Chapter~\ref{c.abmnformfixed}, $s(\phi_i) = \phi_{i+1}$ for each $i \in \Z$: for any given \abmnkpmacspace solution, $s$ acts as the unit left-shift on the just specified $\phi$-sequence. Lemma~\ref{l.battlefield} will show that the $s$-orbit from any positive real passes exactly once through the central domain as it is next defined. As such, this lemma furnishes the existence and uniqueness claims on which the next definition depends.

\begin{definition}\label{d.battlefield}
For $(\kappa,\rho) \in W$,
let $(a,b,m,n) \in$ \abmnkpmac.
The {\em central domain}\hfff{domain} $D$ is $\big( \frac{2 - \kappa \rho}{2+\kappa \rho}, \frac{2 + \kappa \rho}{2-\kappa \rho} \big]$.
The {\em battlefield index}\hfff{battlefield} is set equal to $k \in \Z$ such that~$\phi_k \in D$.
\end{definition}
This definition extends the (1,1)-case in~\cite{LostPennies}, where $D = (1/3,3]$.

Here is our result offering asymptotics for each component in $(a,b,m,n) \in$ \abmnkpmacspace in terms of distance of the index from the battlefield. The result extends the case $(\kappa,\rho) = (1,1)$
treated in \cite[Theorem~$2.14$]{LostPennies}. There are three regimes in $(0,1]^2 \setminus \{ (1,1 \})$: the interior, and the upper and right sides.
\begin{theorem}\label{t.abmn}
Let $(a,b,m,n)$ be an element of \abmnkpmacspace with battlefield index zero.
\begin{enumerate}
\item Suppose that $\kappa \in (0,1)$ and $\rho \in (0,1)$. Then, for $i >  0$,
\begin{eqnarray*}
   m_{-i-1,-i}  & = & m_{-1,0}  \cdot \sigma  \cdot   i^{\frac{1-\rho}{2\rho^2}}
 \left( \frac{1 - \kappa}{1 + \kappa} \right)^{i} \left( 1 + O\left(i^{-1}\right) \right) \\
 a _{-i} & = & m_{-1,0}  \cdot  \sigma  \,
  \tfrac{1 + \kappa}{4\rho} \cdot  i^{\frac{1-\rho}{2\rho^2} - 1}
  \left( \frac{1-\kappa}{1+\kappa} \right)^{i}  \left( 1 + O\left(i^{-1}\right) \right) \, .
 \end{eqnarray*}
 The ratios  $n_{-i,-i-1}/m_{-i-1,-i}$ and $b_{-i}/a_{-i}$  take the form  $\left( \frac{8 \rho^2 \kappa}{1 - \kappa^2} \right)^{1/\rho} \, i^{1/\rho}  \left( 1 + O\left(i^{-1}\right) \right)$.

 Here (and in the following part), \( \sigma = \sigma(\phi_0; \rho, \kappa) \) is a unit-order constant depending on \( \phi_0 \), remaining bounded above and below as \( \phi_0 \) ranges over~$D$, uniformly for $\rho$ and $\kappa$ valued in compact subsets of $(0,1]$.
\item  For $\kappa \in (0,1)$, $\rho =1$ and $i > 0$,
\begin{eqnarray*}
m_{-i-1,-i}   & = & m_{-1,0}    \cdot \sigma \cdot
       \left( \frac{1-\kappa}{1+\kappa} \right)^{i}  \Big( 1 + O(i^{-1}) \Big) \\
a_{-i} & = &  m_{-1,0} \cdot \sigma  \cdot \frac{1+\kappa}{4}  \cdot  i^{-1}
       \left( \frac{1-\kappa}{1+\kappa} \right)^{i}
        \Bigl( 1 + O\big(  i^{-1} \log i \big) \Bigr) \,  ,
        \end{eqnarray*}
And   $n_{-i,-i-1}/m_{-i-1,-i}$ and $b_{-i}/a_{-i}$ equal
$\frac{8\kappa}{1 - \kappa^2}i +  O\left(\log i \right)$.
       \item
Now suppose that $\kappa=1$ and $\rho \in (0,1)$. For $i > 0$,
the quantities
$m_{-i-1,-i}$ and $a_{-i}$ take the form
$$
m_{-1,0}\left(\frac{1-\rho}{1+\rho}\right)^{\rho i^2/2} e^{\chi i}  \cdot e^{o(i)} \, ,
$$
 and the ratios  $n_{-i,-i-1}/m_{-i-1,-i}$ and $b_{-i}/a_{-i}$ equal
$\left( \tfrac{1+\rho}{1-\rho} \right)^i O(1)$.
The constant $\chi = \chi(\phi_0,\rho)$ is bounded away from zero and infinity for $\phi_0$ of battlefield index zero provided that $\rho$ lies in a compact subset of $(0,1)$.
\item For all the  statements above, the components of $(a_i,b_i,m_i,n_i)$ for $i \geq 0$ satisfy the same asymptotics  as the respective elements of $(b_{-i},a_{-i},n_{-i},m_{-i})$.
\item Suppose now that~$(a,b,m,n)$ has battlefield index $k$. Then all statements remain valid after $i$ is replaced by $i-k$ in the conditions $i > 0$ and $i < 0$ and in every right-hand side, and $\phi_0$ is replaced by~$\phi_k$.
 \end{enumerate}
\end{theorem}
Fixed-parameter asymptotics in the region $(\kappa,\rho) \in W$ above $\rho =1$ may also be obtained, but this regime has been omitted since it lies outside the purview of Theorem~\ref{t.nashabmn}, leaving unclear its relevance to the trail game.

Consider battlefield zero and negative territory.
When $\kappa < 1$, in the first two parts of the theorem, the dominant decay (of~$b_{-i}$ say) is exponential, with $a_{-i}/b_{-i}$ decaying as $i^{-1/\rho}$. The exponential decay, with factor $\frac{1-\kappa}{1+\kappa}$, becomes rapid in the low-noise $\kappa \nearrow 1$ limit. Along the right boundary $\kappa =1$, $b_{-i}$ has more rapid $e^{-\Theta_\rho(1) i^2}$ decay, with $\Theta_\rho(1)$ exploding as the point~$(1,1)$ is approached from below; the ratio $a_{-i}/b_{-i}$ has exponential decay.

These results suggest that the point~$(1,1)$ may have singular behaviour, with the most rapid decay. This is borne out by~\cite[Theorem~$2.14$]{LostPennies}: $b_{-i}$ has doubly exponential leading-order decay, of the form $\exp \big\{ - 2 \cdot 2^i A \}$ for some $A > 0$, while $a_{-i}/b_{-i}$ also decays doubly exponentially, having the form~$\exp \big\{ - 2^i A \}$ to leading order.
So the premise~$a_{-i} \ll b_{-i} \ll 1$ and $b_{-i} \ll a_{-i}/b_{-i}$
that we presented in~(\ref{e.abstake}) via the left boundary $\kappa = 0^+$ Brownian Boost case is supported in all four regimes of $(\kappa,\rho) \in (0,1]^2$.

\subsection{Incentive Inch, Outcome Mile}\label{s.iiom}

We may set
\begin{equation}\label{e.lambdamaxkappazero}
\lambdamax(0,\rho) \, = \, \sup \, \left\{ \, \frac{\int_\R g_\rho(\xmac,r) \dd r}{\int_\R f_\rho(\xmac,r) \dd r} : \xmac \in (0,\infty) \, \right\}
\end{equation}
to specify a Brownian Boost counterpart to $\lambdamax(\kappa,\rho)$ from Definition~\ref{d.lambdamax}.
Indeed, by Theorem~\ref{t.fg}, the supremum is over all default solutions of the ODE pair, so that $\lambdamax(0,\rho)$ measures the maximum ratio of prize for Mina relative to Maxine compatible with equilibrium existence.


We have already noted a result that is formally stated as Proposition~\ref{p.sfacts}(5): the $f$- and $g$-integrals are always equal, so $\lambdamax(0,\rho) = 1$. This holds for any $\rho \in (0,\infty)$, though the interpretation via Nash equilibria is known for \bbrhospace only when $\rho \in (0,1]$, via \tlpkpspace and Theorem~\ref{t.nashabmn}. 
This is the conclusion of the second preface premise, validated:
 the non-existence of equilibria in the imbalanced $\lambda \neq 1$ game. (That said, the argument in the preface suggests that the more incentivized player will win \bbrhospace at almost no running cost. But the absence of equilibria in the game gives neither this player nor her opponent any guidance as to how to play it.)

The heuristic may be applied to the Trail of Lost Pennies, where it predicts $\lambdamaxkp = 1$ for any $(\kappa,\rho) \in (0,1]^2$. The premise was examined for {\rm TLP}($1,1$) in~\cite{LostPennies}, which concluded, rigorously and by numerical evidence for the respective bounds\footnote{As will be reported in a forthcoming article, U.C. Berkeley undergraduates Neo Lee and Adam Ousterovitch have obtained a computer-assisted proof of the upper bound.},
\begin{equation}\label{e.lambdamaxbounds}
1.000096 \leq    \lambdamax(1,1) \leq 1.000098 \, .
\end{equation}
So while the heuristic when literally interpreted is false, equilibria are fragile under asymmetric perturbation of incentive, with a ratio of relative incentive of order $10^{-4}$ being enough to disrupt their existence, the sense of which the phrase `Incentive Inch, Outcome Mile' seeks to capture.

Investigating the function $\lambdamax:(0,1]^2 \to [1,\infty)$ offers a way of testing the strength and robustness of the discouragement effect. We will prove the next result, which quantifies the conclusion that $\lambdamax(0,\rho) = 1$ by bounding above the rate of convergence of $\lambdamax(\kappa,\rho)$ to one as $\kappa \searrow 0$.
\begin{theorem}\label{t.lowlambdamax}
There exist $C > 0$ and $\eta,\kappa_0 \in (0,1)$ such that, for $\kappa \in (0,\kappa_0)$ and $\rho \in (0,1]$, $\big\vert \lambdamax(\kappa,\rho) - 1  \big\vert \leq C \kappa^\eta$.
\end{theorem}
{\em Remark.} The result may be extended to the regime $\rho > 1$ when $(\kappa,\rho) \in W$ (that is, $\kappa^2 \rho \leq 1$), with $\eta = \eta(\rho)$ decaying to zero in the high-$\rho$ limit.

In the \hyperref[c.directions]{Epilogue}, we report on $\lambdamax:(0,1]^2 \to [1,\infty)$ numerically, finding this function to have several remarkable features. The numerics also prompt the following conjecture.
\begin{conjecture}\label{c.lambdamax}
\leavevmode
\begin{enumerate}
\item 
The maximum value of $\lambdamax:(0,1]^2 \to [1,\infty)$ is at most $1 + 10^{-4}$.
\item
The maximum value of $(0,1] \to [1,\infty): (1,\rho) \mapsto \lambdamax(1,\rho)$ is attained at $\rho =1$.
\end{enumerate}
\end{conjecture}
This conjecture replaces a counterpart in an earlier released version, which claimed that
 the maximum of the map 
 $\lambdamax:(0,1]^2 \to [1,\infty)$ is attained at the point $(1,1)$. As we will discuss in the Epilogue, numerics do not in fact support this claim, though the maximizer does appear to be close to~$(1,1)$.
 
In {\rm TLP}($1,1$), an imbalance of incentive of order~$10^{-4}$ is enough to prevent equilibria from existing. Tiny as this amount is, the counterpart imbalance in any of the games~\tlpkpspace for $(\kappa,\rho) \in (0,1]^2 \setminus  \{ (1,1) \}$   appears to be greater by at most only by an even smaller amount, with this imbalance never exceeding~$10^{-4}$. The above conjecture reflects an unexpected aspect of the discouragement effect and asymmetric stake decay.
The trailing player is discouraged, cuts her losses, and thereby contributes to stake-decay asymmetry. As reviewed in the preceding subsection, $(\kappa,\rho) = (1,1)$ is the point where this asymmetry is greatest. Paradoxically, our conjecture implies that this is also a site of comparatively low discouragement---the weakest possible if only $\rho$ varies; and close to the weakest when both parameters do---since equilibria indexed by~$(1,1)$ would exist with relative incentives nearing the most asymmetric in the broader parameter space.

\subsection{Structure of the article}

The article 
is divided into 
two Parts.

Part~\ref{p.one} is devoted to the discrete Trail of Lost Pennies games.
It develops an analytic framework for studying the \abmnkpmacspace system, and classifies the games' time-invariant Nash equilibria.

Part~\ref{p.two} returns to Brownian Boost.
It presents a heuristic derivation and analytic study of the \bbrhospace ODE pair governing equilibrium behaviour, and then realises Brownian Boost as a fine-mesh, high-noise limit of the discrete games.

The \hyperref[c.directions]{Epilogue} reports on numerical investigations, outlines several directions for further inquiry and 
is followed by
 a~\hyperref[c.glossary]{glossary of notation}.

\subsection{Acknowledgements}
The author warmly thanks a referee of~\cite{LostPennies} whose pivotal role has been indicated in Section~\ref{s.tugofwar}. 
He thanks several contributors in peer review whose helpful input has informed the structure and literature review of the article.
He further thanks Daniil Kardava for numerical investigations including the outputs shown in Figures~\ref{f.doublelogheatmap} and~\ref{f.blackandwhite}.
The 
Berkeley Research Computing Cluster provided resources to compute these numerics.
The author is supported by NSF-DMS grants~$2153359$ and~$2450726$.

\part{Discrete Tug-of-War Games\label{p.one}}

\chapter{Some basic symmetries and tools}\label{c.symmetrytools}

The discrete Trail of Lost Pennies games are the object of study in this Part.
Their analysis is organised around the $\Z$-indexed \abmnkpmacspace system from Definition~\ref{d.abmn}, which encodes equilibrium stakes and values.
Central to our approach are structural features of this system, including the  map
$s:(0,\infty)\to(0,\infty)$ sending $\phi_0$ to $\phi_1$ introduced in Definition~\ref{d.scd}.

The aim of this three-chapter Part is to develop a rigorous analytic framework for the discrete games and to relate their equilibria to solutions of \abmnkpmac.
In the present chapter, we assemble the basic tools needed for this programme.
In Chapter~\ref{c.abmnformfixed}, we show that $s$ acts as a shift on the $\phi$-sequence associated to any  \abmnkpmacspace element by Definition~\ref{d.deltai}; we prove the explicit-form Theorem~\ref{t.defexplicit} for \abmnkpmac; and, by developing $s$-orbit asymptotics, we establish the fixed-parameter \abmnkpmacspace Theorem~\ref{t.abmn}.
In Chapter~\ref{c.nashabmn}, we complete the picture by proving the Nash--\abmnmacspace equivalence Theorem~\ref{t.nashabmn}.

In this chapter, we focus on the foundational symmetries and constructions.
In consecutive sections, we prove basic properties of $\phi_0$, $\phi_1$ and the resulting map $s(\phi_0) = \phi_1$ and describe the escape behaviour of forward and backward $s$-orbits;
introduce a role-reversal symmetry that inverts $s$;  define the battlefield index of an element of \abmnkpmacspace via its $s$-orbit; and introduce the Mina margin map, a key device for analysing equilibrium structure.
In a final section, we analyse Penny Forfeit, the one-step sub-game of \tlpkp.

\section{The bijections $\phi_0$ and $\phi_1$, and the orbit of~$s$}

The functions $\phi_0$ and $\phi_1$ 
are algebraically specified, in Definition~\ref{d.fourfunctions}.  Our first 
result records  some of their basic properties,  
which have permitted the specification of the map $s:(0,\infty) \to (0,\infty)$ that sends $\phi_0$ to $\phi_1$ in Definition~\ref{d.scd}.
\begin{lemma}\label{l.incphi}
Suppose that $(\kappa,\rhomac) \in (0,1] \times (0,\infty)$.  
\begin{enumerate}
 \item Each of $\phi_0$ and $\phi_1$ satisfies
$$
 \lim_{\beta \searrow 0} \phi(\kappa,\rhomac,\beta) = 0 \, \, \, \, \textrm{and} \, \, \, \, 
 \lim_{\beta \nearrow \infty} \phi(\kappa,\rhomac,\beta) = \infty  \, ,
$$
where in the case $\kappa =1$, we also suppose that $\rho \leq 1$. 
 \item  If $\rhomac^2 \kappa \leq 1$, then  $(0,\infty) \to (0,\infty):\beta \to \phi_i(\kappa,\rhomac,\beta)$ is an increasing bijection  for $i \in \{0,1\}$. 
 \item  If  $\kappa \rho < 1 + \sqrt{1-\kappa^2}$,  then $\phi_0(\kappa,\rhomac,\beta) >  \phi_1(\kappa,\rhomac,\beta)$. In particular, this holds when $\rho^2 \kappa \leq 1$.
\end{enumerate}
\end{lemma}

\begin{figure}[htbp]
\centering
\includegraphics[width=0.25\textwidth]{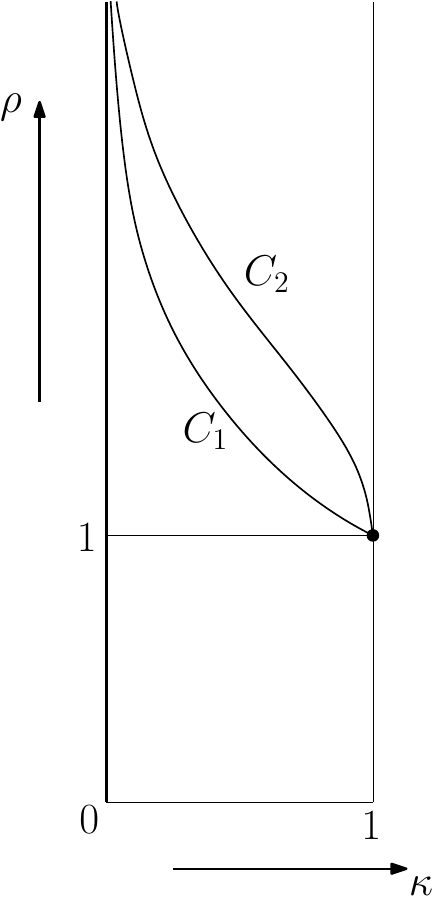}
\caption{The main contours on the $(\kappa,\rho)$-map. The curve $C_1$ is the locus of $\kappa^2 \rho =1$, which is the upper boundary of the region~$W$ in~(\ref{e.weakregion}); $C_2$ is the locus of $\kappa \rho = 1 + \sqrt{1-\kappa^2}$. The curves emanate from $(1,1)$, the point indexing the game studied in~\cite{LostPennies} and the upper-right corner of the unit box in which Nash-\abmnmacspace equivalence is established by Theorem~\ref{t.nashabmn}. The \bbrho-line $\kappa =0$ lies below $C_1$, which is indicative of how the ODE-pair Theorem~\ref{t.fg} is valid for all~$\rho \in (0,\infty)$.}\label{f.kapparho}
\end{figure}

Figure~\ref{f.kapparho} shows how  the contours specified by the conditions in the lemma lie in 
 the $(\kappa,\rho)$-strip.

{\bf Proof of Lemma~\ref{l.incphi}(1).} 
 When $\kappa \in (0,1)$,
 then $\phi_0(\kappa,\rho,\beta)$ and   $\phi_1(\kappa,\rho,\beta)$ are asymptotic to $\beta$, for $\beta$ both high and low, whatever the value of $\rho \in (0,\infty)$. 
When $\kappa =1$, we suppose $\rho \in (0,1]$. Hence, the asymptotics
$$
\phi_0(1,\rhomac,\beta)  \, \stackrel{\beta \nearrow \infty}{\sim} \, \begin{cases} \, \tfrac{1+\rho}{1-\rho}\beta \, \, \, \, & \textrm{if $\rho \in (0,1)$}  \\ 
  \, 2\beta^2 \, \, & \textrm{if $\rho = 1$} \end{cases}  \, \, \, \, \, \textrm{and} \, \, \, \, \, 
\phi_1(1,\rhomac,\beta)  \, \stackrel{\beta \searrow 0}{\sim} \, \begin{cases} \, \tfrac{1-\rho}{1+\rho}\beta \, \, \, \, & \textrm{if $\rho \in (0,1)$} \\ 
  \, \beta^2/2 \, \, & \textrm{if $\rho = 1$} \end{cases}  \, ,
$$
suffice to treat the remaining cases.

{\bf (2).} In view of the preceding part, it is enough to argue that each function is increasing.

It is useful to permit negative $\kappa$ and apply the 
 symmetry $\phi_1(\kappa,\rho,\beta) = \phi_0(-\kappa,\rho,\beta)$. Indeed it then suffices to show that $(0,\infty) \to (0,\infty):\beta \to \phi_0(\kappa,\rhomac,\beta)$ is increasing 
whenever non-zero $\kappa \in [-1,1]$ satisfies  $\rhomac^2 \vert \kappa \vert \leq 1$. 

Writing $q$ for the right-hand denominator in~(\ref{e.phizero.rho}), $ \frac{\partial \, \phi_0(\kappa,\rhomac,\beta)}{\partial \beta} \, = \, P(\beta^\rhomac) q^{-2}$
where of the coefficients of the quartic $P(x) = \sum_{i=0}^4 h_i x^i$, 	$h_0 = (1+\kappa)^2$ and $h_4 = (1-\kappa)^2$ are evidently non-negative.
That $h_2 = 2\big(3 - \kappa^2(1+2\rhomac^2)\big)$ is likewise follows from the conditions, which are weaker than our hypothesis,  that $\vert \kappa \vert \leq 1$ and $\rhomac \vert \kappa \vert \leq 1$.
The coefficients $h_1 = 4(1-\kappa)(1- \kappa \rhomac^2)$ 
 and $h_3 = 4(\kappa+1)(\rhomac^2 \kappa +1)$
 are also non-negative: in one case trivially; in the other,  
 as our hypothesis is tailored to show; and with the sign of $\kappa$ determining which case applies.

{\bf (3).} First consider $\kappa \in (0,1)$. In this case, $\rho^2 \kappa \leq 1$ is evidently a stronger hypothesis, so we suppose $\kappa \rho < 1 + \sqrt{1-\kappa^2}$. Note that
\[
\phi_0(\kappa,\rho,\beta) - \phi_1(\kappa,\rho,\beta) \, = \,  
\frac{8 \kappa \rho\, \beta^{1 + \rho} (1 + \beta^\rho)^2}
{g(\kappa,\beta^\rho) \cdot g(-\kappa,\beta^\rho)},
\]
where $g(\kappa,x) =  (1 - \kappa)x^2 + 2(1 - \rho\kappa)x + 1 + \kappa$.
When $\kappa \in (0,1)$ and $\rho > 0$, the quadratic $g(-\kappa,\bullet):\R \to \R$ is always positive, while for such $\kappa$-values,  $g(\kappa,\bullet):\R \to \R$ is as well, because the discriminant sign condition $\rho < \kappa^{-1}(1 + \sqrt{1-\kappa^2})$ has been hypothesised. Thus, $\phi_0$ exceeds $\phi_1$ when $\kappa \in (0,1)$.

For $\kappa = 1$, we have $\rho \in (0,1]$. (We include $\rho =1$ because it meets the condition $\rho^2 \kappa \leq 1$.)
Then $g > 0$ is readily checked, so $\phi_0 > \phi_1$ in this case also. \qed

\begin{definition}\label{d.subdiagonal}
A map $f:(0,\infty) \to (0,\infty)$ is {\em sub-diagonal} if $f(x) < x$. 
\end{definition}
By Lemma~\ref{l.incphi}(3), $s$ meets this definition. 
Write $\incbij$ for the space of increasing bijections of $(0,\infty)$ and note that any element of $\incbij$ is continuous. 
Since $\phi_0$ and $\phi_1$ belong to $\incbij$ by Lemma~\ref{l.incphi}(1,2), and $s(\phi_0) = \phi_1$,
we see that $s \in \incbij$.  Hence the next result implies the following corollary.
\begin{lemma}\label{l.siterate}
Let $f:(0,\infty) \to (0,\infty)$ be a continuous sub-diagonal bijection. 
For $x \in (0,\infty)$, set $x_0 = x$ and iteratively define the forward and backward orbits $x_i = f(x_{i-1})$ and $x_{-i} = f^{-1}(x_{1-i})$ for $i \in \nwozero$.
Then $\big\{ ( x_i,x_{i-1} ]: i \in \Z \big\}$ is (in decreasing order) a partition of $(0,\infty)$.
\end{lemma}
\begin{corollary}\label{c.orbitescape}
For any $x \in (0,\infty)$, $s_{-i}(x) \to \infty$ and $s_i(x) \to 0$ as $i \to \infty$.
\end{corollary}
{\bf Proof of Lemma~\ref{l.siterate}.} The orbit sequence $\big\{ x_i: i \in \Z \big\}$ is decreasing because $f$ is sub-diagonal. If its left limit~$x_\infty$ were positive, this limit point would lie in the domain of the continuous map~$f$, so that, absurdly, $x_i$ would converge in high~$i$ both to~$x_\infty$ and to the smaller value $f(x_\infty)$. Hence, $x_\infty = 0$. With a similar notation and argument, $x_{-\infty} = \infty$. Thus every positive real lies in  $( x_i,x_{i-1} ]$ for precisely one integer~$i$. \qed

\section{Role-reversal symmetry and the inverse of~$s$}\label{s.roleshift}

Recall the region~$W$ of $(\kappa,\rho)$-parameter space from~(\ref{e.weakregion}).
 \begin{proposition}\label{p.sminusone}
 Let $(\kappa,\rho) \in W$. The function $s:(0,\infty) \to (0,\infty)$ from Definition~\ref{d.scd}(1) is invertible, with $s^{-1}(x) \, = \,  1/s(1/x)$ for ~$x \in (0,\infty)$. 
 \end{proposition}
 {\bf Proof.} By definition, $s$ sends $\phi_0$ to $\phi_1$. Since both maps $\phi_i(\beta) = \phi_i(\kappa,\rho,\beta)$ are bijections $(0,\infty) \to (0,\infty)$ by Lemma~\ref{l.incphi}(2), the inverse map $s^{-1}$
 sending $\phi_1$ to $\phi_0$ is well defined. The formula $ s^{-1}(x) \, = \,  1/s(1/x)$ amounts to $1/\phi_0 = s(1/\phi_1)$. To see this, note from~(\ref{e.phizero.rho}) and~(\ref{e.phione.rho}) that 
 when $\phi_0 = \phi_0(\kappa,\rho,\beta)$ and  $\phi_1 = \phi_1(\kappa,\rho,\beta)$, we have that $1/\phi_0 =  \phi_1(\kappa,\rho,1/\beta)$ and  $1/\phi_1 =  \phi_0(\kappa,\rho,1/\beta)$. 
So the sought equality  $1/\phi_0 = s(1/\phi_1)$ is then the instance of $s(\phi_0) = \phi_1$ corresponding to $1/\beta$. \qed

The preceding proof establishes the claimed identity by a direct algebraic calculation.
It is conceptually helpful, however, to complement this argument with a game-theoretic point of view,
which also motivates the term \emph{role reversal} for the symmetry just identified.
Some apparatus developed later in Part~\ref{p.one} is needed for this interpretation.

Two ingredients are required.
First, the Nash--ABMN equivalence Theorem~\ref{t.nashabmn} furnishes a game-theoretic interpretation
of any given $(a,b,m,n)\in$ \abmnkpmac.
Second, while the proof above works with the formulaic definitions of $\phi_0$ and $\phi_1$,
we will wish to invoke their interpretation in terms of sequence elements.
In Lemma~\ref{l.gdphiphi}, we will establish the coincidence of these two descriptions,
as discussed shortly after Definition~\ref{d.deltai}:
for $(a,b,m,n)\in$ \abmnkpmac, writing $\beta=\beta_0:=n_{1,-1}/m_{-1,1}$ yields
$\phi_0=n_{0,-1}/m_{-1,0}$ and $\phi_1=n_{1,0}/m_{0,1}$.

With these identifications in place, the symmetry underlying Proposition~\ref{p.sminusone}
admits the following interpretation.
Reflect gameplay governed by the strategy pair $(a,b)$\hfff{notab} through the origin.
The players now occupy opposite ends of the gameboard, and under the reflected dynamics
each would be acting against her own interest. But their play makes sense if they now change ends. 
The resulting gameplay is governed by the strategy pair $(b(-\bullet),a(-\bullet))$,
which extends naturally to
$$
\big(b(-\bullet),a(-\bullet),n(-\bullet),m(-\bullet)\big)\in \text{\abmnkpmac} \, .
$$
Under this transformation, the parameters satisfy
$\phi_0\mapsto 1/\phi_1$, $\phi_1\mapsto 1/\phi_0$, and $\beta_0\mapsto 1/\beta_0$,
which explains both the formula $s^{-1}(x)=1/s(1/x)$ and the reciprocal $\beta$-parametrization
that appears in the proof above.

The solution class \abmnkpmacspace is invariant under $\Z$-shifts of the index,
and the role-reversal symmetry identified here provides a second invariance of the system.

We keep a record of another consequence which we have noted, along with an extension.
\begin{corollary}\label{c.rolereversal}
If $(a,b,m,n)$ is an element of \abmnkpmac, then so is
$$
\big( (b_{k-i}, a_{k-i}, n_{k-i}, m_{k-i}): i  \in \Z \big) \, , \, \, \, \, \textrm{for any $k \in \Z$} \, . 
$$ 
\end{corollary}
We have noted this result for $k=0$, and apply the $\Z$-shift symmetry to obtain the other choices.
Alternatively, note that $\Z$ is reflection-symmetric not only about integers but also about half-integers: for example, we may reflect gameplay about minus one-half instead of zero to obtain the solution with~$k=-1$.

For the reader who prefers a direct algebraic check, the corollary can readily be confirmed by examining the \abmnmacspace equations in Definition~\ref{d.abmn}.

 \section{The battlefield index}
 
 Next we clarify that the battlefield index as specified in Definition~\ref{d.battlefield} is well-defined. Recall that $D = \big( \frac{2 - \kappa \rho}{2+\kappa \rho}, \frac{2 + \kappa \rho}{2-\kappa \rho} \big]$ is the central domain.
 \begin{lemma}\label{l.battlefield}
For $(\kappa,\rho) \in (0,1]^2$,
 let $(a,b,m,n)$ be an ABMN($\kappa,\rho$) solution. 
 \begin{enumerate}
\item  There is a unique value of $x \in (0,\infty)$ such that $x \, s(\kappa,\rho,x)= 1$.
\item This value is given by $x = \frac{2 + \kappa \rho}{2-\kappa \rho}$, with $s(x) =  \frac{2 - \kappa \rho}{2+\kappa \rho}$.
\item We have that  $1 \in D \subset \big\{ x \in \R: \vert x - 1 \vert \leq 2 \kappa \rho \big\}$.
\item There is a unique value $k \in \Z$ for which  $\phi_k \in D$.
 \end{enumerate} 
 \end{lemma}
 {\bf Proof: (1).} The function $s = s(\kappa,\rho,\bullet):(0,\infty) \to (0,\infty)$ belongs to $\incbij$, so it meets the decreasing map $1/x$ at exactly one $x \in (0,\infty)$.\\
 {\bf (2).} Evaluating $\phi_0$ and $\phi_1$ at $\beta=1$ gives $\phi_0(\kappa,\rho,1) =  \frac{2 + \kappa \rho}{2-\kappa \rho}$ and  $\phi_1(\kappa,\rho,1) =  \frac{2 - \kappa \rho}{2+\kappa \rho}$.
 Since $s$ maps $\phi_0$ to $\phi_1$, we identify $x = \frac{2 + \kappa \rho}{2-\kappa \rho}$ as the unique solution of $x s(x) =1$.\\ 
 {\bf (3).} Since $D = (s(x),x]$ with $x s(x) =1$ and $s$ is subdiagonal, $1 \in D$. The endpoints of $D$ lie at distances from one of  $\frac{2 \kappa \rho}{2-\kappa \rho}$ and  $\frac{2 \kappa \rho}{2+\kappa \rho}$, the former expression the larger and bounded above by $2\kappa \rho$ since $(\kappa,\rho) \in (0,1]^2$. \\ 
 {\bf (4).}  Let $p \in (0,\infty)$.
By Lemma~\ref{l.siterate}, the intervals $\big( s_{i+1}(p) ,  s_i(p) \big]$, indexed by~$i \in \Z$, partition~$(0,\infty)$.
The orbit $\phi_i$ visits each interval in the partition exactly once, doing so in decreasing order of index. Taking $p=   \frac{2 + \kappa \rho}{2-\kappa \rho}$ yields what is claimed. \qed

\section{The Mina margin map}

In Definition~\ref{d.lambdamax}, $\lambdamax(\kappa,\rho)$ has been defined to be the supremum of the Mina margin $n_{\infty,-\infty}/m_{-\infty,\infty}$
over all \abmnkpmacspace solutions. It is worth noting that the several symmetries enjoyed by \abmnkpmacspace permit a more restricted supremum to be taken, and
the Mina margin map is a useful device for making this point.
Recall from Section~\ref{s.solvingabmn} that for $x \in (0,\infty)$ there is a unique element of \abmnkpmacspace  with $m_{-\infty}= n_\infty = 0$, $m_\infty =1$ and $\frac{n_{0,-1}}{m_{-1,0}} = x$.
This is the standard solution $\big(\asta_i(x),\bst_i(x),\mst_i(x),\nst_i(x): i \in \Z \big)$.

\begin{figure}[htbp]
\centering
\includegraphics[width=0.8\textwidth]{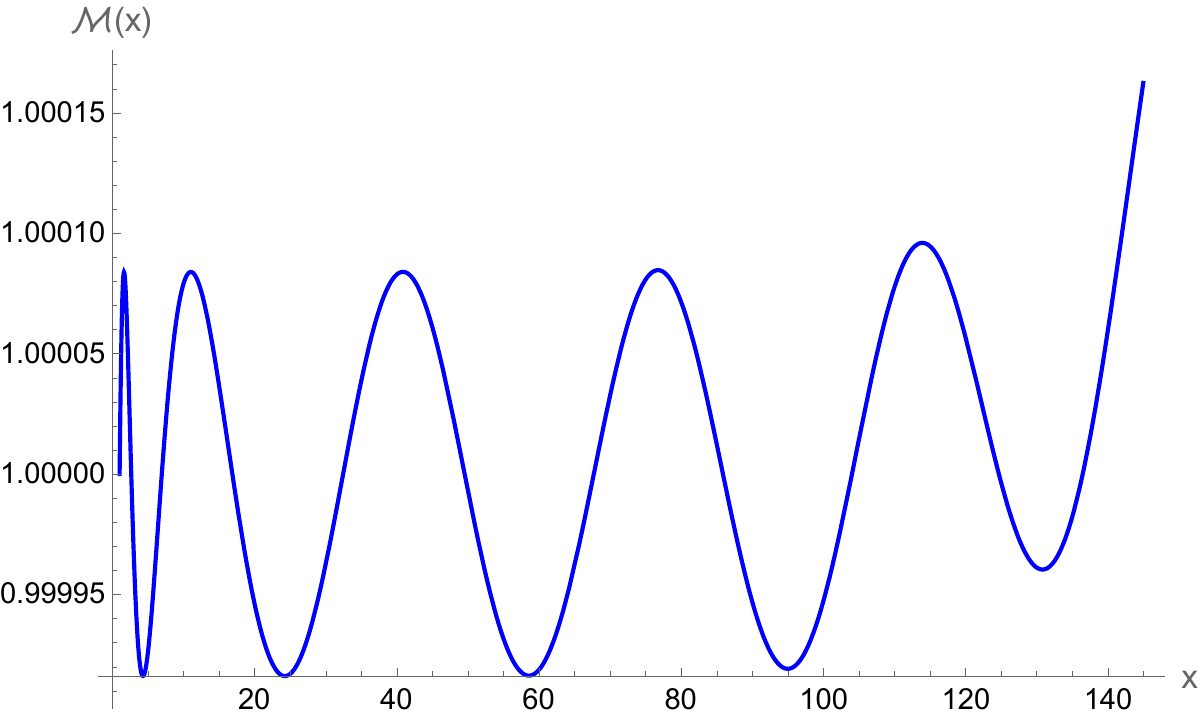}
\caption{The finite-trail Mina margin map $x \to \mc{M}_{0.9,1}^{-9,9}(x)$ is plotted on $(1,145)$. The map rises rapidly to the right of the plotted range, and its values on $(0,1)$ are determined by the symmetry $\mc{M}(x) = \mc{M}(x^{-1})^{-1}$. There are twenty-one roots of $\mc{M}(x) = 1$, given by $x=1$ and ten pairs $(z,z^{-1})$.
}\label{f.mmm}
\end{figure}

\begin{definition}\label{d.mmm}
Let the Mina margin map $\minammkp:(0,\infty) \to (0,\infty)$\hfff{mmm} be given by 
$$
\textrm{$\minammkp(x) =\nst_{-\infty}(\kappa,\rho,x)$, \,      $x \in (0,\infty)$} \, .
$$
 Namely, $\minammkp(x)$ is 
the Mina margin of $\quadstand$. 
\end{definition}
\begin{proposition}\label{p.relativereward}
\leavevmode
\begin{enumerate}
\item The function $\minammkp:(0,\infty) \to (0,\infty)$ satisfies $\minammkp(s(x)) = \minammkp(x)$ for $x \in (0,\infty)$.
\item The map $x \to \minammkp(x)$ is continuous and is given by
$$
 \minammkp(x)   \, \, = \, \, \Bigg(  \sum_{k \in \Z} \, \, \prod_{i=0}^k \big( c_i(x) - 1 \big)  \Bigg)^{-1} \, \cdot \,  \minanum \, .
$$
\item For $x \in (0,\infty)$, $\minammkp(x^{-1}) = \minammkp(x)^{-1}$. In particular, $\minammkp(1) =1$.
\item We have $\minammkp(0,\infty) = \minammkp(D)  = \big[\lambdamaxkp^{-1},\lambdamaxkp\big]$.
\end{enumerate}
\end{proposition}
For the proof, we define  a Mina margin map associated to the finite trail $\llbracket -k , k \rrbracket$
by setting $\minammkp^{-k,k}(x) = \frac{\nst_{k,-k}(x)}{\mst_{-k,k}(x)}$\hfff{mmm.finite}: see Figure~\ref{f.mmm} for a depiction.
 
{\bf Proof  of Proposition~\ref{p.relativereward}(1).} 
Taking the high $k$ limit, $\nst_{k,-k} \to \nst_{\infty,-\infty}$ and $\mst_{-k,k} \to \mst_{-\infty,\infty} = 1$, so that
\begin{equation}\label{e.rkkconvergence}
 \minammkp^{-k,k}(x)  \, \longrightarrow \, 
 \minammkp(x) \, ,
\end{equation}
the limit in~$\R$ by  Theorem~\ref{t.abmnpositive}(2).
Since replacing $x \to s(x)$ in  $\big(\asta(x),\bst(x),\mst(x),\nst(x) \big)$ results in a left shift by one place, 
$$
\minammkp^{-k,k}\big( s(x) \big) \, = \, \frac{\nst_{1+k,1-k}(x)}{\mst_{1-k,1+k}(x)} \, .
$$
As $k \to \infty$, the left-hand side converges to $\minammkp\big( s(x) \big)$, by~(\ref{e.rkkconvergence}) with $x \to s(x)$, while  the right-hand side converges to  $\frac{\nst_{\infty,-\infty}}{\mst_{-\infty,\infty}} = \minammkp(x)$ by~(\ref{e.rkkconvergence}) and  the decay of high-indexed $m$- and $n$-differences in Theorem~\ref{t.abmn}.  Hence $\minammkp(s(x)) = \minammkp(x)$ for $x \in (0,\infty)$.

{\bf (2).} There is a  Mina margin map associated to each finite trail. For the gameboard~$\llbracket -k,k \rrbracket$,  it takes the form
$$
\minammkp^{-k,k}(x) = \frac{n_{k,-k}}{m_{-k,k}}
$$  
for any element $(a,b,m,n) \in$ \abmnkpmacspace with $\phi_0 = x$.
The decay in high $\vert i \vert$ for $m_{i,i+1}$ and $n_{i+1,i}$ is shown in Theorem~\ref{t.abmn} to be at least as rapid as exponential, uniformly in choices $x \in D$ that correspond to battlefield index zero.  Thus 
$\minammkp^{-k,k}$ converges uniformly to $\minammkp$ on $D$. Writing 
$\minammkp^{-k,k}(x) = \frac{\ndefault_{k,-k}}{\mdefault_{-k,k}}$ as a ratio of default values, we may sum the explicit product formulas from Definition~\ref{d.zdef} to show that this prelimit function is continuous on~$D$; whence, so is $\minammkp$.
For~$x \in D$, the claimed formula for  $\minammkp(x)$ emerges by taking the high~$k$ limit of this ratio of explicit expressions.

 The map $s:(0,\infty) \to (0,\infty)$ is invertible and may be iterated, forwards and backwards, so that Proposition~\ref{p.relativereward}(1) yields  $\minammkp(s_i(x)) = \minammkp(x)$ for $(x,i) \in (0,\infty) \times \Z$. From any $x \in (0,\infty)$, the orbit $\big\{ s_i(x): i \in \Z \big\}$ visits~$D$ exactly once, as noted in Lemma~\ref{l.battlefield}(4). 
Hence, $\minammkp(0,\infty) = \minammkp(D)$. Since $s$ is continuous, we learn that $\minammkp$ is continuous on all of~$\R$. 
The stated formula for this function is invariant under the replacement $x \to s_i(x)$ for any $i \in \Z$, so the validity of the formula passes from $D$ to~$\R$.

{\bf (3).} Consider the symmetry \abmnkpmacspace $\to$ \abmnkpmacspace that sends 
$$
\textrm{$(a,b,m,n)$ to $\big(  a(-1 - \bullet),b(-1 - \bullet),m(1 - \bullet),n(-1 - \bullet)  \big)$}
$$
given by taking $k=-1$ in Corollary~\ref{c.rolereversal}. 
By reflecting about minus one-half, the midpoint of $[-1,0]$, it acts as the inversion $x \mapsto x^{-1}$ on the central ratio~$\phi_0 = n_{0,-1}/m_{-1,0}$; and it does likewise on the Mina margin~$n_{\infty,-\infty}/m_{-\infty}$. But the Mina margin map sends the central ratio to the Mina margin. By considering $(a,b,m,n) \in$ \abmnkpmacspace  with $\phi_0 = x$, we confirm that $\minammkp(x^{-1}) = \minammkp(x)^{-1}$. Take $x=1$ and note $\minammkp \geq 0$ to find that $\minammkp(1) =1$.

{\bf (4).} 
As noted in the proof of the second part, $\minammkp(0,\infty) = \minammkp(D)$. By Proposition~\ref{p.relativereward}(1,2), the range of~$\minammkp$ takes the form $[\lambda^{-1},\lambda]$ where $\lambda$ is the supremum of the adopted values. But $\lambda = \lambdamaxkp$ since the supremum in Definition~\ref{d.lambdamax} is unchanged when taken over standard solutions. \qed

{\bf Proof of Theorem~\ref{t.minamarginvalues}.}  The value $\lambdamaxkp$ has been identified as the supremum of the values taken by the continuous map $\minammkp$ on the precompact set~$D \subset (0,\infty)$, so this value is finite. Since $1 \in D$ and $\minammkp(1) = 1$,  $\lambdamaxkp \geq 1$.  
As  noted in the preceding proof, the values of the Mina margin adopted by elements of \abmnkpmacspace are not restricted by considering only standard elements; the resulting set of values is $\minammkp(0,\infty)$, which equals $\big[\lambdamaxkp^{-1},\lambdamaxkp\big]$ by  Proposition~\ref{p.relativereward}(4). This establishes the claims made by  Theorem~\ref{t.minamarginvalues}. \qed

{\em Remark.} By Proposition~\ref{p.relativereward}(1,3), $\mc{M}_{\kappa,\rho}\big( \frac{2 - \kappa \rho}{2+\kappa \rho}\big) =  \mc{M}_{\kappa,\rho}\big( \frac{2 + \kappa \rho}{2-\kappa \rho}\big) =1$.
So the function  $\lambdamax(\kappa,\rho) - 1$ vanishes at the endpoints of~$D = D_{\kappa,\rho}$, and its oscillations thereon determine its range.
The element $(\kappa,\rho) \in (0,1]^2$ for which $D$ is maximal is $(1,1)$, with $D = (1/3,3]$. This offers circumstantial support for Conjecture~\ref{c.lambdamax}.  

\section{Penny Forfeit}\label{s.pennyforfeit}

\leavevmode

We now solve  the one-step sub-game of \tlpkp, which we call $(\kappa,\rho)$-Penny Forfeit or~\pfkp. In doing so, we will see the point of entry of the stronger condition 
 $\rho \leq 1$, which is found in the Nash-\abmnmacspace relationship as stated in Theorem~\ref{t.nashabmn}.

Let $(\kappa,\rho) \in (0,1] \times (0,\infty)$. In \pfkpspace with boundary condition $(m_{-1},m_1,n_{-1},n_1) \in \R^4$ satisfying $m_{-1} < m_1$ and $n_1 < n_{-1}$, 
Maxine and Mina stake $a$ and $b$, and Maxine wins with probability~$\frac{a^\rho}{a^\rho + b^\rho}$. Maxine and Mina's mean winnings  are 
\begin{equation}\label{e.maxineminawinnings}
 \Big( \tfrac{\kappa a^\rho}{a^\rho+b^\rho} + \tfrac{1-\kappa}{2}  \Big) m_1 +   \Big( \tfrac{\kappa b^\rho}{a^\rho+b^\rho} + \tfrac{1-\kappa}{2}  \Big) m_{-1} - a  \, \, \, \, \textrm{and} \, \, \, \, 
 \Big( \tfrac{\kappa b^\rho}{a^\rho+b^\rho} + \tfrac{1-\kappa}{2}  \Big) n_{-1} +   \Big( \tfrac{\kappa a^\rho}{a^\rho+b^\rho} + \tfrac{1-\kappa}{2}  \Big) n_1 - b \, . 
 \end{equation}
 \begin{lemma}\label{l.pennyforfeit}
Suppose that $\rho \in (0,1]$. For $\kappa \in (0,1]$,
there is a unique pair $(a,b) \in [0,\infty)^2$
for which 
  the expressions in~(\ref{e.maxineminawinnings}) are both global maxima as the variables $a$ and $b$ are respectively varied over $[0,\infty)$. It is given by 
\begin{equation}\label{e.absolution}
(a,b) \, = \, \kappa \rho \cdot \bigg(\frac{M^{1+\rho} N^\rho}{(M^\rho+N^\rho)^2},\frac{M^\rho N^{1+\rho}}{(M^\rho+N^\rho)^2} \bigg) \, , \, \, \, \, \textrm{with} \, \, \, \, M = m_{-1,1}  \, \, \, \, \textrm{and} \, \, \, \, N = n_{1,-1}  \, .
\end{equation}
Note that $a$ and $b$ are strictly positive.
\end{lemma}
{\bf Proof.} 
The maximizing pair cannot be $(0,0)$.
Indeed, if for example $a$ equals zero,
then an infinitesimal increase of $b$ from zero will increase Mina's expected payoff from~$\tfrac{n_{-1}+n_1}{2}$ to 
$$
 \big( \tfrac{1-\kappa}{2} + \kappa \big) n_{-1} + \tfrac{1-\kappa}{2} n_1 = \tfrac{1+\kappa}{2}n_{-1} + \tfrac{1-\kappa}{2} n_1 \, . 
$$

A critical point $(a,b)$ is given by setting the respective partial derivatives in $a$ and~$b$ of the two expressions in 
(\ref{e.maxineminawinnings})
equal to zero: the conditions are 
\begin{equation}\label{e.critcond}
\frac{\kappa \rho \, b^\rho a^{\rho-1}}{(a^\rho+b^\rho)^2}M - 1 \, = \, \frac{\kappa \rho \, a^\rho b^{\rho-1}}{(a^\rho+b^\rho)^2}N - 1 \, = \, 0  
\end{equation}
and these imply that $\kappa \rho \, b^\rho a^{\rho-1} M \;=\; \kappa \rho \,  a^\rho b^{\rho-1} N$. Since $ab \neq 0$,
$bM = aN$. Substituting $b = aN/M$ into $\kappa \rho \, b^\rho a^{\rho-1} M = (a^\rho+b^\rho)^2$, dividing by $a^{2\rho-1}$ and rearranging yields the formula for~$a$ in~(\ref{e.absolution}), with the formula for $b$ following from $b = aN/M$. The solution is positive and unique.

That the solution is a global maximum is due to $\rho \leq 1$. Indeed,~$\frac{a^{\rho-1}}{(a^\rho+b^\rho)^2}$ then has numerator that is decreasing in~$a \geq 0$, so that, since the denominator is increasing in this variable, the expression is decreasing. With an analogous property for $\frac{b^{\rho-1}}{(a^\rho+b^\rho)^2}$, this has the implication that the critical point in~(\ref{e.critcond})  is global in the sense of Lemma~\ref{l.pennyforfeit}, completing the proof of this result.  \qed

\medskip

\subsubsection{Other parameter regimes}
When $\rho > 1$, the argument above continues to identify the pair~$(a,b)$ in~(\ref{e.critcond}) as a critical point. However, the numerator in~$\frac{a^{\rho-1}}{(a^\rho+b^\rho)^2}$ is now increasing, and this sets up $a=0$ as a rival for the global maximizer of the first function in~(\ref{e.maxineminawinnings}). The condition
$M/N \geq (\rho - 1)^{1/\rho}$
characterises when the rival $a=0$ falls short and when the putative critical point is global.  
Switching $M$ and~$N$ in the last bound yields the applicable condition in regard to the second  function in~(\ref{e.maxineminawinnings}). 
 In summary, a global maximum in the sense of Lemma~\ref{l.pennyforfeit} never exists when $\rho > 2$; when $\rho \in (1,2]$, it entails that $M/N$ be suitably close to one, by lying in   $\big[(\rho - 1)^{1/\rho},(\rho-1)^{-1/\rho}\big]$. 

The game ${\rm PF}(1,\rho)$ is a two-player  Tullock contest whose   equilibrium analysis  has been addressed in all cases.
 The global maximum when it exists was found in~\cite{Nti1999}. When the global maximum fails to exist, mixed equilibria have been shown to exist~\cite{Wang2010} and to be unique~\cite{Ewerhart2017b,FengLu2017} when $\rho \in (1,2)$ and also to exist uniquely~\cite{Ewerhart2025} when $\rho \geq 2$. 
 
 \medskip

 \subsubsection{`Hand It Over!': A strategically equivalent zero-sum game when terminal rewards are equal}

Moulin and Vial~\cite{MoulinVial1978} developed a theory reducing certain two-person games to strategically equivalent
zero-sum games, in the sense that equilibrium strategies coincide up to constant payoff shifts.
This framework is applicable to ${\rm PF}(\kappa,\rho)$ in the equal-reward case
$m_{-1,1} = n_{1,-1}$.

In the strategically equivalent variant `Hand It Over!' game, stakes are not surrendered to the bank; rather, the
stake of a given player is handed to the opponent. After suitable constant payments are added to
each player's payoff, a zero-sum game results whose equilibria coincide with those of
${\rm PF}(\kappa,\rho)$. Indeed, the respective mean winnings in~(\ref{e.maxineminawinnings})
are shifted by amounts that are independent of $a$ and of $b$, leaving the solutions
in~(\ref{e.absolution}) unaltered.

This strategic equivalence breaks as soon as the gameboard has at least two sites in open play,
because the equilibrium values $m_i$ and $n_i$ then differ from those of the one-step game by
terms that depend nontrivially on equilibrium stakes from both players.

\chapter{\abmnmacspace solutions: explicit forms and fixed-parameter asymptotics}\label{c.abmnformfixed}
Here we solve the \abmnkpmacspace equations explicitly and deduce consequences. 
After giving the straightforward proof of the strict monotonicity of $m$- and $n$-differences recorded in Theorem~\ref{t.abmnpositive}(1), we re-express in the first section the \abmnmacspace system via the two-variable-per-site \mnmacspace  equations. This system permits the iterative computation of consecutive $m$- and $n$-differences, leading to the explicit sum-of-products representation in Theorem~\ref{t.defexplicit}.
The fixed-parameter asymptotics Theorem~\ref{t.abmn} will be obtained by analysing this representation. In the next two sections, we offer elements needed for that analysis: first, the asymptotics of the map~$s$; and then the resulting asymptotics for the $s$-orbit. Obtaining also needed asymptotics for the $c$ and $d$ maps that appear in the products in the representation, we give the proof of Theorem~\ref{t.abmn} in the fourth section. The chapter
 ends with the proof of  Theorem~\ref{t.abmnpositive}(2) on the finiteness of boundary data for elements of \abmnkpmac, which is a quick  corollary of Theorem~\ref{t.abmn}.

{\bf Proof of Theorem~\ref{t.abmnpositive}(1).} Since $a_i$ and $b_i$ are positive,  \abmnmac$(3)$ implies that $m_{i+1} \geq m_{i-1}$
(and in fact the bound is strict).
  Rearrange \abmnmac$(1)$ in the form 
$$
m_i   =  \Big( \kappa \tfrac{a_i^\rho}{a_i^\rho + b_i^\rho} + \tfrac{1-\kappa}{2} \Big) m_{i+1} +  \Big( \kappa \tfrac{b_i^\rho}{a_i^\rho + b_i^\rho} + \tfrac{1-\kappa}{2} \Big)  m_{i-1} - a_i \, .
$$
The right-hand $m$-coefficients are non-negative with unit sum. 
From  $m_{i-1} \leq m_{i+1}$, we thus find that
  $m_i \leq m_{i+1} - a_i$. Since $a_i > 0$,
$m_i < m_{i+1}$. That $n_{i+1} < n_i$ follows similarly. \qed

Note that Theorem~\ref{t.abmnpositive}(1) yields that the boundary vector $(m_{-\infty},m_\infty,n_{-\infty},n_\infty)$ exists, with $m_{-\infty} < m_\infty$ and $n_{-\infty} > n_\infty$, which is part of the inference stated in Theorem~\ref{t.abmnpositive}(2). However, in principle  $m_{-\infty}$ or $n_\infty$ may be $-\infty$
and $m_\infty$ or $n_{-\infty}$,~$\infty$. These possibilities will be excluded (and the proof of Theorem~\ref{t.abmnpositive}(2) completed) at the end of Section~\ref{s.fixedparam}, on the basis of the \abmnkpmacspace asymptotics Theorem~\ref{t.abmn}.

\section{Explicit ABMN solutions}\label{s.explicitabmn}
  The  \rvv $\big\{ m_i,n_i: i \in \Z \big\}$ satisfy the \mnmacspace system on $\Z$ if, for $i \in \Z$,
\begin{eqnarray*}
 m_{i-1,i}  \big( M_i^\rhomac + N_i^\rhomac \big)^2 & = & \kappa M_i^{2\rhomac+1} \, + \, \tfrac{1-\kappa}{2} \cdot M_i  \big( M_i^\rhomac + N_i^\rhomac \big)^2   +\kappa(1-\rhomac)M_i^{1+\rhomac} N_i^\rhomac   \\ 
 n_{i+1,i}  \big( M_i^\rhomac + N_i^\rhomac \big)^2  & = & \kappa N_i^{2\rhomac+1} \, + \, \tfrac{1-\kappa}{2} \cdot N_i  \big( M_i^\rhomac + N_i^\rhomac \big)^2  +\kappa(1-\rhomac)M_i^\rhomac N_i^{1+\rhomac} \, ,
\end{eqnarray*}
where $M_i := m_{i-1,i+1} = m_{i+1} - m_{i-1} $ and $N_i := n_{i+1,i-1} = n_{i-1} - n_{i+1}$. 
We will call these equations $\textrm{MN}(1)$ and  $\textrm{MN}(2)$.
  \begin{proposition}\label{p.abmnsolvesmn}
 Let $(a,b,m,n) \in $ \abmnkpmac. 
 The $(m,n)$-components 
 solve the \mnmacspace system on~$\Z$. We have that
  \begin{equation}\label{e.abclaim}
   a_i = \frac{\rhomac \kappa M_i^{1+\rhomac} N_i^\rhomac}{(M_i^\rhomac+N_i^\rhomac)^2}  \, \, \, \, , \, \, \, \,  b_i = \frac{\rhomac  \kappa  M_i^\rhomac N_i^{1+\rhomac}}{(M_i^\rhomac+N_i^\rhomac)^2}  \, \, \, \, \textrm{and} \, \, \, \,
    \frac{a_i^\rhomac}{a_i^\rhomac+b_i^\rhomac} =     \frac{M_i^\rhomac}{M_i^\rhomac+N_i^\rhomac}
  \end{equation}
for each $i \in \Z$.
  \end{proposition}
  {\bf Proof.} 
  From \abmnmac$(3,4)$ follows~(\ref{e.abclaim}). Since $a_i^\rhomac + b_i^\rhomac > 0$, we may express \abmnmac$(1)$  in the form
$$
 m_i \, = \, \Big( \kappa \tfrac{a_i^\rhomac}{a_i^\rhomac + b_i^\rhomac} +  \tfrac{1-\kappa}{2} \Big) m_{i+1} + \Big( \kappa \tfrac{b_i^\rhomac}{a_i^\rhomac + b_i^\rhomac} +  \tfrac{1-\kappa}{2} \Big) m_{i-1} \, - \, a_i \, ,
$$
we find from~(\ref{e.abclaim}) that 
  $$
  m_{i-1,i}   \, = \,  \bigg( \frac{M_i^\rhomac}{M_i^\rhomac+N_i^\rhomac} + \frac{1-\kappa}{2}\bigg)  M_i \, - \, \frac{\rhomac \kappa M_i^{1+\rhomac} N_i^\rhomac}{(M_i^\rhomac+N_i^\rhomac)^2} \, ,
  $$
   whence \mnmac$(1)$ holds. Equation \mnmac$(2)$ is obtained similarly, from \abmnmac$(2)$. \qed

\begin{definition}\label{d.fourq}
Let $(a,b,m,n) \in$ \abmnkpmac.
Define $\Z$-indexed sequences $\delta$, $\beta$, $\gamma$ and $\phi$ so that
$$
\delta_i = \frac{n_{i,i-1}}{n_{i+1,i-1}} \, \, , \, \, \beta_i =  \frac{n_{i+1,i-1}}{m_{i-1,i+1}} \, \, , \, \, \gamma_i = \frac{m_{i-1,i}}{m_{i-1,i+1}} \, \, \, \textrm{and} \, \, \phi_i =  \frac{n_{i,i-1}}{m_{i-1,i}} \, .
$$
\end{definition}
Two useful relations that result 
are $\phi_i \gamma_i =  \beta_i \delta_i$ and $\phi_{i+1} (1-\gamma_i) = \beta_i  (1-\delta_i)$.

From Definition~\ref{d.fourfunctions},
recall the four basic functions $\gamma$, $\delta$, $\phi_0$ and $\phi_1$ that map $(\kappa,\rho,\beta) \in W \times (0,\infty)$ to~$(0,\infty)$.
\begin{lemma}\label{l.gdphiphi}
For $i \in \Z$,
 $\gamma_i = \gamma(\kappa,\rho,\beta_i)$, 
 $\delta_i = \delta(\kappa,\rho,\beta_i)$,
 $\phi_i = \phi_0(\kappa,\rho,\beta_i)$ and  $\phi_{i+1} = \phi_1(\kappa,\rho,\beta_i)$. 
\end{lemma}
{\bf Proof.}
In MN($1$), write $m_{i-1,i} = M_i \gamma_i$.  Then divide by $M_i^{2\rho+1}$, use $\beta_i = N_i/M_i$, and rearrange to obtain
$$
   \big( \gamma_i -  \tfrac{1-\kappa}{2}  \big) \big(1+\beta_i^\rho\big)^2   \, = \,   \kappa  \big( 1+  (1-\rho)  \beta_i^\rho \big) \, .
$$
In MN($2$), write $n_{i+1,i} = M_i \beta_i  (1-\delta_i)$,
divide by $M_i^{2\rho+1}$, use $\beta_i = N_i/M_i$, cancel a factor of $\beta_i$ and rearrange,  to arrive at
$$
 \big( \tfrac{1+\kappa}{2} - \delta_i \big) \big(1+\beta_i^\rho\big)^2  = \kappa \beta_i^{\rho} \big( \beta_i^\rho  + 1-\rho \big)  \, .
$$
Rearranging the preceding two displays yields 
\begin{eqnarray}
 \gamma_i  & = & \frac{(1-\kappa)\beta_i^{2\rho} + 2(1-\rho\kappa)\beta_i^\rho+  1+\kappa}{2 (1+\beta_i^\rho)^2} \nonumber \\
\delta_i & = & \frac{(1-\kappa)\beta_i^{2\rho} + 2(1+\rho\kappa) \beta_i^\rho + 1 + \kappa}{2(1+\beta_i^\rho)^2} \, , \nonumber 
\end{eqnarray}
or  $\gamma_i = \gamma(\kappa,\rho,\beta_i)$ and 
 $\delta_i = \delta(\kappa,\rho,\beta_i)$ in view of the form of the functions $\gamma$ and $\delta$ presented in Definition~\ref{d.fourfunctions}.
 
Using the first relation noted after Definition~\ref{d.fourq} yields  $\phi_i = \phi_0(\kappa,\rho,\beta_i)$; the second,  $\phi_{i+1} = \phi_1(\kappa,\rho,\beta_i)$. \qed

{\bf Proof of Theorem~\ref{t.defexplicit}.}
In Definition~\ref{d.scd}(2), the function $c:(0,\infty) \to (0,\infty)$, $c(\bullet) = c(\kappa,\rho,\bullet)$, is specified so that $c(x) = 1/\gamma(\kappa,\rho,\beta)$ where $\beta \in (0,\infty)$ satisfies $\phi_0(\kappa,\rho,\beta) = x$.
By Definition~\ref{d.scd}(1), the map $s:(0,\infty) \to (0,\infty)$  sends any value adopted by the function $\phi_0$, for some choice of $\beta \in (0,\infty)$, to the value of $\phi_1$ assumed for that same~$\beta$. In view of the relations for $\phi_i$ and $\phi_{i+1}$ identified in Lemma~\ref{l.gdphiphi}, the action of $s$ on elements of the sequence $\big\{ \phi_i: i \in \Z \}$ specified by an \abmnkpmacspace element is simply the shift: $s(\phi_i) = \phi_{i+1}$ for $i \in \Z$.

For $x \in (0,\infty)$ given, consider then an \abmnkpmacspace element $(a,b,m,n)$
for which $\phi_0 = n_{0,-1}/m_{-1,0}$ equals~$x$.
 Recalling Definition~\ref{d.stabc}, we have $c_i(x) = c(s_i(x))$ where $s_i(x)$ equals $\phi_i$ due to $x=\phi_0$ and iteration of the shift action of~$s$.
 In light of the preceding paragraph then,  $c_i(x) = 1/\gamma(\kappa,\rho,\beta_i)$\hfff{gammaess} since $\phi_0(\kappa,\rho,\beta_i) = \phi_i$ by Lemma~\ref{l.gdphiphi}.
 Hence, we obtain the first equality as we write
 \begin{equation}\label{e.ciformula}
  c_i(x) - 1 \, = \, \frac{1 - \gamma(\kappa,\rho,\beta_i)}{\gamma(\kappa,\rho,\beta_i)}  \, = \,  \frac{1 - \gamma_i}{\gamma_i}   \, = \,  \frac{m_{i,i+1}}{m_{i,i-1}} \, ,
\end{equation}
the second equality\footnote{Usages such as  Lemma~\ref{l.gdphiphi}($\gamma$) and  Proposition~\ref{p.sfacts}(2,$f$) refer to the statement made the result in question about the object~$\gamma$ or~$f$.} due to Lemma~\ref{l.gdphiphi}($\gamma$) and the third to Definition~\ref{d.fourq}($\gamma$).
With the product notation from Definition~\ref{d.zdef} applying when the index $j$ is negative, 
we find that 
\begin{equation}\label{e.mdifferenceratio}
\frac{m_{j,j+1}}{m_{-1,0}} \, = \, \prod_{i=0}^j \big( c_i(x) - 1 \big)
\end{equation}
for any $j \in \Z$. 
The \abmnkpmacspace element $(a,b,m,n)$ under consideration may be dilated by varying $m_{-1,0} \in (0,\infty)$,
in correspondence with the dilation factor $\mu$ that appears in Theorem~\ref{t.defexplicit}.
By setting $m_{-1,0} =1$, we reduce the task of proving the theorem to checking that 
$(a,b,m,n)$ equals the default quadruple 
$\big(\adefault(x),\bdefault(x),\mdefault(x),\ndefault(x)\big)$ (so $\mu=1$).
And indeed we have proved the $m$-component projection of the desired identity, because the right-hand side in~(\ref{e.mdifferenceratio}) is $\mdefault_{k+1}(x) - \mdefault_k(x)$ from 
Definition~\ref{d.zdef}.

 Evident variations of the argument leading to~(\ref{e.ciformula}) yield
 $$
  d_i(x) - 1  \, = \,  \frac{1 - \delta(\kappa,\rho,\beta_i)}{\delta(\kappa,\rho,\beta_i)}  \, = \,  \frac{1 - \delta_i}{\delta_i}   \, = \,   \frac{n_{i+1,i}}{n_{i,i-1}} \, ,  \, \, \, \,  \textrm{so that} \, \, \, \, \frac{n_{j+1,j}}{n_{0,-1}} \, = \, \prod_{i=0}^j \big( d_i(x) - 1 \big)
  $$
for $j \in \Z$. Since $n_{0,-1} = x \, m_{-1,0} = x$ by our normalization, we obtain
the $n$-component claim made in  Theorem~\ref{t.defexplicit}.
The expressions for the sequence $a_i$ and $b_i$ in Proposition~\ref{p.abmnsolvesmn} coincide with the formulaic counterparts in Definition~\ref{d.zdef}.
This completes the proof of Theorem~\ref{t.defexplicit}. 
\qed

\section{$\phi_0$, $\phi_1$ and $s$ asymptotics}

We record large~$\beta$ asymptotics of the functions $\phi_0(\bullet) = \phi_0(\kappa,\rho,\bullet)$ and $\phi_1(\bullet) = \phi_1(\kappa,\rho,\bullet)$ from Definition~\ref{d.fourfunctions}, and of the mapping $s: \phi_0 \mapsto \phi_1$.

\begin{lemma}\label{l.twophiands}
\leavevmode
\begin{enumerate}
\item
For  $\rho \in (0,1]$ and $\kappa \in (0,1)$,
\begin{align*}
\phi_0 &= \beta + \frac{4\rho\kappa}{1 - \kappa} \, \beta^{1 - \rho} + O\left( \beta^{1 - 2\rho} \right) \, , \\
\phi_1 &= \beta - \frac{4\rho\kappa}{1 + \kappa} \, \beta^{1 - \rho} + O\left( \beta^{1 - 2\rho} \right) \, \, \, \, \textrm{as $\beta \to \infty$} \, ,
\end{align*}
and
\[
s(x) = x - \frac{8\rho\kappa}{1 - \kappa^2} \, x^{1 - \rho} + O\left( x^{1 - 2\rho} \right)
\quad \text{as } x \to \infty \, .
\] 
\item For $\kappa = 1$ and $\rho \in (0,1)$,
$$
\phi_0 = \tfrac{1+\rho}{1-\rho} \beta + O (\beta^{1-\rho}) \, \, \,, \, \, \, 
\phi_1 = \beta + O (\beta^{1-\rho}) \, \, \, \, \,  \textrm{and} \, \, \, \, \, 
s(x) = \tfrac{1-\rho}{1+\rho} x + O \big( x^{1-\rho} \big) \, . 
$$
\end{enumerate}
\end{lemma}
{\bf Proof.} In either case, the weaker condition we consider for a $(\kappa,\rho)$ pair, namely membership of $W$ as specified in~(\ref{e.weakregion}), is met; this enables the use of Lemma~\ref{l.incphi}, so $s$ is well defined.
 The $\phi$-asymptotics are computed by working with the formulas in Definition~\ref{d.fourfunctions}. 
For example, in the latter case $\kappa=1$ and $\rho\in(0,1)$, the formulas
(\ref{e.phizero.rho}) and (\ref{e.phione.rho}) reduce to
\[
\phi_0(1,\rho,\beta)
= \beta\cdot \frac{2(1+\rho)\beta^\rho+2}{2(1-\rho)\beta^\rho+2}
= \beta \cdot \frac{(1+\rho)\beta^\rho+1}{(1-\rho)\beta^\rho+1}
\]
and
\[
\phi_1(1,\rho,\beta)
= \beta\cdot \frac{2\beta^{2\rho}+2(1-\rho)\beta^\rho}{2\beta^{2\rho}+2(1+\rho)\beta^\rho}
= \beta\cdot \frac{\beta^\rho+(1-\rho)}{\beta^\rho+(1+\rho)} \, .
\]
Thus, as $\beta\to\infty$,
\[
\phi_0(1,\rho,\beta)
= \tfrac{1+\rho}{1-\rho}\,\beta + O(\beta^{1-\rho})
\qquad\text{and}\qquad
\phi_1(1,\rho,\beta)
= \beta + O(\beta^{1-\rho}) \, .
\]
Both $s$-estimates follow straightforwardly from the $\phi$-asymptotics in view of $s: \phi_0 \mapsto \phi_1$. \qed

\section{Asymptotics for the backward orbit of $s$}

Recall from Definition~\ref{d.stabc} that $s_{-i}$ denotes the $i$-fold backward iterate $s$.
For $x$ close to one, the asymptotics of $s_{-i}(x)$ differ according to whether $\rho$ lies  in $(0,1)$ or equals $1$.

\begin{lemma}\label{l.phibeta}
Let $(a,b,m,n)$ be an element of \abmnkpmacspace of battlefield index zero (so that $\phi_0 \in D$). 
\begin{enumerate}
\item For $\kappa \in (0,1)$ and $i > 0$,
\[
\phi_{-i} \, = \, 
\begin{cases}
\displaystyle \, 
\frac{8\kappa}{1 - \kappa^2} \, i +  O(\log i), & \text{for } \rho = 1 \, , \\[1ex]
\displaystyle \,
\left( \frac{8 \rho^2 \kappa}{1 - \kappa^2} \right)^{1/\rho} \, i^{1/\rho} + O\left( i^{(1 -\rho)/\rho} \right), & \text{for } \rho \in (0,1) \, .
\end{cases}
\]
These asymptotics are equally valid for $\beta_{-i}$.
\item Now let $\kappa = 1$
 and $\rho \in (0,1)$.
For  $i > 0$, 
$\phi_{-i} = \big( \tfrac{1+\rho}{1-\rho} \big)^{i + \sigma + o(1)}$, 
where $\sigma = \sigma(\phi_0)$ is bounded in absolute value.
And $\beta_{-i}$ is likewise, with $\sigma -1$ in place of $\sigma$.
\end{enumerate}
\end{lemma}
{\bf Proof: (1).} Note that $\phi_{-i}$ equals the $i$\textsuperscript{th}  element $s_{-i}(x)$ on the backward orbit of $s$ 
whose starting point  $x = \phi_0$ lies in the central domain~$D$,
 since the battlefield index equals zero.

Whenever~$\rho \in (0,1]$, the $s$-asymptotic in Lemma~\ref{l.twophiands} implies that the inverse map $s_{-1}$ satisfies 
\begin{equation}\label{e.sminusone}
s_{-1}(x) = x + \frac{8\rho\kappa}{1 - \kappa^2} \, x^{1 - \rho} + O\left( x^{1 - 2\rho} \right)
\quad \text{as } x \to \infty,
\end{equation}
We explain  how to obtain an asymptotic for $s_{-i}(x)$ from this input, doing so first in outline.
 Writing $C = \frac{8\rho\kappa}{1 - \kappa^2}$, set $x_0 = x$ and iterate the recursion $x_{n+1} = x_n + C x_n^{1-\rho} + O\! \big( x_n^{1-2\rho} \big)$. Neglecting the $O(\cdot)$ term permits us to interpret $x$ as an approximate solution to the differential equation 
$\frac{dx}{dn} = C x^{1-\rho}$, whence $x_n$ is seen to grow as $ A \, n^{1/\rho}$, with $A = (C \rho)^{1/\rho}$. Reintroducing the neglected terms introduces a perturbation $\sum_{k=1}^n x_k^{1-2\rho}$ 
to the value of~$x_n$. Since $x_n^{1-2\rho} = O(n^{\rho^{-1}-2})$, this perturbation is  $O\!\left( n^{\rho^{-1}- 1} \right)$ when $\rho \in (0,1)$; a factor of $\log n$ is required when $\rho=1$. From $\phi_{-i} = s_{-i}(x)$, the $\phi$-asymptotics claimed in the lemma are obtained, at least heuristically.

We give a rigorous argument for $\rho \in (0,1)$; the $\rho =1 $ involves introducing suitable logarithmic factors.  Find large positive constants $C_0$ and $D$ such that for $x \geq C_0$ the implied constant in the big-$O$ term in~(\ref{e.sminusone}) is at most~$D$.
By orbit escape Corollary~\ref{c.orbitescape}, select $n_0$ such that $\phi_{-n} \geq C_0$  for $n \geq n_0$. 
Set $t_n = (\phi_{-n}/A)^{\rho} - n$, so that
$\phi_{-n} = A (n+t_n)^{1/\rho}$.
Let $l = l(n)\in\Z$ be the nearest integer to $t_n$.
Since $\phi_{-n} \geq C_0$, for $n \geq n_0$, the unit neighbourhood of $n+t_n$ lies in $(0,\infty)$, so that we may apply the mean value theorem to the map 
$u\mapsto A(n+u)^{1/\rho}$ to learn that 
\[
\big| \phi_{-n} - A(n+l(n))^{1/\rho} \big|
= A\big|(n+t_n)^{1/\rho}-(n+l(n))^{1/\rho}\big|
\le 2 A \rho^{-1} (n+l)^{\rho^{-1} - 1} 
\]
where recall that $A$ equals $(C \rho)^{1/\rho}$. 
For a positive constant $K$ suitably determined by~$D$, we select $n_1 \geq \max \{ n_0, K^2\}$, and set $\ell = l(n_1)$ (so that the offset value $\ell \in \N$ is now fixed, independently of~$n$).
Setting $x_n = \phi_{-(n-\ell)}$ and $e_n = x_n - A n^{1/\rho}$, we will argue by induction on $n \geq n_1 +\ell$ that  
$\big|  e_n \big| \le K\, n^{\rho^{-1} -1}$.  Expanding the power of $i + \ell$ in the resulting upper bound on $\vert e_{i+\ell} \vert$ yields the claim asymptotic on~$\phi_{-i}$.

The last display assures the inductive base case $n = n_1 + \ell$. Suppose then that the inductive hypothesis holds for some $n \geq n_1 + \ell$.
The $x_n$-sequence satisfies 
$x_{n+1} = x_n + C x_n^{1-\rho} + O(x_n^{1-2\rho})$.
 Substitute  $x_n = A n^{1/\rho} + e_n$ into this recursion
to find that 
\begin{equation}\label{e.back}
A (n+1)^{1/\rho} + e_{n+1} = A n^{1/\rho} + e_n + C (A n^{1/\rho} + e_n)^{1-\rho} + O(n^{\rho^{-1}-2}).
\end{equation}
By 
Taylor expansion, $(A n^{1/\rho} + e_n)^{1-\rho}$ equals 
\begin{eqnarray*}
 & &  A^{1-\rho} n^{\rho^{-1} - 1} + (1-\rho) A^{-\rho} n^{-1} e_n - \frac{\rho(1-\rho)}{2} A^{-1-\rho} n^{-\rho^{-1} -1 } e_n^2 \big( 1 + O( n^{-\rho^{-1}}e_n) \big) \\
 & = &  A^{1-\rho} n^{\rho^{-1} - 1} + (1-\rho) A^{-\rho} n^{-1} e_n  + O(1) K^2 n^{\rho^{-1} - 3} \, , 
 \end{eqnarray*}
where the displayed equality is due the inductive hypothesis in the guise $n^{-1/\rho - 1} e_n^2 \leq K^2 n^{\rho^{-1} - 3}$ and $n^{-\rho^{-1}}e_n \leq K n^{-1}$.
The final displayed term may be written  $O(1) n^{\rho^{-1} - 2}$ since $n \geq n_1 + \ell \geq n_1  \geq K^2$.
Substituting back into~(\ref{e.back}), and noting that the resulting right-hand $C A^{1 - \rho}  n^{\rho^{-1} - 1}$ term equals 
$$
 A \big( (n+1)^{1/\rho} - n^{1/\rho} \big) + O(n^{\rho^{-1}-2})
 $$ in view of $A = (C \rho)^{1/\rho}$,
the 
$A (n+1)^{1/\rho} $ terms cancel and we obtain
\[
|e_{n+1}| \le |e_n| \Big( 1 + \frac{1-\rho}{\rho n} + \frac{C_1}{n^2} \Big) + C_2 n^{\rho^{-1} - 2},
\]
for suitable constants $C_1, C_2 > 0$.  
We obtain  $\big|  e_{n+1} \big| \le K\, (n+1)^{\rho^{-1} -1}$, with the $C_2 n^{\rho^{-1} - 2}$ term
being absorbed into the right-hand side since the value of $C_2$ is determined by~$D$ and we may specify $K = K(D)$ suitably. The induction thus closes, implying that $|e_n| \le K n^{\rho^{-1}-1}$ holds for all $n \geq n_1 + \ell$.

In regard to~$\beta$-asymptotics, note that  $\beta_i = \beta(\phi_i)$, so that
$\phi_i  = \beta_i \big( 1 +  O ( \beta_i^{- \rho}) \big)$ whenever $\kappa \in (0,1)$ and $\rho \in (0,1]$, by Lemma~\ref{l.twophiands}(1); consequently $\beta_i  = \phi_i \big( 1 +  O ( \phi_i^{- \rho}) \big)$ .
So the $\phi$-asymptotics pass to the $\beta$-sequence.

{\bf (2).} By orbit escape, the sequence of  inverse-$s$ iterates $\phi_{-i}(x) = s_{-1} \big(  \phi_{-i+1}(x) \big)$, with $x = \phi_0 \in (0,\infty)$ given,
  tends to infinity in high~$i$. It is straightforward from Lemma~\ref{l.twophiands}(2) that $s_{-1}(x) = \tfrac{1+\rho}{1-\rho} x \, \big( 1 + O ( x^{-\rho}) \big)$. As such, $\phi_{-i}$ grows exponentially in $i$, and, if we write $\phi_{-i}$ in the form $\left( \tfrac{1+\rho}{1-\rho} \right)^i \psi_i$, the correction factors are seen to satisfy $\psi_{i+1} = \psi_1 \left( 1 + O(e^{-ci}) \right)$.
  Thus $\psi$-sequence is bounded away from zero and infinity, uniformly for $x$ in the central domain~$D$. In this way, we obtain the claimed $\phi_{-i}$-asymptotics.   By Lemma~\ref{l.twophiands}(2),
$\phi_{-i} = \tfrac{1+\rho}{1-\rho} \beta_{-i} + O (\beta_{-i}^{1-\rho})$, whence $\beta_{-i} = \tfrac{1-\rho}{1+\rho} \phi_{-i} \big( 1 + O (\phi_{-i}^{-\rho}) \big)$, yielding the $\beta_{-i}$-asymptotics. \qed

\section{Fixed-parameter \abmnmacspace asymptotics}\label{s.fixedparam} \leavevmode
 
 The obtained control on the $s$-orbit equips us for the next derivation.

{\bf Proof of Theorem~\ref{t.abmn}.}
The role-reversal  and shift symmetries for \abmnmacspace solutions noted in Section~\ref{s.roleshift} reduce Theorem~\ref{t.abmn}(4,5) to Theorem~\ref{t.abmn}(1,2,3). And since $n_{-i,-i-1}/m_{-i-1,-i}$ equals~$\phi_{-i}$,
and  $b_{-i}$ is $a_{-i} N_{-i}/M_{-i} = a_{-i} \beta_{-i}$, the $\phi_{-i}$- and $\beta_{-i}$-asymptotics offered by Lemma~\ref{l.phibeta} reduce Theorem~\ref{t.abmn}(1,2,3) to the claims made there regarding $m_{-i-1,-i}$  and $a_{-i}$. In addressing the first and second parts, and then the third, we will thus be concerned only with the $m_{-i-1,-i}$  and $a_{-i}$ estimates.

{\bf (1,2).}
 From $c(x) =1/\gamma(\kappa,\rho,\beta)$ and~(\ref{e.cbeta.rho}), we find that  
\[
c(x) =
\begin{cases}
\displaystyle
\frac{2}{1 - \kappa}
- \frac{4\kappa}{(1 - \kappa)^2} \cdot (1 + \beta)^{-2}
+ O\left( (1 + \beta)^{-4} \right),
& \text{for } \rho = 1 \, , \\[1.5ex]
\displaystyle
\frac{2}{1 - \kappa}
- \frac{4(1 - \rho)\kappa}{(1 - \kappa)^2} \cdot \beta^{-\rho}
+ O\left( \beta^{-2\rho} \right),
& \text{for } \rho \in (0,1)
  \, .
\end{cases}
\]
We have that $c_{-j}(x) = c \big( s_{-j}(x) \big) = 1/\gamma(\kappa,\rho,\beta_{-j})$, with $\beta_{-j}$-asymptotics offered by Lemma~\ref{l.phibeta}.
When $\rho = 1$,
 $$
  c_{-j}-1 \, = \, \frac{1+\kappa}{1-\kappa} -  \frac{(1+\kappa)^2}{16\kappa j^2} \left( 1 + O \left( j^{-1}\log j \right) \right) 
  =  \frac{1+\kappa}{1-\kappa} \left( 1 -  \frac{1-\kappa^2}{16\kappa j^2} \left( 1 + O \left( j^{-1} \log j \right) \right) \right) \, ,
$$ 
so that
\begin{eqnarray*}
  \frac{m_{-(i+1),-i}}{ m_{-1,0}}   & = &  \prod_{j=1}^i \big( c_{-j} - 1 \big)^{-1}  =   \prod_{j=1}^i 
  \left( \frac{1 - \kappa}{1 + \kappa} \right) 
   \left( 1 +  \frac{1-\kappa^2}{16\kappa j^2} \left( 1 + O \left( j^{-1} \log j \right) \right) \right)  \\
   & = & \sigma
\left( \frac{1 - \kappa}{1 + \kappa} \right)^i  \exp \big\{ {O(1) (\kappa^{-1} - \kappa)} \big\}  \Big( 1 + \tfrac{1-\kappa^2}{\kappa} O(i^{-1}) \Big)  \, ;
\end{eqnarray*}
we obtain the sought $m_{-(i+1),-i}$-asymptotic for  $\rho=1$ by absorbing  the  $\exp \big\{ {O(1) (\kappa^{-1} - \kappa)} \big\}$ factor into~$\sigma$. 

For $\rho \in (0,1)$,
  $$
  c_{-j}-1 \, = \, \frac{1+\kappa}{1-\kappa} - \frac{4(1-\rho)\kappa}{(1-\kappa)^2} \frac{1-\kappa^2}{8\rho^2 \kappa} j^{-1} \Big( 1 + O(j^{-1}) \Big)
  =  \frac{1+\kappa}{1-\kappa} \left( 1 - \frac{1-\rho}{2\rho^2} j^{-1} \Big( 1 + O(j^{-1}) \Big) \right) \, ,
  $$
 leading to
\begin{eqnarray*}
  \frac{m_{-(i+1),-i}}{ m_{-1,0}} 
 & = &   \prod_{j=1}^i \left( c_{-j} - 1 \right)^{-1}  
 =  \prod_{j=1}^i \left( \frac{1 - \kappa}{1 + \kappa} \right) 
\left( 1 + \frac{1 - \rho}{2\rho^2} j^{-1} \Big( 1 + O(j^{-1}) \Big) \right) \\
& = & \sigma
\left( \frac{1 - \kappa}{1 + \kappa} \right)^i 
  i^{\frac{1-\rho}{2\rho^2}}  \Big( 1 + O(i^{-1}) \Big)  \, ,
 \end{eqnarray*}
whence the claimed asymptotics for $m_{-i-1,-i}$. 

It remains to compute $a_{-i}$-asymptotics. The formula given in Proposition~\ref{p.abmnfacts}(1) expresses $a_{-i}$ in terms of $M_{-i} = m_{-i-1,-i} + m_{-i,-i+1}$ and $N_{-i} = M_{-i} \beta_{-i}$.
When $\rho \in (0,1)$,  we apply the derived $m$-asymptotics to find that  
$$ 
M_{-i}     =  m_{-1,0}  \cdot \sigma \cdot  \frac{2}{1-\kappa}  \cdot i^{{\frac{1-\rho}{2\rho^2}}}  \left( \frac{1 - \kappa}{1 + \kappa} \right)^{i} \left( 1 + O\left(i^{-1}\right) \right)
$$
with $\beta_{-i}$-asymptotics from Lemma~\ref{l.phibeta} then yielding
$$ 
N_{-i}     =  m_{-1,0}  \cdot \sigma \cdot  \frac{2}{1-\kappa} \left( \frac{8\rho^2\kappa}{1 - \kappa^2} \right)^{1/\rho} \cdot i^{{\frac{1+\rho}{2\rho^2}}}  \left( \frac{1 - \kappa}{1 + \kappa} \right)^{i} \left( 1 + O\left(i^{-1}\right) \right) \, .
$$
Noting $\big( M_{-i}/N_{-i} \big)^\rho = \beta_{-i}^{-\rho} \leq C/i$ offers a simplified asymptotic formula for $a_{-i}$, namely
$$
a_{-i} = \frac{\kappa \rho M_{-i}^{1+\rho} N_{-i}^\rho}{(M_{-i}^\rho + N_{-i}^\rho)^2} = \kappa \rho M_{-i}\big( M_{-i}/N_{-i} \big)^\rho \big( 1 + O(i^{-1}) \big) \, .
$$
Thus, when $\rho \in (0,1)$ $\big( M_{-i}/N_{-i} \big)^\rho = \frac{1-\kappa^2}{8 \rho^2 \kappa} i^{-1} \big( 1 + O(i^{-1}) \big)$ and the above $M_{-i}$-asymptotic yield the claimed $a_{-i}$-asymptotic.

For $\rho =1$, we adopt the same approach, and merely need to note the accurate form of $\beta_{-i}$-asymptotics from Lemma~\ref{l.phibeta}. We find that
$$
M_{-i}   \,  = \,  m_{-1,0}    \cdot \sigma \cdot \frac{2}{1-\kappa} 
       \left( \frac{1-\kappa}{1+\kappa} \right)^{i}  \Big( 1 + O(i^{-1}) \Big) \
$$
and 
$$
N_{-i}  \,  = \,  m_{-1,0}    \cdot \sigma \cdot \frac{16\kappa}{(1-\kappa)^2(1+\kappa)} \cdot i 
       \left( \frac{1-\kappa}{1+\kappa} \right)^{i} 
       \Bigl( 1 + O\big(  i^{-1} \log i \big) \Bigr) \,  ;
$$
from $M_{-i}/N_{-i} \leq C/i$, we note $a_{-i} = \kappa M_{-i}^2 N_{-i}^{-1} \big( 1 + O(i^{-1}) \big)$.  Substituting into this formula gives the sought $a_{-i}$-asymptotics for $\rho=1$.

{\bf (3).} 
When $\kappa =1$ and $\rho \in (0,1)$, we have $c = \frac{(1+\beta^\rho)^2}{1 + (1-\rho) \beta^\rho}$.
Thus, $c =  \tfrac{1}{1-\rho}  \beta^\rho + \tfrac{1-2\rho}{(1-\rho)^2} + O(\beta^{-\rho})$.
Applying~$\beta_{-j}$-asymptotics from Lemma~\ref{l.phibeta}(2) to  $c(\phi_{-j}) = \tfrac{1}{1-\rho}  \beta_{-j}^\rho + O(1)$,
we obtain
$$
 c(\phi_{-j} ) - 1 \, = \, \frac{1}{1-\rho} \left( \frac{1 + \rho}{1 - \rho} \right)^{\rho ( j  +  \sigma - 1) + o(1)} +  O(1) \, .
$$
We find that
$$
 \big( c(\phi_{-j} ) - 1 \big)^{-1} = \left(  {1 - \rho} \right)^{\rho j + 1} (1+\rho)^{-\rho j} \Big( 1 + O \big( \tfrac{1-\rho}{1+\rho} \big)^j  \Big) \cdot   \big( \tfrac{1-\rho}{1+\rho} \big)^{\rho(\sigma - 1) + o(1)}
$$
Using 
$m_{-(i+1),-i}  =  m_{-1,0}  \prod_{j=1}^i \big( c_{-j} - 1 \big)^{-1}$,
we have
$$
m_{-i-1,-i}
= m_{-1,0}\,(1-\rho)^i
\left(\frac{1-\rho}{1+\rho}\right)^{\rho\big(\frac{i(i+1)}{2} + i(\sigma - 1)\big) + o(i)} \, .
$$
Equivalently, 
$$
m_{-i-1,-i}
\, = \, m_{-1,0}\left(\frac{1-\rho}{1+\rho}\right)^{\rho i^{2}/2} e^{\chi i + o(i)} \, ,
$$
for a suitable constant $\chi = \chi(\rho,\sigma)$. 
From $M_{-i} = 
m_{-(i+1),-i}+
m_{-i,-(i-1)}$ and 
$m_{-(i+1),-i} \ll 
m_{-i,-(i-1)}$, we see that
$$
 M_{-i} =  m_{-1,0}\left(\frac{1-\rho}{1+\rho}\right)^{\rho (i-1)^2/2} e^{\chi i + o(i)} \, .
$$
Now $N_{-i}= M_{-i} \beta_{-i} = M_{-i} \left( \tfrac{1+\rho}{1-\rho} \right)^{i+\sigma - 1 + o(1)}$ via Lemma~\ref{l.phibeta}(2). The smallness of $M_{-i}$ relative to~$N_{-i}$ permits the same simplified asymptotic formula for $a_{-i}$ as seen earlier:  
$$
a_{-i} = \frac{\kappa \rho M_{-i}^{1+\rho} N_{-i}^\rho}{(M_{-i}^\rho + N_{-i}^\rho)^2} = \kappa \rho M_{-i}\big( M_{-i}/N_{-i} \big)^\rho \left( 1 + O(1) \left(\tfrac{1-\rho}{1+\rho}\right)^i \right) \, .
$$
We have that  $\big( M_{-i}/N_{-i} \big)^\rho = \beta_{-i}^{-\rho} =  \left(\frac{1-\rho}{1+\rho}\right)^{\rho(i+\sigma - 1 +o(1))}$, so that
$$
 a_{-i} =  m_{-1,0}\left(\frac{1-\rho}{1+\rho}\right)^{\rho i^2/2} e^{\chi i}  e^{o(i)} \, ,
$$
which is the $a_{-i}$-asymptotic asserted in Theorem~\ref{t.abmn}(3). 
This completes the proof of the theorem. \qed

The obtained estimates permit us to note the finiteness of \abmnkpmacspace boundary data.

{\bf Proof of Theorem~\ref{t.abmnpositive}(2).}
As noted after the proof of Theorem~\ref{t.abmnpositive}(1), $m_\infty$, $m_{-\infty}$, $n_\infty$ and $n_{-\infty}$ exist as elements of $\R \cup \{ \infty \} \cup \{ - \infty \}$. Since  $m_0,n_0 \in \R$, it is enough, in order  to infer that the four quantities are finite real numbers, to show that
the non-negative differences $m_{0,i}$, $m_{-i,0}$,  $n_{i,0}$ and  $n_{0,-i}$ are bounded above as $i$ varies over~$\N$. 
 These bounds  may be obtained by summing the estimates on consecutive differences $m_{j-1,j}$ and $n_{j,j-1}$ provided by Theorem~\ref{t.abmn}. \qed

 \chapter{Nash equilibria and the ABMN equations}\label{c.nashabmn}

Here we prove  Theorem~\ref{t.nashabmn} on  Nash-\abmnmacspace equivalence. 
The forward implication $(1)\!\!\!\!\implies\!\!\!\!(2)$  is proved in the first four sections, the reverse in the fifth.
The derivations follow the template given by the proof of the counterpart \cite[Theorem~$2.6$]{LostPennies} in~\cite[Chapter~$4$]{LostPennies}, with some substantial changes.

In the forward-implication proof, some arguments are new and others closely follow  counterparts  in~\cite[Chapter~$4$]{LostPennies}. To make our presentation self-contained while indicating where the overlap lies,  
the first three sections use the convention that {\bf Proof} denotes the start of an argument with substantial new elements, while {\bf Derivation} indicates one that is close to one in~\cite{LostPennies}. No lack of rigour should be inferred from use of the latter label, though we have sometimes opted for a more verbal style of presentation of such arguments.   
A different approach has been adopted for the reverse implication, as we explain in Section~\ref{s.reverse}.

\section{Escape is almost certain at a time-invariant Nash equilibrium}

To prove the forward implication, we consider $(S_-,S_+) \in \mc{N}_{\kappa,\rho} \cap \tis^2$.  As in Definition~\ref{d.quadruple}, write\footnote{The order $(S_-,S_+)$ is governed by the convention $- < +$ in which Mina precedes Maxine. Since Maxine stakes $a$ and Mina~$b$,
this results in the identification of $(S_-,S_+)$ with $(b,a)$.}  $b_i$ and $a_i$ for the stakes dictated by $S_-$ and $S_+$ when the counter is at~$i$, and also specify 
 $m_i$ and $n_i$ by the same definition. Our task is to show~$(a,b,m,n) \in$ \abmnkpmac.

Here we prove a useful property of~$(S_-,S_+)$: under gameplay governed by this pair, $\vert X_n \vert \to \infty$ is almost certain.
  
 \begin{proposition}\label{p.nashescape}
 For $(S_-,S_+) \in \mc{N}_{\kappa,\rho} \cap \tis^2$ and $i \in \Z$, $\pgameplay{S_-}{S_+}{i}(E) = 1$. 
 \end{proposition}

Recall the payoff notation~(\ref{e.receipt}).   A strategy pair $(S_-,S_+) \in \mc{S}^2$ is said to have {\em \fmc} if neither $\E^k_{S_-,S_+}[P_-]$ nor $\E^k_{S_-,S_+}[P_+]$ equals minus infinity, for any $k \in \Z$. 

Let $(S_-,S_+) \in \tis^2$.  Denote by $b_i$ and $a_i$ the stakes offered by $S_-$ and $S_+$, respectively, when the counter is at site $i \in \Z$
(without supposing $(S_-,S_+) \in \mc{N}_{\kappa,\rho}$).
The {\em idle zone} $\mc{I}$\hfff{idle} is set equal to
$$
\mc{I} \, = \, \big\{ j \in \Z : a_j = b_j = 0 \big\} .
$$
\begin{lemma}\label{l.idlezone}
Suppose that  $(S_-,S_+) \in \tis^2$ is such that $\mc{I}$ is non-empty.   For $k \in \Z$, consider the counter evolution $X:\N \to \Z$ under 
$\pgameplay{S_-}{S_+}{k}$. For  given $i \in \nwithzero$,
condition on  $X_i$ being a given element of $\mc{I}$.
(If $i$ equals zero, suppose that $k \in \mc{I}$.) Let $j$ be the first time after~$i$ for which $X_j  \not\in \mc{I}$. Then the conditional law of $X: \llbracket i,j\rrbracket \to \mathbb{Z}$ equals simple random walk begun at the given value $X_i$ and stopped on leaving~$\mc{I}$.
\end{lemma}
{\bf Derivation.} At each turn whose index lies in $\llbracket i,j-1 \rrbracket$, the counter lies in the idle zone and no stakes are offered.
The counter thus evolves as a symmetric simple random walk: on flip moves, by definition; on stake moves, by the zero-stake rule given in Section~\ref{s.tlp}. \qed

 An element of $\mc{S}_0^2$ is non-zero when at least one of its components is not zero at some vertex. 
\begin{proposition}\label{p.finitemeancosts}
Let $(S_-,S_+) \in \tis^2$ be non-zero, with \fmc.
Then escape is almost  certain: 
$\pgameplay{S_-}{S_+}{k}(E)=1$
 for $k \in \Z$.
\end{proposition}
{\bf Derivation.} Suppose on the contrary that $\pgameplay{S_-}{S_+}{k}(E^c) > 0$ for some $k \in \Z$.
Find $\ell \in \Z$ such that  it is with positive probability that the process $X$ under the law 
$\pgameplay{S_-}{S_+}{k}$
visits $\ell$ infinitely often. If $a_\ell + b_\ell >0$, then one or other of the players will incur mean infinite running cost due to stakes offered at site~$\ell$.
If $a_\ell = b_\ell =0$, let $I$ be an interval that is maximal under inclusion among those contained in the idle zone~$\mc{I}$ and containing~$\ell$. Since  $(S_-,S_+)$ is non-zero, we may select $j \in \Z \setminus I$ to be adjacent to an element of~$I$. By Lemma~\ref{l.idlezone}, each visit by~$X$ to $\ell$ leads with probability $2^{-\vert \ell - j \vert}$ to a visit to~$j$ after a further $\vert \ell - j  \vert$ turns of the game. So the mean number of visits to~$j \not\in \mc{I}$ is infinite. At least one player incurs infinite running cost as a result of these visits, contrary to hypothesis. \qed

For $S \in \tis$, we write  $\macleft(S) \in \Z \cup \{ -\infty\} \cup \{\infty\}$ 
and  $\macright(S) \in \Z \cup \{ -\infty\} \cup \{\infty\}$ 
for the infimum and supremum of the set  $\{ i \in \Z : S(i,1) > 0 \}$.  
The strategy $S$ is said to be {\em wide}\hfff{wide} if $\macleft(S) = -\infty$ and $\macright(S) = \infty$; if $S$ is not wide, it is {\em narrow}.

When a pair of narrow strategies is used, a player may secure victory by adding small stakes on the side where she leads. And if a wide strategy is played against a narrow one, the wide-staking player may harmlessly cut costs by lowering stakes in the infinite region where she offers a positive stake unopposed. We now specify {\em rocket} and {\em drag} stake-changing operations that act as tools for players with these respective needs.

For $\psi \in (0,1)$, the right $\psi$-rocket $\rocket_\psi^{i\rightarrow}$ at $i \in \Z$  is the element of $\tis$ given by
$$
 \rocket_\psi^{i\rightarrow}(j) \, = \, \psi^{j-i+1} {\bf 1}_{j \geq i} \, \, \, , \, \, \, j \in \Z \, ,
$$
while the left  $\psi$-rocket $\rocket_\psi^{\leftarrow i}$ at $i \in \Z$  is the element of $\tis$ given by
$$
 \rocket_\psi^{\leftarrow i}(j) \, = \, \psi^{i-j+1} {\bf 1}_{j \leq i} \, \, \, , \, \, \, j \in \Z \, .
$$

The right drag  at $i \in \Z$  is  the map  $\dragmapright$   that sends $q  \in \tis$ to 
$$
\Z \to (0,\infty ): j \to \, \,  \begin{cases}
  \, \qhalf  &  \text{if $j \geq i$} \\
 \, q_j  &  \text{if $j < i$} \, ,
\end{cases}
$$
and the  left drag  $\dragmapleft$ sends $q  \in \tis$ to 
$$
\Z \to (0,\infty ): j \to \, \,  \begin{cases}
  \,  \qhalf &  \text{if $j \leq i$} \\
 \, q_j  &  \text{if $j > i$} \, .
\end{cases}
$$

\begin{lemma}\label{l.boostdrag}
Let $(S_-,S_+) \in \tis^2$. 
\begin{enumerate}
\item Suppose that the quantities $\macright(S_-)$ and $\macright(S_+)$ are finite. Let $i \in \Z$ exceed both, and let $\psi \in \big( \tfrac{1-\kappa}{1+\kappa},1 \big)$. 
Choose $k \in \N$ so that
\begin{equation}\label{e.kappahyp}
 \Big( \frac{1-\kappa}{1+\kappa} \Big)^{k+1} \big( m_\infty - m_* \big)
 + \bigg( \psi^k +  \Big( \frac{1-\kappa}{1+\kappa} \Big)^{k+1} \bigg)  \kappa^{-1} \bigg( \frac{\psi}{1 - \psi} + \frac{1-\kappa}{(1+\kappa)\big(1 - \frac{1-\kappa}{(1+\kappa)\psi}\big)} \bigg)
\end{equation}
 is strictly less than $m_{-\infty,\infty}$.
Then $\egameplay{S_-}{\rocket_\psi^{i\rightarrow}}{i+k}[P_+] > \egameplay{S_-}{S_+}{i+k}[P_+]$. 
\item  Suppose that $\macright(S_+) = \infty$ and $\macright(S_-) < \infty$. Let $i \in \Z$ satisfy  $i > \macright(S_-)$ and  $S_+(i,1) > 0$. Then  $\egameplay{S_-}{\drag^{i\rightarrow}(S_+)}{i}[P_+] > \egameplay{S_-}{S_+}{i}[P_+]$.
\item  If $\macleft(S_-)$ and $\macleft(S_+)$ exceed $-\infty$ and $i \in \Z$ is \lessmin, then, provided that the quantity given by replacing $m_\infty - m_*$ by $n_{-\infty} - n_*$ in~(\ref{e.kappahyp}) is strictly less than $n_{\infty,-\infty}$, we have that $\egameplay{\rocket_\psi^{\leftarrow i}(S_-)}{S_+}{i-k}[P_-] > \egameplay{S_-}{S_+}{i-k}[P_-]$.
\item  If $\macleft(S_-) = -\infty$, $\macleft(S_+) > - \infty$ and $i \in \Z$ is such that $i < \macleft(S_+)$ and $S_-(i,1) > 0$, then  $\dragbound$.
\end{enumerate}
\end{lemma}
{\bf Proof: (1,3).} We prove only (1), since (3) has the same proof in essence. Let $Z:\N \to \Z$ denote simple random walk {\rm SRW}$\big(\tfrac{1+\kappa}{2}\big)$ with $Z(0)= i \in \Z$ (and the indicated right-move probability) under the law $\pimac$. 
Let $\#_j(Z)$ denote the cardinality of the set of visits made by $Z$ to $j \in \Z$. It is readily seen that
$$
 \E_i \big[ \#_j(Z) \big] \, \, = \, \,   \begin{cases} \, \kappa^{-1} & \textrm{for $i \geq j$} \, , \\ \, \kappa^{-1} \big( \tfrac{1 - \kappa}{1+\kappa} \big)^{j-i} & \textrm{for $i < j$} \, . \end{cases}
$$
Under the strategy pair  $\big(S_-,\rocket_\psi^{i\rightarrow}\big)$, Mina offers no stake at sites at or to the right of $i$, while Maxine always offers some positive stake at such locations. 
The counter trajectory under $\pgameplay{S_-}{\rocket_\psi^{i\rightarrow}}{i+k}$ stopped at $\tau_{i-1}$ thus has the law of ${\rm SRW}\big(\tfrac{1+\kappa}{2}\big)$ begun at~$i+k$ and stopped on arrival at $i-1$ (at a time that may be infinite).

Note that
\begin{eqnarray}
\egameplay{S_-}{\rocket_\psi^{i\rightarrow}}{i+k} \big[ C_+[0,\tau_{i-1}) \big] & = & \sum_{j=i}^\infty \E_{i+k} \big[ \#^{[0,\tau_{i-1})}_j(Z) \big] \psi^{j-i+1} \nonumber \\
 & \leq & \sum_{j \in \Z} \E_{i+k} \big[ \#_j(Z) \big] \psi^{j-i+1}  \nonumber \ \\
 & = & \sum_{j = i+k}^\infty \kappa^{-1} \psi^{j-i+1} \, + \,     \sum_{j = -\infty}^{ i+k-1} \kappa^{-1} \big( \tfrac{1-\kappa}{1+\kappa} \big)^{i+k-j} \psi^{j-i+1}  \nonumber  \\
  & = & \psi^k \kappa^{-1} \bigg(  \frac{\psi}{1-\psi} +       \frac{1-\kappa}{(1+\kappa)\big(1 - \frac{1-\kappa}{(1+\kappa)\psi} \big)} \bigg) \, ,  \label{e.rocketearly}
\end{eqnarray}
and that
\begin{eqnarray*}
\egameplay{S_-}{\rocket_\psi^{i\rightarrow}}{i+k} \big[ C_+[\tau_{i-1},\infty) \big] & = & 
\pgameplay{S_-}{\rocket_\psi^{i\rightarrow}}{i+k} \big( \tau_{i-1} < \infty \big) 
\egameplay{S_-}{\rocket_\psi^{i\rightarrow}}{i-1} \big[ C_+ \big] \, .
\end{eqnarray*}
Note that $\pgameplay{S_-}{\rocket_\psi^{i\rightarrow}}{i+k} \big( \tau_{i-1} < \infty \big) = \big( \tfrac{1-\kappa}{1+\kappa} \big)^{k+1}$ and that
\begin{eqnarray*}
\egameplay{S_-}{\rocket_\psi^{i\rightarrow}}{i-1} \big[ C_+ \big] & \leq & \egameplay{S_-}{\rocket_\psi^{i\rightarrow}}{i} \big[ C_+ \big] \\
& = & \sum_{j = i}^\infty \kappa^{-1} \psi^{j-i+1} \, + \,     \sum_{j = -\infty}^{ i-1} \kappa^{-1} \big( \tfrac{1-\kappa}{1+\kappa} \big)^{i+k-j} \psi^{j-i+1} \\
  & = &  \kappa^{-1} \bigg(  \frac{\psi}{1-\psi} +       \frac{1-\kappa}{(1+\kappa)\big(1 - \frac{1-\kappa}{(1+\kappa)\psi} \big)} \bigg) \, ,
\end{eqnarray*}
so that 
\begin{equation}\label{e.rocketlate}
\egameplay{S_-}{\rocket_\psi^{i\rightarrow}}{i+k} \big[ C_+[\tau_{i-1},\infty) \big] \, \leq \, \big( \tfrac{1-\kappa}{1+\kappa} \big)^{k+1} \kappa^{-1} \bigg(  \frac{\psi}{1-\psi} +       \frac{1-\kappa}{(1+\kappa)\big(1 - \frac{1-\kappa}{(1+\kappa)\psi} \big)} \bigg) \, .
\end{equation}
From~(\ref{e.rocketearly}) and~(\ref{e.rocketlate}),   we find that
\begin{equation}\label{e.rocketub}
\egameplay{S_-}{\rocket_\psi^{i\rightarrow}}{i+k} \big[ C_+ \big] \, \leq \,  \Big( \psi^k + \big( \tfrac{1-\kappa}{1+\kappa} \big)^{k+1} \Big) \kappa^{-1} \bigg(  \frac{\psi}{1-\psi} +       \frac{1-\kappa}{(1+\kappa)\big(1 - \frac{1-\kappa}{(1+\kappa)\psi} \big)} \bigg) \, .
\end{equation}
Since 
$\pgameplay{S_-}{\rocket_\psi^{i\rightarrow}}{i+k} \big( \tau_{i-1} = \infty \big) = 1 -  \big( \tfrac{1-\kappa}{1+\kappa} \big)^{k+1}$
and $m_* \leq m_{-\infty}$,
$$
\egameplay{S_-}{\rocket_\psi^{i\rightarrow}}{i+k} \big[ T_+ \big] \, \geq \,  \Big( 1 -  \big( \tfrac{1-\kappa}{1+\kappa} \big)^{k+1} \Big) m_\infty + \big( \tfrac{1-\kappa}{1+\kappa} \big)^{k+1}m_* \, .
$$
Since 
$\pgameplay{S_-}{S_+}{i+k}(E_+) = 0$ and $m_* \leq m_{-\infty}$, we have that $\egameplay{S_-}{S_+}{i+k}[T_+] \leq m_{-\infty}$.
We write
\begin{eqnarray*}
 & & \egameplay{S_-}{\rocket_\psi^{i\rightarrow}}{i+k} \big[ P_+ \big] -  \egameplay{S_-}{S_+}{i+k} \big[ P_+ \big] \\
 & = &
\Big( \egameplay{S_-}{\rocket_\psi^{i\rightarrow}}{i+k} \big[ T_+ \big] - 
\egameplay{S_-}{S_+}{i+k} \big[ T_+ \big] \Big) - 
\Big(\egameplay{S_-}{\rocket_\psi^{i\rightarrow}}{i+k} \big[ C_+ \big]  - \egameplay{S_-}{S_+}{i+k} \big[ C_+ \big] \Big)
\end{eqnarray*}
and note that first bracketed right-hand term is at least 
$$
 \Big( 1 -  \big( \tfrac{1-\kappa}{1+\kappa} \big)^{k+1} \Big) m_\infty + \big( \tfrac{1-\kappa}{1+\kappa} \big)^{k+1}m_* - m_{-\infty} \, ,
$$ 
while the second is at most the right-hand side of~(\ref{e.rocketub}). Hence, the hypothesis on~$k$ expressed in terms of~(\ref{e.kappahyp}) implies that 
$\egameplay{S_-}{\rocket_\psi^{i\rightarrow}}{i+k} \big[ P_+ \big] -  \egameplay{S_-}{S_+}{i+k} \big[ P_+ \big]$ is strictly positive, as we seek to show in proving Lemma~\ref{l.boostdrag}(1). 

{\bf Derivation: (2,4).} We derive (2), (4) being symmetrically obtained. The switch from $(S_-,S_+)$ to $(S_-,\drag^{i\rightarrow}(S_+))$ does not change the law of gameplay, because it merely
causes Maxine to decrease, by a factor of one-half, certain positive stakes on occasions when Mina offers no stake.  The switch thus saves on running cost for Maxine while leaving unchanged her terminal receipt. \qed

\begin{definition}
To $(S_-,S_+) \in \tis^2$, associate $(b,a):\Z \to [0,\infty)$ as usual. 
\begin{enumerate}
\item Let $S_-' \in \tis$ be associated to $b':Z \to [0,\infty)$. If $b'_i \geq b_i$ for all $i \in Z$, then $(S_-',S_+)$ is called a left strengthening of $(S_-,S_+)$.
\item Now let  $S_+' \in \tis$ be associated to $a':Z \to [0,\infty)$. If $a'_i \geq a_i$ for all $i \in Z$, then $(S_-,S_+')$ is called a right strengthening of $(S_-,S_+)$.  
\end{enumerate} 
When the assumed bounds are reversed, we speak of a left or right weakening.
\end{definition}
The straightforward proof of the next fact is omitted.
\begin{lemma}\label{l.monotonicity}
Let $(S_-',S_+)$  be a left strengthening of~$(S_-,S_+)$.
For $i \in \Z$, there is a coupling of gameplays $X,X':\N \to \Z$ under $\pgameplay{S_-}{S_+}{i}$ such that $X'(j) \leq X(j)$ for $j \in \N$ almost surely. Couplings with the evidently needed direction for the bounds exist for each of the three other variations. 
\end{lemma}
Recall that $(b_i,a_i)$ denotes the stake-pair offered at $i \in \Z$ under given  $(S_-,S_+) \in \mc{S}_0^2$.
\begin{lemma}\label{l.znn}
\leavevmode
\begin{enumerate}
\item Any element of $\mc{N}_{\kappa,\rho}$ has \fmc.
\item 
  If $(S_-,S_+) \in \mc{S}_0^2$ satisfies $\macleft(S_-) > -\infty$ and $\macleft(S_+) = -\infty$, let
 $i \in \Z$ satisfy $a_i > 0$ and $b_j =0$ for $j \in (-\infty,i-1 \rrbracket$.
  Then $\PP^i_{S_-,S_+}(E_-)$ equals zero.
\item[$(3)$] 
If $(S_-,S_+) \in \tis^2$ is an element of~$\mc{N}_{\kappa,\rho}$ then  $S_-$ and $S_+$ are wide.
\end{enumerate}
\end{lemma}
In the ensuing proof and later, the identically zero strategy is denoted by $0$. 

{\bf Derivation of Lemma~\ref{l.znn}(1).} 
 For $(S_-,S_+) \in \mc{N}_{\kappa,\rho}$ and $i \in \Z$,
 $\egameplay{S_-}{S_+}{i} [P_+] \geq  \egameplay{S_-}{0}{i}  [P_+] \geq  \min \{ m_{-\infty},m_\infty,m_* \} = m_* > -\infty$, the respective bounds due to  $(S_-,S_+) \in \mc{N}_{\kappa,\rho}$; absence of running cost for Maxine implying that $P_-$ is some average of the possible terminal receipt values $m_{-\infty}$, $m_\infty$ and $m_*$; and assumption on $m_*$.   
  Likewise,   $\egameplay{S_-}{S_+}{i}[  P_-] > -\infty$.

{\bf Proof: (2).} It is enough to argue that if $X$ under $\PP^i_{S_-,S_+}$ visits $i-1$, then its return to $i$ is assured. Consider $X$ under~$\PP^i_{S_-,S_+}$  from the time of a first visit to $i-1$ until such a return is made (if at all).
Since $S_-$ is zero on $j \in (-\infty,i-1 \rrbracket$, this subtrajectory of $X$ has the law of $X$ under $\PP^{i-1}_{0,S_+}$ stopped at $i$. Since $(0,S_+)$ is a right strengthening of $(0,0)$, and $X$ under $\PP^{i-1}_{0,0}$, being a symmetric simple random walk, necessarily visits $i$, Lemma~\ref{l.monotonicity} implies that the subtrajectory will reach $i$. This confirms the sought statement. \qed

 {\bf Derivation: (3).} We argue by contradiction and suppose without loss of generality that $S_-$ is narrow.
  Either $\macleft(S_-) > -\infty$ or $\macright(S_-) < \infty$.
 
 Suppose that $\macright(S_-) < \infty$. If $\macright(S_+) < \infty$, then Lemma~\ref{l.boostdrag}(1) provides $\hat{S}_+$ and $i \in \Z$ such that
 $\egameplay{S_-}{\hat{S}_+}{i}[P_+] > \egameplay{S_-}{S_+}{i}[P_+]$. 
 If $\macright(S_+) = \infty$, then Lemma~\ref{l.boostdrag}(2) does so.
 Suppose instead that  $\macleft(S_-) > -\infty$. If $\macleft(S_+) > -\infty$, then Lemma~\ref{l.boostdrag}(3) furnishes  $\hat{S}_-$ for Mina and $i \in \Z$ 
 for which $\egameplay{\hat{S}_-}{S_+}{i}[P_-] > \egameplay{S_-}{S_+}{i}[P_-]$ holds.
 
 In the remaining case,  $\macleft(S_-) > -\infty$ and $\macleft(S_+) = -\infty$. 
 The pair $(S_-,S_+) \in \tis^2 \cap \mc{N}_{\kappa,\rho}$ is non-zero, because $S_+$ is; it has \fmcspace by Lemma~\ref{l.znn}(1).
 Thus $\pgameplay{S_-}{S_+}{i}(E^c) = 0$ by Proposition~\ref{p.finitemeancosts}. Select $i \in \Z$ such that $a_i > 0$ and $b_j = 0$ for $j \in (-\infty, i \rrbracket$.
 Lemma~\ref{l.znn}(2) implies that
  $\pgameplay{S_-}{S_+}{i}(E_-) = 0$. Thus,  $\pgameplay{S_-}{S_+}{i}(E_+) = 1$, so that $T_+$ equals $m_\infty$ almost surely. If Maxine plays a strategy $\hat{S}_+$ formed from $S_+$ by reducing the stake she offers at $i$ by a factor of one-half, then gameplay $X:\N \to \Z$ is equal in law  under $\pgameplay{S_-}{S_+}{i}$ and $\pgameplay{S_-}{\hat{S}_+}{i}$; $T_+ = m_\infty$ almost surely under both laws; but Maxine's running cost is almost surely less under $\pgameplay{S_-}{\hat{S}_+}{i}$ than it is under $\pgameplay{S_-}{S_+}{i}$, because the first cost, incurred at site~$i$, is lower. 
  Thus,
   $\egameplay{S_-}{\hat{S}_+}{i}[P_+] > \egameplay{S_-}{S_+}{i}[P_+]$.

 We have obtained a contradiction to  $(S_-,S_+) \in \mc{N}_{\kappa,\rho}$ in each case we considered. This completes the proof of Lemma~\ref{l.znn}(3). \qed

{\bf Proof of Proposition~\ref{p.nashescape}.} This result follows from Proposition~\ref{p.finitemeancosts} and Lemma~\ref{l.znn}(1,2). \qed

 \section{A Nash component wins against zero}

Suppose that  Mina plays a time-invariant strategy $S_- \in \mc{S}_0$ that forms part of a Nash equilibrium $(S_-,S_+) \in \mc{N}_{\kappa,\rho}$,
  in a game in which Maxine offers no opposition, playing the zero-stake strategy.
Here we prove the next result, which  asserts,  plausibly enough, that Mina wins in the sense that
 $\PP^i_{(S_-,0)}(E_-) = 1$, no matter the value of the starting location~$X(0) = i \in \Z$.

 \begin{proposition}\label{p.leftescape}
 Let $(S_-,S_+) \in \mc{N}_{\kappa,\rho}$ with $S_- \in \mc{S}_0$. Then $\PP^i_{(S_-,0)}(E_-) = 1$ holds for all $i \in \Z$. 
 \end{proposition}

 The presence of flip moves, when $\kappa \in (0,1)$, makes the proposition 
  non-trivial, as we now explain.
 In the setup in question, $S_-$ is known to be wide by Lemma~\ref{l.znn}(3); so Mina offers positive stakes at an infinite set~$K$ of integer sites. When $\kappa$ equals one (as it is in~\cite{LostPennies}), so that every move is stake, this is enough to reach the desired conclusion that left escape~$E_-$ is almost certain starting from given $i \in \Z$. Indeed, when $X$ visits $K$, a left move is assured; while at sites in $\Z \setminus K$, no stakes are offered by either player, and the next move has equal chance of being left or right, according to the rule for zero stakes given in Section~\ref{s.tlp}. It is easily seen that this dynamics forces the counter leftward, through a sequence of one-way locks. However, when $\kappa \in (0,1)$, flip moves occur with probability $1-\kappa$; so, when $X$ visits $K$, the next move is left with probability $(1+\kappa)/2$.  The counter thus evolves as a symmetric simple random walk on~$\Z \setminus K$, with moves biased to the left by a uniform amount on visits to~$K$. Although~$K$ is infinite (since $S_-$ is wide), this set could in principle be arbitrarily sparse; in which case, this dynamics will not realize left escape $E_-$
 for some (or indeed all) starting points.  
 
  We see then that, to derive Proposition~\ref{p.leftescape}, we must harness the hypothesis $(S_-,S_+) \in \mc{N}_{\kappa,\rho}$  in a stronger form than the mere inference that $S_-$ is wide. 
 To survey the proof, we first mention that it is enough to reach the weaker conclusion that $\PP^i_{(S_-,0)}(E_-) \longrightarrow 1$ as $i \to -\infty$ because, as we will see in proving the next stated Lemma~\ref{l.leftescape}(1), it is simple to conclude as desired from this inference. We will then suppose that this weaker conclusion is false and contradict the hypothesis of Proposition~\ref{p.leftescape}. Lemma~\ref{l.leftescape}(2)
 shows that  $\PP^i_{(S_-,0)}(E_-) \centernot\longrightarrow 1$ as $i \to -\infty$ in fact implies that left escape~$E_-$ never occurs. This information will enable an argument that $(S_-,S_+)$ is not a Nash equilibrium, so that the desired contradiction to the hypotheses of  Proposition~\ref{p.leftescape} may be obtained. 
 \begin{lemma}\label{l.leftescape}
 Let $(S_-,S_+) \in \mc{N}_{\kappa,\rho}$ with $S_- \in \mc{S}_0$.
 \begin{enumerate}
 \item 
  If the sequence $\big\{ \PP^i_{(S_-,0)}(E_-): i \in \Z \big\}$ converges to the value one in the limit $i \to -\infty$, then $\PP^i_{(S_-,0)}(E_-) = 1$ for all $i \in \Z$. 
 \item If this convergence does not hold,  then  $\PP^i_{(S_-,0)}(E_-)$ equals zero  for all $i \in \Z$. 
 \end{enumerate}
 \end{lemma}
 {\bf Proof: (1).} Let $\tau_j = \min \big\{ k \in \N: X_k = j \big\}$. Since $(S_-,0)$ is a left strengthening of $(0,0)$, and~$X$ under $\PP^i_{(0,0)}$ is symmetric simple random walk, Lemma~\ref{l.monotonicity}  implies that
 $\tau_j < \infty$ occurs almost surely under $\PP^i_{(S_-,0)}$ whenever $j \leq i$. By hypothesis, we may find for any $\e > 0$ a sequence $j_k \to -\infty$ as $k \to \infty$
 such that $\PP^{j_k}_{(S_-,0)}(E_-) \geq 1 - \e$. Since $\tau_{j_k} < \infty$ is assured to occur under $\PP^i_{(S_-,0)}$, and~$X$ viewed from time $\tau_{j_k}$ onwards realizes $E_-$
 with probability at least $1-\e$ by the strong Markov property,  $\PP^i_{(S_-,0)}(E_-) \geq 1 - \e$. Since $\e > 0$ is arbitrary, we obtain Lemma~\ref{l.leftescape}(1).
 
 {\bf (2).} Let $i \in \N$ be given. The hypothesised lack of convergence permits us to find $\e > 0$ and a strictly decreasing sequence $\big\{ v_j : j \in \N_+ \big\}$\hfff{withandwithout} such that $v_1 < i$ and $\PP^{v_j}_{(S_-,0)}(E_-) \leq 1 - \e$. By the definition of~$E_-$, we may choose $u_j < v_j$ such that  $\PP^{v_j}_{(S_-,0)}(\tau_{u_j} = \infty) \geq \e/2$. By thinning the sequence of $v_j$ as needed, we may further suppose that $v_{j+1} \leq u_j$. We also set $v_0 = i$.
 
 View the evolving trajectory $X:\N \to \Z$ under $\PP^i_{(S_-,0)}$. Think of an experiment in which
 time passes discretely: $0,1,2,\cdots$. 
 If $X$ reaches $v_i$ but not $v_{i+1}$, shout `stop!' between times $i$ and $i+1$. If time $i \geq 1$ arrives without `stop!' being shouted, 
 then it will be shouted between times $i$ and $i+1$ with conditional probability at least $\e/2$: indeed, since `stop!' has not been shouted by time~$i$, $X$ has reached $v_i$; if it does not then reach $u_i$, `stop!' will be shouted between times~$i$ and~$i+1$; but if $X$ reaches $u_i$, it will, by the strong Markov property, fail to reach $v_{i+1}$ with conditional probability at least $\e/2$, in which event, `stop!' will be shouted between times $i$ and $i+1$. In this way, the index~$I$ such that `stop!' is shouted between times $I$ and $I+1$ under $\PP^i_{(S_-,0)}$ is stochastically dominated by a geometric random variable $G \geq 1$ of success parameter~$\e/2$. If left escape~$E_-$ occurs, `stop!' is never shouted. This event forces the random index~$I$ to be infinite, which is a singular event. Thus,  $\PP^i_{(S_-,0)}(E_-) = 0$. \qed

{ \bf Proof of Proposition~\ref{p.leftescape}.}
We will argue that, when $(S_-,S_+) \in \mc{N}_{\kappa,\rho}$ with $S_- \in \mc{S}_0$ satisfies  $\PP^i_{(S_-,0)}(E_-) = 0$
for all $i \in \Z$, then  $(S_-,S_+) \not\in \mc{N}_{\kappa,\rho}$.  In light of Lemma~\ref{l.leftescape},
 this is enough to prove the proposition by contradicition.
 
 We will in fact prove the stronger assertion that, when $(S_-,S_+) \in \mc{N}_{\kappa,\rho}$ with $S_- \in \mc{S}_0$ satisfies  $\PP^i_{(S_-,0)}(E_-) = 0$ for some $i \in \Z$ such that Mina's first-turn stake~$S_-(0)$ is strictly positive when $X_0 = i$, then $(S_-,S_+) \not\in \mc{N}_{\kappa,\rho}$.
 Fixing such an $i$, we will show that
 \begin{equation}\label{e.minainterest}
  \E^i_{(0,S_+)} [P_-] >  \E^i_{(S_-,S_+)} [P_-]  \, :
 \end{equation}
  it is in Mina's interests to play the zero strategy, rather than $S_-$, against Maxine's~$S_+$, when play starts at~$i$.
  Naturally,~(\ref{e.minainterest}) implies that $(S_-,S_+) \not\in \mc{N}_{\kappa,\rho}$, so proving~(\ref{e.minainterest}) is enough.
  
Preparing  to show~(\ref{e.minainterest}), note that
  \begin{equation}\label{e.monoinference}
  \PP^i_{(S_-,S_+)}(E_-) =  \PP^i_{(0,S_+)}(E_-) = 0 \, .
  \end{equation}
  Indeed, $(S_-,0) \to (S_-,S_+)$ is a right strengthening and  $(S_-,S_+) \to (0,S_+)$ is a left weakening, so~(\ref{e.monoinference}) follows from~$\PP^i_{(S_-,0)}(E_-) = 0$ and Lemma~\ref{l.monotonicity}.
  
  Why may we expect~(\ref{e.minainterest}) to hold? In other words, 
  why would Mina switch from $S_-$ to $0$ against~$S_+$? That  $\PP^i_{(S_-,S_+)}(E_-)$ is zero makes Mina's motivation simple: $S_-$ is not working out for her, because her victory $E_-$ never happens. By switching to $0$, she will save on running costs. As for terminal receipts, these are split between non-escape~$E^c$ and right escape~$E_+$ when she plays~$S_-$. By playing~$0$ instead, Mina will cease to exert any left pressure, so, in an instance of right strengthening and monotonicity.  any change to this split will take the form of a rightward move of probability mass from $E^c$ to $E_+$.
   But that would help Mina, because $E^c$ is the worse outcome for her in the sense that $n_* \leq n_\infty$.
  To record these inferences symbolically,
 \begin{eqnarray*}
   \E^i_{(0,S_+)}[P_-] =  \E^i_{(0,S_+)}[T_-] & = &  \PP^i_{(0,S_+)}(E^c) n_* +  \PP^i_{(0,S_+)}(E_+) n_\infty \\ 
   & \geq & \PP^i_{(S_-,S_+)}(E^c) n_* +  \PP^i_{(S_-,S_+)}(E_+) n_\infty  
   =  \E^i_{(S_-,S_+)}[T_-] >  \E^i_{(S_-,S_+)}[P_-] \, ,
 \end{eqnarray*}
  where the first equality is due to absence of running cost for Mina when she plays zero; the second equality crucially invokes  $\PP^i_{(S_-,0)}(E_-) = 0$;
  the first inequality is due to the~(\ref{e.monoinference})-consequences  
  $$
  \PP^i_{(S_-,S_+)}(E^c) +  \PP^i_{(S_-,S_+)}(E_+) =  \PP^i_{(0,S_+)}(E^c) +  \PP^i_{(0,S_+)}(E_+) =1 \, ,
  $$
   and the monotonicity deduction
   $\PP^i_{(0,S_+)}(E_+) \geq \PP^i_{(S_-,S_+)}(E_+)$; the next equality depends on~(\ref{e.monoinference}) for $(S_-,S_+)$;
   and the strict inequality is due to the running cost $C_-$ in~(\ref{e.runningcosts}) being a sum of non-negative terms whose first, $S_-(0)$, is positive under $\PP^i_{(S_-,S_+)}$. We have proved~(\ref{e.minainterest}) and with it Proposition~\ref{p.leftescape}. \qed

\section{Positive stakes at Nash equilibrium}

Recall that to $(S_-,S_+) \in \tis^2$ Definition~\ref{d.quadruple}
associates  $\big\{ (a_i,b_i,m_i,n_i): i \in \Z \big\}$.
Here we show that when $(S_-,S_+)$ is Nash, stakes and $m$- and $n$-increments are positive.
\begin{proposition}\label{p.allpos}
Let $(S_-,S_+) \in \mc{S}_0^2 \cap \mc{N}_{\kappa,\rho}$. For all~$i \in \Z$,  $a_i > 0$, $b_i > 0$, $m_{i+1} > m_i$ and $n_i > n_{i+1}$.
\end{proposition}

Four lemmas lead to the proof.

\begin{lemma}\label{l.mnincdec}
Suppose that $(S_-,S_+) \in \mc{N}_{\kappa,\rho} \cap \tis^2$. Then $m_i \leq m_{i+1}$ and $n_{i+1} \leq n_i$ for $i \in \Z$.
\end{lemma}
{\bf Derivation.} Under
$\pgameplay{S_-}{S_+}{i}$, 
let $\sigma_{i+1} \in \nwozero \cup \{ \infty \}$ denote \stoppingtime. 
In the specification~(\ref{e.receipt}) of Maxine's net receipt $P_+$ as $T_+ - C_+$, the running cost $C_+$ may be written $C_+\llbracket 1,t \rrbracket$ and $C_+\llbracket t+1,\infty \rrbracket$
where Maxine's stakes up to the $t$\textsuperscript{th} turn enter as summands in the first term. 
Taking $t = \sigma_{i+1}$, we use the strong Markov property at time~$\sigma_{i+x1}$ and drop $C_+\llbracket 1,\sigma_{i+1} \rrbracket \geq 0$ to obtain 
$$
\egameplay{S_-}{S_+}{i} [P_+] \, \leq \,  \erhs \, .
$$
Here, the left-hand side equals $m_i$ by definition while the right-hand side is
$$
 m_{i+1} \pgameplay{S_-}{S_+}{i} \big( \sigma_{i+1} < \infty \big) + \macmid
 + m_* \pgameplay{S_-}{S_+}{i} \big( \sigma_{i+1} = \infty, E^c \big) \, .
 $$
 The third right-hand term vanishes by Proposition~\ref{p.nashescape}, so that $m_i$ is bounded above by a weighted average of $m_{-\infty}$ and $m_{i+1}$.
 We will find as desired that $m_i \leq m_{i+1}$ by showing $m_{-\infty} \leq m_{i+1}$.
 In this regard, we first {\em claim} that $\egameplay{S_-}{0}{i+1} [P_+] = m_{-\infty}$. To check this, note that Lemma~\ref{l.znn}(3) implies  that~$S_-$ is wide.
 We may now make use of Proposition~\ref{p.leftescape} to learn that  
 $E_-$, and thus also $T_+ = m_{-\infty}$, are $\pgameplay{S_-}{0}{i+1}$-almost certain. 
The absence of running costs for Maxine means that $P_+ = T_+$ under $\pgameplay{S_-}{0}{i+1}$. The claim obtained, we use it and $(S_-,S_+) \in \mc{N}_{\kappa,\rho}$ to find that
$m_{i+1} \, = \, \egameplay{S_-}{S_+}{i+1}[ P_+] \, \geq \, \egameplay{S_-}{0}{i+1} [P_+] = m_{-\infty}$, thereby confirming
  $m_i \leq m_{i+1}$. Omitting the similar proof that  $n_{i+1} \leq n_i$, we obtain Lemma~\ref{l.mnincdec}. \qed

\begin{lemma}\label{l.firstrearranged}
Let  $\big\{ (b_i,a_i): i \in \Z \big\} \in \mc{N}_{\kappa,\rho} \cap \tis^2$.  
Recall from Definition~\ref{d.quadruple} that $m_i$ equals Maxine's mean receipt when the counter starts at $i \in \Z$.
Suppose that $a_i + b_i > 0$. Then
\begin{equation}\label{e.firstrearranged}
 m_i \, = \, \Big( \kappa \tfrac{a_i^\rhomac}{a_i^\rhomac + b_i^\rhomac} +  \tfrac{1-\kappa}{2} \Big) m_{i+1} + \Big( \kappa \tfrac{b_i^\rhomac}{a_i^\rhomac + b_i^\rhomac} +  \tfrac{1-\kappa}{2} \Big) m_{i-1} \, - \, a_i \, .
\end{equation}
\end{lemma}
{\bf Proof.}
 Maxine will spend $a_i$ at the first turn; the move will be stake with probability~$\kappa$ and then she win it with conditional probability $\tfrac{a_i^\rhomac}{a_i^\rhomac  + b_i^\rhomac}$; if she does so, the counter will reach $i+1$,
  and her resulting conditional mean receipt will be $m_{i+1}$; 
  and this circumstance will equally arise if a fair coin lands heads on a flip move, with probability $(1-\kappa)/2$. Otherwise, Maxine's receipt will be $m_{i-1}$. Note that \tworatios of~(\ref{e.firstrearranged}) are well defined, because  $a_i + b_i > 0$. \qed

\begin{lemma}\label{l.condpositive}
Let $(S_-,S_+) \in \mc{N}_{\kappa,\rho} \cap \tis^2$, and let $i \in \Z$. Then $a_i > 0$ implies that $m_{i+1} > m_i$. And $b_i > 0$ implies that $n_{i-1} > n_i$.
\end{lemma}
{\bf Proof.} Lemma~\ref{l.firstrearranged} and $a_i > 0$ imply that $m_i < \max \{ m_{i-1},m_{i+1} \}$. But the maximum is attained by $m_{i+1}$ in view of Lemma~\ref{l.mnincdec}. 
The second assertion in the lemma is similarly obtained.
\qed

\begin{lemma}\label{l.abfacts}
Let $(S_-,S_+) \in \mc{S}_0^2 \cap \mc{N}_{\kappa,\rho}$. Then
\begin{enumerate}
\item $a_i > 0$ implies that $a_{i+1} + b_{i+1} > 0$.
\item $a_i > 0$ implies that $b_i > 0$. 
\item  $b_i > 0$ implies that $a_i > 0$. 
\end{enumerate}
\end{lemma}
{\bf Proof: (1).} If $a_{i+1} = b_{i+1} = 0$, then $m_i = (m_{i-1} + m_{i+1})/2$ by the zero-stakes fair-coin rule. But $a_i > 0$ implies that $m_{i+1} > m_i$ by Lemma~\ref{l.condpositive}.
A one-turn variation for Maxine, in which she stakes $0^+$ rather than $0$ with the counter at $i+1$,  would result in her mean receipt equalling $\tfrac{1-\kappa}{2}m_{i-1} + \tfrac{1+\kappa}{2} m_{i+1}$.
Since this strictly exceeds $m_i$, we learn that $(S_-,S_+) \not\in \mc{N}_{\kappa,\rho}$. Thus $a_i > 0$ is inconsistent with $a_{i+1} + b_{i+1} = 0$.

{\bf (2).} Suppose that $a_i > 0$ and $b_i = 0$. Let $S'_i$ denote the strategy for Mina formed from $S_-$ by replacing her stake at site~$i$ by $a_i/2$, so that it is reduced but remains positive. Gameplay under $(S_-,S_+)$ and under $(S'_-,S_+)$ are equal in law, because Mina will win every stake turn at site~$i$ in either case. Mina will save a positive amount on running cost whenever $X$ visits~$i$.
Thus,
$\E^i_{(S'_-,S_+)}[P_-] > \E^i_{(S_-,S_+)}[P_-]$, so that $(S_-,S_+) \not\in \mc{N}_{\kappa,\rho}$. This contradiction shows that $a_i > 0$ implies $b_i > 0$.

{\bf (3).} This argument is in essence identical to the preceding one. \qed

{\bf Proof of Proposition~\ref{p.allpos}.} By Lemma~\ref{l.znn}(3), $S_+$ is wide. By Lemma~\ref{l.abfacts}, $a_i > 0$ implies that $a_{i+1} > 0$. Hence, all coefficients $a_i$ are positive; by Lemma~\ref{l.abfacts}(2), so are all the $b_i$. By Lemma~\ref{l.condpositive}, the differences $m_{i,i+1}$ and  $n_{i+1,i}$ are also found to be positive.\qed

\section{The forward implication}

We are ready for the next derivation.  The argument follows the lines of the proof of~\cite[Theorem~$2.6(1)$]{LostPennies}, with a different approach used at the end to handle flip moves.
 
{\bf Proof of Theorem~\ref{t.nashabmn}(1).}  Suppose  that 
$(S_-,S_+) \in \mc{N}_{\kappa,\rho} \cap \tis^2$ for \tlpkpspace with boundary data $(m_{-\infty},m_\infty,n_{-\infty},n_\infty)$.  
Note that, in view of Proposition~\ref{p.allpos}, each $a_i$ and $b_i$, and each difference $m_{i,i+1}$ and $n_{i+1,i}$, is positive.

Equation \abmnmac$(1)$ is a rearrangement of the formula in Lemma~\ref{l.firstrearranged}, and \abmnmac$(2)$ is obtained similarly.

To derive \abmnmac$(3,4)$, recall that $S_-(i,j) = b_i$ and $S_+(i,j)=a_i$ for each $(i,j) \in \Z \times \nwozero$. 
For given $i \in \Z$, we will consider a perturbed strategy $\hat{S}_+ \in \mc{S}$ for Maxine in which only her first-turn stake is altered, and only then if the counter is at $i$. In this way,  $\hat{S}_+(j,k) = a_j$ for $j \in \Z$ and $k \geq 2$; and also for $k=1$ and $j \in \Z$, $j \not= i$. We let $\eta > -a_i$ be small in absolute value, and set $\hat{S}_+(1,i) = a_i + \eta$.

The {\em original} scenario refers to $\pgameplay{S_-}{S_+}{i}$, the law governing~$X:\nwithzero \to \Z$ given the initial condition $X_0 = i$ under the strategy pair $(S_-,S_+)$. The {\em altered} scenario refers to the same law, but now governed by the pair $(S_-,\hat{S}_+)$. Write $O_+ = \egameplay{S_-}{S_+}{i} [P_+]$ and $A_+ = \egameplay{S_-}{\hat{S}_+}{i} [P_+]$ for the mean payoffs to Maxine in the original and altered scenarios.    
Then\[
\begin{aligned}
O_+ &= 
\left( \kappa \frac{a_i^\rho}{a_i^\rho+b_i^\rho} + \frac{1-\kappa}{2} \right) m_{i+1} 
+ \left( \kappa \frac{b_i^\rho}{a_i^\rho+b_i^\rho} + \frac{1-\kappa}{2} \right) m_{i-1} 
- a_i \,, \, \, \, \, \textrm{and} \\[6pt]
A_+ &= 
\left( \kappa \frac{(a_i+\eta)^\rho}{(a_i+\eta)^\rho + b_i^\rho} + \frac{1-\kappa}{2} \right) m_{i+1} 
+ \left( \kappa \frac{b_i^\rho}{(a_i+\eta)^\rho + b_i^\rho} + \frac{1-\kappa}{2} \right) m_{i-1} 
- (a_i+\eta) \, .
\end{aligned}
\]
Hence,
\begin{equation}\label{e.aodifference}
A_+ - O_+ \, = \, \bigg( \frac{\rho \, a_i^{\rho -1} b_i^\rho}{(a_i^\rho+b_i^\rho)^2} \kappa \, m_{i-1,i+1}  - 1 \bigg) \cdot \eta \cdot \big( 1 + o(1) \big) \, ,
\end{equation}
where  $\vert \eta \vert \to 0$ for the $o(1)$ term.
Since $(S_-,S_+) \in \mc{N}_{\kappa,\rho}$, $A_+$ is at most $O_+$, for any value $\eta > - a_i$. Hence, the derivative in $\eta$ of $A_+ - O_+$  vanishes at zero, so that 
$\tfrac{\rho \, a_i^{\rho -1} b_i^\rho}{(a_i^\rho+b_i^\rho)^2} \kappa m_{i-1,i+1}   - 1 = 0$
or equivalently
\begin{equation}\label{e.bma}
 \rho \, a_i^{\rho -1} b_i^\rho  \kappa \, m_{i-1,i+1}  \,  = \,  \big( a_i^\rho + b_i^\rho \big)^2 \, .
\end{equation}
Now consider the same original scenario alongside a \nas in which it is Mina who employs a perturbed strategy $\hat{S}_-$ (as a function of given  $i \in \Z$). Similarly  as we have done, we choose $\eta > - b_i$,
and set $\hat{S}_-(j,k)$ equal to $b_j$ for  $j \in \Z$ and $k \geq 2$  or when $k=1$ and $j \in \Z \setminus \{ i \}$; and then we set $\hatseta$. 
Denote $O_-  = \egameplay{S_-}{S_+}{i} [P_-]$ and $A_-  =  \egameplay{\hat{S}_-}{S_+}{i} [P_-]$. We find  that
\[
\begin{aligned}
O_- &= 
\left( \kappa \frac{b_i^\rho}{a_i^\rho+b_i^\rho} + \frac{1-\kappa}{2} \right) n_{i-1} 
+ \left( \kappa \frac{a_i^\rho}{a_i^\rho+b_i^\rho} + \frac{1-\kappa}{2} \right) n_{i+1} 
- b_i \,, \, \, \, \, \textrm{and} \\[6pt]
A_- &= 
\left( \kappa \frac{(b_i+\eta)^\rho}{a_i^\rho+(b_i+\eta)^\rho} + \frac{1-\kappa}{2} \right) n_{i-1} 
+ \left( \kappa \frac{a_i^\rho}{a_i^\rho+(b_i+\eta)^\rho} + \frac{1-\kappa}{2} \right) n_{i+1} 
- (b_i+\eta) \, .
\end{aligned}
\]
Thus, similarly to~(\ref{e.aodifference}), 
\[
A_- - O_- \;=\; \Bigg( \frac{\rho \, a_i^\rho b_i^{\rho-1}}{(a_i^\rho+b_i^\rho)^2} \, \kappa \, n_{i+1,i-1}  - 1 \Bigg) \cdot \eta \cdot \big( 1 + o(1) \big) \, .
\]
The condition that $(S_-,S_+) \in \mc{N}_{\kappa,\rho}$ gives $O_- \geq A_-$, for any  $\eta > - b_i$. Thus,
\begin{equation}\label{e.anb}
\rho \, a_i^\rho b_i^{\rho-1} \, \kappa \, n_{i+1,i-1} \;=\; \big( a_i^\rho + b_i^\rho \big)^2 \, .
\end{equation}
The obtained equations~(\ref{e.bma}) and~(\ref{e.anb}) are \abmnmac$(3,4)$ with index $i$. 

We have shown that  $\big\{ (a_i,b_i,m_i,n_i): i \in \Z \big\}$
is an element of  \abmnkpmac. 
To finish the proof of Theorem~\ref{t.nashabmn}(1), it remains to confirm that the boundary values~(\ref{e.boundarydata}) are achieved.
We will prove that $\lim_{i \to \infty} m_{-i} = m_{-\infty}$; the three other limits are similarly obtained. 
The sequence $\big\{ m_{-i}: i \in \nwithzero \big\}$ decreases by Proposition~\ref{p.allpos} to a limit that we call $\mathfrak{m}_{-\infty}$.

By Definition~\ref{d.quadruple} and $(S_-,S_+) \in \mc{N}_{\kappa,\rho}$, $m_i = \pgameplay{S_-}{S_+}{i} [P_+] \geq \pgameplay{S_-}{0}{i} [P_+]$.
It is Proposition~\ref{p.leftescape} that now permits us to identify the right-hand term as being equal to~$m_{-\infty}$.
Hence, $\mathfrak{m}_{-\infty} \geq m_{-\infty}$; we wish to obtain the opposing inequality.
We 
take the mean of the equality $P_+ = T_+ - C_+$ in~(\ref{e.receipt}) and
 remove non-negative running costs~$C_+$ to find that $m_i \leq \mupperbound$  where we invoked Proposition~\ref{p.nashescape} to eliminate a non-escape~$E^c$ term. Thus $\mathfrak{m}_{-\infty} \leq m_{-\infty}$ provided that we show that  $\lim_{i \to -\infty} \pgameplay{S_-}{S_+}{i}(E_+)$ equals zero: far to the left,  Mina's victory is close to assured.  

 Let $k \in \Z$ denote the battlefield index of $(a,b,m,n) \in$ \abmnkpmacspace as specified in Definition~\ref{d.battlefield}. 
 Here we turn to the fixed-parameter asymptotic Theorem~\ref{t.abmn}.
 It would be of interest to harness this theorem\footnote{An application of Proposition~\ref{p.allpos} is technically needed to permit this use of Theorem~\ref{t.abmn}, because this proposition tells us that the right limit~$\mathfrak{m}_\infty$ strictly exceeds $\mathfrak{m}_{-\infty}$, so that the trivial zero \abmnmacspace  solution is eliminated from consideration, and $(a,b,m,n) \in$  \abmnkpmacspace is established.} to prove say a  $\sim$-asymptotic for the  decay of the probability $\pgameplay{S_-}{S_+}{i}(E_+)$ of `escape across the battlefield', but a rough leading-order estimate suffices for our purposes. From Theorem~\ref{t.abmn}, we need the simple inference, valid in each of the three treated $(\kappa,\rho)$-regimes, that $b_i \gg a_i$ as $i \to -\infty$, at a rate determined by $k-i$. Far to the left of the battlefield, Mina dominates the stakes and wins asymptotically all stake moves. Her turn victory probability tends to $\kappa + \tfrac{1}{2}(1-\kappa) = \tfrac{1}{2}(1+\kappa)$. Simple random walk with this left-move probability hits the point~$\ell$ steps to the right of its starting location with probability $\big( \tfrac{1-\kappa}{1+\kappa} \big)^\ell$ for $\ell \in \N$. Crudely absorbing the effect of discrepancy from the limiting move probability into a factor in the exponent, we infer that   $\pgameplay{S_-}{S_+}{i}(E_+) \leq \big( \tfrac{1-\kappa}{1+\kappa} \big)^{(k-i) (1-o(1))}$ where $o(1) \geq 0$ vanishes as $i \to -\infty$.
 Hence holds the bound  $\mathfrak{m}_{-\infty} \leq m_{-\infty}$ to which we reduced the proof of 
   Theorem~\ref{t.nashabmn}(1).
  \qed

 \section{The reverse implication}\label{s.reverse}

Here we prove Theorem~\ref{t.nashabmn}(2), the step at which the infinite-turn game is controlled by comparison with finite-trail counterparts.
Throughout,~$\big\{ (a_i,b_i,m_i,n_i) : i \in \Z \big\}$ denotes an element of \abmnkpmac,
with boundary data $(m_{-\infty},m_\infty,n_{-\infty},n_\infty)$ that satisfies~(\ref{e.quadruple}).
We define strategies $S_-,S_+ \in \mc{S}$ that offer $b$- and $a$-stake compatibly with the rule~(\ref{e.ba}). 

Since all counter moves are $\pm 1$, counter location is constrained by parity. First we denote   the set of space-time sites that are thus in principle accessible for gameplay $X:\nwithzero \to \Z$ under $\pgameplay{S_1}{S_2}{\imac}$ for some strategy pair $(S_1,S_2) \in \mc{S}^2$.
 Here and throughout the section, the symbol~$\imac$ is reserved for  denoting a choice of initial counter location. 
  
\begin{definition}\label{d.fpc}
For $\imac \in \Z$, the forward play-cone $F_\imac$ of $\imac$ is set equal to 
$$
F_\imac \, = \, \Big\{ \, (k,\ell) \in \Z \times \nwozero: \fpcrhs \, \Big\} \, .
$$

Let $S \in \mc{S}$ (and recall the formulation of the strategy space~$\mc{S}$ from Section~\ref{s.tlp}). A  {\em Mina deviation point} is an element $(q,\ell) \in F_\imac$ for which 
there exists a trajectory $\psi: \llbracket 0, \ell \rrbracket \to \Z$ with $\psi(0) = \imac$ and $\psi(\ell) = q$
such that $S(\psi) \neq b_q$.
Write  $\mathsf{D}_-(S,\imac) \subseteq F_\imac$ for the set of Mina deviation points. The strategy $S$ is  {\em deviating for Mina} if  $\mathsf{D}_-(S,\imac) \neq \emptyset$.
 A {\em Maxine deviation point} is an element $(q,\ell) \in F_\imac$ such that  $S(\psi) \not= a_q$ for some path~$\psi$ as above. Write $\mathsf{D}_+(S,\imac)$ for the set of these points; if  $\mathsf{D}_+(S,\imac) \not= \emptyset$, then $S$ is deviating for Maxine.  
\end{definition}

Mina deviation points~$(u,\ell)$ are instances in space-time at which at least one counter history leading to the point would prompt her to stake an amount other than $b_u$ against Maxine's~$a_u$. Such choices by Mina may be viewed as mistakes; to substantiate this notion, we wish to argue that Mina will receive a penalty in the sense of mean total receipt as a consequence of offering deviant stakes. The next two propositions offer results to this effect. 
The first concerns finite trail games and asserts that Mina will receive a penalty by playing the given deviating strategy $\sdev$ in any such game whose gameboard is broad enough to encompass a deviating move under~$\sdev$; moreover, the penalty is uniformly bounded below over such gameboards.

Write $P_-^{j,k}$ for Mina's total receipt in playing the trail game on $\llbracket -j-1,k+1\rrbracket$, the counter stopping on arrival at $-j-1$ or $k+1$ with terminal payments given by $(m_{-j-1},m_{k+1},n_{-j-1},n_{k+1})$.
\begin{proposition}\label{p.supjk}
Let 
$\imac \in
 \Z$ be given, and let $\sdev \in \mc{S}$ be deviating for Mina.  Suppose that $\PP_{\sdev,S_+}^\imac
(E) = 1$.
For any given  $(u,\ell) \in \mathsf{D}_-(\sdev,\imac)$,
$$
\sup \, 
 \E_{\sdev,S_+}^\imac
 [P_-^{j,k}]  \, < \, 
 \E_{S_-,S_+}^\imac
 [P_-]  \, ,
$$
with the supremum  taken over those $j,k \in \nwozero$  for which $u \in \llbracket -j+\ell,k-\ell \rrbracket$.
\end{proposition}
The second result expresses that a penalty is also suffered in the infinite trail game. In essence, this result captures the notion that $(S_-,S_+)$ is a Nash equilibrium and thus the content of Theorem~\ref{t.nashabmn}(2).
\begin{proposition}\label{p.comp.sminus}
Let $\imac \in
 \Z$, and let $\sdev \in \mc{S}$ be \devMina. Then
$$
 \E_{\sdev,S_+}^\imac
 [P_-]  <  \sminusrhs \, .
$$
\end{proposition}
This pair of propositions forms the backbone of the proof of Theorem~\ref{t.nashabmn}(2). They are simply the assertions made by~\cite[Propositions~$4.11$ and~$4.12$]{LostPennies}  in regard to Mina's deviation. 
(For a reason to be explained shortly,  Proposition~\ref{p.supjk} is phrased a little differently than~\cite[Proposition~$4.11(1)$]{LostPennies} and includes a new hypothesis.) 
Alongside symmetric assertions regarding Maxine's deviant play made in these results from~\cite{LostPennies} but omitted here, Theorem~\ref{t.nashabmn}(2) follows directly. Indeed, Mina's replacement of $S_-$ by another strategy~$S$ when playing against $S_+$ will either affect no change in her mean outcome---namely, $\egameplay{S}{S_+}{\imac} [P_-] = \egameplay{S_-}{S_+}{\imac} [P_-]$---if $S$ is not deviating; or a negative change, $\egameplay{S}{S_+}{\imac} [P_-] < \egameplay{S_-}{S_+}{\imac} [P_-]$, by Proposition~\ref{p.comp.sminus}. And of course likewise if Maxine is the one to deviate. 

The derivation of  Theorem~\ref{t.nashabmn}(2) thus substantially coincides with that of the counterpart Theorem~$2.6(2)$ in~\cite{LostPennies}. But one significant change is needed.

Our presentation of the proof of Theorem~\ref{t.nashabmn}(2)  is intended to be comprehensive in describing changes to the counterpart in~\cite[Section~$4.2$]{LostPennies}, and to offer a substantially complete conceptual guide to the proof while avoiding excessive repetition of~\cite{LostPennies}. We will begin by describing the more major change, which concerns the proof of Proposition~\ref{p.comp.sminus} and will entail presenting a further result, Proposition~\ref{p.altered}. We will describe why this result is needed and state it. An overview of the derivation at large will then be offered, in which some more minor changes to the proof in~\cite{LostPennies} will be noted. Then we will prove Proposition~\ref{p.altered}.

\subsection{The substantial new element, which handles possible non-escape} 
In the proof of~\cite[Proposition 4.11]{LostPennies}, counterpart 
to Proposition~\ref{p.comp.sminus}, 
the case~$\PP^\imac
_{\sdev,S_+}(E^c) > 0$ of possible non-escape is treated separately, by a simple argument asserting that, in this case,  $\E_{\sdev,S_+}^\imac
 [P_-]=-\infty$ while\footnote{That  $\E_{S_-,S_+}^\imac
 [P_+] = m_\imac$ and $\E_{S_-,S_+}^\imac
 [P_-] = n_\imac$ is proved in~\cite[Lemma~$3.11(2)$]{LostPennies}   which is contingent on~\cite[Lemma~$3.7$]{LostPennies}. The latter result has an invalid proof for the present context (where $\kappa$ may be less than one), but in the application in question, the pair $(S_-,S_+)$ lies in~$\tis^2$ with the stake amounts $a_i$ and $b_i$ all being positive; and, in this case,~\cite[Lemma~$3.7$]{LostPennies} is readily obtained for~$\kappa \in (0,1)$.} $\E_{S_-,S_+}^\imac [P_-] = n_\imac > -\infty$. The conclusion that  $\E_{\sdev,S_+}^\imac [P_-]=-\infty$ is easy to reach in the pure stake $\kappa =1$ case: since $\PP^\imac_{\sdev,S_+}(E^c) > 0$, an edge $[i,i+1]$ indexed by some $i \in \Z$ may be found that  is traversed from right-to-left infinitely often with positive probability. When the counter is at~$i+1$, Mina consistently faces a stake of $a_{i+1} >0$, so that, in order to win infinitely many of the moves from site~$i+1$, she has to expend infinitely in stake payments. In the present case, where $\kappa \in (0,1)$, this reasoning is flawed, because each move from site~$i+1$ is flip with probability~$1-\kappa > 0$, so that the edge $[i,i+1]$ may in principle be traversed from right to left by the counter on infinitely many occasions without Mina spending a dime when the counter is at~$i+1$.

We will circumvent this difficulty: rather than establishing that $\E_{\sdev,S_+}^\imac[P_-]=-\infty$ when  $\PP^\imac_{\sdev,S_+}(E^c)$ is positive, we will invoke the next result. We write $\trap$ for the complement of the escape event~$E$.
\begin{proposition}\label{p.altered}
Suppose that $\PP_{\sdev,S_+}^\imac (\trap) > 0$. 
There exists an altered strategy for Mina $\sdevalt \in \mc{S}$ such that 
\begin{equation}\label{e.altprop}
\textrm{$\pidevalt( \trap ) = 0$ and $\eidev[P_-] \leq \eidevalt[P_-]$} \, .
\end{equation}
\end{proposition}
 Mina will be willing to use the altered strategy in place of the original deviating one, and her doing so permits us to reduce the proof of Proposition~\ref{p.comp.sminus} to the case where escape is almost certain under  $(\sdev,S_+)$.
  The argument needed to treat the case of certain escape is identical to the corresponding one in~\cite{LostPennies},
  and our discussion of it is subsumed in the overview to which we now turn.
 
\subsection{Overview of the proof at large}
 
Given the reduction of the proof of Theorem~\ref{t.nashabmn}(2) that we summarised verbally after Propositions~\ref{p.supjk} and~\ref{p.comp.sminus}, which is recorded more formally in~\cite[Section~$4.2$]{LostPennies}, the substantial elements for this overview are the proofs of this pair of results.  We discuss them in turn.

{\em Deriving Proposition~\ref{p.supjk}.} There may be infinitely many  Mina deviation points for $\sdev$ whose spatial coordinate lies in~$\llbracket -j,k \rrbracket$. We begin by reducing to a finite number by eliminating late deviating moves. For $\hmac \in \N$, let the strategy $\sdev[\hmac]$ be formed from $\sdev$ 
by removing every deviating move after time~$\hmac$: thus, Mina will stake $b_u$ at $(u,t) \in \Z \times \N$ when $t \geq \hmac$.  Since $\PP_{\sdev,S_+}^\imac
(E) = 1$, the strategy pair $(\sdev[\hmac],S_+)$ when played from $\imac$ on gameboard~$\llbracket -j,k \rrbracket$ results in termination at a random finite time; so if Mina plays $\sdev[\hmac]$ in place of $\sdev$ for high~$\hmac$, there will be merely an arbitrarily small shift in the mean outcomes. 

Restricting to such finitely deviating strategies permits the fundamental game-theoretic technique of backward induction to be applied. We first describe the basic plan. Take a given strategy $\sdev$ with finitely many Mina deviation points whose spatial coordinate lies in  $\llbracket -j,k \rrbracket$. Let $g$ be the earliest time of one of the deviating points. Form a strategy~$S'$ by correcting all deviating play for Mina at time~$g$. Since there are fewer deviating points, an inductive hypothesis may be invoked to conclude that Mina's mean total receipt at any space-time $(v,g+1)$ is no higher than the value $b_v$ obtainable under non-deviant play via $(S_-,S_+)$. Now undo the time-$g$ corrections $S' \to \sdev$ and consider a location $(w,g)$ of deviating play for Mina. The inductive step is completed by arguing that Mina's outcome is strictly worse than it would be under non-deviating play from $(w,g)$. As we have seen, the boundary condition at time $g+1$ is not better; the argument analyses the one-step game played from~$(w,g)$ with these boundary conditions. It is at this point that one of the variations of the proof from~\cite{LostPennies} is made. The needed input is the analysis of the one-step game $(\kappa,\rho)$-Penny Forfeit from Lemma~\ref{l.pennyforfeit}: for $\kappa \in (0,1)$ and $\rho \in (0,1]$ in the present context, but with $(\kappa,\rho)=(1,1)$ in~\cite{LostPennies}. 

To state the formal change needed: 
we write $\PP_{S,S_+}^{u,\ell}$ for the law of a delayed-start game (governed by the indicated strategy pair), in accordance with notation set out in  \cite[Section~$3.3$]{LostPennies}.
Under 
this law,
the counter starts at~$u$ on  the turn with index~$\ell$. We record the stake offered under strategy~$S$ at this starting turn as $S(u,\ell)$, with the history that brought the counter to this space-time point being implicit.
The two displayed equations in the proof of~\cite[Lemma 4.16(2)]{LostPennies} will now read
\begin{eqnarray*}
\E^{u,\ell}_{S,S_+} \big[ P^{j,k}_- \big] & = & \bigg(  \frac{\kappa \, S(u,\ell)^\rho}{a_u^\rho+S(u,\ell)^\rho}  + \frac{1-\kappa}{2} \bigg) \, \E^{u-1,\ell+1}_{S,S_+} \big[ P^{j,k}_- \big]  \\
& & \qquad \qquad \qquad 
+ \,  \bigg( \frac{\kappa \, a_u^\rho}{a_u^\rho+S(u,\ell)^\rho}+ \frac{1-\kappa}{2} \bigg)  \, \E^{u+1,\ell+1}_{S,S_+} \big[ P^{j,k}_- \big] \, - \, S(u,\ell)  \, ,
\end{eqnarray*}
and
 $$
\E^{u,\ell}_{S,S_+} \big[ P^{j,k}_- \big] \, \leq \,  \bigg( \frac{\kappa \, S(u,\ell)^\rho}{a_u^\rho+S(u,\ell)^\rho}  \, + \, \frac{1-\kappa}{2} \bigg)  n_{u-1}
+  \bigg(  \frac{\kappa \, a_u^\rho}{a_u^\rho+S(u,\ell)^\rho}+ \frac{1-\kappa}{2} \bigg)    n_{u+1} \, - \, S(u,\ell) \, .
$$
We then invoke Lemma~\ref{l.pennyforfeit} to find that the preceding right-hand side has a unique maximum in $b$ at $b = b_u$, when it assumes the value $n_u$.

A second  variation addresses a point that has been elided in the above summary. There is a difference in strategy definition between~\cite{LostPennies} and the present article. While~\cite{LostPennies} specifies strategies simply as functions of space-time, we permit them to depend on the counter history to the present moment. This has led us to a definition of Mina deviation point whereby there must exist at least one history leading to the point in question which would cause her to place a deviant stake in playing from there. In order that the proof of Proposition~\ref{p.supjk} leads to a strict inequality in its conclusion, it is enough to argue that, for at least one Mina deviation point~$(u,\ell)$ with $u \in \llbracket -j,k \rrbracket$, 
every element in the path space~$\Lambda$ that begins at~$(\imac,0)$ and ends at~$(u,\ell)$ lies in the rectangle~$\llbracket -j,k \rrbracket \times \llbracket 0, \ell \rrbracket$.
Indeed, for such a point~$(u,\ell)$, there exists a history~$(\imac,0) \to (u,\ell)$---choose one and call it~$\hmacprime$!---which induces Mina to play a deviant move at the $(\ell + 1)$\textsuperscript{st} turn. The counter may follow this path without the game played on~$\llbracket -j,k \rrbracket$ ending. The trajectory follows this history with positive probability (if $\kappa \in (0,1)$, via a sequence of flip moves; if $\kappa=1$, by an argument in~\cite{LostPennies}). Consequently, the introduction of this Mina deviation point in the iterative procedure discussed above leads to a positive loss in her mean payoff, as in the proof we are adapting.  The loss is determined by~$(u,\ell)$ and~$\hmacprime$. The introduction of other deviation points has a non-positive effect on her payoff, so the cumulative effect is bounded above by the said loss. In Proposition~\ref{p.supjk}, a given deviation point $(u,\ell)$ is considered, and the hypothesis $u \in \llbracket -j+\ell,k-\ell \rrbracket$ is imposed on $j$ and $k$. It is this hypothesis that ensures that $(u,\ell)$ meets the condition on path-space inclusion. The values of $j$ and $k$ may be chosen to exceed some large constant specified by the given~$(u,\ell)$, so the resulting loss is independent of such~$(j,k)$; this leads to the uniformity asserted in Proposition~\ref{p.supjk}.

In summary, an inductive argument based on noting that deviating play is punished in the one-step game leads to the inference that the above discussed finite-deviating strategies $\sdev[\hmac]$ are uniformly punished on finite trails. By choosing the finite trails to be broad enough, the condition $\PP_{\sdev,S_+}^\imac(E) = 1$ implies that the error arising from the use of $\sdev[\hmac]$  in place of $\sdev$ is for high~$\hmac$ smaller than the incurred penalty. In this way, Proposition~\ref{p.supjk} is derived.

{\em Obtaining Proposition~\ref{p.comp.sminus} from Proposition~\ref{p.supjk}.}
Proposition~\ref{p.altered} permits us to reduce to the case where $\PP^\imac_{\sdev,S_+}(E) = 1$. 
The certainty of escape means that the counter will leave a broad enough board on the side on which it escapes globally. 
This permits us to truncate to a broad finite board incurring an arbitrarily small discrepancy in mean terminal payment. Removing non-negative running costs incurred beyond departure from the finite board then yields  $\E_{\sdev,S_+}^\imac [P_-] \leq \E_{\sdev,S_+}^\imac [P_-^{j,k}]$ up to the same small error. Proposition~\ref{p.supjk} may then be invoked, with the uniform penalty there identified overcoming the small opposing error, yielding the sought bound   $\E_{\sdev,S_+}^\imac [P_-] < \E_{S_-,S_+}^\imac [P_-]$.

\subsection{Obtaining Proposition~\ref{p.altered}}

Our discussion of the proof of Theorem~\ref{t.nashabmn}(2) concludes with the following derivation. 

{\bf Proof of Proposition~\ref{p.altered}.}
The trap event $\trap$ is a costly one for Mina because her terminal receipt in this event will be~$n_*$, which is by assumption strictly lower than her losing receipt~$n_\infty$; and, moreover, she may have running costs to pay. She would be happier with an altered strategy in which she instead consistently stakes zero in the trapping event, leading to an improved terminal receipt of $n_\infty$
alongside zero running cost. The problem with this idea is that the proposed alteration is not a well-defined strategy, because the proposed change is contingent on the occurrence of~$E^c = \trap$, an event undetermined by any finite-step evaluation of gameplay.  We will resolve this difficulty by introducing an event $\proxytrap$ determined by an initial move-sequence that nearly coincides with $\trap$, and defining Mina's altered strategy~$\sdevalt$  by asking her to stake zero after the moment at which $\proxytrap$ has been determined to occur.  The definition will yield an admissible strategy because the specification of the strategy space~$\mc{S}$ in Section~\ref{s.tlp} permits a player to consult counter history in deciding how to stake.  
   
Before elaborating this construction, we first address a simpler case, in which it is not needed: this is when trapping is not merely possible, but certain. That is, if $\PP_{\sdev,S_+}^\imac (\trap) = 1$,
then 
 we may simply take $\sdevalt$ equal to the zero strategy. Doing so results in Maxine winning every stake move under $\pidevalt$, with counter evolution $X:\N \to \Z$ 
given by {\rm SRW}$\big( (1+\kappa)/2 \big)$ begun at $X_0 = \imac$, this entailing the occurrence of~$E_+$; since Mina has no running costs, we see then that $P_-  = T_- = n_\infty$ holds $\pidevalt$-almost surely. 
The desired properties~(\ref{e.altprop}) hold, the inequality due to  $\eidev[P_-] \leq n_* < n_\infty = \eidevalt[P_-]$.
  
Now assume that  $\PP_{\sdev,S_+}^\imac (\trap) \in (0,1)$. In constructing and analysing the altered strategy $\sdevalt$, we will couple the gameplays under $(\sdev,S_+)$
and~$(\sdevalt,S_+)$. We will write $\PP^\imac$ for the law governing  these two gameplays, and will distinguish between them by indicating the strategy pair associated to a given random variable. For example, $T_-(\sdev,S_+)$ under $\PP^\imac$ denotes Mina's terminal payoff for gameplay governed by the strategy pair $(\sdev,S_+)$ under the coupling; in law, this random variable is equal to $T_-$ under $\pidev$. 
  
 Under $\PP^\imac$, we will define $\proxytrap$ in terms of a parameter~$\e > 0$ measuring the approximation of~$\trap$. 
 We will set $\proxytrap = \big\{ \tau_\e < \infty \big\}$ for an $\N$-valued stopping time $\tau_\e$ in such a way that
  $$
   \textrm{$\trap \subseteq \proxytrap$ holds up to a $\PP^\imac$-null set, and} \, \, \, \, \PP^\imac \big( \proxytrap \setminus \trap \big) \leq \e \, .
 $$
To construct $\tau_\e$, let $\#_j$  for $j \in \Z$ denote the total number of visits made by $X:\N \to \Z$ to the site~$j$, for the copy of counter evolution under $(\sdev,S_+)$
  offered by~$\PP^\imac$. It follows readily from the meaning of absence of escape that for any $i \in \N_+$, we may find a non-random finite subset $J_i \subset \Z$
  such that
  \begin{equation}\label{e.maxprobtrap}
   \PP^\imac \Big( \max_{j \in J_i} \#_j = \infty \, \Big\vert \, \trap \Big) \, \geq \, 1 - \e/2^i \, ,
  \end{equation}
  where here it is understood that $\trap$ is specified in terms of  counter evolution under $(\sdev,S_+)$. We may further select $N_i \in \N$ for which
  \begin{equation}\label{e.maxtrapcomplement}
   \PP^\imac  \Big( \max_{j \in J_i} \#_j \geq N_i \, \Big\vert \,  \trap^c \Big) \, < \, \e/2^i  \, .
  \end{equation}
  Writing $\#_j(n)$ for the cardinality of the set of times at most $n$ at which the counter visits $j$ (so that $\#_j(\infty) = \#_j$), we set
  $$
  \phi_i \, = \, \min \Big\{ n \in \N: \max_{j \in J_i} \#_j(n) \geq N_i \Big\} \, .
  $$
  Now we set $\tau_\e = \min_{i \in \N} \phi_i$. To define~$\sdevalt$, recall the path spaces $\Lambda_k$ used to specify~$\mc{S}$ in Section~\ref{s.tlp}. 
  We set $\sdevalt(\psi) = \sdev(\psi)$ whenever $\psi \in \Lambda_k$  for some $k \in \N$ such that  $\tau_\e > k$ for any counter evolution~$X$ with $X \big\vert_{\llbracket 0, k \rrbracket} = \psi$.
  The value of  $\sdevalt(\psi)$ is set to zero for any $\psi$ in the remainder of~$\Lambda$. It is a straightforward check that $\tau_\e$ is a stopping time and $\sdevalt$ an element of~$\mc{S}$
(which coincides with $\sdev$ before $\tau_\e$ and is zero at that time and thereafter).  
   The coupling $\PP^\imac$ is constructed so that counter evolutions under~$(\sdev,S_+)$ and~$(\sdevalt,S_+)$ almost surely coincide until~ $\tau_\e$.   
  
  Next we verify the desired properties that $\trap \subset \proxytrap$ and $\PP^\imac(\proxytrap \setminus \trap ) < \e$. To do so, note that~(\ref{e.maxprobtrap}) implies that, conditionally on $\trap$,   $\max_{j \in J_i} \#_j$ equals infinity for all but finitely many~$i$ almost surely, so that $\phi_i$ is finite with the exception of at most finitely many $i$; this implies that $\tau_\e < \infty$, so that $\proxytrap$ occurs.
 Note further that
 $$
 \PP^\imac(\proxytrap \setminus \trap)  =  \PP^\imac \Big(  \trap^c \cap \big\{ \exists \, i \in \N_+, j \in J_i: \#_j \geq N_i \big\} \Big) <  \sum_{i=1}^\infty\ \e/2^i = \e \, ,
 $$ 
 the bound due to~(\ref{e.maxtrapcomplement}).

 We now use the constructed $\tau_\e$ to prove that the desired~(\ref{e.altprop}) holds.  Note that
 \begin{equation}\label{e.eiproxy}
  \E^\imac
 \Big[ P_-(\sdevalt,S_+) {\bf 1}_{\proxytrap^c} \Big] =   \E^\imac
 \Big[ P_-(\sdev,S_+) {\bf 1}_{\proxytrap^c} \Big]  \, .
 \end{equation}
 Running costs under $\sdevalt$ and $\sdev$ coincide until~$\tau_\e$, after which they are cancelled under~$\sdevalt$. The switch to the altered strategy when $\trap$ occurs leads to a terminal receipt of $n_\infty$ in place of~$n_*$ for Mina (she loses the game, but at least it finishes). 
Thus the bound holds in the following:  
 \begin{eqnarray*}  
 \E^\imac
 \Big[ P_-(\sdevalt,S_+) {\bf 1}_{\trap} \Big] & = & \E^\imac
 \, \bigg[ \Big(  T_-(\sdevalt,S_+) - \sum_{i=0}^\infty \sdevalt(i) \Big)  {\bf 1}_{\trap}  \bigg] \\
  & \geq & \big(n_\infty - n_* \big) \PP^\imac(\trap) \, + \, \E^\imac
 \, \bigg[ \Big( T_-(\sdev,S_+)  - \sum_{i=0}^\infty \sdev(i) \Big) {\bf 1}_{\trap} \bigg]  \\
   & = & \big( n_\infty - n_* \big) \PP^\imac(\trap) 
   + \E^\imac \big[ P_- (\sdev,S_+) {\bf 1}_\trap \big] \, .
 \end{eqnarray*} 
(Note that in the summands a standard usage is made to refer to stakes offered at the $(i+1)$\textsuperscript{st} turn by the strategy in question.) 
 
The same inequality on running costs implies the first bound as we write
 \begin{eqnarray*}
 & &  \E^\imac
 \Big[ P_-(\sdevalt,S_+) {\bf 1}_{\proxytrap \setminus \trap} \Big] -   \E^\imac
 \Big[ P_-(\sdev,S_+) {\bf 1}_{\proxytrap \setminus \trap} \Big] \\
  & \geq &   \E^\imac
 \Big[ T_-(\sdevalt,S_+) {\bf 1}_{\proxytrap \setminus \trap} \Big] -   \E^\imac
 \Big[ T_-(\sdev,S_+) {\bf 1}_{\proxytrap \setminus \trap} \Big] \\
  & \geq &  \Big(  n_\infty - \big( \mu n_\infty + (1-\mu) n_{-\infty} \big) \Big) \PP^\imac \big( \proxytrap \setminus \trap \big) \geq - \big( n_{\infty,-\infty} \big) \e \, ,
 \end{eqnarray*}
where the convex combination $\mu n_\infty + (1-\mu) n_{-\infty}$ appears because $T_-$ equals either $n_\infty$ or $n_{-\infty}$ on~$E$ under $(\sdev,S_+)$.
 
 Since $\trap$, $\proxytrap \setminus \trap$ and $\proxytrap^c$ partition the space of outcomes, we may add~(\ref{e.eiproxy}) and the two bounds that follow it to obtain
 $$
  \E^\imac
 \Big[ P_-(\sdevalt,S_+) \Big]  \geq  \E^\imac
 \Big[ P_-(\sdev,S_+) \Big] + \big( n_\infty - n_* \big) \PP^\imac(\trap) -     \big( n_{\infty,-\infty} \big) \e \, .
 $$
 Choosing $\e$ to be less than $\tfrac{n_\infty - n_*}{n_{\infty,-\infty}} \PP^\imac(\trap)$, we find that $\E^\imac \big[  P_-(\sdevalt,S_+)  \big] > \E^\imac \big[  P_-(\sdev,S_+)  \big]$, as claimed in~(\ref{e.altprop}). This completes the proof of Proposition~\ref{p.altered}. \qed

\part{Back to Brownian Boost\label{p.two}}

\chapter{The Brownian Boost ODE pair}\label{c.brownianboost}

This part returns to Brownian Boost and develops its analysis in two chapters.
In the present chapter, we begin with a heuristic explanation for why time-invariant equilibria of the game are governed by solutions of the \bbrhospace ODE pair in Definition~\ref{d.depair}, and we then undertake a detailed analytic study of these solutions.
In the following chapter, Brownian Boost is realised as a fine-mesh, high-noise limit of the Trail of Lost Pennies, leading to the low-$\kappa$ $\lambdamax$ Theorem~\ref{t.lowlambdamax} and the asymptotic stakes-and-gameplay Theorem~\ref{t.lowkappasde}.

The present chapter is divided into two sections.
The first offers a formal but heuristic derivation of the \bbrhospace ODE pair, obtained by analysing \bbrhospace directly in continuous time.
The second solves this ODE system explicitly, establishing Theorem~\ref{t.fg}, describing the qualitative behaviour of solutions in Proposition~\ref{p.fgab}, and recording further analytic properties that will be used in Chapter~\ref{c.highnoise} to analyse the low-$\kappa$ regime of \tlpkp.

\section{Coupled HJB equations for Brownian Boost}\label{s.bbrho} 
The Hamilton-Jacobi equation arises from the Euler-Lagrange equation in a reformulation of Newtonian mechanics. Bellman~\cite{Bellman} generalized the context to control theory (with one agent)
and Isaacs~\cite{Isaacs65} to zero-sum differential game theory (with two or more players). 
In our non-zero-sum context, there is a system of HJB equations indexed by the players.
For conceptual clarity, here we give a formal argument exhibiting the ODE pair as coupled HJB equations for \bbrho.

In our rigorous treatment, $\rho$-Brownian Boost is regularized as \tlpkpspace  for low~$\kappa$. For the present purpose, we disregard niceties concerning how feedback loops interfere with specifying gameplay in~\bbrho, and study the game directly.

Suppose then that \bbrhospace is played at a Nash equilibrium, with  Maxine and Mina adopting stake profiles $a,b:\R \to [0,\infty)$
from which neither would unilaterally choose to deviate. 
 For $x \in \R$, let 
 $m(x)$ and $n(x)$
 denote the mean total receipt accruing respectively to Maxine and Mina when $X_0 = x$ and the stake profile pair $(a,b)$ is adopted.
 
 Consider Maxine's point of view in the first $[0,{\rm d}t]$ of time. In this duration, she will spend $a(x) {\rm d}t$, where the error due to taking $a(X_s) = a(x)$ for $s \in [0,{\rm d}t]$
 is negligible. 
 Writing  $N(0,r)$ for a centred Gaussian of variance~$r$, note that $X({\rm d}t)$ equals  $x + \tfrac{a^\rho-b^\rho}{a^\rho+b^\rho} \dd t + N(0,\dd t)$ in law. 
 Maxine's mean net receipt equals  her mean receipt subsequent to time~${\rm d}t$ minus the running cost that she accrues on~$[0,{\rm d}t]$: that is,
   $$
   m(x) = \E \, m(X_{{\rm d}t}) \, - \, a(x) {\rm d}t \, ,
   $$
or
$$
 m(x) \, = \,-a(x) \dd t \, + \, \E \, m \Big( x + \tfrac{a(x)^\rho-b(x)^\rho}{a(x)^\rho+b(x)^\rho} \dd t + N(0,\dd t)  \Big) \, .
 $$
With
  $\mu(y,r)$ denoting the law of $N(0,r)$, the latter expected value is $m(x) + \frac{a-b}{a+b} m'(x) \dd t + I$ where $I = \int_\R \big(  m(x+y) - m(x) \big) \dd \mu(y,\dd t)$ equals $\tfrac{1}{2}m''(x) \dd t$. Cancelling $m(x)$, and omitting to denote the argument~$x$,  
\begin{equation}\label{e.am}
 \frac{a^\rho-b^\rho}{a^\rho+b^\rho} m' + \frac{m''}{2} - a \, = \, 0  \, .
\end{equation}
Mina's point of view offers the analogous
\begin{equation}\label{e.bn}
 \frac{a^\rho-b^\rho}{a^\rho+b^\rho} n' + \frac{n''}{2} - b \, = \, 0 \, .
\end{equation}
This is a pair of Markovian dynamic programming equations, valid for any stake profile pair $(a,b)$. As we show next, a further equation pair arises from the consideration that $(a,b)$ is a Nash equilibrium. The stability under unilateral deviation manifest in this concept is gauged in terms of mean total net receipt, with the class of perturbed strategies being broader than time-invariant ones. Indeed, let $z \in [0,\infty)$, and suppose that Maxine stakes at rate~$z$ during $[0,{\rm d}t]$, after which she reverts to the dictates of the stake profile $a:\R \to [0,\infty)$. Writing $m(x,z)$ for her mean net receipt when she plays against Mina's stake profile~$b$, we have that
$$
 m(x,z) \, = \, - z \, {\rm d} t \, + \, \E \, m \Big( x + \tfrac{z^\rho-b^\rho}{z^\rho+b^\rho} \dd t + N(0,\dd t)  \Big) \, ,
$$
whence $m(x,z) = m(x) +  \big( \tfrac{z^\rho-b^\rho}{z^\rho+b^\rho} m'(x) - z + m''(x)/2 \big) {\rm d}t$.
Since $(a,b)$ is a Nash equilibrium, the $z$-indexed variant strategy does not tempt Maxine, and $z \to m(x,z)$ has a maximum at $z = a(x)$, so that the partial derivative in $z$ of the just recorded ${\rm d}t$-coefficient vanishes at $z = a(x)$. Rearranging,
\begin{equation}\label{e.amprime}
 2\rho b^\rho a^{\rho - 1} m' = \big(a^\rho+b^\rho \big)^2 \, .
\end{equation}
Mina's counterpart variation completes the second equation pair:
\begin{equation}\label{e.bmprime}
 - 2\rho a^\rho b^{\rho -1} n' = \big(a^\rho+b^\rho \big)^2 \, ,
\end{equation}
where note that $n' < 0$, since Mina plays left.

Supposing that $a$ and $b$ are positive (and we omit to justify this in these heuristics), 
the just obtained pair implies
$m' b = - n' a$. Returning to the same equation pair with this fact, and introducing the notation 
$f = m' > 0$ and $g = -n' > 0$, we obtain 
\begin{equation}\label{e.ab}
 a =  \frac{2\rho f^{1+\rho}g^\rho}{\big( f^\rho+ g^\rho \big)^2} \, \, \, \, \textrm{and} \, \, \, \,  b =  \frac{2\rho f^\rho g^{1+\rho}}{\big( f^\rho+  g^\rho \big)^2} \, .
\end{equation} 
Hence,
$$
 \frac{a^\rho}{a^\rho+b^\rho} =  \frac{f^\rho}{f^\rho+ g^\rho}  \, \, \, \, \textrm{and} \, \, \, \,   \frac{b^\rho}{a^\rho+b^\rho} =  \frac{g^\rho}{f^\rho+ g^\rho} \, .  
$$
Revisiting~(\ref{e.am}) and~(\ref{e.bn}) with these inferences and notation, 
$$
 a \, = \,  \frac{f^\rho  -  g^\rho}{f^\rho+ g^\rho}\, f  \, + \, \frac{f'}{2}  \, \, \, \, \textrm{and} \, \, \, \,  b \, = \, - \frac{f^\rho  -  g^\rho}{f^\rho+ g^\rho}\, g \, - \, \frac{g'}{2} \, . 
$$ 
Substituting these stake profile formulas into~(\ref{e.ab}) yields
\begin{eqnarray}
 2 \rho f^{1+\rho}  g^\rho & = & \, \, \, \, \, \big( f^{2\rho}- g^{2\rho} \big)f + \tfrac{1}{2} f' \big(f^\rho +  g^\rho \big)^2 \nonumber \\ 
 2 \rho f^\rho  g^{1+\rho} & = & - \big( f^{2\rho}- g^{2\rho} \big)g - \tfrac{1}{2} g' \big(f^\rho +  g^\rho \big)^2 \, , \nonumber 
\end{eqnarray}
so that $(f,g)$ solves the $\rho$-Brownian Boost ODE pair specified in Definition~\ref{d.depair}.

{\em Remark.} This formal derivation of the \bbrho -ODE pair is valid for $\rho \in (0,\infty)$. On the discrete side, the region~$W$ of $(\kappa,\rho)$-space specified in~(\ref{e.weakregion}) and utilized in computing \abmnkpmacspace elements contains in its boundary the whole `Brownian Boost line' $\{ 0 \} \times (0,\infty)$. These indications may suggest that discrete solution impediments  lose their force in the low~$\kappa$ limit, and that \bbrhospace equilibria may be computed for all $\rho \in (0,\infty)$.

These impressions are misleading: the derivation above identifies
critical points of the players' optimization problems for all
$\rho\in(0,\infty)$, but it does not guarantee their global optimality.
As noted in Section~\ref{s.pennyforfeit}, the one-step game Penny
Forfeit admits a global maximizer in pure strategies throughout the
regime $\rho\le1$, while for $\rho>1$ global optimality may fail,
independently of the value of $\kappa\in(0,1)$. Accordingly, our
rigorous results for Brownian Boost are confined to $\rho\in(0,1]$,
and the analysis of higher values of~$\rho$ remains open.

\medskip

\subsubsection{Moscarini and Smith's game}\label{s.ms}
In~\cite{MoscariniSmith2007} a continuous-time tug-of-war game in which the counter's drift is given by
$a(X_t)^{1/2}-b(X_t)^{1/2}$ is studied---in our terms, the game is Brownian Boost with $\rho=1/2$ and the drift denominator omitted, so that the drift term in~(\ref{e.countersde}) becomes  homogeneous  in the stake pair~$(a,b)$ with degree one-half rather than zero.
Moscarini and Smith focus on Markov-perfect equilibria for games played on finite intervals centred at the origin, restricting attention to symmetric equilibria, in which reflection through the origin corresponds to interchange of the players' continuation values; they obtain explicit characterisations of such equilibria.

Their analysis may be recast in a form closely parallel to the heuristic HJB derivation given above for Brownian Boost.
Introducing value derivatives $f=m'>0$ and $g=-n'>0$, the equilibrium conditions lead to the coupled ODE system
\[
f' = f\big(g- f/2 \big),
\qquad
g' = g\big(g/2-f\big),
\]
with accompanying stake rates
$a=f^2/4$ and $b=g^2/4$.
By translation invariance, one may take $f(0)=g(0)>0$.

Analysing the MS-ODE pair shows that
for positive arguments $g<f$ and $g/f \to 0$,
so 
a discouragement effect is present. Moreover,
$f'(u)\sim -f(u)^2/2$ and  $f(u)\sim 2/u$ for large~$u$.
As a consequence,
$\int_{\R} f(r)\,{\rm d}r
=
\int_{\R} g(r)\,{\rm d}r
=
\infty$,
and thus no MS equilibria on the real line exist with finite terminal prizes.
This necessitates the restriction, adopted in~\cite{MoscariniSmith2007}, to finite gameboards. 

Although Moscarini and Smith recover a discouragement phenomenon akin to that identified by Harris and Vickers,
there does not appear to be a counterpart of the delicate $\lambdamax=1$ and $\lambdamax - 1 \in [0,10^{-4}]$ effects that we find rigorously and numerically in \bbrhospace and \tlpkp.
Indeed, $\lambdamax$ may naturally be interpreted as the supremum, over
finite gameboards~$I$ and solutions $(f,g)$ of the above ODE pair, of the
statistic $\int_I g(r)\,{\rm d}r / \int_I f(r)\,{\rm d}r$.
Gameboards must be finite intervals in order for these quantities to be
well-defined positive real values, but variation over all such intervals
simply yields $\lambdamax=\infty$, since $g\gg f$ for $u\ll0$.

Nor does a sharply localised battlefield region emerge: on long intervals, equilibrium stakes and continuation values vary relatively slowly across the interior.
It would be interesting to extend Moscarini and Smith's classification beyond the symmetric case.
From the ODE perspective, one natural approach is to fix a solution $(f,g)$ with $f(0)=g(0)$ and associate to any finite interval $I\subset\R$ the triple
$\big( \int_I f , \, \int_I g , \, |I| \big)$.
Determining whether this correspondence is bijective is an open question whose resolution would yield a concrete parametrisation of finite-interval MS equilibria
in terms of terminal payoffs and gameboard size.

\section{Proving properties of the ODE pair}\label{s.odeproofs}

Here we prove  Theorem~\ref{t.fg} and Proposition~\ref{p.fgab}. We also derive further information on \bbrhospace ODE pair solutions in Proposition~\ref{p.sfacts}. This includes the key identity $\int_{\R} f_\rho(\xmac,u) \, {\rm d} u =  \int_{\R} g_\rho(\xmac,u) \, {\rm d} u$,
valid for any choice of the flow index $\xmac \in (0,\infty)$: in the notation of~ (\ref{e.lambdamaxkappazero}), $\lambdamax(0,\rho) =1$, so that no incentive asymmetry may exist at equilibrium. This information will be central to deriving Theorem~\ref{t.lowlambdamax} on low-$\kappa$ $\lambdamax(\kappa,\rho)$ in Chapter~\ref{c.highnoise}.


Theorem~\ref{t.fg} classifies default solutions of the \bbrhospace ODE pair and, as we now note in passing, it thereby also classifies all {\em strictly positive} solutions. 
Indeed, the dilation $(f,g) \mapsto (af,ag)$ preserves the positive solution space for any $a>0$, so that any such solution is equivalent, up to dilation, to a default one. 
It follows immediately from Theorem~\ref{t.fg} that the strictly positive solutions are precisely the positive dilations of default solutions. 
Note also that $\big(A \, e^{-2u}, 0 \big)$ and $\big( 0, A \, e^{2u} \big)$, for $A \in [0,\infty)$, are non-negative solutions.

 We begin the analytic derivations by recasting the ODE satisfied by $S_\rho(\xmac,\bullet)$ in Definition~\ref{d.fg} by means of the $\rho$\textsuperscript{th} power of this function. 
\begin{lemma}\label{l.de}
For $\rho,\xmac \in (0,\infty)$, set
$J(u) := S_\rho(\xmac,u)^\rho$, where $S_\rho(\xmac,\bullet)$ is specified in Definition~\ref{d.fg}.
Then $J$ is the unique solution of the differential equation
\[
\frac{{\rm d} J(u)}{{\rm d} u} \, = \,  -8 \rho^2 \, \frac{J(u)^2}{\big(1 + J(u)\big)^2} \, \, \, \, \textrm{with} \, \, \, \,  J(0) = \xmac^\rho \, .
\]
\end{lemma}
{\em Remarks: (1).} 
This result has the consequence that $S_\rho(\xmac,u)^\rho= S_1(\xmac^\rho,\rho^2 u)$, since the right-hand expression is also a solution of the  equation.

{\em (2).}
Integrating the equation, we find that $J(u)^2 e^{J(u) - J(u)^{-1}} = e^{-8 \rho^2 u}$ when $\xmac=1$.
In view then of what we just noted, 
$$
 S_\rho(1,u)^\rho = S_1(1,\rho^2 u)   \, \sim \, \begin{cases}  \,  8 \rho^2 \vert u \vert  \qquad \quad \, \,  u \ll 0 \, ,  \\ \,  \big( 8\rho^2 u \big)^{-1} \qquad u \gg 0 \, .
  \end{cases} 
$$

{\bf Proof of Lemma~\ref{l.de}.} 
The initial condition $J(0) = S_\rho(\xmac,0)^\rho = \xmac^\rho$ holds. Differentiating $J(u) = S_\rho(\xmac,u)^\rho$ gives
\[
J'(u) \, = \, \rho S_\rho(\xmac,u)^{\rho-1} S_\rho'(\xmac,u) 
\, = \, \rho S_\rho^{\rho-1} \times \frac{- 8 \rho \, S_\rho^{1+\rho}}{(1+S_\rho^\rho)^2} 
\, = \, -8 \rho^2 \frac{J(u)^2}{(1 + J(u))^2} \, ,
\]
as desired.

We now argue that the solution $J$ of the initial-value problem is unique. Let $K$ be another, and set $A = \big\{ x \in \R: J(x) = K(x) \big\}$.
Then $0 \in A$ by assumption. Since $J$ is everywhere positive, and $J$ and $K$ are continuous, the right-hand side of the differential equation being Lipschitz in $J$ implies that $A$ is open. Since $J$ and $K$ are continuous, $A$ is closed. Thus $A = \R$ and $J =K$. \qed

A pair of logarithmic derivatives offers a convenient reformulation of the \bbrhospace ODE pair.
\begin{lemma}\label{l.fgcharacterization}
Let $\rho \in (0,\infty)$. For a pair of differentiable functions $f,g:\R \to (0,\infty)$, set $\phi_f = \frac{f'}{2f}$, $\phi_g = - \frac{g'}{2g}$ and $j = (g/f)^\rho$.
The pair $(f,g)$ is a solution of~(\ref{e.fg.rho}) if and only if the pair of equations 
$$
 \big( \phi_f, \phi_g \big) \, = \, \bigg( \frac{2\rho j - (1-j^2)}{(1+j)^2} , \frac{2\rho j + (1-j^2)}{(1+j)^2} \bigg) 
$$
is satisfied.
\end{lemma}
{\bf Proof.} 
Divide the first equation in the pair~(\ref{e.fg.rho})
  by $f$ and write in terms of $F := f^\rho$ and $G := g^\rho$ to obtain
\[
2\rho F G = F^2 - G^2 + \frac12 \frac{f'}{f} (F+G)^2
= F^2 - G^2 + \phi_f (F+G)^2   \, ,
\]
where we use $\phi_f = \frac{f'}{2f}$.  Thus the pair~$(f,g)$ satisfies the first equation in~(\ref{e.fg.rho}) if and only if
$$
\phi_f = \frac{2\rho F G - (F^2-G^2)}{(F+G)^2} \, ,
$$
or 
\begin{equation}\label{e.phif}
\phi_f  =  \frac{2\rho j - (1-j^2)}{(1+j)^2} \, ,
\end{equation}
where we have introduced the function $j = G/F$ after dividing by the positive $F^2$.

Now divide the second equation~(\ref{e.fg.rho}) by $g$ to find that 
$2\rho F G$ equals $G^2 - F^2 + \phi_g (F+G)^2$.  Again dividing by $F^2 > 0$, we see that $(f,g)$ satisfies this second equation precisely when
$\phi_g = \frac{2\rho j + (1-j^2)}{(1+j)^2}$.

By intersecting the pair of established equivalences, we complete 
 the proof of Lemma~\ref{l.fgcharacterization}. \qed

{\bf Proof of Theorem~\ref{t.fg}.} Let $\xmac \in (0,\infty)$. In shorthand, we will denote  $f_\rho(\bullet) = f_\rho(\xmac,\bullet)$,  $g_\rho(\bullet) = g_\rho(\xmac,\bullet)$ and $S(\bullet) = S_\rho(\xmac,\bullet)$.
Write $F_\rho = f_\rho(\bullet)^\rho$ and $G_\rho = g_\rho(\bullet)^\rho$, and note that these functions are everywhere positive.

We will show that $(f_\rho,g_\rho)$ solves~(\ref{e.fg.rho}). To this end, note that, by Definition~\ref{d.fg}, 
\begin{equation}\label{e.fhrep}
f_\rho(\rindep) = \exp\!\Bigg\{ 2 \int_0^{\rindep} \Phi_f(\udummy)\, {\rm d}\udummy \Bigg\}, 
\qquad 
g_\rho(\rindep) = \xmac \cdot \exp\!\Bigg\{ -2 \int_0^{\rindep} \Phi_g(\udummy)\, {\rm d}\udummy \Bigg\} \, .
\end{equation}
where 
\begin{equation}\label{e.phifg}
\Phi_f = 1 - \frac{2\big(1+(1-\rho)S^\rho\big)}{(1+S^\rho)^2}, 
\qquad 
\Phi_g = 1 - \frac{2\big((1-\rho)S^\rho + S^{2\rho}\big)}{(1+S^\rho)^2}.
\end{equation}
Moreover, from these expressions for $f_\rho$ and $g_\rho$, we see that 
$\Phi_f = \frac{f'_\rho}{2f_\rho}$ and $\Phi_g = - \frac{g_\rho'}{2g_\rho}$.
Differentiating $F_\rho =f_\rho^\rho$ and $G_\rho=g_\rho^\rho$, we also have that
\begin{equation}\label{e.FH}
F_\rho' = 2\rho \, F_\rho \Phi_f \, \, \, \,  \textrm{and}  \, \, \, \, G_\rho' = -2\rho \, G_\rho \Phi_g \, .
\end{equation}

In order to argue that $(f,g) = (f_\rho,g_\rho)$ solves~\eqref{e.fg.rho},
we write $j_\rho = G_\rho/F_\rho$ for the $j$-function attached to the pair~$(f_\rho,g_\rho)$.

\begin{lemma}\label{l.js}
\leavevmode
\begin{enumerate}
\item We have that
$\Phi_f +\Phi_g = \frac{4\rho\,S^\rho}{(1+S^\rho )^2}$. 
\item And that
$$
\frac{ {\rm d}j_\rho(u)}{{\rm d}u} \;=\; -8\rho^2 \frac{j_\rho(u)\, S_\rho(\xmac,u)^\rho}{\big(1+S_\rho(\xmac,u)^\rho\big)^2} \, .
$$
\item 
And also that  \(j_\rho(u) = S_\rho(\xmac,u)^\rho\) for all \(u\in\R\).
\end{enumerate}
\end{lemma}
{\bf Proof: (1,2).} Since $g_\rho(0) = \xmac$ and $f_\rho(0) =1$, $j_\rho(0) = \xmac^\rho$. By~(\ref{e.FH}),
\begin{equation}\label{e.jfact}
\frac{ {\rm d}j_\rho(u)}{{\rm d}u} \;=\;  
 \frac{G_\rho'}{F_\rho} - \frac{G_\rho F_\rho'}{F_\rho^2}
= -2\rho j_\rho\big(\Phi_f+\Phi_g\big) \, .
\end{equation}

Writing $J = S^\rho >0$ as in  Lemma~\ref{l.de}, we find from~(\ref{e.phifg}) that 
\[
\Phi_f+\Phi_g
\, = \, 2 - \frac{2\big(1+2(1-\rho)J+J^2\big)}{(1+J)^2}
\, = \,  2 - \frac{2\big((1+J)^2 - 2\rho J\big)}{(1+J)^2} \, ,
\]
whence Lemma~\ref{l.js}(1) holds. Returning to~(\ref{e.jfact}), we obtain Lemma~\ref{l.js}(2).

{\bf (3).} By  Lemma~\ref{l.de}, $J= S^\rho$ from satisfies the differential equation in that result and may be compared to the solution~$j_\rho$ of the related differential equation in the preceding part. Consequently, $j_\rho' J = j_\rho J'$. Consider the ratio $q(u) = j_\rho(u)/J(u)$. The derivative is a fraction whose denominator is $J^2 > 0$ and whose numerator vanishes by the just obtained identity. So $q' = 0$ identically. Thus $q =1$ since $q(0) = j_\rho(0)/J(0) = \xmac^\rho/\xmac^\rho = 1$. Hence  $j_\rho = J = S^\rho$ and we obtain 
 Lemma~\ref{l.js}(3). \qed

By Lemma~\ref{l.fgcharacterization}, we may prove that $(f_\rho,g_\rho)$ solves~(\ref{e.fg.rho}) by showing that 
$$
 \big( \Phi_f ,\Phi_g \big) \, = \,  \bigg( \frac{2\rho j_\rho - 1+ j_\rho^2}{(1+j_\rho)^2} , \frac{2\rho j_\rho + 1-j_\rho^2}{(1+j_\rho)^2} \bigg) \, ,
$$
where $j_\rho = (g_\rho/f_\rho)^\rho$. But $j_\rho$ equals $S_\rho^\rho$ by Lemma~\ref{l.js}(3), so that this pair of conditions results from~(\ref{e.phifg}) by a simple rearrangement. 
Note further that, by taking $\rho$\textsuperscript{th} roots, we obtain $g_\rho(\xmac,\bullet)=  f_\rho(\xmac,\bullet) S_\rho(\xmac,\bullet)$, as claimed in Theorem~\ref{t.fg}.

To prove the converse direction in this theorem, 
let $(f,g)$ be a default solution of~\eqref{e.fg.rho}, so that $f(0) = 1$ and $g(0) > 0$.
By Lemma~\ref{l.fgcharacterization}, the pair $(\phi_f,-\phi_g)$ of one-half logarithmic derivatives satisfies
\begin{equation}\label{e.phifphih}
\phi_f = \frac{2\rho j - (1-j^2)}{(1+j)^2} \, \, \, \, \textrm{and} \, \, \, \,
\phi_g = \frac{2\rho j + (1-j^2)}{(1+j)^2} \, \, \, \, \textrm{with} \, \, \, \, j = (g/f)^\rho \, ,
\end{equation}
whence $\phi_f + \phi_g = 4 \rho j/(1+j)^2$. But $j'/j = -2 \rho \big( \phi_f + \phi_g \big)$, so that $j'  = -8 \rho^2 j^2/(1+j)^2$. Note that $j(0) = \xmac^\rho$ where we set $\xmac = g(0) > 0$.
Thus~$j$ solves the initial value problem satisfied by $J(u) = S_\rho(\xmac,u)^\rho$ in Lemma~\ref{l.de}. By the uniqueness claim in this lemma, $j = J$.  Hence, $j(u) = J(u) = S_\rho(\xmac,u)^\rho$ for all $u\in\R$.
Since $\phi_f$ and $-\phi_g$ are one-half logarithmic derivatives, we have
\[
f(\rindep) = f(0) \cdot \exp\Big\{ 2 \int_0^{\rindep} \phi_f(\udummy) \, {\rm d}\udummy \Big\} \, \, \, \, \textrm{and} \, \, \, \,
g(\rindep)
= g(0) \cdot \exp\Big\{ -2 \int_0^{\rindep} \phi_g(\udummy) \, {\rm d}\udummy \Big\} 
\]
with $\big(f(0),g(0)\big) = (1,\xmac)$.
Now~(\ref{e.phifphih}) alongside $j = S_\rho^\rho$ exhibits the pair $(f,g)$ in the desired form $\big(f_\rho(\xmac,\bullet), g_\rho(\xmac,\bullet)\big)$.
The converse direction thus treated, this completes the proof of  Theorem~\ref{t.fg}. \qed

We now gather analytic facts needed to prove Proposition~\ref{p.fgab} next and Theorem~\ref{t.lowlambdamax} later.

\begin{proposition}\label{p.sfacts}
Let $\rho,\xmac \in (0,\infty)$ and $\rindep \in \R$.
\begin{enumerate}
\item  $S_\rho(\xmac,-\rindep) = S_\rho(\xmac^{-1},\rindep)^{-1}$. In particular,
 $S_\rho(1,-\rindep) = S_\rho(1,\rindep)^{-1}$.
\item Now let $\vmac = \vmac(\xmac)$ denote the unique real number such that $8 \rho \vmac = 2\log \xmac + \rho^{-1} \big( \xmac^\rho - \xmac^{-\rho} \big)$. 
The map $(0,\infty) \to \R: \xmac \to \vmac$ is an increasing bijection. The value $\vmac = \vmac(\xmac)$ is the unique solution of $S_\rho(\xmac,\vmac) =1$. We have
$$
f_\rho(1,\rindep)=\frac{f_\rho(\xmac,\vmac+\rindep)}{f_\rho(\xmac,\vmac)} \, \, \, \,  \textrm{and} \, \, \, \, 
g_\rho(1,\rindep)=\frac{g_\rho(\xmac,\vmac+\rindep)}{g_\rho(\xmac,\vmac)} \, .
$$
\item $f_\rho(1,\rindep) = g_\rho(1,-\rindep)$.
\item $f_\rho(1,\rindep)$ and $g_\rho(1,\rindep)$ are bounded above by $e^{-2 \vert \rindep \vert+o(\vert \rindep \vert)}$ as $\vert \rindep \vert  \to \infty$.
\item $\int_{\R} f_\rho(\xmac,\udummy) \, {\rm d} \udummy =  \int_{\R} g_\rho(\xmac,\udummy) \, {\rm d} \udummy$.
\end{enumerate}  
\end{proposition}{\em Remark.}  In signposting Brownian Boost in Subsection~\ref{s.solvingbb},  
we indicated that TINE are indexed by a battlefield value $\battlefield\in\R$,
the equilibrium profile for this general value being obtained from the zero-indexed profile by a shift of~$v$. It is Proposition~\ref{p.sfacts}(2)
that identifies the value $\battlefield = \battlefield(\xmac)$. Since $\xmac \mapsto \battlefield$ is a bijection~$(0,\infty) \to \R$, the flow index and battlefield parameterisations are equivalent.
However, translation invariance of equilibrium profiles is more naturally expressed via 
the game theoretic interpretation, using the standard normalization for solutions specified in Definition~\ref{d.bbstandard}, rather than by means of the default normalization we are presently using in this analytic discussion. We will elaborate at the end of the chapter.


{\bf Proof of Proposition~\ref{p.sfacts}(1,2).}  Write $J(\rindep) = S_\rho(\xmac,\rindep)^\rho$ and integrate the differential equation in Lemma~\ref{l.de} on $[0,\rindep]$.
Since $J(0) = \xmac^\rho$, we find that 
$$
\rindep =-\frac{1}{8\rho^2}\Big( H\big(J(\rindep)\big) - H(\xmac^\rho) \Big) \, , 
$$
for  $H(z):=z+2\log z - z^{-1}$. 

The function   $H$
is an increasing bijection from $(0,\infty)$ to $\R$;
thus, so is $\xmac\mapsto H(\xmac^\rho)$, and also $\xmac\mapsto\vmac(\xmac)$, as we claim in the second part.

Now we  treat the special case in the first part, by taking $\xmac = 1$. Since $H(1) = 0$, we have 
$-8\rho^2\,\rindep =H(J(\rindep))$. 
The bijection~$H$  satisfies $H(1/z)=-H(z)$. 
We learn that $8\rho^2 \rindep$ equals both $-H(J(\rindep))$ and $H(J(-\rindep))$. So $ H(J(-\rindep)) = -H(J(\rindep)) = H\big(J(\rindep)^{-1}\big)$ whence $J(-\rindep) = J(\rindep)^{-1}$ since $H$ is invertible. Taking the $\rho$\textsuperscript{th} root yields $S_\rho(1,-\rindep) = S_\rho(1,\rindep)^{-1}$.

Rewriting the last display, $- 8 \rho^2 \rindep = H\big( S_\rho(\xmac,\rindep)^\rho \big) - H(\xmac^\rho)$. Since $H(1) = 0$, the unique solution~$\vmac = \vmac(\xmac)$ of $8 \rho^2 \vmac = H(\xmac^\rho)$
(which is the value identified in the second part of the proposition)
 is that time for which $S_\rho(\xmac,\vmac) =1$ (as we seek to prove in that part). 
Now $\rindep \to S_\rho(1,\rindep)$ and $\rindep \to S_\rho(\xmac,\vmac+\rindep)$ solve the initial-value problem stated in Lemma~\ref{l.de}. The uniqueness of the solution to this problem implied by this lemma shows that these two functions mapping $\R$ to $(0,\infty)$ are equal. 

We may now complete the proof of the first part by noting that
$$
 S_\rho(\xmac,-\rindep) = S_\rho(1,-\rindep - \vmac(\xmac)) = S_\rho(1,\rindep+\vmac(\xmac))^{-1} = S_\rho(1,\rindep - \vmac(\xmac^{-1}))^{-1} = S_\rho(\xmac^{-1},\rindep)^{-1} \, ,
$$
where the first and last equalities arise from the just obtained equality of functions applied for $\xmac$ and $\xmac^{-1}$. The second equality is an instance of  $S_\rho(1,-\rindep) = S_\rho(1,\rindep)^{-1}$, while the third is due to $\vmac(\xmac^{-1}) = -\vmac(\xmac)$, a fact seen from $8 \rho^2 \vmac(\xmac) = H(\xmac^\rho) = - H(\xmac^{-\rho}) = -8 \rho^2 \vmac(\xmac^{-1})$.

We use the representation~(\ref{e.fhrep}) and~(\ref{e.phifg}) in the first and last equalities as we write, with $\vmac = \vmac(\xmac)$,
\begin{eqnarray*}
\frac{f_\rho(\xmac,\vmac+\rindep)}{f_\rho(\xmac,\vmac)}
& = & \exp  \left\{2\int_\vmac^{\vmac+\rindep}\Phi_f\big(S_\rho(\xmac,\udummy)\big)\, {\rm d} \udummy \right\}
 \, = \, \exp \left\{2\int_0^{\rindep}\Phi_f\big(S_\rho(\xmac,\vmac+\udummy)\big)\, {\rm d}  \udummy \right\} \\
 & = & 
 \exp \left\{2\int_0^{\rindep}\Phi_f\big(S_\rho(1,\udummy)\big) \dd \udummy \right\}
= f_\rho(1,\rindep) \, ,
\end{eqnarray*}
the penultimate equality due to  $S_\rho(1,\rdummy) = S_\rho(\xmac,\vmac+\rdummy)$. Thus we obtain the second part of the proposition in regard to~$f$; the very similar argument for~$g$ is omitted.

{\bf (3).} Regarding $\Phi_f$ and $\Phi_g$ as functions on $(0,\infty)$, we have 
\[
\Phi_f(s)=1-\frac{2}{(1+s^\rho)^2}\Big(1+(1-\rho)s^\rho\Big),\qquad
\Phi_g(s)=1-\frac{2}{(1+s^\rho)^2}\Big((1-\rho)s^\rho+s^{2\rho}\Big) \, ,
\]
which satisfy
\begin{equation}\label{e.phisymmetry}
\Phi_g(1/s)=\Phi_f(s) \, ,
\end{equation}
since
$\Phi_g(1/s)
=1-\frac{2\big((1-\rho)s^{-\rho}+s^{-2\rho}\big)}{(1+s^{-\rho})^2}
=1-\frac{2\big((1-\rho)s^\rho+1\big)}{(1+s^\rho)^2}
=\Phi_f(s)$.
Using again the expressions~(\ref{e.fhrep}) and~(\ref{e.phifg}), 
\begin{eqnarray*}
\log g_\rho(1,-\rindep)
 & = &  -2\int_0^{-\rindep}\Phi_g\big(S_\rho(1,\udummy)\big)\, {\rm d} \udummy
=2\int_0^{\rindep}\Phi_g\big(S_\rho(1,-\udummy)\big)\, {\rm d} \udummy \\
&  = &  2\int_0^{\rindep}\Phi_f\big(S_\rho(1,\udummy)\big)\,{\rm d} \udummy =\log f_\rho(1,\rindep) \, ,
\end{eqnarray*}
where $\Phi_g\big(S_\rho(1,-\udummy)\big)=\Phi_f\big(S_\rho(1,\udummy)\big)$ is due to Proposition~\ref{p.sfacts}(1) and~(\ref{e.phisymmetry}). Exponentiating, we obtain the sought statement. 

{\bf (4).} As the solution to the differential equation in Definition~\ref{d.fg}, $S_\rho(\xmac,u) > 0$ is readily seen to converge to zero and infinity in the respective limits of large positive and negative~$u$.
So 
$$
\Phi_f(S_\rho(\xmac,u)) \to \begin{cases} 1 &  \textrm{as} \, \, \, u \to -\infty \\ -1   &  \textrm{as}  \, \, \,  u \to \infty   \end{cases}
\, \, \, \, \textrm{and} \, \, \, \, 
\Phi_g(S_\rho(\xmac,u)) \to \begin{cases} -1 &  \textrm{as}  \, \, \,  u \to -\infty \\ 1   &  \textrm{as}  \, \, \,  u \to \infty   \end{cases} \, .
$$ 
Note that the convention $\int_a^b f = - \int_b^a f$ is in force as we interpret~(\ref{e.fhrep}) and~(\ref{e.phifg}).
We see that $- \vert \rindep \vert^{-1} \log f_\rho(1,\rindep)$ and $- \vert \rindep \vert^{-1}  \log g_\rho(1,\rindep)$ converge to $2$, as $\rindep$ tends to both minus and plus infinity. This yields the sought statement. 

{\bf (5).} First note the special case when $\xmac=1$: $\int_{\R} f_\rho(1,\udummy) \, {\rm d} \udummy = \int_{\R} g_\rho(1,\udummy) \, {\rm d} \udummy$. This is due to the symmetry and integrability offered by the preceding two parts.

We make use of the special case in asserting the middle equality as we write 
\[
 \frac{1}{f_\rho(\xmac,\vmac)} \int_{\R} f_\rho(\xmac,\udummy)\, {\rm d} \udummy  = \int_{\R} f_\rho(1,\udummy)\, {\rm d} \udummy =
\int_{\R} g_\rho(1,\udummy)\, {\rm d} \udummy = \frac{1}{g_\rho(\xmac,\vmac)} \int_{\R} g_\rho(\xmac,\udummy)\, {\rm d} \udummy \, .
\]
Here,  Proposition~\ref{p.sfacts}(2) furnishes $\vmac = \vmac(\xmac)$, and the other displayed 
equalities are obtained by integrating the identities in this result over~$\R$.

As noted in Theorem~\ref{t.fg}, $g_\rho(\xmac,\bullet)=  f_\rho(\xmac,\bullet) S_\rho(\xmac,\bullet)$. Hence,
\[
\frac{\int_{\R} g_\rho(\xmac,\udummy)\, {\rm d} \udummy}{\int_{\R} f_\rho(\xmac,\udummy)\, {\rm d} \udummy} \;=\; \frac{g_\rho(\xmac,\vmac)}{f_\rho(\xmac,\vmac)} = S_\rho(\xmac,\vmac) \, .
\]
But by Proposition~\ref{p.sfacts}(2),  $S_\rho(\xmac,\vmac) =1$: so the integrals are equal.  \qed

We turn to proving the high~$\vert u \vert$ asymptotics of $f_\rho(1,u)$ and~$g_\rho(1,u)$, thereby refining Proposition~\ref{p.sfacts}(4), and concomitant results for the stake functions $a_\rho(1,u)$ and~$b_\rho(1,u)$.

{\bf Proof of Proposition~\ref{p.fgab}.} As remarked after Lemma~\ref{l.de}, $J(u)=S_\rho(1,u)^\rho$ satisfies $J(u) \sim \frac{1}{8 \rho^2 u}$ for $u \gg 0$.
In the representations~(\ref{e.fhrep}) and~(\ref{e.phifg}), $J$ enters in the role of~$S$ in the functions~$\Phi_f$ and~$\Phi_g$; for small~$J$,
$\Phi_f = -1 + 2(1+\rho)J + O(J^2)$ and $\Phi_g = 1 - 2(1-\rho)J + O(J^2)$. 
Using the high-$\rindep$ asymptotics \(\int_0^u J(\rdummy) \dd \rdummy \sim \frac{1}{8\rho^2}\log u\) and \(\int_0^u J(\rdummy)^2 \dd \rdummy = O(1)\) in these representations, we obtain
\[
\log f_\rho(1,u) = -2u + \frac{1+\rho}{2\rho^2}\log u + O_\rho(1),
\quad
\log g_\rho(1,u) = -2u + \frac{1-\rho}{2\rho^2}\log u + O_\rho(1) \, ,
\]
with continuous dependence on $\rho \in (0,\infty)$ for the implied constants in the $O_\rho(1)$-terms;
whence the claimed asymptotics for $f_\rho(1,u)$ and $g_\rho(1,u)$ as $u \to \infty$.

By Definition~\ref{d.arhobrho}, $S_\rho(\xmac,u) = g_\rho(\xmac,u)/f_\rho(\xmac,u)$ (from Theorem~\ref{t.fg}), and $J(u)=S_\rho(1,u)^\rho$, we see that
$$
a_\rho(1,u) \, = \, 2\rho \, f_\rho(1,u) \frac{J(u)}{(1+ J(u))^2}
$$
and
$$
b_\rho(1,u) \, = \,  2\rho \, g_\rho(1,u)  \frac{J(u)}{(1+ J(u))^2}  \, .
$$
Since $J(u) \to 0$ as $u \to \infty$, 
$a_\rho(1,u) \sim 2\rho \, f_\rho(1,u) J(u)$ and 
$b_\rho(1,u) \sim 2\rho \, g_\rho(1,u) J(u)$. 
So  $J(u) \sim \frac{1}{8 \rho^2 u}$ yields the high-$u$ $a_\rho$- and $b_\rho$-asymptotics, with $\zeta_a = \zeta_f -1$ and $\zeta_b = \zeta_g - 1$, as claimed.

Consider now negative~$u$.
By Proposition~\ref{p.sfacts}(3), we may replace $u$ by $\vert u \vert$ in 
the expressions $f_\rho(1,u)$ and  $g_\rho(1,u)$ provided that we exchange their roles. In this way, the asymptotics as $u \to -\infty$ reduce to what we have proved, after the stated interchanges are made. \qed

We end the chapter by offering the explanation, promised after Proposition~\ref{p.sfacts}, that concerns the translation invariance of \bbrho-equilibrium profiles
and the indexing role for the battlefield value $\battlefield = \battlefield(\xmac)$.

Noting that  $\battlefield(1)=0$, we see from the 
the relation $S_\rho(\xmac,\battlefield(\xmac))=1$ that default solutions satisfy 
\[
\big(f_\rho(\xmac,\bullet),g_\rho(\xmac,\bullet)\big)
=
C(\xmac)\,
\big(f_\rho(1,\bullet-\battlefield(\xmac)),\,
     g_\rho(1,\bullet-\battlefield(\xmac))\big)
\]
for some positive factor $C(\xmac)$.  
In the present context, however, the suitable normalization is the standard one,
 from Definition~\ref{d.bbstandard}, in accordance with our convention that Maxine's winning receipt in \bbrhospace equals one. The counterpart to the preceding display for standard solutions has $C(\xmac)$ identically equal to one, for it is this choice that ensures the property that $\int_\R f_\rho(\xmac,r) {\rm d} r$ is independent of $\xmac \in (0,\infty)$ which is required in the standard normalization.
 Note that the formulas in Definition~\ref{d.arhobrho} for equilibrium profiles are homogeneous of degree one in the inputs $(f,g)$.
 Writing $(a,b)$ 
 for the profile thus associated to the value $\xmac =1$, 
  we thus confirm that the profile indexed by $\xmac \in (0,\infty)$
 takes the form $\big(a(\bullet-\battlefield),\,b(\bullet-\battlefield)\big)$, with $\battlefield = \battlefield(\xmac)$.

  \chapter{The fine-mesh high-noise limit}\label{c.highnoise}

When \tlpkpspace is viewed in the fine-mesh, high-noise limit $\kappa \searrow 0$,
space is squeezed (by a factor of~$\kappa$), time is sped up (by a factor of~$\kappa^{-2}$),
and stakes are revalued (again by a factor of~$\kappa^{-2}$).
In our treatment, Brownian Boost is given a rigorous interpretation as the scaling limit of the Trail of Lost Pennies, with its time-invariant equilibria  governed by the $\rho$-Brownian Boost ODE pair.
This chapter is devoted to obtaining the scaled description and to drawing out some of its consequences.

\section{Two routes to Brownian Boost}

In this section, we present a four-part proposition concerning \abmnmacspace elements whose first two parts offer simple and useful stake formulas
and whose latter parts permit us to discuss competing routes to our analysis of $\rho$-Brownian Boost. After the discussion and proof, we will signpost the structure of Chapter~\ref{c.highnoise}.
  
 Recall that $M_i$ equals $m_{i-1,i+1} = m_{i+1} - m_{i-1}$ and  $N_i$ equals $n_{i+1,i-1} = n_{i-1} - n_{i+1}$.  
\begin{proposition}\label{p.abmnfacts}
Let $(a,b,m,n) \in$ \abmnkpmacspace and let $i \in \Z$.
\begin{enumerate}
\item  $a_i = \frac{\kappa \rho M_i^{1+\rho}N_i^\rho}{(M_i^\rho+N_i^\rho)^2}  \qquad \textrm{and} \qquad  b_i =  \frac{\kappa \rho M_i^\rho N_i^{1+\rho}}{(M_i^\rho+N_i^\rho)^2}$. 
\item  $\frac{a_i^\rho}{a_i^\rho + b_i^\rho} = \frac{M_i^\rho}{M_i^\rho+N_i^\rho}  \qquad \textrm{and} \qquad   \frac{b_i^\rho}{a_i^\rho + b_i^\rho} = \frac{N_i^\rho}{M_i^\rho+N_i^\rho}$. 
\end{enumerate} 
Write $\Delta_i m = m_{i+1} + m_{i-1} - 2m_i$ and $\Delta_i n = n_{i-1} + n_{i+1} - 2n_i$.
\begin{enumerate} 
 \setcounter{enumi}{2} 
\item 
$\kappa \rho  M_i^{1+\rho} N_i^\rho \, = \,  \tfrac{\kappa}{2} \cdot M_i \big( M_i^{2\rho} -  N_i^{2\rho} \big)
\, + \,  \tfrac{1}{2}\cdot  \big( M_i^\rho + N_i^\rho \big)^2 \Delta_i m$. 
\item $\kappa \rho  M_i^\rho N_i^{1+\rho} \, = \,  \tfrac{\kappa}{2} \cdot N_i \big( N_i^{2\rho} -  M_i^{2\rho} \big)
 \, + \,  \tfrac{1}{2}\cdot  \big( M_i^\rho + N_i^\rho \big)^2 \Delta_i n$. 
\end{enumerate} 
\end{proposition}
 Proposition~\ref{p.abmnfacts}(1,2) recasts the \abmnkpmacspace formulas to give explicit expressions for stakes, and records formulas for the players' win probabilities on stake turns.  
 
The equations in Proposition~\ref{p.abmnfacts}(3,4) are discrete counterparts to the $\rho$-Brownian Boost ODE pair~(\ref{e.fg.rho}), the pairs' respective elements 
identified under the correspondence of $m_{i,i+1}$ with $m' = f$ and $n_{i,i-1}$ with $n' = - g'$. Indeed, suppose that we permit the comparisons 
$\kappa^{-1}  m_{\kappa^{-1}u -1,\kappa^{-1}u}  = f(u) + O(\kappa)$
and 
$\kappa^{-1}  n_{\kappa^{-1}u,\kappa^{-1}u -1}  = -g(u) + O(\kappa)$ and their corollaries 
$\deltam_{\kappa^{-1}u}  = \kappa^2 f'(u) + O(\kappa^3)$, 
 $\deltan_{\kappa^{-1}u} = -\kappa^2 g'(x) + O(\kappa^3)$, 
 $M_{\kappa^{-1}u} = 2 \, m_{\kappa^{-1}u -1,\kappa^{-1}u}  + O(\kappa^2)$
and  $N_{\kappa^{-1}u} =  2 \, n_{\kappa^{-1}u,\kappa^{-1}u -1} + O(\kappa^2)$. 
Then on dividing the Proposition~\ref{p.abmnfacts}(3,4) equations by~$2^{2\rho}\kappa^{2(1+\rho)}$,
 we would learn that $f$ and $g$ satisfy the ODE pair~(\ref{e.fg.rho}) up to an $O(\kappa)$ error that must vanish since $\kappa > 0$ may tend to zero.
 Suitably elaborated, such an approach would lead to a rigorous discrete counterpart to the analysis of Brownian Boost offered in  Section~\ref{s.bbrho} wherein~(\ref{e.fg.rho}) was heuristically derived.
 
 So Proposition~\ref{p.abmnfacts}(3,4)  could be used on a route to showing that low-$\kappa$ \abmnkpmacspace solutions are governed by equations solving the Brownian Boost ODE pair. If we took this route, we might then exploit the record of solutions to the ODE pair in Theorem~\ref{t.fg} to describe explicitly \abmnkpmacspace solution asymptotics  as $\kappa \searrow 0$.
 
 However, we prefer to reach such conclusions by following a slightly longer path that we hope offers a more satisfying prospect on the conceptual relationship between low-$\kappa$ \abmnmacspace and Brownian Boost. We will show in Proposition~\ref{p.slowkappa}  how $S_\rho$, the solution of the ODE in Lemma~\ref{l.de}, gives a scaled description of suitably speeded iterates of the positive-$\kappa$ $s$-map that sends $\phi_0$ to $\phi_1$.  Our representation of the components of \abmnmacspace solutions as sums of products in Theorem~\ref{t.defexplicit}
 will then respond to the rapid-time scaling of $s_i$ iterates to the $S_\rho$-flow, with the product of many terms nearly equal to one leading to an integral of exponentials. In this way, the representations of $f_\rho$ and $g_\rho$ in Definition~\ref{d.fg} will emerge directly, in Proposition~\ref{p.mnlowkappa}, which is a detailed version of the stake-function asymptotic Theorem~\ref{t.lowkappasde}(1).

So in proofs we will make no use of  Proposition~\ref{p.abmnfacts}(3,4). These results offer comparison to Brownian Boost 
at the level of equations; our proofs will do so in the sense of solutions, by 
 monitoring the explicit positive-$\kappa$ solutions and showing how they track their Brownian Boost counterparts. 

{\bf Proof of Proposition~\ref{p.abmnfacts}(1,2).}
Use of the shorthand $*_{i,j} = *_j - *_i$ for $* \in \{m,n\}$ continues.
Analogous to $a = \frac{a^\rho - b^\rho}{a^\rho + b^\rho} m' + m''/2$ and to  $b = \frac{a^\rho - b^\rho}{a^\rho + b^\rho} n' + n''/2$ in~(\ref{e.am}) and~(\ref{e.bn}) are the equations
$$
 a_i \, = \, - \frac{a_i^\rho}{a_i^\rho+b_i^\rho}\kappa \cdot m_{i,i+1} \, - \, \frac{b_i^\rho}{a_i^\rho+b_i^\rho}  \kappa  \cdot  m_{i-1,i}   \, + \, \frac{1-\kappa}{2} \Delta_i m
$$
and
$$
 b_i \, = \,  \frac{b_i^\rho}{a_i^\rho+b_i^\rho}\kappa  \cdot  n_{i,i-1}  \, - \, \frac{a_i^\rho}{a_i^\rho+b_i^\rho}  \kappa  \cdot  n_{i+1,i}   \, + \, \frac{1-\kappa}{2} \Delta_i n \, ,
$$
given by  rearranging $\abmnmac(1)$ and  $\abmnmac(2)$ with index $i$.

We seek a counterpart to~(\ref{e.am}). Rearranging the above gives
\begin{equation}\label{e.akappa.rho}
 a_i \, = \, - \frac{b_i^\rho}{a_i^\rho+b_i^\rho}\kappa \cdot m_{i-1,i+1} \, + \,  \kappa  \cdot  m_{i,i+1}   \, + \, \frac{1-\kappa}{2} \Delta_i m
 \end{equation}
and
\begin{equation}\label{e.bkappa.rho}
 b_i \, = \,  - \frac{a_i^\rho}{a_i^\rho+b_i^\rho}\kappa  \cdot  n_{i+1,i-1}  \, - \,  \kappa  \cdot  n_{i,i-1}   \, + \, \frac{1-\kappa}{2} \Delta_i n \, . 
 \end{equation}
Differentiating these respective identities partially with respect to $a_i$ and $b_i$ and rearranging,
$$
 \kappa \rho \cdot b_i^\rho a_i^{\rho -1} m_{i-1,i+1} \, = \,  \kappa \rho \cdot a_i^\rho b_i^{\rho -1} m_{i-1,i+1} \, = \, \big(a_i^\rho+b_i^\rho\big)^2 \, .  
$$
Recall that $M_i = m_{i-1,i+1}$ and  $N_i = n_{i+1,i-1}$. We find then that $M_i/N_i = a_i/b_i$. 
Abbreviating, this yields 
\begin{equation}\label{e.abformulas.rho}
  a_i = \frac{\kappa \rho M_i^{1+\rho}N_i^\rho}{(M_i^\rho+N_i^\rho)^2}  \qquad \textrm{and} \qquad  b_i =  \frac{\kappa \rho M_i^\rho N_i^{1+\rho}}{(M_i^\rho+N_i^\rho)^2} \, , 
\end{equation}
which is Proposition~\ref{p.abmnfacts}(1) and
from which we learn that
$$
 \frac{a_i^\rho}{a_i^\rho + b_i^\rho} = \frac{M_i^\rho}{M_i^\rho+N_i^\rho}  \qquad \textrm{and} \qquad   \frac{b_i^\rho}{a_i^\rho + b_i^\rho} = \frac{N_i^\rho}{M_i^\rho+N_i^\rho} 
$$
or Proposition~\ref{p.abmnfacts}(2).
 
 {\bf (3,4).}
 Returning to~(\ref{e.akappa.rho}) and~(\ref{e.bkappa.rho}) with the expressions~(\ref{e.abformulas.rho}), and multiplying both of the resulting equations  by $\big( M_i^\rho + N_i^\rho \big)^2$, we obtain 
\begin{equation}\label{e.kappanmm.rho}
\kappa \rho  M_i^{1+\rho} N_i^\rho \, = \, \kappa M_i^\rho \big( M_i^\rho + N_i^\rho \big)  m_{i,i+1} \, - \,  \kappa N_i^\rho \big( M_i^\rho + N_i^\rho \big)  m_{i-1,i}
 \, + \, \tfrac{1-\kappa}{2}\cdot  \big( M_i^\rho + N_i^\rho \big)^2 \Delta_i m
\end{equation}
and
\begin{equation}\label{e.kappannm.rho}
\kappa \rho M_i^\rho N_i^{1+\rho}  \, = \, \kappa N_i^\rho \big( M_i^\rho + N_i^\rho \big)  n_{i,i-1} \, - \,  \kappa M_i^\rho \big( M_i^\rho + N_i^\rho \big)  n_{i+1,i}
 \, + \, \tfrac{1-\kappa}{2} \cdot \big( M_i^\rho + N_i^\rho \big) ^2  \Delta_i n \, .
\end{equation}
The facts $2m_{i,i+1} = M_i + \Delta_i m$ and  $2m_{i-1,i} = M_i - \Delta_i m$ respectively imply that
$$
\kappa M_i^\rho \big( M_i^\rho + N_i^\rho \big)  m_{i,i+1} \, = \, \frac{\kappa}{2} M_i^{1+\rho} \big( M_i^\rho + N_i^\rho \big) \,  + \,  \frac{\kappa \Delta_i m}{2}
 \big( M_i^\rho + N_i^\rho \big) M_i^\rho
$$
and
$$
\kappa N_i^\rho \big( M_i^\rho + N_i^\rho \big)  m_{i-1,i} \, = \, \frac{\kappa}{2} M_i N_i^\rho \big( M_i^\rho + N_i^\rho \big)  \, - \, \frac{\kappa \Delta_i m}{2}
 \big( M_i^\rho + N_i^\rho \big) N_i^\rho \, .
$$
Taking the difference of these equations, we may substitute the outcome into~(\ref{e.kappanmm.rho}), thereby finding that
\begin{equation}\label{e.kappanmm.rho.new}
\kappa \rho  M_i^{1+\rho} N_i^\rho \, = \,  \tfrac{\kappa}{2} \cdot M_i \big( M_i^{2\rho} -  N_i^{2\rho} \big)
\, + \,  \tfrac{1}{2}\cdot  \big( M_i^\rho + N_i^\rho \big)^2 \Delta_i m \, .
\end{equation}
where a cancellation $\alpha - \alpha = 0$ with $\alpha =  \tfrac{\kappa}{2} \cdot \Delta_i m \big( M_i^\rho +  N_i^\rho \big)^2$ has simplified the right-hand side. 
Thus we obtain Proposition~\ref{p.abmnfacts}(3).

Similarly, $2n_{i,i-1} = N_i + \Delta_i n$ and  $2n_{i+1,i} = N_i - \Delta_i n$ imply that
$$
\kappa N_i^\rho \big( M_i^\rho + N_i^\rho \big)  n_{i,i-1} \, = \, \frac{\kappa}{2} N_i^{1+\rho} \big( M_i^\rho + N_i^\rho \big)  \, + \,  \frac{\kappa \Delta_i n}{2}
 \big( M_i^\rho + N_i^\rho \big) N_i^\rho
$$
and
$$
\kappa M_i^\rho \big( M_i^\rho + N_i^\rho \big)  n_{i+1,i} \, = \, \frac{\kappa}{2} M_i^\rho N_i \big( M_i^\rho + N_i^\rho \big)  \, - \, \frac{\kappa \Delta_i n}{2}
 \big( M_i^\rho + N_i^\rho \big) M_i^\rho \, ,
$$
which substituted into~(\ref{e.kappannm.rho}) yield
\begin{equation}\label{e.kappannm.rho.new}
\kappa \rho  M_i^\rho N_i^{1+\rho} \, = \,  \tfrac{\kappa}{2} \cdot N_i \big( N_i^{2\rho} -  M_i^{2\rho} \big)
 \, + \,  \tfrac{1}{2}\cdot  \big( M_i^\rho + N_i^\rho \big)^2 \Delta_i n \, ,
\end{equation}
where  the cancellation $\zeta - \zeta = 0$ with $\zeta =  \tfrac{\kappa}{2} \cdot \Delta_i n \big( M_i^\rho +  N_i^\rho \big)^2$ has been made. This proves Proposition~\ref{p.abmnfacts}(4). \qed

There are five further sections in this chapter.
In the next three, we compare one-step of the application $s:\phi_0 \mapsto \phi_1$ to a suitably short passage along the $S_\rho$-flow; infer how a $\kappa^{-1}$-speeding of time leads to a description via this flow; and prove as a consequence the stake-asymptotic Theorem~\ref{t.lowkappasde}(1).
In the two further sections, we prove Theorem~\ref{t.lowlambdamax} on the approach of $\lambdamax$ to one; and derive the scaled gameplay Theorem~\ref{t.lowkappasde}(2).

\section{The action of $s$ mimics a $\kappa$-length ride on the $S_\rho$-flow}

In the discrete context, the map $s:\phi_0 \mapsto \phi_1$ sends the central ratio of any
\abmnkpmacspace element to the counterpart quantity associated with the left shift of this
element by one unit. The integer lattice is scaled by a factor of~$\kappa$ in the transition
from the discrete to the continuum. Counterpart to the $s$-argument $x = \phi_0 \in (0,\infty)$ in the scaled setting is the
flow index $\xmac \in (0,\infty)$, as noted after Definition~\ref{d.deltai}. Recalling the flow
$S_\rho(\xmac,\bullet)$, specified as the solution of an ODE in Definition~\ref{d.fg}, we regard
the  time-$\kappa$ evolution   $\xmac = S_\rho(\xmac,0) \mapsto S_\rho(\xmac,\kappa)$ as counterpart
to a single application of the map~$s$.

The purpose of this section is to make this correspondence quantitative.
Lemma~\ref{l.smallkappa.phi} provides a sharp single-step comparison between the discrete map
$s$ and the   time-$\kappa$ evolution  of the flow $S_\rho$, showing that these two operations are close
to each other.

More precisely, the lemma determines the linear coefficient in the small-$\kappa$ expansions
of $\phi_0$ and $\phi_1$ and controls the $O(\kappa^2)$ terms. It thereby yields such an
estimate for the map $s:\phi_0 \mapsto \phi_1$; and it furnishes a counterpart for the   time-$\kappa$ evolution  of $S_\rho$. This prepares for a comparison of iterates of
$s$ with progression along the $S_\rho$-flow.

The choice of domain variables $x$ and $\xmac$ in Lemma~\ref{l.smallkappa.phi}(2,3) is
characteristic of the discrete and continuous contexts. Setting these variables
equal serves to make explicit the approximation discussed above.

\begin{lemma}\label{l.smallkappa.phi}
Let $(\rho,\beta) \in (0,\infty)^2$ and suppose that $\kappa > 0$ is less than $\min \big\{ \rho^{-2} , (1+\rho)^{-1}/2  \big\}$.
\begin{enumerate}
\item
We have that
\begin{eqnarray*}
  \phi_0(\kappa,\rho,\beta) &= & \beta \;+\; 
     \frac{4\rho\,\beta^{1+\rho}}{(1+\beta^\rho)^2}\,\kappa \;+\; \beta \, \Theta_1(\kappa,\rho,\beta) \kappa^2   \qquad \textrm{and} \\
  \phi_1(\kappa,\rho,\beta) &= & \beta \;-\; 
     \frac{4\rho\,\beta^{1+\rho}}{(1+\beta^\rho)^2}\,\kappa \;+\; \beta \, \Theta_2(\kappa,\rho,\beta) \kappa^2 \, ,
\end{eqnarray*}
where  $\vert \Theta_i(\kappa,\rho,\beta) \vert \, \leq \, 2\rho(1+\rho)$ for $i \in \{1,2\}$.
\item 
We also have that
$$
 s(x) \, = \, x - \frac{8\rho\, x^{1+\rho}}{(1+x^\rho)^2}\,\kappa  \, + \,  x \,  \Theta_3(\kappa,\rho,x) \kappa^2 
$$
with $\vert \Theta_3(\kappa,\rho,x) \vert \leq 52 \rho (1+\rho)^3$. 
\end{enumerate}
Now suppose only that $(\rho,\beta) \in (0,\infty)^2$. 
\begin{enumerate}
    \setcounter{enumi}{2} 
\item For $\xmac \in (0,\infty)$ and $\kappa > 0$,
\[
 S_\rho(\xmac,\kappa) \,=\, \xmac - \frac{8\rho\,\xmac^{1+\rho}}{(1+\xmac^\rho)^2}\,\kappa \,+\, \xmac\,\Theta_4(\kappa,\rho,\xmac)\,\kappa^2
\]
with $\vert \Theta_4(\kappa,\rho,\xmac)\vert \leq 64 \rho^2 (1+\rho)$.
\end{enumerate}
\end{lemma}
{\bf Proof.} Recall that, for the given value of the pair $(\kappa,\rho)$, $\phi_0$ and $\phi_1$ are functions of $\beta \in (0,\infty)$ specified by the formulas in Definition~\ref{d.fourfunctions}.
The variable $x = \phi_0$ is the input for the map~$s$. Our hypotheses ensure that the $(\kappa,\rho)$-pair lies in the region~$W$ specified by $\rho^2 \kappa \leq 1$  in~(\ref{e.weakregion}), so that Lemma~\ref{l.incphi}(2)
implies that the relationship $x \longleftrightarrow \beta$ is bijective on~$(0,\infty)$. 

{\bf  (1).}
We may express $\phi_0=\beta N/D$, with
\begin{equation}\label{e.ndnot}
N=(1-\kappa)\beta^{2\rho}+2(1+\rho\kappa)\beta^\rho+1+\kappa,
\quad
D=(1-\kappa)\beta^{2\rho}+2(1-\rho\kappa)\beta^\rho+1+\kappa.
\end{equation}
Writing $N = N_0 + \kappa N_1$ and $D = D_0 + \kappa D_1$, we have $N_0=D_0=(1+\beta^\rho)^2$, 
$$
N_1 = - \beta^{2\rho} + 2\rho \beta^\rho +1 \, \, \, \, 
\textrm{and}  \, \, \, \, D_1 = - \beta^{2\rho} - 2\rho \beta^\rho +1 \, .
$$
 Note that 
$$
\frac{N}{D} \;=\; 1 + \kappa \,\frac{N_1-D_1}{D_0} \;+\; \kappa^2 R(\kappa,\beta) \, ,
$$
where $R(\kappa,\beta) \;=\; \frac{D_1(D_1 - N_1)}{D_0(D_0+\kappa D_1)}$.
We have 
$$
\frac{N_1-D_1}{D_0} \, = \, \frac{4\rho \beta^\rho}{(1+\beta^\rho)^2} \, .
$$
It remains to control the remainder $R$. Write $t = \beta^\rho$. From the forms displayed above, we see that  $|N_1|$ and $|D_1|$ are at most $(1+\rho)(1+t^2)$, for any $\rho \in (0,\infty)$;
 while $D_0=(1+t)^2\ge 1+t^2$.
Hence $|N_1|/D_0 \leq 1+\rho$. We find then that, when $\kappa \leq (1+\rho)^{-1}/2$, $\vert D_0 + \kappa D_1 \vert$ is at least $\vert D_0 \vert/2$. We also have $\vert N_1 - D_1 \vert/D_0 \leq 4\rho \frac{t}{(1+t)^2} \leq \rho$. Consequently, 
$$
R(\kappa,\beta) \, \leq \,  2 \frac{\vert D_1 \vert \cdot \vert D_1   -  N_1 \vert}{D_0^2} \, \leq \, 2\rho (1+\rho) 
$$
under this circumstance.
This completes the proof of the assertion made in regard to $\phi_0$. For $\phi_1$, the same decomposition applies. The coefficients $N_0$ and $D_0$ remain unchanged while the first-order coefficients are negated: $N_1 \to - N_1$ and $D_1 \to - D_1$.  Consequently,  the linear term in $\kappa$ in the obtained formula for $\phi_1$ flips sign. All the bounds on absolute value remain valid,  so the uniform estimate on $\Theta_2$ follows.
This completes the proof of Lemma~\ref{l.smallkappa.phi}(1).

{\bf (2).} 
As we noted at the outset of the proof,
$\phi_0$ and $\phi_1$ are increasing bijections of $(0,\infty)$ under our hypothesis; thus, the map $s:(0,\infty) \to (0,\infty)$, which by definition sends $\phi_0$ to $\phi_1$, is well defined. 

Note that 
$$
\phi_1(\kappa,\rho,\beta) = \phi_0(\kappa,\rho,\beta) \, - \,   
     \frac{8\rho\,\beta^{1+\rho}}{(1+\beta^\rho)^2}\,\kappa \;+\; \beta \,  \Theta(\kappa,\rho,\beta) \kappa^2 \, ,
$$
where $\Theta = \Theta_2 - \Theta_1$ satisfies $\vert \Theta \vert \leq 8(1+\rho)^2$. Writing $x =  \phi_0(\kappa,\rho,\beta)$, we find that
$$
s(x) = x \, - \,   
     \frac{8\rho\,\beta^{1+\rho}}{(1+\beta^\rho)^2}\,\kappa \;+\; \beta \,  \Theta(\kappa,\rho,\beta) \kappa^2 \, .
$$
But  $\frac{8\rho\,\beta^{1+\rho}}{(1+\beta^\rho)^2} =  \frac{8\rho x^{1+\rho}}{(1+x^\rho)^2} + R$ where
$\vert R \vert \leq 8 \rho \vert \beta - x \vert \Dslope$, with 
$\Dslope = \sup \big\{ \big\vert \tfrac{{\rm d}}{{\rm d}z} \tfrac{z^{1+\rho}}{(1+z^\rho)^2} \big\vert : z \in (0,\infty) \big\}$.
By Lemma~\ref{l.smallkappa.phi}(1),
$$
  \vert \beta - x \vert  \leq   \frac{4\rho\,\beta^{1+\rho}}{(1+\beta^\rho)^2}\,\kappa \;+\; \beta \, \vert \Theta_1(\kappa,\rho,\beta) \vert \kappa^2 \, .
$$
So 
\begin{eqnarray*}
s(x) & = & x \, - \,   
     \frac{8\rho\, x^{1+\rho}}{(1+x^\rho)^2}\,\kappa \;+\;  8\rho \Dslope \bigg(  \frac{4\rho\,\beta^{1+\rho}}{(1+\beta^\rho)^2}\,\kappa \;+\; \beta \, \vert \Theta_1(\kappa,\rho,\beta) \vert \kappa^2 \bigg) \kappa \, + \,  \beta \,  \Theta(\kappa,\rho,\beta) \kappa^2 \\
 & = & x - \frac{8\rho\, x^{1+\rho}}{(1+x^\rho)^2}\,\kappa  \, + \,  \beta \,  \Theta_3(\kappa,\rho,x) \kappa^2 
\end{eqnarray*}
where we are able to take $\Theta_3$ to be a function of $x$, rather than of $\beta$, because of the bijective relationship between the two variables;
since $\kappa \leq 1$, we have
$$
\vert \Theta_3 \vert \leq \vert \Theta \vert + 32 \rho^2 \Dslope \Dsup + 8\rho \Dslope \vert \Theta_1 \vert 
$$
with  $\Dsup = \sup \big\{ \big\vert  \tfrac{z^\rho}{(1+z^\rho)^2} \big\vert : z \in (0,\infty) \big\}$.
Since $\Dsup \leq 1/4$ and $\Dslope \leq 2(1+\rho)$, the latter right-hand side is at most
$$
 4\rho(1+\rho) + 32 \rho^2 \Dslope \Dsup + 16 \Dslope \rho^2(1+\rho) \leq  4\rho(1+\rho) + 16 \rho^2 (1+\rho) + 32\rho^2 (1+\rho)^2   \leq 52 \rho (1+\rho)^3 \, .
$$
In the notation from~(\ref{e.ndnot}) $x=\phi_0=\beta N/D$ with $D>0$ and $N-D=4\rho\kappa\beta^\rho\ge0$. Thus, $x\ge\beta$.
This bound permits us to replace the prefactor of $\beta$ by $x$ in the last right-hand term in the above expression for $s(x)$. Thus we obtain Lemma~\ref{l.smallkappa.phi}(2). 

{\bf (3).} Writing
$f(\xaux) = - \frac{8 \rho \xaux^{1+\rho}}{(1+\xaux^\rho)^2}$, recall that $S_\rho(\xmac,u)$ solves
$\frac{{\rm d}}{{\rm d}u} S_\rho(\xmac,u) = f(S_\rho(\xmac,u))$ with $S_\rho(\xmac,0) = \xmac$.
Observe that $f(\xaux) < 0$ for all $\xaux>0$, so $S_\rho(\xmac,u)$ is decreasing in $u$ and satisfies
$0 < S_\rho(\xmac,u) \le \xmac$ for all $u \ge 0$.
Since $f$ is Lipschitz, the Picard-Lindel\"of theorem~\cite[Theorem~$I.3.1$]{CoddingtonLevinson} implies that
 $S_\rho(\xmac,u)$ has the integral form
\[
S_\rho(\xmac,\kappa) = \xmac + \int_0^\kappa f(S_\rho(\xmac,\udummy))\,d\udummy \, .
\]
Since $S_\rho(\xmac,u) \le \xmac$, we have for $u \in [0,\kappa]$,
\[
|S_\rho(\xmac,u)-\xmac| = \left| \int_0^u f(S_\rho(\xmac,v)) \dd v \right| \le u \, \sup_{z \in [0,\xmac]} |f(z)| \le \kappa \cdot 8 \rho \xmac \, .
\]
By the mean-value theorem, $f(S_\rho(\xmac,u)) = f(\xmac) + f'(\ximac_u) \, (S_\rho(\xmac,u)-\xmac)$
for some \(\ximac_u \in [S_\rho(\xmac,u), \xmac] \subset [0, \xmac]\). Integrating, the remainder is
\[
R(\xmac,\kappa) := S_\rho(\xmac,\kappa) - \xmac - \kappa f(\xmac) = \int_0^\kappa f'(\ximac_u) \big(S_\rho(\xmac,\udummy) - \xmac\big)  \dd \udummy \, .
\]
Differentiating $f$, we readily see that $\vert f'(z) \vert \leq 8 \rho (1+\rho)/(1+z^\rho)^2$.
Using the preceding bound on \(|S_\rho(\xmac,u)-\xmac|\) alongside \(|f'(z)| \le 8 \rho (1+\rho)\), we find that
\[
|R(\xmac,\kappa)| \le \int_0^\kappa |f'(\ximac_\udummy)| \, |S_\rho(\xmac,\udummy)-\xmac|  \dd \udummy \le \kappa \cdot 8 \rho (1+\rho) \cdot (\kappa \cdot 8 \rho \xmac) = 64 \rho^2 (1+\rho) \xmac \, \kappa^2 \, .
\]
Then setting \(\Theta_4(\kappa,\rho,\xmac) = R(\xmac,\kappa)/(\xmac \kappa^2)\), we obtain
\[
S_\rho(\xmac,\kappa) = \xmac + \kappa f(\xmac) + \xmac \, \Theta_4(\kappa,\rho,\xmac) \, \kappa^2 \, ,
\]
with $\Theta_4(\kappa,\rho,\xmac)| \le 64 \rho^2 (1+\rho)$,
which completes the proof. \qed

\section{The scaled $s$-orbit tracks the $S_\rho$-flow}   
We presented precise hypotheses on $(\kappa,\rho)$-pairs  in Lemma~\ref{l.smallkappa.phi}.
However, in the conclusions we seek in this section, $\rho \in (0,1]$. Expressions such as error bounds are a little simpler when this condition is in force, and we apply it henceforth, occasionally remarking on how it may be relaxed.

A compact notation is useful to present our proposition linking the orbit and the flow.
For $\kappa \in (0,1]$, functions $h,h':\R \to (0,\infty)$ satisfy $h \stackrel{\kappa}{\simeq} h'$ provided that, for $z \in \R$  (and some positive $C_0$ and $C_1$), 
$$
\vert h(z) - h'(z) \vert \, \leq \, C_0 \, e^{C_1 \vert z \vert} \kappa \cdot \max \big\{ \vert h(z) \vert , \vert h'(z) \vert \big\} \, .
$$ 
\begin{proposition}\label{p.slowkappa}
Let $\rho \in (0,1]$ and $\xmac \in D$. As functions of the argument $\bullet \in \R$, we have that
$$
 s_{\lfloor \kappa^{-1}\bullet \rfloor}(\xmac) \, \stackrel{\kappa}{\simeq} \, S_\rho(\xmac,\bullet) \, .
$$
\end{proposition}
Our usual notation for the variable in the map $s$ and its iterates is~$x$. It is taken to be $\xmac$ here because the scaled $s$-iterate is being directly compared to the flow~$S_\rho(\xmac,\bullet)$ 
whose index~$\xmac$ governs the choice of notation.
We also mention that  $\rho \geq 1$ may be taken in the above result, provided that $(\kappa,\rho)$ lies in~$W$ as specified in~(\ref{e.weakregion}) 
and the $\stackrel{\kappa}{\simeq}$ notation is modified to permit $\rho$-dependent constants.

{\bf Proof of Proposition~\ref{p.slowkappa}.} For $C > 0$, let $\mc{I}_\rho(\kappa,C)$ denote the set of functions $h:(0,\infty) \to (0,\infty)$
such that
$$
 h(x) \, = \, x - \frac{8\rho \, x^{1+\rho}}{(1+x^\rho)^2}\kappa \, + \,  x  \, O(1)  \kappa^2 \, ,
$$
where $\vert O(1) \vert \leq C$ for all $x \in (0,\infty)$.

Developing a concept from Section~\ref{s.roleshift}, we say that 
a bijection $h:(0,\infty) \to (0,\infty)$ is {\em role-reversal symmetric} if its inverse satisfies $h^{-1}(x) = 1/h(1/x)$.
Indeed, Proposition~\ref{p.sminusone} shows this property for~$s$. The flow satisfies it also: by Proposition~\ref{p.sfacts}(1), $\xmac \to S_\rho(\xmac,\rindep)$
is role-reversal symmetric for any $\rindep > 0$.

By Lemma~\ref{l.smallkappa.phi}(2,3) (and $\rho \leq 1$), 
the maps from $(0,\infty)$ to $(0,\infty)$ given by  $x \to  s(x)$ and $\xmac \to S_\rho(\xmac,\kappa)$ belong to~$\mc{I}_\rho(\kappa,C)$ with  $C = 500$.

As such, the next result in essence delivers the proposition. A subscript~$i$ denotes the $i$\textsuperscript{th} iterate.
\begin{lemma}\label{l.handhprime}
Let $h$ and $h'$ belong to $\mc{I}_\rho(\kappa,C)$.
\begin{enumerate}
\item The iterate difference sequence satisfies the recursion
$$
 \big\vert h_{i+1}(x) - h_{i+1}'(x) \big\vert \leq (1 + 8C' \rho \kappa) 
 \big\vert h_i(x) - h_i'(x) \big\vert + 1000 \, \kappa^2  \max \big\{ \vert h_i(x) \vert , \vert h'_i(x) \vert \big\} 
$$
for $x \in (0,\infty)$.
Here $C'$ denotes
the supremum (which is readily seen to be finite) of the absolute value of the derivative of $\frac{z^{1+\rho}}{(1+z^\rho)^2}$ over~$(0,\infty)$.
\end{enumerate}
Suppose that $\kappa$ is at most a small universal positive constant. 
\begin{enumerate}
    \setcounter{enumi}{1} 
\item  For $x \in D$ and $i \in \N$,
 $\big\vert h_i(x) - h_i'(x) \big\vert \leq 2C   \kappa \, \exp \{ C_2 \kappa i \} \max \big\{ \vert h_i(x) \vert , \vert h'_i(x) \vert \big\}$.
 
\item Suppose further that $h$ and $h'$ are role-reversal symmetric. Then for $i \in \N$
$$
\big\vert h_{-i}(x) - h_{-i}'(x) \big\vert \leq 2C   \kappa \, \exp \{ C_2 \kappa i \} \Big( \min \big\{ \vert h_i(1/x) \vert , \vert h'_i(1/x) \vert \big\} \Big)^{-1} \, .
$$
 \end{enumerate}
\end{lemma}
To confirm that the proposition follows from the lemma, take 
$h(\bullet) =  s(\bullet)$ and $h'(\bullet) = S_\rho(\bullet,\kappa)$.
 For
$\rindep \in \R$, set $i = \lfloor \kappa^{-1} \rindep \rfloor$. Since $\xmac$ lies in the bounded central domain~$D$ and $\rho \leq 1$, Lemma~\ref{l.handhprime}(2) implies that  
$$
\lvert h_i(\xmac) - h_i'(\xmac) \rvert 
   \,\leq\, 2C \kappa \,  e^{C_2 \rindep} \, \max \big\{ \vert h_i(\xmac) \vert , \vert h'_i(\xmac) \vert \big\}  
   $$ 
   holds for positive integers~$i$, and  Lemma~\ref{l.handhprime}(3) delivers the same conclusion for negative~$i$. 
   Hence the desired $\stackrel{\kappa}{\simeq}$ relation
  holds 
with $C_0 = 2C$ and $C_1 = C_2$. \qed


{\bf Proof of Lemma~\ref{l.handhprime}(1).}
Set $\alpha_i =  h_{i+1}(x) - h_{i+1}'(x) - \big( h_i(x) - h_i'(x) \big)$, so that 
$$
 \alpha_i = h_{i+1}(x) - h_i(x) - \big(  h'_{i+1}(x) - h_i'(x)  \big) \, .
$$
Since $h,h' \in \mc{I}_\rho(\kappa,C)$, 
$$
 \alpha_i \, = \,  -8\kappa \rho \Bigg( \frac{h_i(x)^{1+\rho}}{\big( 1 +  h_i(x)^\rho \big)^2} -  \frac{h'_i(x)^{1+\rho}}{\big( 1 +  h'_i(x)^\rho \big)^2} \Bigg) \, 
 + \,    \big( \vert h_i(x) \vert + \vert h'_i(x) \vert \big) O(1) \kappa^2 \, , 
$$
where $\vert O(1) \vert \leq 450$. Hence,
$$
  \vert \alpha_i \vert \, \leq \, 8 \kappa \rho C' \big\vert h_i(x) - h_i'(x) \big\vert + O(1) \kappa^2  \max \big\{ \vert h_i(x) \vert , \vert h'_i(x) \vert \big\} \, ,
$$
where $C'$ is the stated derivative supremum. 
Using $\big\vert h_{i+1}(x) - h_{i+1}'(x)  \big\vert \leq \big\vert h_i(x) - h_i'(x)  \big\vert + \vert \alpha_i \vert$, we obtain the sought statement.

{\bf (2).} Set $\zeta_i = \big\vert h_i(x) - h'_i(x) \big\vert$, and note that $\zeta_0 = 0$. For  $C_2 > 0$ whose value will be later specified, 
we will induct on $i \in \N$ to show that  $\zeta_i \leq C e^{x \kappa C_2 C' i} \kappa \omega_i(x)$ where $\omega_i(x) = \max \big\{ \vert h_i(x) \vert , \vert h'_i(x) \vert \big\}$. By the inductive hypothesis {\rm IH}($i$) indexed by $i \in \N$, we find from the preceding part of the lemma that
$$
 \zeta_{i+1} \leq (1+ 8C' \rho \kappa) C e^{C_2 C' \kappa  x  i} \kappa \, \omega_i(x) + C_1  \kappa^2  \omega_i(x) 
$$
where $C_1$ is suitably high.  Here 
the right-hand side takes the form 
\[
C e^{x \kappa C_2 C'(i+1)} \kappa \, \omega_{i+1}(x) + \psi \, ,  \, \, \, \textrm{where}
\]
\[
\psi = C e^{C_2 C' \kappa  x  i}  \kappa \, \omega_i(x)\,\Big( 1+8C' \rho \kappa - e^{C_2 C' \kappa  x}  \tfrac{\omega_{i+1}(x)}{\omega_i(x)} \Big) 
   + C_1   \kappa^2 \omega_i(x) \, .
\]
Since $\psi \leq 0$ establishes \rm IH($i+1$), the inductive argument will be complete provided that we show the above right-hand side is at most zero, for which it suffices to prove a bound of the form 
$ \tfrac{\omega_{i+1}(x)}{\omega_i(x)}  \geq 1- B \kappa \rho$ with $B = B(C)$ alongside
\begin{equation}\label{e.needednegative}
C \big(   1+8C' \rho \kappa - (1 + C_2 C' \kappa  x )(1 - B \kappa) \big) + C_1 \kappa \, e^{-C_2 C' \kappa  x  i}  \leq 0 \, .
\end{equation}
We justify the lower bound on the ratio $\omega_{i+1}(x)/\omega_i(x)$ as follows.  Since $h\in\mc{I}_\rho(\kappa,C)$, each iterate satisfies
\[
h_{i+1}(x)
= h_i(x) - \frac{8 \rho\, h_i(x)^{1+\rho}}{(1 + h_i(x)^\rho)^2}\,\kappa
+ h_i(x)\,O(1)\,\kappa^2
\]
with $|O(1)|\le C$; and similarly for $h'_i(x)$. 
Since $\frac{z^\rho}{(1+z^\rho)^2}$ is bounded on $(0,\infty)$ and $\kappa$ is taken sufficiently small, it follows that
\[
|h_{i+1}(x)| \ge |h_i(x)|\big(1 - B\rho\kappa\big),
\]
and likewise for the $h'$-sequence. Taking the maximum yields the claimed lower bound on
$\omega_{i+1}(x)/\omega_i(x)$.

To obtain~(\ref{e.needednegative}), note that its left-hand side is at most
$$
-CC_2 C' \kappa \, x   +  \Big( 8C'C \rho  + BC + C_1 \Big) \kappa + C' C_2 C    B \kappa^2 x  \, .
$$
The displayed expression becomes negative with a suitably high choice of the constant~$C_2$.
To confirm this, note that $x$ lies in the central domain~$D$, so $x \geq d : = \inf D > 0$. Supposing (as we may) that $\kappa$ is at most $(2B)^{-1}$,
a choice of $C_2$ high enough that $C C_2 C' d/2 \geq 8 C' C \rho + BC + C_1$ works for our purpose. 
 In this way, we justify  the bound $\psi \le 0$, 
 and thus complete the inductive step.
 Since $x \in D$, we absorb the factor~$C'x$ in the argument of the exponential with an increase in the value of~$C_2$, and so obtain
   Lemma~\ref{l.handhprime}(2).

{\bf (3).} Noting that  
$h_{-i}(x) - h_{-i}'(x) \,  = \, \frac{h_i'(1/x) - h_i(1/x)}{h_i(1/x)h_i'(1/x)}$ by role-reversal symmetry, the result follows from Lemma~\ref{l.handhprime}(2) given the invariance of $D$ under the inversion $x \mapsto x^{-1}$. \qed

\section{Equilibria converge to the putative Brownian Boost counterparts as $\kappa$ vanishes}
 
 With the $s$-orbit run rapidly tracking the $S_\rho$-flow, we are ready to see how the product expressions leading to the explicit \abmnmacspace solutions in Theorem~\ref{t.defexplicit}
 may be recast as integrals of exponential functions. We  need to understand low-$\kappa$ asymptotics for the basic functions $c$ and $d$ from Definition~\ref{d.scd}(2) that enter into these products. 

\begin{lemma}\label{l.cdasymptotics}
When  $\kappa$ is supposed to be at most a universal positive constant,
\begin{equation}\label{e.ceqn}
 c(x)  \, = \, 2 + \kappa \cdot  2 \bigg( 1 - \frac{2\big(1+(1-\rho)x^\rho \big)}{(1+x^\rho)^2} \bigg) \, + \, O(\kappa^2)
\end{equation}
 and
\begin{equation}\label{e.deqn}
d(x)  \, =   \, 2 - \kappa \cdot 2 \bigg( 1 - \frac{2}{(1+x^\rho)^2} \Big( (1-\rho)x^\rho + x^{2\rho} \Big) \bigg)  \, + \,  O(\kappa^2) \, .
\end{equation}
\end{lemma}
{\bf Proof.} Recall that $c(x) = 1/\gamma(\kappa,\rho,\beta)$ where $x = \phi_0(\kappa,\rho,\beta)$.
From Definition~\ref{d.fourfunctions}, we thus have 
$$
c(x) \, = \,  \frac{2 (1+\beta^\rhomac)^2}{(1-\kappa)\beta^{2\rhomac}  + 2(1-\rhomac\kappa)\beta^\rhomac + 1 +\kappa} \, , \, \, \, \textrm{and}
$$
$$
d(x) \, = \, \frac{2(1+\beta^\rhomac)^2}{(1-\kappa)\beta^{2\rhomac} + 2(1+\rhomac\kappa) \beta^\rhomac + 1 + \kappa} \, .
$$

First we argue that  $\vert \beta - x \vert \leq O(1) \rho^2 \, x^{1-\rho} \kappa$.
 By Lemma~\ref{l.smallkappa.phi}(1),
$$
  \vert \beta - x \vert  \leq   \frac{4\rho\,\beta^{1+\rho}}{(1+\beta^\rho)^2}\,\kappa \;+\; \beta \, \vert \Theta_1(\kappa,\rho,\beta) \vert \kappa^2 \, .
$$
But  $\frac{4\rho\,\beta^{1+\rho}}{(1+\beta^\rho)^2} =  \frac{4\rho \, x^{1+\rho}}{(1+x^\rho)^2} + R$ where
$\vert R \vert \leq 4 \rho \vert \beta - x \vert \macD$, with $\macD = \sup \big\{ \big\vert \tfrac{{\rm d}}{{\rm d}z} \tfrac{z^{1+\rho}}{(1+z^\rho)^2} \big\vert : z \in (0,\infty) \big\}$.
Hence,
$$
  \vert \beta - x \vert  \leq   \frac{4\rho\,x^{1+\rho}}{(1+x^\rho)^2}\,\kappa \;+\;  
  4\rho \vert \beta - x \vert \macD \kappa +
  \beta \, \vert \Theta_1(\kappa,\rho,\beta) \vert \kappa^2 \, ,
$$
where recall that $\vert \Theta_1 \vert \leq 2\rho(1+\rho)$,
so that 
$$
  \vert \beta - x \vert  \leq  \bigg( \frac{4\rho\,x^{1+\rho}}{(1+x^\rho)^2}\,\kappa \;+\;  
  \beta \, \vert \Theta_1(\kappa,\rho,\beta) \vert \kappa^2 \bigg) \big( 1 - 4\rho \macD \kappa \big)^{-1} \, .
$$
Since $\beta  \leq x$,
we find  that, provided that $\kappa  \leq (8\macD)^{-1}$, 
$$
\vert \beta - x \vert \leq  \frac{4\rho\,x^{1+\rho}}{(1+x^\rho)^2} \kappa + D_0 \, x^{1-\rho} \rho^2 (1+\rho) \kappa^2
$$
for a universal constant $D_0 > 0$. We obtain  $\vert \beta - x \vert \leq O(1) \rho^2 \, x^{1-\rho} \kappa$ as sought.

Next note that $c(x) = H(\beta^\rho)$ where $H(u) = \frac{2(1+u)^2}{(1-\kappa)u^2 + 2(1-\rho \kappa) u + 1 + \kappa}$ satisfies
$$
 H(u) = \tfrac{2}{1 + \frac{1-2\rho u - u^2}{(u+1)^2}\kappa} = 2 + \kappa \cdot 2 \bigg( 1 - \tfrac{2 \big( (1-\rho)u+1\big)}{(u+1)^2}\bigg) + O \big( (\rho+1) \kappa^2\big) \, .
$$ 
Writing $v = u+1 \geq 1$ and  $\Delta  =(1-\kappa) v^2 + 2 \kappa (1-\rho) v + 2 \kappa \rho$, we have 
$H'(u) = 4 \kappa v \big( (1-\rho) v + 2 \rho  \big)\Delta^{-2}$.
Since $\kappa$ is at most a small positive constant, and $\rho \leq 1$, we see that $\vert H'(u) \vert = O(\kappa)$ for $u \in [0,\infty)$.

Note that
$c(x) = H(x^\rho) + \big( H(\beta^\rho) - H(x^\rho)\big)$ and
$$
\big\vert H(\beta^\rho) - H(x^\rho) \big\vert \leq \sup \vert H' \vert \cdot \vert \beta^\rho - x^\rho \vert \leq O(1) \kappa  \vert \beta - x \vert \rho x^{\rho-1} \leq
 O(1) \rho^2 \kappa^2 \, ,
$$
by $\beta \leq x$ and the $\vert H' \vert$ bound in the second inequality, and the $\vert \beta - x \vert$ bound  in the third. Thus,
$$
 c(x) = 2 + \kappa  \cdot 2 \bigg( 1 - \tfrac{2 \big( (1-\rho)x^\rho+1\big)}{(x^\rho+1)^2}\bigg) + O(1) \big( \rho^2 + 1 \big) \kappa^2 \, ,  
$$
which since $\rho \leq 1$ is the desired asymptotic for $c$. The formula~(\ref{e.deqn}) for~$d$ 
differs from that for $c$ in~(\ref{e.ceqn}) only in a change $\rho \to - \rho$ in the linear-in-$\kappa$ coefficient in the denominator. 
The form of the estimates in the resulting proof are unaffected by this change, and the claimed $d$-asymptotic results. \qed

We may now formulate and prove a technical development of the stake-asymptotics Theorem~\ref{t.lowkappasde}(1).  
\begin{proposition}\label{p.mnlowkappa}
Let $(\kappa,\rho) \in (0,1]^2$.
Write $m_i(x) = \mdefault_i(\kappa,\rho,x)$ (and use other like abbreviations) for the default solution. For $\xmac \in D$ and $\rindep \in \R$, we have that
\begin{eqnarray}
\kappa^{-1} m_{\lfloor \kappa^{-1 }\rindep\rfloor , \lfloor \kappa^{-1 }\rindep\rfloor +1}( \xmac ) & = &  f_\rho(\xmac,\rindep) \big( 1 + \kappa E_{\rindep}\big) \, \, \, \, \textrm{and} \label{e.massertion} \\
\kappa^{-1} n_{\lfloor \kappa^{-1 }\rindep\rfloor + 1, \lfloor \kappa^{-1 }\rindep\rfloor}( \xmac ) & = &   g_\rho(\xmac,\rindep) \big( 1 + \kappa E_{\rindep}\big)    \, , \nonumber
\end{eqnarray}
where in each case the error $E_{\rindep}$ is $O(1+ \xmac^\rho) e^{C_1 \vert \rindep \vert}$;
here, the constant $C_1$ is inherited from Proposition~\ref{p.slowkappa} via the relation $\stackrel{\kappa}{\simeq}$.
The quantities 
$$
\kappa^{-1} M_{\lfloor \kappa^{-1 }\rindep\rfloor} = 
\kappa^{-1} m_{\lfloor \kappa^{-1 }\rindep\rfloor -1 , \lfloor \kappa^{-1 }\rindep\rfloor +1}( \xmac )
$$ 
and  $$
\kappa^{-1} N_{\lfloor \kappa^{-1 }\rindep\rfloor} = 
\kappa^{-1} n_{\lfloor \kappa^{-1 }\rindep\rfloor -1 , \lfloor \kappa^{-1 }\rindep\rfloor +1}( \xmac )
$$ 
satisfy these respective estimates after the insertion of right-hand factors of two.

Further,
\begin{align}
\kappa^{-2} a_{\lfloor \kappa^{-1} \rindep\rfloor}(\xmac)
&= 2\rho \;\frac{ f_\rho(\xmac,\rindep)^{\,1+\rho}\; g_\rho(\xmac,\rindep)^{\,\rho} }{\big( f_\rho(\xmac,\rindep)^{\rho}+g_\rho(\xmac,\rindep)^{\rho}\big)^2}
\;\big(1+\kappa E_{\rindep}\big),
\label{e.aasymp} \\
\kappa^{-2} b_{\lfloor \kappa^{-1} \rindep\rfloor}(\xmac)
&= 2\rho \;\frac{ f_\rho(\xmac,\rindep)^{\,\rho}\; g_\rho(\xmac,\rindep)^{\,1+\rho} }{\big( f_\rho(\xmac,\rindep)^{\rho}+g_\rho(\xmac,\rindep)^{\rho}\big)^2}
\;\big(1+\kappa  E_{\rindep}\big),
\label{e.basymp}
\end{align}
where the errors satisfy the same bounds as above.
\end{proposition}
{\bf Proof of Theorem~\ref{t.lowkappasde}(1).} This is due to the estimates~(\ref{e.aasymp}) and~(\ref{e.basymp}). 

{\bf Proof of  Proposition~\ref{p.mnlowkappa}.}
Since $m_{-1,0}(x) = \kappa$, we have that
\begin{eqnarray}
  n_{k+1,k}(x) 
  & = &  \kappa \, x  \cdot  \prod_{i=0}^k \big( d_i (x)  - 1 \big) \, , \, \, \, \textrm{and} \label{e.nproduct} \\
 m_{k,k+1}(x)   & = & \,     \, \, \, \kappa \cdot \prod_{i=0}^k \big( c_i(x)  - 1 \big) \, . \nonumber
\end{eqnarray}
We now set $x$ equal to the element $\xmac \in D$ given in the proposition; we also adopt the shorthand $S(u) = S_\rho(\xmac,u)$.
It is straightforward that Proposition~\ref{p.slowkappa} implies that
\begin{equation}\label{e.sdouble}
 s_{\lfloor \kappa^{-1}u \rfloor}(\xmac) \, = \, S(u) \, \Big( 1  + O \big( e^{C_1 \vert u \vert} \big) \kappa \Big) \, ;
\end{equation}
the exponential growth has a uniform rate since $\xmac$ lies in the central domain~$D$, where it is bounded away from zero and infinity.
Apply the map $c$ to this relation, use\footnote{Recall the convention for result references stated after~(\ref{e.ciformula})} Lemma~\ref{l.cdasymptotics}(c) and that $\frac{\xaux^\rho}{(1+\xaux^\rho)^2}$ has a derivative that is uniformly bounded in absolute value to find that
$$
c \big( s_{\lfloor \kappa^{-1}u \rfloor}(\xmac) \big) = c \big( S(u) \big) + O \big( (1+ S(u)^{\rho}) e^{C_1 \vert u \vert} \big) \kappa \, .
$$
The function $J(u) = S(u)^\rho$ solves the initial-value problem in Lemma~\ref{l.de}. Since $\vert J' \vert \leq 8 \rho^2$ with $J(0) = \xmac^\rho$,
we have $J(u) \leq \xmac^\rho e^{8\rho^2 \vert u \vert}$ for $u \in \R$. Since $\rho \leq 1$, we obtain the naive upper bound on $S(u)^\rho$ of $\xmac^\rho e^{8 \vert u \vert}$; 
absorbing the additive eight by an increase in the value of $C_1$,
the coefficient of~$\kappa$ in the preceding display may thus be written $O(1+\xmac^\rho) e^{C_1 \vert u \vert}$ for $u \in \R$.
 
By Lemma~\ref{l.cdasymptotics}(c) again, we find that, for $u \in \R$,
\begin{equation}\label{e.ceu}
 c \big( s_{\lfloor \kappa^{-1}u \rfloor}(\xmac) \big) - 1 \, = \,  1 + \kappa \cdot 2 \bigg( 
  1 - \frac{2\big(1+(1-\rho)S(u)^\rho \big)}{(1+S(u)^\rho)^2}  \bigg) \, + \,   E_u \, \kappa^2 \, ,
\end{equation}
where  $E_u =   (1+ \xmac^{\rho}) e^{C_1 \vert u \vert}  O(1)$.  
We see then that, for $\rindep > 0$,
\begin{eqnarray*}
  \prod_{i=0}^{\lfloor \kappa^{-1}\rindep\rfloor} \big( c_i (\xmac)  - 1 \big) & = & \prod_{\udummy \in \kappa \Z \cap [0,\rindep]} \Bigg( 
  1 + \kappa \cdot 2 \bigg( 
  1 - \frac{2\big(1+(1-\rho)S(\udummy)^\rho \big)}{(1+S(\udummy)^\rho)^2}  \bigg)  + E_\udummy \, \kappa^2
  \Bigg) \\
   & = &  \exp \Bigg\{ 2\kappa \sum_{\udummy \in \kappa \Z \cap [0,\rindep]} \bigg( 
  1 - \frac{2\big(1+(1-\rho)S(\udummy)^\rho \big)}{(1+S(\udummy)^\rho)^2}  \bigg) \, + \, \kappa \, E_{\rindep}  \Bigg\} \\
   & = & \exp \Bigg\{ 2 \int_0^{\rindep}  \bigg( 
  1 - \frac{2\big(1+(1-\rho)S(\udummy)^\rho \big)}{(1+S(\udummy)^\rho)^2}  \bigg) \dd \udummy  \Bigg\} \Big( 1 + \kappa E_{\rindep} \Big) \, ,
\end{eqnarray*}
the second equality due to the standard estimate $\prod (1+t_i) = \exp \big\{ \sum t_i + O(\sum t_i^2)\big\}$.
The error terms $E$ may differ from line to line, subject to the condition given when they were introduced above. Since the exponential expression in the final line equals $f_\rho(\xmac,\rindep)$,
we obtain the sought bound on $\kappa^{-1} m_{\lfloor \kappa^{-1 }\rindep\rfloor , \lfloor \kappa^{-1 }\rindep\rfloor +1}(\xmac )$ for $\rindep > 0$. And also when $\rindep < 0$, provided that the product and sum expressions in the  preceding display are interpreted compatibly with the convention for negatively indexed products in~(\ref{d.zdef}). 

Instead applying Lemma~\ref{l.cdasymptotics}(d), we have in  counterpart for $d$,
$$
 d \big( s_{\lfloor \kappa^{-1}u \rfloor}(\xmac) \big) - 1 \, = \,  1 - \kappa \cdot 2 \bigg( 
  1 - \frac{2\big((1-\rho)S(u)^\rho + S(u)^{2\rho}  \big)}{(1+S(u)^\rho)^2}  \bigg) \, + \,   E_u \, \kappa^2 \, , 
$$
where the bound satisfied by $E_u$ is unchanged.
Hence,
\begin{equation}\label{e.drel}
  \prod_{i=0}^{\lfloor \kappa^{-1 }\rindep\rfloor} \big( d_i (\xmac)  - 1 \big) \, =  \,   \exp \Bigg\{ -2 \int_0^{\rindep}   \bigg( 1 -
  \frac{2\big((1-\rho)S(\udummy)^\rho + S(\udummy)^{2\rho}  \big)}{(1+S(\udummy)^\rho)^2}  \bigg) \, {\rm d} \udummy  \Bigg\}  \Big( 1 + \kappa E_{\rindep} \Big)   \, .
\end{equation}
Noting the factor of $\xmac$ on the right-hand side of~(\ref{e.nproduct}), we multiply~(\ref{e.drel}) by~$\xmac$ and find that the resulting right-hand term 
$\xmac \cdot \exp \{ -2I \}$ equals $g_\rho(\xmac,\rindep)$.  The bound on $\kappa^{-1} n_{\lfloor \kappa^{-1 }\rindep\rfloor + 1, \lfloor \kappa^{-1 }\rindep\rfloor}( \xmac )$ follows.

To obtain the assertion made in regard to 
$\kappa^{-1} M_{\lfloor \kappa^{-1 }\rindep\rfloor}$,
sum~(\ref{e.massertion}) for values $\rindep - \kappa$ and $\rindep$, and use the differentiability of $f_\rho(\xmac,\rindep)$ at $\rindep$ to absorb via a factor of $1 + O(\kappa)$
 the error arising from the microscopic unit index displacement. (The derivative in~$\rindep$ is readily seen to be bounded on compact subsets of~$\R$;  in fact, decay at infinity means that this is true on all of~$\R$. So the implied constant in the $O(\kappa)$ term  may be chosen independently of $\rindep \in \R$.)
 
 Likewise for 
$\kappa^{-1} N_{\lfloor \kappa^{-1 }\rindep\rfloor}$. Applying these estimates to the formulas for $a_i$ and $b_i$ in terms of $M_i$ and~$N_i$ in Proposition~\ref{p.abmnfacts}(1), we obtain the stated asymptotics for $a_i$ and $b_i$ and thus 
 complete the proof of Proposition~\ref{p.mnlowkappa}. \qed

\section{The low $\kappa$ limit of $\lambda_{\rm max}$}

Here we prove Theorem~\ref{t.lowlambdamax} concerning the approach of $\lambdamaxkp$ to one in the limit of low $\kappa$.

{\bf Proof of  Theorem~\ref{t.lowlambdamax}.} In light of reduction to standard solutions by basic symmetries, Definition~\ref{d.mmm} and Proposition~\ref{p.relativereward}(4), it suffices to show that there exist positive $C$ and $\eta$ such that, for $\kappa$ small enough, 
\begin{equation}\label{e.ratiobound}
 \left\vert \, \frac{n_{\infty,-\infty}}{m_{-\infty,\infty}} -   1 \, \right\vert \, \leq \,  C  \kappa^\eta 
\end{equation}
for any element of \abmnkpmacspace  with $\phi_0$ in the central domain~$D$. 

Write $x = \phi_0$.
The plan is to argue that 
$m_{-\infty,\infty}$ equals  $m_{-1,0} \int_\R f_\rho(x,\udummy)\, {\rm d} \udummy$
up to an error that vanishes as~$\kappa \searrow 0$,
and that 
$n_{\infty,-\infty}$ similarly approximates  $m_{-1,0} \int_\R g_\rho(x,\udummy)\, {\rm d} \udummy$. The integrals are equal by Proposition~\ref{p.sfacts}(5), as desired.
To implement this approach, we will use the approximations of $\kappa$-scaled $m$- and $n$-differences by  $f_\rho(x,u)$ and  $g_\rho(x,u)$ found in Proposition~\ref{p.mnlowkappa}.
 These approximations worsen for indices that are high multiples of $\kappa^{-1}$ because the mimicry of the $s$-orbit by the $S_\rho$-flow (as gauged by Proposition~\ref{p.slowkappa}) may have deteriorated.
 So we will attempt the comparison only on a short scale, delimited by a continuous-time parameter~$z$ of the form $\eta \log \kappa^{-1}$. We handle the longer scale via the next result.

\begin{lemma}\label{l.mntail}
There exist positive $\eta$, $c_0$ and~$C$ such that, for $(a,b,m,n) \in$ \abmnkpmacspace with $\phi_0 \in D$,
$$
 m_{-\infty, - \lfloor z \kappa^{-1} \rfloor} +   m_{\lfloor z \kappa^{-1} \rfloor, \infty} \leq C \kappa^{c_0} m_{-\infty,\infty} 
$$
and 
$$
 n_{- \lfloor z \kappa^{-1} \rfloor,-\infty} +   n_{\infty,\lfloor z \kappa^{-1} \rfloor} \leq C \kappa^{c_0} n_{\infty,-\infty} \, ,
$$
where $z = \eta \log \kappa^{-1}$.
\end{lemma}
In the argument for  Theorem~\ref{t.lowlambdamax}, the flow $S_\rho$ acts as an analytic proxy for the discrete dynamics.
So we set the index $\xmac$ that specifies the flow $S_\rho(\xmac,\bullet):\R \to (0,\infty)$ equal to the given initial ratio~$x=\phi_0$, and write
$S(u)=S_\rho(x,u)$ for the resulting flow. We also denote the scaled $s$-iterates that the flow tracks by
$\trackS_\kappa(u):=s_{\lfloor \kappa^{-1}u \rfloor}(x)$.

{\bf Proof of Lemma~\ref{l.mntail}.} 
 Since $z = \eta \log \kappa^{-1}$,
 the bound  $\trackS_\kappa(u) \leq 2 S(z)$ holds for $u = z$ 
 by~(\ref{e.sdouble}), provided that we make a suitably small choice of $\eta > 0$. And since $s$ is sub-diagonal by Lemma~\ref{l.incphi}(3), this bound also holds for all~$u \geq z$.

By Lemma~\ref{l.cdasymptotics}(c),
 $$
 c\big(\trackS_\kappa(u) \big) - 1 = 1 + 2 \, \Phi_f\big(\trackS_\kappa(u)\big) \kappa  + E_u \kappa^2
 $$ where $\Phi_f(y) =   
  1 - 2\big(1+(1-\rho)y^\rho \big)(1+y^\rho)^{-2}$ and $E_u = O(1)$. 
  
  Note that $S$ solves the differential equation in Definition~\ref{d.fg}, so $S(y) \to 0$ as $y \to \infty$.
  Note also that $\lim_{y \searrow 0} \Phi_f(y) = -1$.
   For a suitably small choice of $\kappa_0 = \kappa_0(\eta)$, the condition $\kappa \in (0,\kappa_0)$ thus ensures that $S(z) < \e$ where $\e > 0$ is such that $\Phi_f(y) \leq -3/4$ for $y \in (0,2\e)$. Since $\trackS_\kappa(u) \leq 2 S(z)$ for $u \geq z$, the linear coefficient in the last display is at most $-3/2$. Since $E_u \kappa^2 \leq \kappa/2$ by choosing $\kappa_0$ suitably, we see that  $c(\trackS_\kappa(u)) -1 \leq 1 - \kappa$ for $u \geq z$.

By taking a ratio of equalities of the form~(\ref{e.mdef}), we obtain 
$$
 m_{\lfloor z \kappa^{-1} \rfloor +i,\lfloor z \kappa^{-1} \rfloor +i+1} m_{\lfloor z \kappa^{-1} \rfloor,\lfloor z \kappa^{-1} \rfloor +1}^{-1} = 
  \prod_{j = 0}^{i-1}
\Big( c\big( s_j\big(\trackS_\kappa(z)\big)  \big) - 1 \Big)
 $$
 since this ratio of $m$-differences coincides with this ratio for the default solution with the same value of~$\phi_0$.
 Certainly the right-hand side is a product of positive quantities: see the discussion leading up to~(\ref{e.ciformula}) in Section~\ref{s.explicitabmn}.
Noting also that  $s_j\big(\trackS_\kappa(z)\big)   = \trackS_\kappa(z + \kappa j)$, 
we find from 
$m_{\lfloor z \kappa^{-1} \rfloor,\infty} = \sum_{j= \lfloor z \kappa^{-1} \rfloor}^\infty m_{j,j+1}$ that 
\begin{equation}\label{e.mupperbound}
m_{\lfloor z \kappa^{-1} \rfloor,\infty} \, = \, m_{\lfloor z \kappa^{-1} \rfloor,\lfloor z \kappa^{-1} \rfloor +1}  \sum_{i=0}^\infty \prod_{j = 0}^{i-1}
\big( c( \trackS_\kappa(z + \kappa j) ) - 1 \big)
\, \leq \, m_{\lfloor z \kappa^{-1} \rfloor,\lfloor z \kappa^{-1} \rfloor +1} \cdot \kappa^{-1} \, ,
\end{equation}
where in the final bound, we applied the just obtained upper bound of $1-\kappa$ on   $c(\trackS_\kappa(u)) -1$.
 
 But by Proposition~\ref{p.mnlowkappa}($m$) and $\mdefault_0 - \mdefault_{-1} = \kappa$, we have that
 $$
  m_{\lfloor z \kappa^{-1} \rfloor, \lfloor z \kappa^{-1} \rfloor + 1}(x) \, = \,  m_{-1,0} \,  f_\rho(x,z) \big( 1 + \kappa E_z \big)
 $$
 where the default solution is understood, $x = \phi_0$, and the error $E_z = O(1+x^\rho) e^{ C_1 \vert z \vert}$ is simply $O(1)e^{C_1 z}$ since $x \in D$ (and $z > 0$).

 We now decrease if need be the value of $\eta$  in  $z = \eta \log \kappa^{-1}$ so that $\vert E_z \vert \leq \kappa^{-1/2}$ (with $\eta < (2C_1)^{-1}$). 
 Using Proposition~\ref{p.sfacts}(4) to bound  $f_\rho(x,z)$ above, we thus find that  
 $$
 m_{\lfloor z \kappa^{-1} \rfloor, \lfloor z \kappa^{-1} \rfloor + 1} \leq m_{-1,0} \cdot e^{-2z(1-\e)}\big( 1+ O(1) \kappa^{1/2} \big) = m_{-1,0} \cdot \kappa^{2\eta(1-\e)}\big( 1+ O(1) \kappa^{1/2} \big) \, ,
 $$
  where $\vert O(1) \vert \leq 1$. From~(\ref{e.mupperbound}), we obtain 
 \begin{equation}\label{e.mconclusion}
 m_{\lfloor z \kappa^{-1} \rfloor, \infty} \leq m_{-1,0} \, O(1) \kappa^{c_0 - 1} 
 \end{equation}
 with $c_0 = 2\eta(1-\e)$.
 Note that
 $$
 m_{-\infty,-\lfloor z \kappa^{-1} \rfloor}  \, = \, \sum_{i=\lfloor z \kappa^{-1} \rfloor}^\infty m_{-i-1,-i} \, .
 $$
  Since $s$ is sub-diagonal, we have $s_{-1}(y) > y$, and so  $\phi_{-i}$ is bounded below uniformly in~$(i,x) \in \N \times D$ where $x = \phi_0$
 (and $D$ is the central domain). Hence,  $m_{-i-1,-i} = n_{-i-1,-i} \phi_{-i}^{-1} \leq O(1) n_{-i-1,-i}$.

  By Corollary~\ref{c.rolereversal}, $n_{-i,-i-1}$
 is equal to $m_{i,i+1}$ for the role-reversed \abmnmacspace solution $(b_{-i},a_{-i},n_{-i},m_{-i})$. Thus, from~(\ref{e.mconclusion}), we infer that  $m_{-\infty,-\lfloor z \kappa^{-1} \rfloor} \leq m_{-1,0} \, O(1) \kappa^{c_0 - 1}$, whence also 
 $$
m_{-\infty,\lfloor z \kappa^{-1} \rfloor}  +  m_{\lfloor z \kappa^{-1} \rfloor, \infty} \, \leq \, m_{-1,0} \, O(1) \kappa^{c_0 - 1} \, .
 $$
 Since $m$-increments are non-negative, and $f_\rho(x,z)$ is bounded away from zero for $(x,z)$ in the precompact  $D \times [-1,1]$, we find by summing Proposition~\ref{p.mnlowkappa}($m$)  that  $m_{-\infty,\infty} \geq m_{-\lfloor \kappa^{-1} \rfloor, \lfloor \kappa^{-1} \rfloor} \geq c_1 \kappa^{-1} m_{-1,0}$ for some small positive~$c_1$. Hence,
 $$
m_{-\infty,\lfloor z \kappa^{-1} \rfloor}  +  m_{\lfloor z \kappa^{-1} \rfloor, \infty} \leq m_{-\infty,\infty} O(1) \kappa^{c_0 } \, ,
 $$ 
as we sought to show in proving Lemma~\ref{l.mntail}(m).

Lemma~\ref{l.mntail}(n) may be obtained by role-reversal symmetry. Indeed, applying reflection about minus one-half yields
$n_{-i,-\infty}(x) = m_{i-1,\infty}(x^{-1})$ and $n_{\infty,i}(x) = m_{-\infty,-i-1}(x^{-1})$. Since the central domain is invariant under $x \mapsto x^{-1}$, 
we may take $i = \lfloor z \kappa^{-1} \rfloor$ and obtain Lemma~\ref{l.mntail}(n) from Lemma~\ref{l.mntail}(m); technically, there is a mismatch of one unit in the indexing, because reflection has been about $-1/2$ rather than zero, but the discrepancy is absorbed by increasing the value of $C > 0$. \qed

Lemma~\ref{l.mntail}(m), and Proposition~\ref{p.mnlowkappa} summed, imply that
 \begin{equation}\label{e.mintegral}
   m_{-\infty,\infty}  \, =  \,
m_{-\lfloor z \kappa^{-1} \rfloor,\lfloor z \kappa^{-1} \rfloor}  \Big( 1 - O(1) \kappa^{c_0}\Big)
  \, = \,  m_{-1,0} \int_{-z}^z f_\rho(x,\udummy)\, {\rm d} \udummy  \cdot  \big( 1 + \kappa E_z \big) (1 - O(1) \kappa^{c_0})
   \end{equation}
where $E_z =O(1) e^{C_1 \vert z \vert}$ (since $x^\rho = O(1)$, from $x \in D$).
Given the selection of $\eta>0$ in the preceding proof, the choice $z = \eta \log \kappa^{-1}$ leads to $\kappa E_z = O(\kappa^{1/2})$.

Our plan calls for integration over~$\R$ in place of $[-z,z]$, so we wish to estimate the discrepancy between these integrals.
\begin{lemma}\label{l.integral}
For $\xmac \in (0,\infty)$, let $\vmac \in \R$ be the value associated to $\xmac$ by Proposition~\ref{p.sfacts}(2). Then
$$
\int_{\R \setminus [-z,z]}  f_\rho(\xmac,\rdummy)\, {\rm d} \rdummy  \, = \, \frac{1}{f_\rho(1,-\vmac)} \int_{\R \setminus [-z-\vmac,z-\vmac]}  f_\rho(1,\rdummy)\, {\rm d} \rdummy \, .
$$
\end{lemma}
{\bf Proof.} By a change of variable and Proposition~\ref{p.sfacts}(2,f),
$$
\int_{[-z,z]^c}  f_\rho(\xmac,\rdummy)\, {\rm d} \rdummy \, = \, 
\int_{[-z-\vmac,z-\vmac]^c}  f_\rho(\xmac,\vmac+\rdummy)\, {\rm d} \rdummy  \, = \,  f_\rho(\xmac,\vmac) \int_{[-z-\vmac,z-\vmac]^c}  f_\rho(1,\rdummy)\, {\rm d} \rdummy  \, .
$$
Take $\rindep =-\vmac$ in Proposition~\ref{p.sfacts}(2,f) and use $f_\rho(\xmac,0)=1$ (which is immediate from Definition~\ref{d.fg})
to find that $f_\rho(\xmac,\vmac) = 1/f_\rho(1,-\vmac)$. \qed

As $\xord$ varies over the precompact~$D$, the quantity $\vmac = \vmac(\xord) \in \R$ remains bounded. So~$\frac{1}{f_\rho(1,-\vmac)} = O(1)$.

We may thus apply Lemma~\ref{l.integral} and Proposition~\ref{p.sfacts}(4) to find that, for any $\e > 0$, and $z > 0$ large enough, 
$$
\int_{\R \setminus [-z,z]}  f_\rho(\xord,\rdummy)\, {\rm d} \rdummy  \leq C \exp \{ -2z(1-\e) \} \, , 
$$
with the constant~$C$ absorbing the influence of the bounded offset~$\vmac$.

The integral  $\int_{\R}  f_\rho(\xord,\udummy)\, {\rm d} \udummy$ is positive and finite,  so  
$$
\int_{-z}^z   f_\rho(\xord,\udummy)\, {\rm d} \udummy \, = \, 
\int_{\R}   f_\rho(\xord,\udummy)\, {\rm d} \udummy \,  \big( 1 - O(1) e^{-z}\big) \, ,
$$
where we took $\e \in (0,1/2)$.
Since $z = \eta \log \kappa^{-1}$, we have $e^{-z} = \kappa^\eta$, so that~(\ref{e.mintegral}) yields
$$
 m_{-\infty,\infty} \,  = \,  m_{-1,0} \int_{\R}   f_\rho(\xord,\udummy)\, {\rm d} \udummy \cdot  \big( 1 + O(1)\kappa^{1/2} \big)  \big( 1 - O(1) \kappa^{c_0} \big)  \big( 1 - O(1) \kappa^\eta \big)
$$
or simply 
$m_{-\infty,\infty} \, = \,  m_{-1,0} \int_{\R}   f_\rho(\xord,\udummy)\, {\rm d} \udummy \, \cdot  \big( 1 - O(1) \kappa^\eta \big)$
by decreasing the value of $\eta$ if need be.

A counterpart argument harnessing 
Lemma~\ref{l.mntail}(n) yields 
$$
 n_{\infty,-\infty} \, = \, m_{-1,0} \int_{\R}   g_\rho(\xord,\udummy)\, {\rm d} \udummy \, \cdot  \big( 1 - O(1) \kappa^\eta \big) \, .
$$

Hence,
$$
 \frac{n_{\infty,-\infty}}{m_{-\infty,\infty}} =  \frac{\int_{\R}   g_\rho(\xord,\udummy)\, {\rm d} \udummy  \,  \big( 1 - O(1) \kappa^\eta \big)}{\int_{\R}   f_\rho(\xord,\udummy)\, {\rm d} \udummy \, \big( 1 + O(1) \kappa^\eta \big)} \, .
$$
As planned, we may note that the two integrals are equal, by Proposition~\ref{p.sfacts}(5). Hence, 
$$
 \frac{n_{\infty,-\infty}}{m_{-\infty,\infty}} =   1 + O(1) \kappa^\eta 
$$
and we obtain~(\ref{e.ratiobound}) as desired. 
This completes the proof of Theorem~\ref{t.lowlambdamax}. \qed

{\em Remark.} If we take $\rho \geq 1$, a more general error estimate (roughly $E_z = \exp \big\{ (1+\rho)^4 O(\vert z \vert) \big\}$) in Proposition~\ref{p.mnlowkappa} will lead to $\eta =\eta(\rho) \searrow 0$ as $\rho \to \infty$ in Theorem~\ref{t.lowlambdamax}. The hypothesis $(\kappa,\rho) \in W$ is also needed, to enable  $S_\rho$-tracking of the $s$-orbit, as in the remark that follows Proposition~\ref{p.slowkappa}.

\section{Scaled gameplay in the low-$\kappa$ limit} 

Here we prove Theorem~\ref{t.lowkappasde}(2). 

\begin{proposition}\label{p.driftscaled}
Consider $\tlp(\kappa,\rho)$ played at a time-invariant Nash equilibrium of battlefield index zero. Let $p(i)$ denote the probability of a rightward move at location~$i$. Then
$$
\kappa^{-1} \Big( 2   p \big(  \lfloor \kappa^{-1} u \rfloor \big) -  1 \Big) \, \lora \, \frac{1- S_\rho(1,u)^\rho}{1+  S_\rho(1,u)^\rho}
$$ uniformly for $u$ lying in compact subsets of~$\R$.
\end{proposition}
{\bf Proof.} 
By Theorem~\ref{t.nashabmn}, gameplay is governed by the stake-profile components of an element $(a,b,m,n) \in$ \abmnkpmac. 
The probability~$p(i)$ is a sum of contributions according to whether the turn is flip or stake:
$$
 p(i) \, = \, \frac{1-\kappa}{2} +  \frac{\kappa \, a_i^\rho}{a_i^\rho + b_i^\rho}
$$
so that
$$
 \kappa^{-1} \big( 2   p (  i ) -  1 \big) \, = \, \frac{a_i^\rho - b_i^\rho}{a_i^\rho + b_i^\rho} \, .
$$
By Proposition~\ref{p.abmnfacts}(1),  
$$
a_i = \frac{\kappa \rho M_i^{1+\rho}N_i^\rho}{(M_i^\rho+N_i^\rho)^2}  \qquad \textrm{and} \qquad  b_i =  \frac{\kappa \rho M_i^\rho N_i^{1+\rho}}{(M_i^\rho+N_i^\rho)^2} 
$$
so that $b_i/a_i$ equals $N_i/M_i$. 
Thus, $\frac{a_i^\rho - b_i^\rho}{a_i^\rho + b_i^\rho}  \, = \,  \frac{1 - \beta_i^\rho}{1 + \beta_i^\rho}$ with $\beta_i = N_i/M_i$.
With an error $E_u$ satisfying the bound in Proposition~\ref{p.mnlowkappa}, this result implies that
$$
 \beta_{\lfloor \kappa^{-1} u \rfloor} = \frac{g_\rho(x_\kappa,u)}{f_\rho(x_\kappa,u)} \big( 1 + \kappa E_u \big) \, .
$$
The proposition is applied with the parameter $\xmac$ set equal to the $\phi_0$-value (which we denote here by~$x_\kappa$) associated to the given element  $(a,b,m,n) \in$ \abmnkpmac. 
The result is applicable since~$x_\kappa$ lies in the $(\kappa,\rho)$-central domain~$D$ in view of the concerned  Nash equilibrium having battlefield zero. 
And $x_\kappa \in D$ implies that $x_\kappa - 1 = O(\kappa)$ given the form of~$D$ in Definition~\ref{d.battlefield}.
Since $g_\rho(x,u)$ and $f_\rho(x,u)$ are smooth positive functions, $g_\rho(x_\kappa,u)/f_\rho(x_\kappa,u) = g_\rho(1,u)/f_\rho(1,u) \big( 1 + O(\kappa)\big)$.
 
But  $\frac{g_\rho(1,u)}{f_\rho(1,u)} = S_\rho(1,u)$ by Theorem~\ref{t.fg}, so that 
$$
 \kappa^{-1} \Big( 2   p \big(  \lfloor \kappa^{-1} u \rfloor \big) -  1 \Big) \, = \, \frac{1- S_\rho(1,u)^\rho}{1+  S_\rho(1,u)^\rho} \big( 1 + \kappa E_u \big) \, , 
  $$
  where the error $\vert E_u \vert$ is bounded on compact subsets. This completes the proof of Proposition~\ref{p.driftscaled}. \qed 

{\bf Proof of Theorem~\ref{t.lowkappasde}(2).} Ethier and Kurtz's~\cite[Corollary~$7.4.2$]{EthierKurtz}
 provides a framework for proving the convergence of discrete Markov chains to diffusion processes. 
For the framework to apply to a  sequence of Markov chains \( \{ Y^n \} \) with transition kernels \( p_n(x, \bullet) \), it is sufficient that the following conditions are met.
\begin{itemize}
\item The scaled drift coefficients
   $b_n(x) := n^2 \int (y-x)\, p_n(x,dy)$ 
    converge uniformly on compact sets to a continuous function $b(x)$.
\item
    The scaled diffusion coefficients
    $a_n(x) := n^2 \int (y-x)^2\, p_n(x,dy)$ converge in the same sense to one (the variance of the limiting diffusion).
    \item
    The jumps of $Y^n$ are uniformly bounded by order~$n^{-1}$. 
    \item The martingale problem for the limiting generator
    \[
        L f(x) =  \tfrac12 f''(x) + b(x) f'(x)  \, , \quad f \in C_c^\infty(\mathbb{R}),
    \]
    is well-posed.
\end{itemize}
The chains  $Y^n$ may be specified on $[0,\infty)$ rather than~$\N$, by linear interpolation. When the above conditions are met, these chains are 
continuous real-valued  processes on $[0,\infty)$ 
whose speeded versions 
$$
[0,\infty) \to \R: u \to Y^n(n^2 u)
$$ 
converge in distribution to the unique solution of the SDE
\[
  {\rm d}Z_t =  b(Z_t)\, {\rm d}t + {\rm d}W_t  \, ,
\]
where $W_t$ is standard Brownian motion. (Convergence occurs in the compact-uniform topology on the space~$\mc{C}$ of continuous functions mapping $[0,\infty)$ to~$\R$, because our interpolated prelimiting processes are continuous, and $\mc{C}$ is a closed subspace of the space of c\`adl\`ag paths with the Skorokhod topology---the $J_1$-topology in Billingsley's~\cite{Billingsley} terminology---employed by Ethier and Kurtz.)

Recall that, for $y \in \R$, $X_{\kappa,\rho}(y,\bullet):\N \to \Z$ denotes gameplay under the time-invariant Nash equilibrium of battlefield index zero, with initial condition $X_{\kappa,\rho}(0)=\lfloor y \rfloor$.
We apply the framework with $Y^n(k) = n^{-1}X_{n^{-1},\rho}(nz,k)$
for $k \in \N$ (and with the domain being extended to $[0,\infty)$ via interpolation). 
In this way,
 $n \in \N$ corresponds to~$\kappa$ in Theorem~\ref{t.lowkappasde} via $n = \kappa^{-1}$.
(It would seem that $\kappa$ must tend to zero through integer reciprocals. But in fact we may equally apply the framework with $n \to \infty$ in an arbitrary fashion.)
To check that the framework is applicable, note that the scaled drift hypothesis is granted by Proposition~\ref{p.driftscaled} with $b(u) = \frac{1- S_\rho(1,u)^\rho}{1+  S_\rho(1,u)^\rho}$.
Since the counter jump is $\pm 1$ at each step in \tlpkp, every jump of  $Y^n$ has magnitude~$n^{-1}$ and $a_n(x) = 1$ identically. By~\cite[Corollary~$6.3.3$]{StroockVaradhan}, 
the martingale problem for  ${\rm d}Z_t =  a(Z_t){\rm d}W_t  + b(Z_t)\, {\rm d}t$   is well-posed 
when $a$ and~$b$ are bounded with bounded continuous derivatives (in our case, $a =1$, and $b$ is smooth with $\vert b \vert \leq 1$).

The outcome is the convergence asserted by Theorem~\ref{t.lowkappasde}(2), where in the SDE-drift $R_\rho(Z_t)$, the function $R_\rho:\R \to (-1,1)$ is given by $R_\rho(u) = \frac{1 - S_\rho(1,u)^\rho}{1 + S_\rho(1,u)^\rho}$. The alternative formula claimed for $R_\rho(u)$ arises from the equality $S_\rho(1,u)^\rho= S_1(1,\rho^2 u)$, which is precisely the identity noted after Lemma~\ref{l.de} with~$\xmac=1$. 
As is also noted there, $S_1(1,\rho^2 u) \sim (8 \rho^2 \vert u \vert )^{{\bf 1}_{u < 0}-{\bf 1}_{u > 0}}$ as $\vert u \vert \to \infty$, which yields the asymptotics claimed for $R_\rho(u)$ when applied to the alternative formula. \qed

\chapter*{Epilogue: Directions}\label{c.directions}

Our treatment has more or less directly posed certain open problems. 
We begin by surveying  several such directions and then turn to some broader themes for further inquiry.

{\em Invariance principles.}
For $\rho \in (0,1]$, we have exhibited Brownian Boost as a fine-mesh high-noise  scaling limit of the Trail of Lost Pennies. How universal is this scaled structure? 
A natural direction is to seek to prove an invariance principle, a promising point of departure being given by variants of \tlpkpspace whose stake and flip moves have symmetric jump distributions no longer supported on plus and minus one.

In another direction, our limit $\kappa \searrow 0$ represents one of
high {\em flip} noise, with the contest rule fixed. Taking $\rho>0$
small in the contest success function $a^\rho/(a^\rho+b^\rho)$ leads to
outcomes close to those of a fair coin, just as making flip moves
dominant does. Thus the limit $\rho \searrow 0$ for fixed
$\kappa\in(0,1]$ yields a second high-noise regime whose scaling
behaviour may be compared with the limit studied here. We have found a
`Brownian Boost line' along the positive $\rho$-axis; the proposed
direction would seek a counterpart along the positive $\kappa$-axis.

A separate issue concerns the formulation of the model itself in continuous time.
The chosen high-flip-noise representation of Brownian Boost has permitted us to exhibit its equilibrium structure
explicitly, but it does so by circumventing the challenges inherent in instantaneous feedback
loops. This may prompt attention to the problem of formulating and solving
$\rho$-Brownian Boost directly via suitable classes of non-anticipatory
strategies (see Section~\ref{s.bbhighnoise}).

{\em Discrete-game variants.}
Other problems concern how the solution theory we have developed breaks down in certain corners of  parameter space. These include proving
non-existence results for solutions of $\abmnmacspace$ when $(\kappa,\rho)$ lies sufficiently
above the region $W$ specified in~(\ref{e.weakregion})---one may begin with $\kappa=1$ and
$\rho>1$---and determining whether non-time-invariant equilibria exist in \tlpkpspace
and \bbrho. Several related directions for ${\rm TLP}(1,1)$ are discussed in
\cite[Section~7]{LostPennies}. 

The `Hand It Over!' one-step game discussed at the end of Chapter~\ref{c.symmetrytools} may act
as the single-step rule for a \tlpkp-type game played on integer intervals or on~$\Z$.
More generally, the Tullock contest and its `hand-over' variant may be interpolated by
introducing a pair of parameters in $[0,1]$ that specify the proportions of the winning and losing
players' stakes at a given turn that are transferred to the opponent, with the remainder
surrendered to the bank.
This yields a family of stake-governed tug-of-war games that interpolate between
 pure rent dissipation, at $(0,0)$,  and pure rent transfer at~$(1,1)$.
It would be of interest to understand how equilibrium structure, incentive asymmetry,
and discouragement effects vary across this family, and whether the high-noise scaling
limits identified for~\tlpkpspace  persist under such perturbations or undergo qualitative change.

{\em Broader themes.}
The explicit ODE equilibria uncovered for Brownian Boost may also serve as a point of entry
for analysts into questions of economic interest. Numerical investigations indicate
unexpected intricacies in the quantified discouragement effect, 
 while natural variants of
Brownian Boost introduce higher-dimensional state variables and lead to multi-component
coupled equations, posing challenges that may be of interest to PDE specialists.

We develop these broader themes through three directions for further study.

{\bf Finite-interval games.}
The Trail of Lost Pennies may be played on a finite interval $\llbracket - j,k \rrbracket$ for $j,k \in \N$. The game ends when $X$ reaches $- j$ or $k$ with terminal payments given by a quadruple $\big(m_{- j},m_k,n_{- j},n_k\big)$. Under an analogue of the Nash-\abmnmacspace Theorem~\ref{t.nashabmn}, time-invariant Nash equilibrium stake profiles would correspond to \abmnkpmacspace elements that extend the boundary data to  $\llbracket - j ,k \rrbracket$. The finite-trail Mina margin map $\minammkp^{-j,k}: \phi_0 \mapsto n_{k,- j}/m_{- j,k}$ satisfies the formula in Proposition~\ref{p.relativereward}(2) with summations over $\llbracket - j,k-1 \rrbracket$ instead of~$\Z$. The level sets of this map index equilibria of given Mina margin (or relative incentive) $n_{k,- j}/m_{- j,k}$.  The finite-interval games were investigated for {\rm TLP}($1,1$) in~\cite{LostPennies}. When the Mina margin is close to one, it appears that there is a unique equilibrium when $k- j  \leq 5$; for $k- j  =6$, there are three, and the number may be expected to grow as $2(k-j) + \Theta(1)$  for longer gameboards: see~\cite[Section~$2.5$]{LostPennies}.

We have not investigated the finite-interval games in this article, but the finite-trail Mina margin map offers a useful perspective on its results, with the low-$\rho$ convergence $\lambdamax(\kappa,\rho) \to 1$ corresponding to $\minammkp^{-j,k} \to 1$ uniformly on compacts. The characteristic zigzag pattern seen in Figure~\ref{f.mmm} takes longer to appear as gameboard length rises when $\kappa$ is smaller: while $\mc{M}^{-9,9}_{1,1} = 1$ has $27$ roots according to~\cite[Equation~$(16)$]{LostPennies}, there are $21$ roots for~$\mc{M}_{0.9,1}^{-9,9} = 1$ as depicted in Figure~\ref{f.mmm}. Likewise,  
the outset gameboard length for  non-unique equilibria at given Mina margin may be expected rise as $\kappa$ drops: longer gameboards are needed at high-flip-noise levels for the effects of stake turns to be felt. 

Nor have we explored $\rho$-Brownian Boost on finite intervals. Given the remark about Penny Forfeit in Section~\ref{s.pennyforfeit}, it seems likely that with suitable boundary conditions the characterization of equilibria in terms of the \bbrhospace ODE pair remains valid when $\rho \in (1,2)$ when the game is played on finite intervals whose length 
satisfies a suitable $\rho$-determined upper bound.

{\bf The map $(\kappa,\rho) \mapsto \lambdamax(\kappa,\rho)$.} In~(\ref{e.lambdamaxkappazero}), we extended the domain of $\lambdamax$ by setting its values on the $\kappa =0$ axis equal to one. 
This accords with the absence of asymmetric equilibria in \bbrhospace due to Proposition~\ref{p.sfacts}(5). 
Regarding the second high-noise limit mooted above, where~$\rho \searrow 0$ with  $\kappa \in (0,1]$,
  it is reasonable to surmise analogously that $\kappa \mapsto \lambdamax(\kappa,0) -1$ on~$(0,1]$ also vanishes identically.

\begin{figure}[htbp]
\centering
\includegraphics[width=0.8\textwidth]{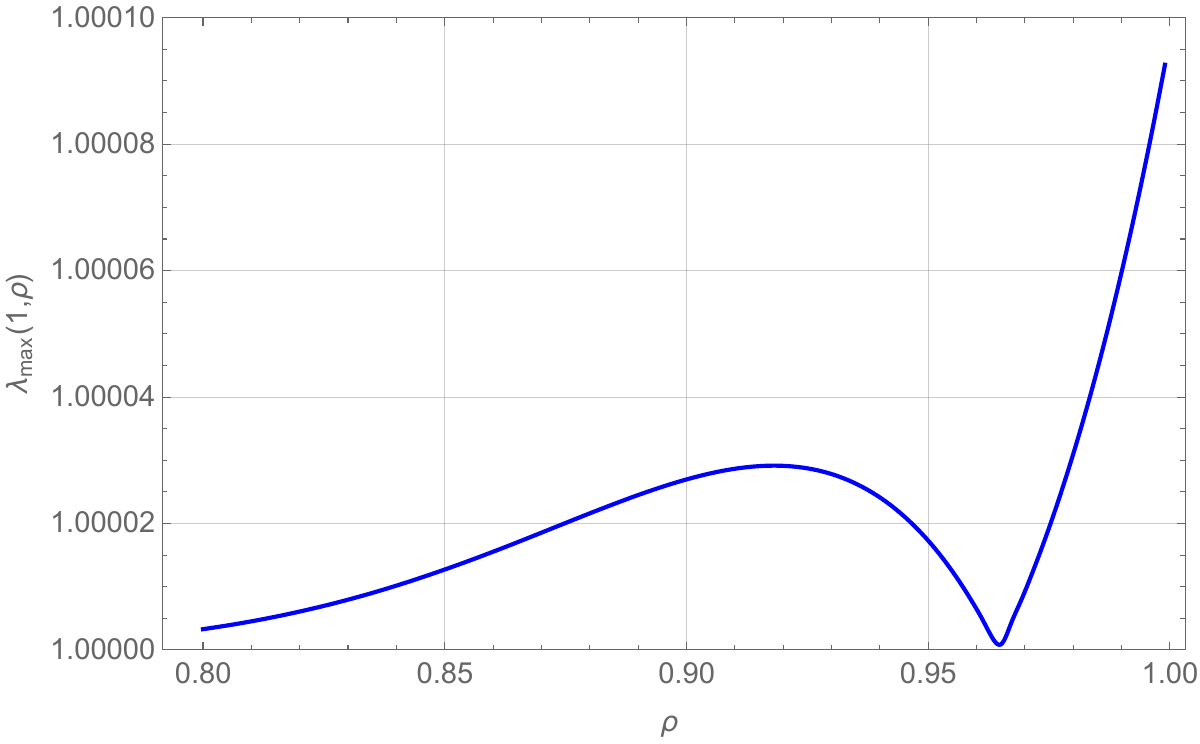}
\caption{A numerical approximation of the curve $\rho \to \lambdamax(1,\rho)$ for $\rho$-values in $(0.8,1)$. The curve shown has been interpolated from a sequence of points $\big(\rho,\lambdamax(1,\rho) \big)$, where each $\lambdamax(1,\rho)$ is approximated by maximizing over a fine mesh the values in the central domain $D = D_{\kappa,\rho}$ of the finite-trail Mina margin map~$\mc{M}_{1,\rho}^{j,k}$ for suitably high $j,k \in \N$.}\label{f.kappaisone}
\end{figure}



As we will elaborate, numerical approximations of the map~$\lambdamax(\kappa,\rho)$ 
suggest that this convergence occurs, and that it does so rapidly.
 They also reveal some surprises.
The function $(0,1] \to [1,\infty): \kappa \mapsto \lambdamax(\kappa,1)$  has left limit equal to one according to Theorem~\ref{t.lowlambdamax}. Although low~$\kappa$-values are intractable numerically, the function appears to be increasing as $\kappa$ rises into a numerically determinable range; it takes values close to  $1 + 10^{-5}$ and  $1 + 8 \cdot 10^{-5}$ at the respective $\kappa$-values $0.65$ and $0.9$. A maximum value close to $1 + 9.92 \times 10^{-5}$
is achieved close to $\kappa=0.976$, with a slight but monotonic decrease witnessed until $1 + 9.68 \times 10^{-5}$ is seen at  $\kappa=1$. (The latter point corresponds to the game in~\cite{LostPennies}: see the bounds recalled in~(\ref{e.lambdamaxbounds}).)
The behaviour appears unimodal, with a peak a little to the left of~$\kappa=1$. 
It would be natural enough to expect variation in~$\rho$ at $\kappa=1$ to behave in a roughly similar way, and in Figure~\ref{f.kappaisone} numerics for the map $(0,1] \to [1,\infty): \rho \mapsto \lambdamax(1,\rho)$ are shown.  This function appears to be maximized at $\rho =1$ (consistently with Conjecture~\ref{c.lambdamax}(2)), and to tend rapidly to one as~$\rho$ falls. But the function is obviously not monotone. Of course its behaviour on $[0.96,0.97]$
compels a higher-digit numerical review there.
 Astonishingly, $\lambdamax(1,\rho) = 1$ appears to have an isolated solution~$p$ that lies in the interval $[0.964556,0.964557]$:  at $p$, the Mina margin map~$\mc{M}_{1,\rho}$ becomes identically equal to one, with its argument maximizer  in the central domain  jumping discontinuously as $\rho$ passes through this value. As $\rho$ varies through the special point,~$\mc{M}_{1,\rho}$ resembles a vibrating string that becomes straight for an instant; as we explain in a moment, this comparison appears to endure as $\rho$ decreases.

\begin{figure}[p]
\centering

\includegraphics[width=0.55\textwidth]{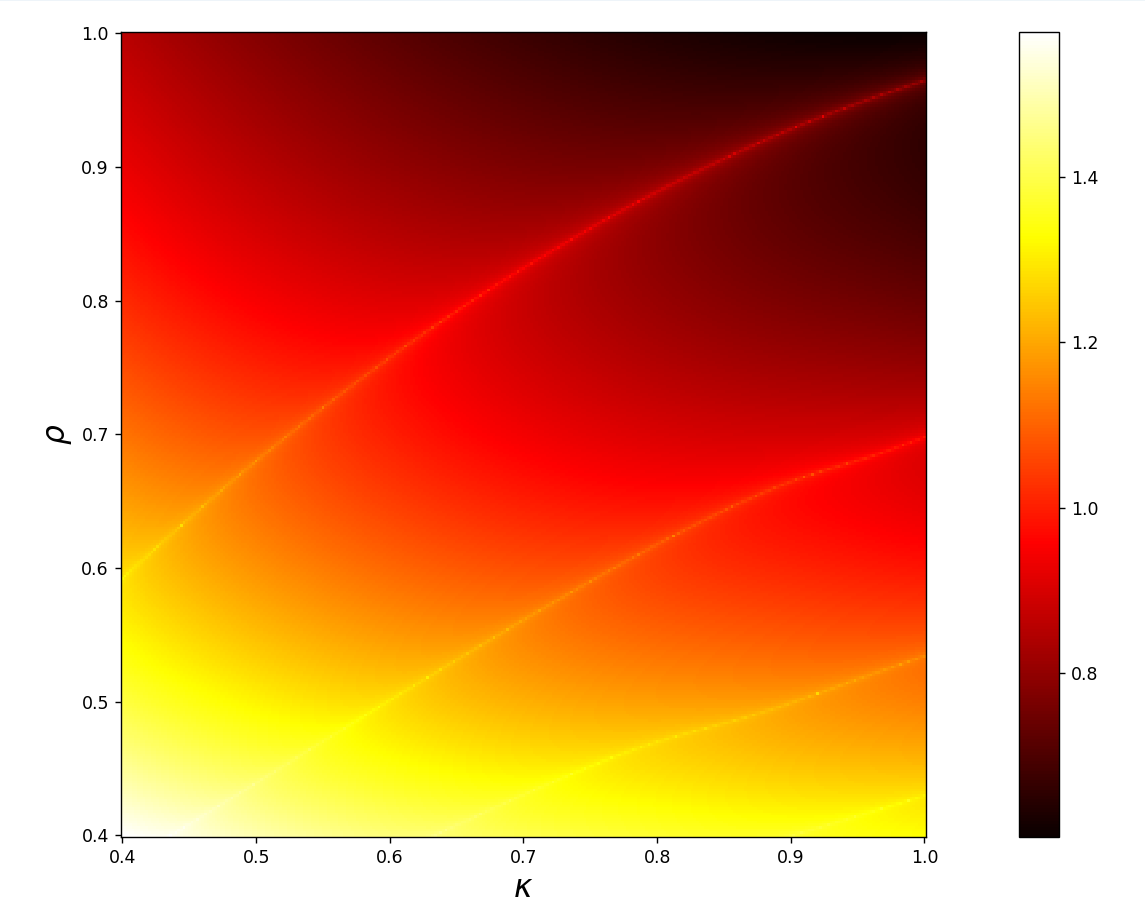}

\caption{A heat map of the function $(\kappa,\rho) \mapsto \lambdamax(\kappa,\rho) -1$
in the box $[0.4,1]^2$ in which the scale for colour aligns with the iterated logarithm of the plotted values.
Specifically, the numbers labelling the colour scale on the right represent the value of $\log_{10} \log_{10} \frac{1}{\lambdamax(\kappa,\rho) -1}$, so that lighter colours correspond to extremely small values:
for example,  $\lambdamax(\kappa,\rho)-1 \approx 10^{-25}$ for $(\kappa,\rho)$-points coloured with the $1.4$-indexed yellow.
From the point~$(1,\rho_1)$, with $\rho_1 \approx 0.96$, which is a root of $\lambdamax - 1$ depicted in Figure~\ref{f.kappaisone},
there emerges a zero-locus path charting a roughly WSW course. At progressively smaller ambient values, three further such paths chart similar routes from $(1,\rho_i)$,
where $\rho_2 \approx 0.70$, $\rho_3 \approx 0.53$, and $\rho_4 \approx 0.43$.}
\label{f.doublelogheatmap}

\vspace{1em}

\includegraphics[width=0.45\textwidth]{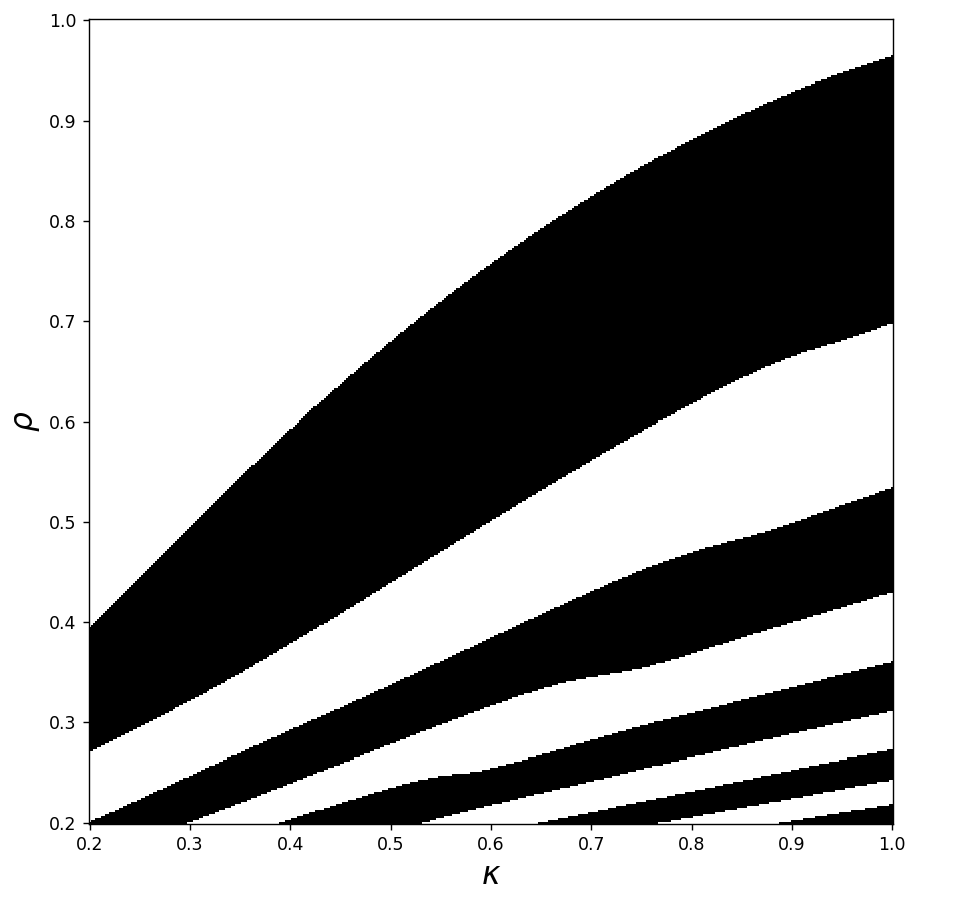}

\caption{For any $(\kappa,\rho) \in (0,1]^2$, the Mina margin map satisfies $\mc{M}_{\kappa,\rho}(1)=1$. In this plot, a point $(\kappa,\rho)$ is marked white if \(1\) is a point of local increase of the map \(x \mapsto \mc{M}_{\kappa,\rho}(x)\), and black if it is a point of local decrease. Along the boundaries of the coloured regions, the map is identically equal to one; these boundaries therefore form part of the zero locus of $\lambdamax(\kappa,\rho)-1$. This plot on $[0.2,1]^2$, which requires roughly three hundred digits of base-ten precision, reveals five further curves in the zero locus beyond those identified in the preceding figure, for a total of nine. The alternating pattern of colour, as \(\rho\) decreases for fixed \(\kappa\), reflects successive phase reversals in the local passage of the graph \(x \mapsto \mc{M}_{\kappa,\rho}(x)\) through~\((1,1)\).}
\label{f.blackandwhite}
\end{figure}

In this way, the locus $\lambdamax(\kappa,\rho) =1$ of parameter pairs where no incentive asymmetry is permitted (so that the discouragement effect is infinitely strong) not only contains one (and perhaps the other) axis; it also appears to contain the point~$(1,p)$, directly south of $(1,1)$ by about four percent. Further numerical investigation, illustrated by Figure~\ref{f.doublelogheatmap}, indicates a path in the locus that starts at~$(1,p)$ and moves roughly west-by-southwest through the $(\kappa,\rho)$-box $(0,1]^2$, passing through $[0.83,0.84]\times \{ 0.9 \}$ and  
$[0.66,0.67]\times \{ 0.8 \}$ and turning gradually to the left onto a more southwesterly course. Moreover,  several further zero-locus paths leave  the $\kappa=1$ line from a succession of $\rho$-values. Recording these points in the form $\rho_i$, with $\rho_1 = p$ as above,  the next entries satisfy  $(\rho_2,\rho_3,\rho_4) \approx (0.70,0.53,0.43)$.
Viewed as a unit-time movie in which the point $(1,\rho)$ travels south along the $\kappa=1$ line, the graph of the Mina margin map $(0,\infty) \to (0,\infty): x \mapsto \mc{M}_{1,\rho}(x)$ resembles a string vibrating through several oscillations.
The spatial frequency of this string rises as $\rho = 0$ is approached, since 
 the central domain $\big( (2-\rho)/(2+\rho), (2+\rho)/(2-\rho) \big]$ contracts to the point~$1$.  The string's amplitude decays rapidly; conclusions are tentative because the range $\rho < 0.2$ remains unchecked,  but numerics suggest that $\lambdamax(1,\rho)$ is approaching one at a rate that is roughly doubly exponential in~$1/\rho$. 
 And the string's temporal frequency also appears to rise, with the special values~$\rho_i$ at which $\mc{M}_{1,\rho}$ is identically equal to one becoming closer together: nine such values have been identified, with the last two satisfying $\rho_8 \approx 0.24$ and $\rho_9 \approx 0.22$.
 
  \FloatBarrier
\FloatBarrier

Figure~\ref{f.blackandwhite} offers a depiction of numerical data that highlights the alternating phases of the Mina margin map, according in essence as to whether the derivative of $x \mapsto \mc{M}_{\kappa,\rho}(x)$ is positive or negative at the value $x=1$ (where this map is always equal to one). 
Nine boundary paths are visible in the depiction, each of which lies in the zero-locus of $\lambdamax(\kappa,\rho) - 1$. In principle, these paths
may bifurcate, join or terminate in their onward journey. The depiction, however, suggests another possibility: that  they continue disjointly, reaching close to or even arriving at the origin~$(0,0)$.
 
 How many roots does the function $(0,1] \to [0,\infty): \rho \mapsto \lambdamax(1,\rho) -1$ have? Apparently, at least nine. Plausibly there are more; perhaps even infinitely many.
 Numerical testing in the $(\kappa,\rho)$-region~$[0.2,1]^2$
 depicted in Figure~\ref{f.blackandwhite}
 is computationally substantial, not only because $\lambdamax(\kappa,\rho)-1$ assumes extremely small values---the value $\lfloor - \log_{10} \big( \lambdamax(0.2,0.2)-1 \big) \rfloor$ is estimated to equal $270$, for example---but also because the approximation of the Mina margin map by finite trail counterparts requires longer interval integers in both the regimes of low $\kappa$ and of low~$\rho$. 
 As such,  there are limits on the degree to which direct numerical testing can further reveal these strange effects.
Naturally, it would be most interesting to explain them theoretically.


{\bf D-TOUR.}
For $d \geq 1$, let
$x:[0,\infty) \to \R^d$ with $x(0) = 0$ satisfy 
$$
\dot{x}(t) \, = \, v(t) \, \, \, \,  \textrm{and} \, \, \, \, \dot{v}(t) \, = \,   F(t) - v    +  \dot{B}_t   \, ,
$$
 with $B$ standard $d$-dimensional Brownian motion.
The trajectory~$x$
  models a small flying vehicle agitated by thermal fluctuations in the ambient air
  and
  subject to both an applied force $F:[0,\infty) \to \R^d$  and aerodynamic drag.
  (This is the Ornstein-Uhlenbeck process in its original physical guise~\cite{OrnsteinUhlenbeck}, where noise acts on the velocity, and with a force applied. Dilations of space and time permit  the diffusivity and linear-drag coefficients to equal one.)
   The
    Dual-Thrust Ornstein-Uhlenbeck Rocket
   comes equipped with two thrusters whose strength and direction may be adjusted independently, under the respective control of two players.
The D-TOUR trajectory~$x$ begins statically at a given point in a domain $D \subset \R^d$.  At time~$t \geq 0$, the applied force is a superposition of thrusts
 $$
 F(t) \, = \, \psi \big( a(t) \big) V_+(t) + \psi \big( b(t) \big) V_-(t) \, , 
 $$
 where, at this time, Maxine\footnote{The names are less apposite unless $d =1$ since the players no longer seek necessarily to maximize or minimize~$x$.} nominates stake rate $a(t) \in [0,\infty)$ and a  direction vector~$V_+(t)$ valued in the Euclidean unit sphere~$S^{d-1}$, while Mina nominates $b(t)$ and $V_-(t)$. 
  The map $\psi:[0,\infty) \to [0,\infty)$ is the magnitude of the thrust offered by a player as a function of her spending rate; it may be supposed to be increasing and convex
 and to vanish at zero, 
   with the choices $\psi(z) = z^\rho$ for $\rho \in (0,1]$ seeming natural.  
  The domain
  boundary comes equipped with functions $f,g: \partial D \to \R$, and the
  game ends when the rocket~$x$ reaches $\partial D$ at time $\tau$, with total net receipt $g(x_\tau) - \int_0^\tau b(t) \, {\rm d}t$ for Mina and $f(x_\tau) - \int_0^\tau a(t) \, {\rm d}t$ for Maxine. 
  
  It would be interesting to study this more physically natural game to see if the conclusions we have  reached for $\rho$-Brownian Boost---the fragility of equilibria under slight changes in relative incentive; the presence of a battlefield zone; the asymmetry in decay away from that zone---are borne out. Such a study could also be contemplated for 
  a variety of discrete-time or stochastic differential games governed by stakes. 
  
  In the Moscarini-Smith game discussed in Subsection~\ref{s.ms}, 
  drift is determined from stake-pair via the rule $a^\rho - b^\rho$ with $\rho =1/2$ (in place of our $(a^\rho - b^\rho )/(a^\rho + b^\rho)$); for a range of~$\rho$-values, this game
   may be examined analytically, or perhaps exhibited as a scaling limit of discrete models (with long-range jumps being one means of producing unbounded drift).   Tug-of-war in discrete and continuous time 
   may be played on more general graphs, or in higher dimensions; these geometrically richer settings likely present probabilistic and analytic subtleties that invite further exploration.

\chapter*{Glossary of notation}\label{c.glossary}

In three tables, we list much of the article's main notation, providing a summarizing phrase for each item, as well as the page number at which the concept is introduced. 
The first table offers principal examples concerning strategy pairs and the associated gameplay in the discrete games. The second also concerns discrete contexts,  where an \abmnkpmacspace element is supposed given.  The third records notation for the ODE apparatus around Brownian Boost.

\vspace{5mm}

\large{\bf The discrete games: strategies, gameplay and equilibria}

\vspace{1mm}

\small{
\bigskip
\def\qq{&}
\begin{center}
\halign{
#\quad\hfill&#\quad\hfill&\quad\hfill#\cr
$\llbracket i,j \rrbracket$ \qq the integer interval $\Z \cap [i,j]$  \fff{intint}
\tlpkp \qq the Trail of Lost Pennies game indexed by $(\kappa,\rho) \in (0,1) \times (0,\infty)$ \fff{tlpkp}
$\kappa \in (0,1]$ \qq the probability that a given turn in \tlpkpspace is stake \fff{kappa}
$\rho \in (0,\infty)$ \qq the Tullock exponent governing the contest at a stake turn \fff{rho}
$\Lambda_k$ \qq the space of nearest-neighbour paths $\psi: \llbracket 0,k \rrbracket \to \Z$ \fff{pathspace}
$\mc{S}$ \qq the space of strategies $S:\Lambda = \bigcup_{k=0}^\infty \Lambda_k \to [0,\infty)$ \fff{strategy}
$\tis$ \qq the space of time-invariant strategies \fff{timeinvariant}
$(S_-,S_+) \in \tis^2$ \qq \typnot \fff{typicalnotation}
$X$ \qq counter evolution $X: \N \to \Z$ \fff{counter}
$P_-$, $T_-$ and $C_-(t)$ \qq Mina's total, and terminal, receipts and time-$t \in \N_+$ running cost \fff{receiptmina}
$P_+$, $T_+$ and $C_+(t)$ \qq Maxine's counterparts \fff{receiptmaxine}
$E$, $E_-$, $E_+$ \qq \escapes \fff{escape}
$\PP^i_{S_-,S_+}$ and $\E^i_{S_-,S_+}$ \qq  \lawandexpectation \fff{lawexpect}
TINE \qq {\bf T}ime-{\bf I}nvariant (or time-homogeneous Markov-perfect) {\bf N}ash {\bf E}quilibrium \fff{tine}
$\mc{N_{\kappa,\rho}}$ \qq the set of Nash equilibria in \tlpkp \fff{nash}
$m_i$ and $n_i$ \qq Maxine and Mina's mean total receipt for play from $i$ under a given TINE
 \fff{mini}
$a_i$ and $b_i$ \qq the players' stakes in this scenario  \fff{aibi}
$S_- = b$ and $S_+ = a$ \qq a notational abuse for $(S_-,S_+) \in \tis^2$; $a_i$ and $b_i$ are stakes at site~$i$ \fff{notab}
 $\mc{I}$ \qq \idlemac \fff{idle}
 wide \qq an element of $\tis$ of infinite support \fff{wide}
}\end{center}}

\newpage

\medskip

\vspace{4mm}

\large{\bf \abmnkpmacspace elements}

\vspace{1mm}

\small{
\bigskip
\def\qq{&}
\begin{center}
\halign{
#\quad\hfill&#\quad\hfill&\quad\hfill#\cr
$(a,b,m,n):\Z \to (0,\infty)^4$ \qq an element of  \abmnkpmac, namely a positive \abmnmacspace solution \fff{positiveabmn}
$\big( m_{-\infty},m_\infty,n_{-\infty},n_\infty \big)$ \qq the boundary data for such a solution \fff{boundarydata}
$m_{i,j},n_{j,i}$ \qq differences $m_j - m_i$ and $n_i - n_j$, for $i < j$ \fff{mndiffernces}
$W \subset (0,1] \times (0,\infty)$ \qq the subset of $(\kappa,\rho)$-space where $\rho^2 \kappa \leq 1$ 
\fff{weakregion}
Mina margin \qq the value $\frac{n_{-\infty} - n_\infty}{m_\infty - m_{-\infty}}$ \fff{minamargin}
$\lambdamax(\kappa,\rho)$ \qq maximum value of Mina's winning payoff for which TINE exist \fff{lambdamax}
 $\phi_0$,$\phi_1$,$\gamma$,$c = 1/\gamma$,$\delta$,$d = 1/\delta$ \qq formulaically specified positive functions of $(\kappa,\rho,\beta)$, with $\beta \in (0,\infty)$ \fff{fourfunctions}
central ratio \qq a name for $\frac{n_{-1} - n_0}{m_0 - m_{-1}}$ \fff{centralratio}
$s:(0,\infty) \to (0,\infty)$ \qq the map that sends $\phi_0$ to $\phi_1$, for any  \abmnkpmacspace element \fff{smap}
standard and default \qq distinguished elements in symmetry-equivalence classes of \abmnkpmacspace    \fff{defaultstandard} 
 $\phi_i$ \qq the quantity $\frac{n_{i-1} - n_i}{m_i - m_{i-1}}$ \fff{phi}
  $D = D_{\kappa,\rho}$ \qq  the central domain $\big( \frac{2 - \kappa \rho}{2+\kappa \rho}, \frac{2 + \kappa \rho}{2-\kappa \rho} \big]$  \fff{domain}
battlefield index \qq the unique index $i \in \Z$ such that $\phi_i \in D$ \fff{battlefield}
$\mc{M}:(0,\infty) \to (0,\infty)$ \qq \macmmm  \fff{mmm}
$\mc{M}_{j+1,k+1}:(0,\infty) \to (0,\infty)$ \qq \counterpart \fff{mmm.finite}
 $M_i$, $N_i$,
$\beta_i$, $\gamma_i$, $c_i$, $\delta_i$, $d_i$ \qq statistics local to site $i$ associated to an \abmnmacspace solution  \fff{gammaess}
}\end{center}}

\large{\bf Brownian Boost: continuous variables and the ODE apparatus}

\vspace{1mm}

\small{
\bigskip
\def\qq{&}
\begin{center}
\halign{
#\quad\hfill&#\quad\hfill&\quad\hfill#\cr
$X_t$ \qq counter location at time~$t$ \fff{X}
\bbrho \qq the Brownian Boost game with Tullock exponent $\rho \in (0,\infty)$ \fff{bbrho}
$\lambda \in (0,\infty)$ \qq Mina's winning receipt in \bbrho \fff{lambda}
the \bbrho-ODE pair \qq the coupled ODE system characterizing equilibrium in \bbrho \fff{bbode}
$(f,g):\R \to (0,\infty)^2$ \qq a solution of the \bbrho-ODE pair  \fff{fg}
$\xmac \in (0,\infty)$ \qq the flow index that parametrises ODE solutions \fff{flowindex}
$S_\rho(\xi,\bullet)$ \qq the ODE flow mapping $\R$ to $(0,\infty)$, with $S_\rho(\xmac,0) = \xmac$  \fff{srho}
$\big(f_\rho(\xmac,\bullet),g_\rho(\xmac,\bullet)\big)$ \qq explicit (default) solution of the \bbrho-ODE pair indexed by $\xi$ \fff{fgrho}
$\big(a_\rho(\xmac,\bullet) ,b_\rho(\xmac,\bullet) \big)$ \qq the accompanying stake-rate profiles  \fff{abrho}
standard \qq normalization of \bbrho-ODE pair solutions for which   $\int_\R f(r) \, {\rm d} r =1$ \fff{bbstandard} 
}\end{center}}

\vspace{-1mm}

\bibliographystyle{plain}

\bibliography{stake}

\end{document}